\documentclass[11pt]{article}

\pdfoutput=1

\usepackage{amscd,amsmath, amssymb, fancyhdr, epsfig,url,color,mathtools}

\usepackage{graphicx}
\usepackage{amsfonts}
\usepackage{indentfirst}
\usepackage{amsthm, amsgen, amstext, amsbsy, amsopn,yfonts,bbm}
\usepackage{xcolor}
\usepackage{mathrsfs}
\numberwithin{equation}{section}

\usepackage{tocloft}

\cftsetindents{section}{5em}{3em}
\cftsetindents{subsection}{5em}{3em}

\usepackage[backref=page]{hyperref}
\renewcommand*{\backrefalt}[4]{%
	\ifcase #1 (Not cited.)%
	\or        (Cited on page~#2.)%
	\else      (Cited on pages~#2.)%
	\fi}
\hypersetup{
	colorlinks   = true,
	citecolor    = magenta,
	linkcolor    =blue,
	urlcolor     =magenta	
}

\newcommand{\version}{v. 2.2, April 29, 2023}

\def\eqref#1{(\ref{#1})}

\newcommand{\Z}{{\mathbb Z}}
\newcommand{\C}{{\mathbb C}}
\newcommand{\F}{{\mathbb F}}

\newcommand{\mP}{{\mathbb P}}
\newcommand{\R}{{\mathbb R}}
\newcommand{\Q}{{\mathbb Q}}

\def\1{\sqrt{-1}\:}

\newcommand{\cntrct}                
{\hspace{2pt}\raisebox{1pt}{\text{$\lrcorner$}}\hspace{2pt}}

\newcommand{\Id}{\operatorname{Id}}

\newcommand{\Sp}{\operatorname{Sp}}

\newcommand{\Spin}{\operatorname{Spin}}

\renewcommand{\Re}{\operatorname{Re}}
\renewcommand{\Im}{\operatorname{Im}}

\newcommand{\Per}{\operatorname{\sf Per}}

\newcommand{\Gal}{\operatorname{Gal}}

\newcounter{Mycounter}[section]
\newcounter{lemma}[section]
\renewcommand{\thelemma}{{Lemma \thesection.\arabic{lemma}}}
\newcommand{\lemma}{%
    \setcounter{lemma}{\value{Mycounter}}
    \refstepcounter{lemma}
    \stepcounter{Mycounter}
    {\noindent \bf \thelemma\ }}

\newcounter{claim}[section]
\newcounter{sublemma}[section]
\newcounter{corollary}[section]
\renewcommand{\thecorollary}{{Corollary \thesection.\arabic{corollary}}}
\newcommand{\corollary}{%
    \setcounter{corollary}{\value{Mycounter}}
    \refstepcounter{corollary}
    \stepcounter{Mycounter}
    {\noindent \bf \thecorollary\ }}

\newcounter{theorem}[section]
\renewcommand{\thetheorem}{{Theorem \thesection.\arabic{theorem}}}
\newcommand{\theorem}{%
    \setcounter{theorem}{\value{Mycounter}}
    \refstepcounter{theorem}
    \stepcounter{Mycounter}
    {\noindent \bf \thetheorem\ }}

\newcounter{conjecture}[section]
\newcounter{proposition}[section]
\renewcommand{\theproposition}
      {{Proposition \thesection.\arabic{proposition}}}
\newcommand{\proposition}{%
    \setcounter{proposition}{\value{Mycounter}}
    \refstepcounter{proposition}
    \stepcounter{Mycounter}
    {\noindent \bf \theproposition\ }}

\newcounter{definition}[section]
\renewcommand{\thedefinition}
      {{Definition~\thesection.\arabic{definition}}}
\newcommand{\definition}{%
    \setcounter{definition}{\value{Mycounter}}
    \refstepcounter{definition}
    \stepcounter{Mycounter}
    {\noindent \bf \thedefinition\ }}

\newcounter{example}[section]
\renewcommand{\theexample}{{Example \thesection.\arabic{example}}}
\newcommand{\example}{%
    \setcounter{example}{\value{Mycounter}}
    \refstepcounter{example}
    \stepcounter{Mycounter}
    {\noindent \bf \theexample\ }}

\newcounter{remark}[section]
\renewcommand{\theremark}{{Remark \thesection.\arabic{remark}}}
\newcommand{\remark}{%
    \setcounter{remark}{\value{Mycounter}}
    \refstepcounter{remark}
    \stepcounter{Mycounter}
    {\noindent \bf \theremark\ }}

\newcounter{problem}[section]
\newcounter{question}[section]
\renewcommand{\thequestion}{{Question \thesection.\arabic{question}}}
\newcommand{\question}{%
    \setcounter{question}{\value{Mycounter}}
    \refstepcounter{question}
    \stepcounter{Mycounter}
    {\noindent \bf \thequestion\ }}

\newcounter{condition}[section]
\makeatletter

\@addtoreset{equation}{section} \@addtoreset{footnote}{section} \makeatother

\def\blacksquare{\hbox{\vrule width 5pt height 5pt depth 0pt}}
\def\endproof{\blacksquare}

\newcommand{\T}{{\mathbb{T}}}
\newcommand{\SP}{{\mathbb{S}}}

\newcommand{\bJ}{{\overline{J}}}

\newcommand{\bomega}{{\overline{\omega}}}

\newcommand{\barf}{{\bar{f}}}

\newcommand{\bl}{{\bar{l}}}

\newcommand{\bbw}{{\bar{\bf w}}}

\newcommand{\bphi}{{\bar{\phi}}}

\newcommand{\bfc}{{\mathbbm c}}

\newcommand{\bfr}{{\bar{r}}}

\newcommand{\bfx}{{\mathbbm x}}
\newcommand{\bfy}{{\mathbbm y}}

\newcommand{\bfU}{{\mathbb U}}
\newcommand{\bfW}{{\mathbb W}}
\newcommand{\bfZ}{{\mathbb Z}}

\newcommand{\cA}{{\mathcal{A}}}
\newcommand{\cB}{{\mathcal{B}}}

\newcommand{\cC}{{\mathcal{C}}}

\newcommand{\cF}{{\mathcal{F}}}

\newcommand{\cH}{{\mathcal{H}}}

\newcommand{\cJ}{{\mathcal{J}}}
\newcommand{\cK}{{\mathcal{K}}}
\newcommand{\cL}{{\mathcal{L}}}
\newcommand{\cM}{{\mathcal{M}}}

\newcommand{\cO}{{\mathcal{O}}}

\newcommand{\cR}{{\mathcal{R}}}
\newcommand{\cS}{{\mathcal{S}}}
\newcommand{\cT}{{\mathcal{T}}}
\newcommand{\cU}{{\mathcal{U}}}
\newcommand{\cV}{{\mathcal{V}}}
\newcommand{\cW}{{\mathcal{W}}}
\newcommand{\cX}{{\mathcal{X}}}
\newcommand{\cY}{{\mathcal{Y}}}

\newcommand{\hJ}{{\hat{J}}}

\newcommand{\hM}{{\hat{M}}}

\newcommand{\goh}{{\mathfrak{h}}}
\newcommand{\gol}{{\mathfrak{l}}}

\newcommand{\goH}{{\mathfrak{H}}}

\newcommand{\goU}{{\mathfrak{U}}}

\newcommand{\Alb}{{{\rm Alb}\,}}
\newcommand{\Span}{{{\rm Span}\,}}

\newcommand{\tB}{{\widetilde{B}}}
\newcommand{\tC}{{\widetilde{C}}}

\newcommand{\tI}{{\widetilde{I}}}
\newcommand{\tJ}{{\widetilde{J}}}

\newcommand{\tM}{{\widetilde{M}}}

\newcommand{\tR}{{\widetilde{R}}}

\newcommand{\tU}{{\widetilde{U}}}
\newcommand{\tV}{{\widetilde{V}}}
\newcommand{\tX}{{\widetilde{X}}}

\newcommand{\ttK}{{\tt{K}}}
\newcommand{\ttU}{{\tt{U}}}

\newcommand{\tildeh}{{\widetilde{h}}}

\newcommand{\tiota}{{\widetilde{\iota}}}

\newcommand{\tomega}{{\widetilde{\omega}}}

\newcommand{\teta}{{\widetilde{\eta}}}
\newcommand{\ttheta}{{\widetilde{\theta}}}
\newcommand{\talpha}{{\widetilde{\alpha}}}

\newcommand{\EmbS}{\textit{Emb}}
\newcommand{\EmbK}{\textit{EK}}

\newcommand{\Pairs}{{\it Pairs}}
\newcommand{\Triples}{{\it Triples}}

\newcommand{\Teich}{{\tt Teich}}
\newcommand{\THodge}{{\tt T}}

\newcommand{\TeichP}{{\tt Pairs}}
\newcommand{\TeichT}{{\tt Triples}}

\newcommand{\Diff}{\mathrm{Diff}}
\newcommand{\Symp}{\mathrm{Symp}}
\newcommand{\SympH}{\mathrm{Symp}_H}
\newcommand{\symp}{\mathit{S}}
\newcommand{\sympK}{\mathit{SK}}

\newcommand{\Lin}{{\it Lin}}
\newcommand{\Lincmpt}{{\it Lin}}
\newcommand{\comp}{\mathit{C}}
\newcommand{\compK}{\mathit{CK}}

\newcommand{\Int}{Int\,}

\newcommand{\cmpt}{{\it Cmpt}}
\newcommand{\cmptteich}{{\tt Teich}}

\newcommand{\Image}{{\rm{Im}\,}}

\newcommand{\Vol}{{\rm{Vol}\,}}

\newcommand{\pr}{{\textit{pr}\,}}

\newcommand{\Pos}{{\rm{Pos}\,}}
\newcommand{\Dom}{{\it D}}

\begin{document}

\begin{center}
{\LARGE\bf
K\"ahler-type embeddings of balls into symplectic manifolds\\
\ \\
}

Michael Entov,
Misha Verbitsky\footnote{
  Partially supported by HSE University Basic Research Program,
FAPERJ E-26/202.912/2018 and CNPq - Process 310952/2021-2. }

\end{center}


\hfill
\hfill

\begin{abstract}
\small

Consider a symplectic embedding of a disjoint union of domains (lying in the standard symplectic $\R^{2n}$) into a symplectic manifold $M$.
 We say that such
 an embedding is K\"ahler-type, or respectively tame, if it is holomorphic with respect to some (not a priori fixed, K\"ahler-type)
 integrable
 complex structure on $M$ compatible with the symplectic form, or respectively tamed by it. Assume that $M$ is either of the following: a complex projective space (with the standard symplectic form);
an even-dimensional torus, or a K3 surface, equipped with an irrational K\"ahler-type symplectic form. Then any two K\"ahler-type embeddings of a disjoint union of balls into $M$ can be mapped into each other by a symplectomorphism acting trivially on homology. If the embeddings are holomorphic with respect to complex structures compatible with the symplectic form and lying in the same connected component of the space of K\"ahler-type complex structures on $M$, then
the symplectomorphism can be chosen to be smoothly isotopic to the identity. For certain $M$ and certain disjoint unions of balls we describe precisely the obstructions to the existence of K\"ahler-type embeddings of the balls into $M$. In particular, symplectic volume is the only obstruction for the existence of K\"ahler-type embeddings of $l^n$ equal balls (for any $l$) into $\C P^n$ with the standard symplectic form and of any number of possibly different balls into a torus or a K3 surface, equipped with an irrational symplectic form.
We also show that symplectic volume is the only obstruction for the existence of tame embeddings of disjoint unions of equal balls, polydisks, or parallelepipeds, into a torus equipped with a generic K\"ahler-type symplectic form. For balls and parallelepipeds the same is true for K3 surfaces.

\end{abstract}


{\small
\tableofcontents
}


\section{Introduction}
\label{_introduction_Section_}


Symplectic embeddings of balls and other domains into symplectic manifolds have been studied since the early days of symplectic topology.
So far, the existence of the embeddings (the so-called symplectic packing problem) got most of attention and much less is known
about the connectedness of the spaces of such embeddings. The notable exception
is the classical result of McDuff \cite{_McD-connectedness_}  about the connectedness of the space of symplectic embeddings of disjoint unions of balls into certain 4-dimensional symplectic manifolds -- in particular, into $\C P^2$, $\C P^1\times \C P^1$ and their blow-ups (see also \cite{_McD-Topology1991_, _McD-Cambridge1990_,_Lalonde-MathAnn-1994_,_Biran-IMRN1996_} for earlier partial results on the subject and
 \cite{_McDuff-ellipsoids-JT2009_},
 \cite{_Borman-Li-Wu_},
\cite{_Crist-Gardiner-JDG2019_}
 for more recent extensions of McDuff's result).

In this paper, in an attempt to get more information about the connectedness question for symplectic manifolds of dimension higher than 4, we consider a more restrictive class of
symplectic embeddings of disjoint unions of domains.

Namely, let $M$ be a
manifold equipped with a
{\bf K\"ahler-type symplectic form} $\omega$ -- i.e., a symplectic form compatible with some ({\sl not a priori fixed})
integrable complex structure on $M$; such complex structures will be called {\bf K\"ahler-type complex structures}.
We will denote by $\compK (M)$ the space of K\"ahler-type complex structures on $M$.

We consider symplectic embeddings (of disjoint unions of domains) into $(M,\omega)$ that are holomorphic with respect to some {\it (not a priori fixed)} integrable (K\"ahler-type) complex structure $I$ on $M$ which is compatible with $\omega$, or tamed by it.
We will say that such a symplectic embedding is of {\bf K\"ahler type}, or, respectively, {\bf tame}.
The existence of K\"ahler-type embeddings, {\sl for a fixed isotopy class of the complex structure $I$ on $M$ compatible with $\omega$}, has been previously studied -- in a different guise --
in \cite{_Eckl2017_}, \cite{_Fleming2021_}, \cite{_Luef-Wang_}, \cite{_Trusiani-AIF2021_}, \cite{_WittNystrom1_}, \cite{_WittNystrom2_}, see \ref{_I-Kahler-type-embs-under-diff-names_Remark_}, \ref{_CPn-l-to-n-balls-flat-Kahler-metric-prev-results_Remark_} and part 3 of \ref{_existence-connectedness-tori-K3-comments_Remark_}.

With this terminology, we have the following hierarchy of properties of the embeddings:
\[
\textrm{K\"ahler-type}\
\Longrightarrow\ \textrm{tame}
\Longrightarrow\ \textrm{symplectic}.
\]
A K\"ahler-type embedding of a domain in $\R^{2n}$ (equipped with the standard symplectic and complex structures) into $(M,\omega)$ which is holomorphic with respect to a complex structure $I$ compatible with $\omega$ can be equivalently described as a K\"ahler-isometric, or just K\"ahler, embedding of the domain into $M$ equipped with the K\"ahler metric $\omega (\cdot, I\cdot) + \sqrt{-1} \omega (\cdot, \cdot)$.
In other words,  an embedding of a domain into $(M,\omega)$ is of K\"ahler type if and only if it is K\"ahler with respect to some K\"ahler metric whose imaginary part is $\omega$. The restriction of the K\"ahler metric to the image of the embedding is the standard flat metric.

There do exist symplectic embeddings of balls that are not of K\"ahler type (see below). We do not know whether these
symplectic embeddings are tame or whether there exist tame embeddings (of balls or any other domains) that are not of K\"ahler-type -- see \ref{_sympl-tame-Kahler-type_Remark_} and
\ref{_existence-of-sympl-not-tame-or-tame-not-K-type-embs_Question_}. Let us note that for a disjoint union of starshaped domains (with a piecewise-smooth boundary)
its K\"ahler-type and tame
embeddings form $\Symp (M)$-invariant open subsets of the space of its symplectic embeddings into $(M,\omega)$ (with respect to the $C^\infty$-topology) -- see Section~\ref{_Diff0-cap-Symp-action-on-embeddings_Subsection_}.

As far as the connectedness of the spaces of K\"ahler-type embeddings is concerned, let us note first that a K\"ahler embedding into $(M,\omega)$ may be holomorphic with respect to different complex structures compatible with $\omega$, and, second, that two K\"ahler-type embeddings lying in the same $\Symp_0 (M,\omega)$-orbit  have to be holomorphic with respect to complex structures lying in the same orbit of the $\Symp_0 (M,\omega)$-action on the space $\compK (M)$ of K\"ahler-type complex structures on $M$,
hence in the same connected component of the set of complex structures on $M$ compatible with $\omega$.

Unfortunately, we know almost nothing about the latter connected components or about the space $\compK (M)/\Symp_0 (M,\omega)$. Because of this, we cannot say much about the orbits of the $\Symp_0 (M,\omega)$-action on the space of K\"ahler-type embeddings of
given domains into $(M,\omega)$. However, we can say more about the actions of $\Symp (M,\omega)\cap\Diff_0 (M)$
and of the symplectic Torelli group $\SympH (M,\omega)$ (the group of symplectomorphisms of $(M,\omega)$ acting trivially on the homology)
on the spaces of K\"ahler-type embeddings of balls into $(M,\omega)$.

Namely, assume that $(M,\omega)$ is either of the following:
$\C P^n$ with the standard symplectic form; $\T^{2n}$ or a K3 surface equipped with a K\"ahler-type symplectic form $\omega$ which is irrational (i.e., the real cohomology class $[\omega]\in H^2 (M;\R)$ is not a real multiple of a rational one). Then
any two
K\"ahler-type embeddings of an arbitrary disjoint union of balls
into $(M,\omega)$ can be mapped into each other by
the action of $\SympH (M,\omega)$. If these embeddings are holomorphic
with respect to complex structures compatible with $\omega$ and lying in the same connected component of
$\compK (M)$, then these embeddings can be mapped into each other by
the action of $\Symp (M,\omega)\cap \Diff_0 (M)$. In particular, these two K\"ahler-type embeddings are holomorphic with respect to isotopic
complex structures compatible with $\omega$.

Note that if $(M,\omega)$ is a rational or ruled closed symplectic 4-manifold, or if $(M,\omega)$ is obtained from a rational or ruled closed symplectic 4-manifold by removing either a closed 2-dimensional symplectic submanifold or a Lagrangian submanifold diffeomorphic to $\SP^2$ or $\R P^2$ (for instance, if $(M,\omega)$ is the standard symplectic ball), then
any two symplectic
embeddings of an arbitrary disjoint union of closed balls
into $(M,\omega)$ can be mapped into each other by
the action of $\Symp^c_0 (M,\omega)$ -- see \cite{_McD-connectedness_}, \cite{_Borman-Li-Wu_} (cf. \cite{_McD-Topology1991_, _McD-Cambridge1990_,_Lalonde-MathAnn-1994_,_Biran-IMRN1996_} for earlier partial results on the subject).
In particular, this is true for any two K\"ahler-type, or tame, embeddings and thus the space of such embeddings is connected (see \ref{_connectedness-rational-ruled-mfds_Theorem_} below).

Let us now discuss our results on the existence of K\"ahler-type and tame embeddings.

Following the terminology in \cite{_EV-JTA_}, we say that K\"ahler-type,
or tame,
embeddings of a disjoint union $\bfW=\bigsqcup_{i=1}^k W_i$ of compact domains with piecewise boundary into $(M,\omega)$ are {\bf unobstructed} if there are no obstructions for the existence of such embeddings
of $\lambda\bfW = \bigsqcup_{i=1}^k \lambda W_i$, $\lambda>0$, into $(M,\omega)$  other than the symplectic volume. (The domains we consider
always lie in the standard symplectic $\R^{2n}$.) Note that if K\"ahler-type embeddings of a disjoint union of domains into $(M,\omega)$ are unobstructed, this means, in particular, that one can choose a complex structure $I$ on $M$ compatible with $\omega$ so that the K\"ahler metric $\omega (\cdot, I\cdot) + \sqrt{-1} \omega (\cdot, \cdot)$ on $M$ is flat outside a set of arbitrarily small volume.

Since the classical work of McDuff-Polterovich \cite{_McD-Polt_}, it has been known that proving the existence of symplectic embeddings of closed balls $\bigsqcup_{i=1}^k B(r_i)$ of radii $r_1,\ldots,r_k$ into a closed symplectic manifold $(M,\omega)$ is equivalent to showing that an appropriate cohomology class of the complex blow-up $\tM$ of $M$ at $k$ points can be represented by a symplectic form with suitable properties.

In this paper we extend the results of McDuff and Polterovich \cite{_McD-Polt_} to K\"ahler-type embeddings of $\bigsqcup_{i=1}^k B(r_i)$. Namely, assume $M$ is closed. For each $\bfx=(x_1,\ldots,x_k)$, $x_1,\ldots,x_k\in M$, and each complex structure $I$ on $M$, denote by
$\tM_{I,\bfx}$ the complex blow-up of $(M,I)$ at $x_1,\ldots,x_k$.
Let $\tI$ be the lift of $I$ to $\tM$.
We show that $\bigsqcup_{i=1}^k B(r_i)$ admits a K\"ahler-type embedding into $(M,\omega)$
if and only if for some complex structure $I$ compatible with $\omega$ and some $\bfx$
the cohomology class $\Pi^* [\omega] - \pi \sum_{i=1}^k r_i^2 e_i\in H^2 (\tM_{I,\bfx};\R)$ is K\"ahler with respect to $\tI$ (i.e., can be represented by a K\"ahler form on $(\tM_{I,\bfx},\tI)$). In fact, if this cohomology class is K\"ahler, then $\bigsqcup_{i=1}^k B(r_i)$ admits a K\"ahler-type embedding into $(M,\omega)$ which is holomorphic with respect to a complex structure on $M$ compatible with $\omega$ and smoothly isotopic to $I$.
Together with information on the K\"ahler cones of $(\tM_{I,\bfx},\tI)$ for various $\tI$ (compatible with $\omega$) and $\bfx$ this allows us to show the existence of K\"ahler-type embeddings of $\bigsqcup_{i=1}^k B(r_i)$ into certain $(M,\omega)$ (see Section~\ref{_embs-of-balls-kahler-classes-on-blowups_Subsection_}).

Here are some of the applications for particular symplectic manifolds (see Section~\ref{_main-results-all_Section_}):

If $(M,\omega)$ is either the standard symplectic $\C P^n$ or a standard symplectic ball of real dimension $2n$, then for any $l\in\Z_{>0}$, K\"ahler-type embeddings of a disjoint union of $l^n$ equal balls are unobstructed. In particular, this implies that for the standard complex structure on $\C P^n$ and
 $l^n$ equal balls of total
 symplectic volume smaller than the volume of $\C P^n$, there exists a K\"ahler form on $\C P^n$, isotopic to the standard Fubini-Study form, so that the resulting K\"ahler manifold admits a {\sl K\"ahler} (i.e., both holomorphic and symplectic) embedding of the $l^n$ equal balls -- meaning, in particular, that the K\"ahler metric is a standard flat metric on the image of the balls.

For $M=\C P^2$ we also describe exactly which disjoint unions of $1\leq k\leq 8$ (possibly different) balls admit K\"ahler-type embeddings into $M$.
 It then follows from \cite{_McD-connectedness_} that any symplectic embedding of such a disjoint union of balls into $\C P^2$ is, in fact, of K\"ahler type.

In the cases where $M$ is $\C P^1\times \C P^1$ or a complex blow-up of $\C P^2$ at one point equipped with
any symplectic form compatible with the orientation, we describe exactly which disjoint unions of $1\leq k\leq 7$ (possibly different) balls admit K\"ahler-type embeddings into $M$ that are holomorphic with respect to a complex structure isotopic to the ``standard" one and compatible with the symplectic form.
 It again follows from \cite{_McD-connectedness_} that any symplectic embedding of such a disjoint union of balls into $M$ is, in fact, of K\"ahler type.

If $M$ is either $\T^{2n}$ or a K3 surface and the K\"ahler-type symplectic form $\omega$ on $M$ is irrational (i.e., the real cohomology class $[\omega]\in H^2 (M;\R)$ is not a real multiple of a rational one), then K\"ahler-type embeddings of any disjoint union of balls into $(M,\omega)$ are unobstructed. If $\omega$ is rational (i.e., $[\omega]\in H^2 (M;\R)$ is a real multiple of a rational class),
then
there may be obstructions for the existence of K\"ahler-type embeddings of disjoint unions of balls into $(M,\omega)$ that are independent of the symplectic volume. For instance, in the case $M=\T^4$ there are obstructions coming from the Seshadri constants  (see \ref{_Seshadri-constants_Remark_}). In particular, this shows that there exist symplectic embeddings of disjoint unions of balls into $(\T^4,\omega)$ that are not K\"ahler-type.

Finally, we show that tame embeddings of disjoint unions of equal balls, or equal polydisks, or equal parallelepipeds, into $(\T^{2n},\omega)$
are unobstructed provided $\omega$ lies in a $\Diff_0 (\T^{2n})$-invariant open dense set of K\"ahler-type symplectic forms on $\T^{2n}$ containing all the irrational forms (the set depends on the embedded domains).
By choosing possibly smaller $\Diff_0 (\T^{2n})$-invariant open dense sets of such $\omega$ we make sure that the unobstructedness holds, in fact, for the tame embeddings into $(\T^{2n},\omega)$ which are holomorphic with respect to complex structures that are not only tamed by $\omega$ but also ``arbitrarily close" (in the sense of Hodge theory) to complex structures compatible with $\omega$.
K\"ahler-type (not just tame) embeddings of disjoint unions of equal balls, or equal polydisks, or equal parallelepipeds, into $(\T^{2n},\omega)$ are unobstructed if $\omega$
is an irrational K\"ahler-type symplectic form lying in an appropriate dense $\Diff_0 (\T^{2n})$-orbit in the space of forms (the orbit depends on the embedded domains). Similar results concerning balls and parallelepipeds hold also for symplectic manifolds underlying K\"ahler K3 surfaces.

We also discuss extensions of the results for K3 surfaces to a larger class of hyperk\"ahler manifolds -- see Section~\ref{_hyperkahler-case_Subsection_}.


\hfill


\section{Preliminaries}
\label{_preliminaries_Section_}

Assume that $M$ is a smooth manifold of real dimension $2n$.

Denote by $\Diff (M)$ the group of diffeomorphisms of $M$ and by $\Diff_0 (M)$ its identity component, and by $\Diff_H (M)$ the subgroup of $\Diff (M)$ formed by the diffeomorphisms acting trivially on the homology of $M$.

For a symplectic form $\omega$ on $M$ denote by $\Symp (M,\omega)$ the group of symplectomorphisms of $(M,\omega)$, by $\Symp_0 (M,\omega)$
its identity component and by $\SympH (M,\omega)$ the {\bf symplectic Torelli group}, i.e. the group of symplectomorphisms of $(M,\omega)$ acting trivially on homology. For a non-compact symplectic manifold $(M,\omega)$ we denote by $\Symp^c (M,\omega)$
the group of compactly supported symplectomorphisms of $(M,\omega)$ and by $\Symp^c_0 (M,\omega)$ its identity component.


\hfill


\subsection{K\"ahler-type complex and symplectic structures}
\label{_Kahler-type-complex-sympl-strs_Subsection_}

Let $\omega$ and $I$ be a symplectic form and a complex structure on $M$. The form $\omega$ {\bf tames} $I$ if $\omega (v,Iv)>0$ for any non-zero $v\in TM$, and is {\bf compatible with $I$} if it tames $I$ and is also $I$-invariant,
meaning that $\omega (\cdot, I\cdot) + i\omega (\cdot, \cdot)$ is a Hermitian metric on $(M,I)$. In the latter case the pair $(\omega,I)$ is called
a {\bf K\"ahler structure} on $M$.

A cohomology class in $H^2 (M;\R)$ is called {\bf K\"ahler (with respect to a complex structure $I$ on $M$)}, if it can be represented by a symplectic form on $M$ compatible with $I$. The cohomology classes that are K\"ahler with respect to $I$ form the {\bf K\"ahler cone of $(M,I)$} in $H^2 (M;\R)$.

A symplectic/complex structure on $M$ is said to be of {\bf K\"ahler type}, if it appears in {\it some} K\"ahler structure on $M$.


\hfill


\example
\label{_Kahler-type-structures-torus_Example_}

Any complex structure on $\C P^2$ is of K\"ahler type \cite{_Yau-PNAS1977_}. Whether this holds also for $\C P^n$, $n>2$, is
a well-known open question -- see e.g. \cite{_Libgober-Wood-JDG1990_,_Tosatti-CPn_}, cf. \cite{_Huckleberry-Kebekus-Peternell_}.

For a closed manifold $M$, $\dim_\C M =
2$, with an even $b_2 (M)$, any complex structure on $M$ is of K\"ahler type \cite{_Buchdahl-AnnInstFourier1999_, _Lamari-AnnInstFourier1999_}. In particular,
any complex structure on $\C P^1\times \C P^1$, $\C P^2\sharp \overline{\C P^2}$ and $\T^4$ is of K\"ahler type.

The K\"ahler-type symplectic/complex structures on $\T^{2n}=\R^{2n}/\Z^{2n}$ are exactly the ones that can be mapped
by a diffeomorphism of $\T^{2n}$ to a linear symplectic/complex structure -- i.e., a symplectic/complex structure whose lift to $\R^{2n}$
has constant coefficients
(see e.g. \cite[Prop. 6.1]{_EV-JTA_}).
There do exist complex structures on $\T^{2n}$, for each $n\geq 3$, that are not of K\"ahler-type \cite[p.212]{_Sommese-Quaternionic-MathAnn-1975_}, \cite{_Borisov-Salamon-Viaclovsky-Duke2011_}.

If $M$ is a smooth manifold underlying a complex K3 surface, then any
complex structure on $M$ is of K\"ahler-type \cite{_Siu_} and
its first Chern class is zero (see e.g. \cite{_Geom-K3-Asterisque_})\footnote{In fact, all smooth manifolds $M$ underlying a complex K3 surface
are diffeomorphic. They all are compact and connected. All complex structures on such an $M$ define the
same orientation on $M$ so that $b_+ (M)=3$ and
$b_- (M)=19$ -- see e.g. \cite{_Geom-K3-Asterisque_}.
Further on, we will always equip a K3 surface with this
standard orientation. Alternatively, to fix the orientation,
one can use the fact that all complex K3 surfaces are deformation equivalent and, in
particular, oriented diffeomorphic -- see e.g. \cite{_Geom-K3-Asterisque_}.}.
It is well-known that any complex
structure on such an $M$ appears in a hyperk\"ahler
structure and so does any K\"ahler-type symplectic form
on $M$  -- see e.g.
\cite[Prop. 5.1]{_EV-K3_}.

It is a well-known open question whether all symplectic forms on $\T^{2n}$ and K3 surfaces (equipped with the standard smooth structures) are of K\"ahler type -- see e.g.
\cite{_Donaldson:ellipt_}. (There do exist exotic smooth structures on topological K3 surfaces that admit symplectic forms -- see \cite{_Fintushel-Stern-InventMath1998_}, cf. \cite{_Chen-ProcAMS-2020_}; it follows from the classification of complex surfaces that these exotic smooth K3 surfaces do not admit complex structures and therefore the symplectic forms on them are not of K\"ahler type). Using the Kodaira-Spencer stability theorem \cite{_Kod-Spen-AnnMath-1960_} and Torelli theorems for tori (see e.g. \cite[Ch. I, Thm. 14.2]{_Barth-Hulek-Peters-vdVen_}) and K3 surfaces (see \cite{_Burns-Rapoport_}, cf. \cite[p.96]{_Geom-K3-Asterisque_}) one can show that the K\"ahler-type forms on these manifolds form an open subset of the space of all symplectic forms (with respect to the $C^\infty$-topology -- see e.g. \cite{_Am-Ver-JGF2015_}, cf. \cite{_Streets-Tian_} for a more general result).


\hfill


Denote by $\sympK (M)$ the space of all K\"ahler-type symplectic forms on $M$ and by $\compK (M)$ the space of all K\"ahler-type complex structures on $M$.

We equip these spaces with the $C^\infty$-topologies (the complex structures are viewed here as integrable almost-complex structures -- i.e., $(1,1)$-tensors).

The space $\compK (M)$ is an open subspace of the space of all complex structures on $M$; both of these spaces are
locally $C^\infty$-path connected -- see \ref{_compK-path-connected-Kahler-coh-classes_Proposition_}.

Given $\omega\in \sympK (M)$, denote by $\cmpt (M,\omega)\subset \compK (M)$ the space of all (K\"ahler-type) complex structures on $M$ compatible with $\omega$.

If $\cC_0$ is a connected component of $\compK (M)$, we say that $\omega$ and
$\cC_0$ are {\bf compatible}, if some complex structure lying in $\cC_0$ is compatible with $\omega$
 -- i.e., if the set $\cC_0\cap \cmpt (M,\omega)$ is non-empty.

With this terminology, every K\"ahler-type $\omega$ is compatible with some connected component of $\compK (M)$.

Note that the space $\compK (M)$ is invariant under the natural action of $\Diff (M)$ (and hence of $\Diff_0 (M)$) on the space of complex structures on $M$, but the space
$\cmpt (M,\omega)$ is not invariant even under the $\Diff_0 (M)$-action.


\hfill


\definition
\label{_Teich_Definition}

Define the {\bf Teichm\"uller space of K\"ahler-type complex structures} on $M$ as
\[
\Teich (M) := \compK (M)/\Diff_0 (M),
\]
and if $\cC_0$ is a connected component of $\compK (M)$, then
\[
\Teich_{\cC_0} (M) := \cC_0/\Diff_0 (M).
\]
The set $\Teich_{\cC_0} (M)$ is a connected component of $\Teich (M)$.

Given $I\in\compK (M)$, denote by $[I]\in\Teich (M)$ the image of $I$ under the projection $\compK (M)\to \Teich (M)$.

Let $\omega\in\sympK (M)$. Denote the image of $\cmpt (M,\omega)$ under
the quotient projection $\compK (M)\to \Teich (M)$
by
\[
\cmptteich (M,\omega) :=\frac{\Diff_0 (M)\cdot \cmpt (M,\omega)}{\Diff_0 (M)}\subset \Teich (M).
\]
In other words, $\cmptteich (M,\omega)$ is the set of points of $\Teich (M)$ that can be represented by complex structures
on $M$ compatible with $\omega$.

Given a connected component $\cC_0$ of $\compK (M)$, the image of $\cmpt (M,\omega)\cap \cC_0$ under
the quotient projection $\cC_0\to \Teich_{\cC_0} (M)$
will be denoted by
\[
\cmptteich_{\cC_0} (M,\omega) :=\frac{\Diff_0 (M)\cdot \left(\cmpt (M,\omega)\cap \cC_0\right)}{\Diff_0 (M)} =
\]
\[
= \cmptteich (M,\omega)\cap\Teich_{\cC_0} (M).
\]


\hfill


Note that the space $\Teich (M)$ may not be Hausdorff -- this is the case, for instance, when $M$ is a smooth manifold underlying a complex K3 surface, see \ref{_Teich-K3-not-Hausdorff_Remark_}.

\subsection{Hodge decomposition}
\label{_Hodge decomposition_Subsection_}

Assume $M$ is closed.

A complex structure $I$ on $M$ defines a $(p,q)$-decomposition of any (comp\-lex-valued) differential $k$-form $\eta$ on $M$, $k=0,1,\ldots, 2n$ -- we will denote it by $\eta = \sum_{p+q=k} \eta_I^{p,q}$.
If $I$ is of K\"ahler type, then the $(p,q)$-decomposition on differential forms, defined by $I$, induces a decomposition, called the Hodge decomposition\footnote{Although the proof of the
existence of the Hodge decomposition involves both the complex structure $I$ and a K\"ahler form on $(M,I)$, one can show (see e.g. \cite[Vol.
1, Prop. 6.11]{_Voisin-book_}) that, in fact, the Hodge
decomposition on the cohomology depends only on $I$.}, on the complex cohomology of $M$: $H^k (M;\C) = \oplus_{p+q=k} H_I^{p,q} (M;\C)$.
Given a cohomology class $\alpha\in H^* (M;\C)$, we will denote its $(p,q)$-component in the Hodge decomposition by $\alpha_I^{p,q}\in H_I^{p,q} (M;\C)$.
We denote $H^{p,q}_I (M;\R) := H^{p,q}_I (M;\C)\cap H^* (M;\R)$.
For a real cohomology class $\alpha\in H^{2p} (M;\R)$ its $(p,p)$-component is also real.
A differential form on $M$ is called a {\bf $(p,q)$-form}, and a cohomology class of $M$ is called a {\bf $(p,q)$-class} (with respect to $I$), if their $(p,q)$-component is their only non-zero component.

If $\omega$ is compatible with $I$, then it is a $(1,1)$-form with respect to $I$, and, accordingly, its cohomology class is a $(1,1)$-class.

If a cohomology class $\alpha\in H^2 (M;\R)$ is K\"ahler with respect to a K\"ahler-type complex structure $I$ on $M$, then
for any K\"ahler-type complex structure $I'$ sufficiently $C^\infty$-close to $I$ the class $\alpha_{I'}^{1,1}$ is K\"ahler with respect to $I'$ -- see \ref{_compK-path-connected-Kahler-coh-classes_Proposition_}.


\hfill


\subsection{K\"ahler-type and tame embeddings}
\label{_partial-Kahler-type-embeddings-preliminaries_Subsection_}

Let $\omega$ be a K\"ahler-type symplectic form on $M$.

The symplectic and complex structures on $\R^{2n}$ are always assumed to be the standard ones and will be denoted, respectively, by $\omega_0$ and $J_0$.

Let $W_1,\ldots, W_k\subset \R^{2n}$  be compact domains with piecewise-smooth boundaries whose interiors contain the origin. Denote
\[
\bfW:= \bigsqcup_{i=1}^k W_i.
\]
Given $\lambda>0$, define
\[
\lambda \bfW:= \bigsqcup_{i=1}^k \lambda W_i.
\]

We say that $f: \bfW\to M$ is a smooth/symplectic/holomorphic embedding if it extends to such an embedding of the disjoint union of some open neighborhoods of $W_i$, $i=1,\ldots,k$. We will denote by $f_i: W_i\to M$, $i=1,\ldots,k$, the restriction of such an extension to
 a neighborhood of $W_i$.

Clearly, if $f: \bfW\to (M,\omega)$ is a symplectic embedding, then
\[
\textstyle \Vol (\bfW,\omega_0) < \Vol (M,\omega) = \langle [\omega]^n, [M]\rangle.
\]
(Recall that $\Vol$ stands for the symplectic volume and $\langle \cdot, \cdot \rangle$ stands for the natural pairing between cohomology and homology classes.)


\hfill


\definition
\label{_Kahler-type-partially-Kahler-type-embs_Definition_}

Let $f: \bfW = \bigsqcup_{i=1}^k W_i\to (M^{2n},\omega)$
be an embedding.

Given a K\"ahler-type complex structure $I\in\cmpt (M,\omega)$, we say that the embedding $f$ is of {\bf $[I]$-K\"ahler type}, if it is symplectic with respect to $\omega$ and holomorphic with respect to some {\it  (not a priori fixed)}  complex structure $J$ on $M$ compatible with $\omega$ and isotopic to $I$, i.e. $[J]=[I]\in \cmptteich (M,\omega)$. (Such a $J$ is tautologically of K\"ahler-type).

We say that the embedding $f$ is of {\bf K\"ahler type}, if it is
symplectic with respect to $\omega$ and holomorphic with respect to some {\it  (not a priori fixed)} complex structure compatible with $\omega$ --
or, in other words, if $f$ is of $[I]$-K\"ahler type for some {\it  (not a priori fixed, and possibly non-unique)} $I\in\cmpt (M,\omega)$.

Given a {\it fixed K\"ahler structure $(\omega, I)$} on $M$, we say that  $f: \bfW\to (M^{2n},\omega, I)$ is a {\bf K\"ahler} embedding if it is K\"ahler-isometric -- i.e. both symplectic with respect to $\omega$ and holomorphic with respect to $I$.

We say that the embedding  $f: \bfW\to (M^{2n},\omega)$ is {\bf tame} if it is symplectic with respect to $\omega$ and holomorphic with respect to some {\it  (not a priori fixed)} K\"ahler-type complex structure $I$ tamed by $\omega$.

Assuming that $M$ is closed
and $\varepsilon>0$, we say that the embedding $f$ is {\bf $\varepsilon$-tame} if it is symplectic with respect to $\omega$ and holomorphic with respect to some {\it  (not a priori fixed)} K\"ahler-type complex structure $I$ so that

\begin{itemize}

\item{} $I$ is tamed by $\omega$,

\item{} the class
$[\omega]_I^{1,1}$ is K\"ahler,

\item{}
$\left|\left\langle \left([\omega]_I^{2,0} + [\omega]_I^{0,2}\right)^n, [M]\right\rangle\right| <\varepsilon$.

\end{itemize}


\hfill


With the terminology as in \ref{_Kahler-type-partially-Kahler-type-embs_Definition_}, for an embedding into a (closed) symplectic manifold
$(M,\omega)$
we have the following equivalences and hierarchy of properties:
\begin{center}

K\"ahler-type $\Longleftrightarrow$ $[I]$-K\"ahler type for some (possibly non-unique!) $I\in\cmpt (M,\omega)$,

\end{center}
\begin{center}

tame $\Longleftrightarrow$ $\varepsilon$-tame for some $\varepsilon$,

\end{center}
\begin{center}

K\"ahler-type $\Longrightarrow$ $\varepsilon$-tame for every $\varepsilon$ $\Longrightarrow$ tame $\Longrightarrow$ symplectic.

\end{center}

In general, there do exist symplectic embeddings of balls that are not of K\"ahler-type (see \ref{_Seshadri-constants_Remark_}),  but we do not know whether there exist
symplectic embeddings of balls or any other domains into any $M$, equipped with a K\"ahler-type symplectic form, that are not tame, or tame embeddings that are not of K\"ahler-type -- see
\ref{_sympl-tame-Kahler-type_Remark_}
and
\ref{_existence-of-sympl-not-tame-or-tame-not-K-type-embs_Question_}.


\hfill


\remark
\label{_existence-of-Kahler-type-embs-of-small-balls_Remark_}

For any $I\in\cmpt (M,\omega)$ any finite disjoint union of sufficiently small balls (and hence of arbitrary sufficiently small domains in $(\R^{2n},\omega_0)$) admits an $[I]$-K\"ahler-type embedding into $(M,\omega)$ -- this follows e.g. from \cite[Prop. 5.5.A]{_McD-Polt_}).


\hfill


\remark

If a K\"ahler manifold $(M,\omega, I)$ admits a K\"ahler embedding of a domain in $\R^{2n}$, its restriction to the image of the embedding is flat. This provides a strong obstruction to the existence of K\"ahler embeddings (if one fixes a K\"ahler structure on the manifold).


\hfill


\remark
\label{_I-Kahler-type-embs-under-diff-names_Remark_}

The existence of $[I]$-K\"ahler-type embeddings has been previously studied -- in a different guise -- in \cite{_Eckl2017_}, \cite{_Fleming2021_}, \cite{_Luef-Wang_}, \cite{_Trusiani-AIF2021_}, \cite{_WittNystrom1_}, \cite{_WittNystrom2_} in the following setup.

Let $(M,\omega,I)$ be a closed K\"ahler manifold (in \cite{_Eckl2017_}, \cite{_Fleming2021_}, \cite{_Trusiani-AIF2021_}, \cite{_WittNystrom1_}, \cite{_WittNystrom2_} it is assumed that the complex manifold $(M,I)$ is projective and the cohomology class $[\omega]\in H^2 (M;\R)$ is the first Chern class of an ample line bundle on $(M,I)$. Consider the following question: Given a disjoint union $\bfW = \bigsqcup_{i=1}^k W_i$ of domains in $\R^{2n}$,
does there exist a K\"ahler form $\eta$ on $(M,I)$, isotopic to $\omega$, such that $(M,\eta,I)$ admits a K\"ahler embedding of $\bfW$?

One easily sees that the answer to the latter question is positive if and only if $(M,\omega)$ admits an $[I]$-K\"ahler-type embedding of $\bfW$.


\hfill


Let us introduce additional terminology related to K\"ahler-type
and tame
embeddings.

Let $\cJ\subset \compK (M)$ be a non-empty set. If a K\"ahler-type, or tame,
embedding $f$ is holomorphic with respect to some complex structure belonging to $\cJ$, we say that $f$ {\bf favors $\cJ$}.

In particular, if $\cJ = \cC_0$ is a connected component $\cC_0$ of $\compK (M)$, we say that such an $f$ {\bf favors $\cC_0$}, or that $\cC_0$ is a {\bf favorite connected component of $f$}.

With this terminology,
tautologically,

\begin{itemize}

\item{} For $I\in\cmpt (M,\omega)$, every $[I]$-K\"ahler-type embedding into $(M,\omega)$ favors the $\Diff_0 (M)$-orbit of $I$ in $\compK (M)$.

\item{} Every K\"ahler-type
embedding into $(M,\omega)$ favors a connected component of $\cmpt (M,\omega)$, and consequently a connected component of $\compK (M)$.

\end{itemize}

Note that, in principle, a K\"ahler-type embedding may be of $[I]$-K\"ahler-type for several $I\in\cmpt (M,\omega)$ and may favor more than one connected component of $\cmpt (M,\omega)$ or of $\compK (M)$.

The group $\Symp (M,\omega)$ acts naturally on the space of K\"ahler-type, or tame, or $\varepsilon$-tame,
embeddings $\bfW\to (M,\omega)$. Consequently, so do the groups $\SympH (M,\omega)$, $\Symp_0 (M,\omega)$ and $\Symp (M,\omega)\cap \Diff_0 (M)$.

By a standard result concerning extensions of symplectic isotopies (proved as in \cite[Thms. 3.3.1, 3.3.2]{_McD-Sal-3_}),
 if the domains $W_1,\ldots,W_k$ are starshaped, then
 any two symplectic embeddings $f: \bfW\to (M,\omega)$, $f': \bfW\to (M,\omega)$
can be connected by a (smooth) path of symplectic embeddings $\bfW\to (M,\omega)$ if and only if $f$ and $f'$
lie in the same $\Symp_0 (M,\omega)$-orbit -- i.e., there exists $\phi\in \Symp_0 (M,\omega)$ such that $\phi \circ f = f'$.

Also note, that if $\phi \circ f = f'$ for some $\phi\in\Diff (M)$ and $f$ is holomorphic with respect to a complex structure $I$ on $M$, then $f'$ is holomorphic with respect to the complex structure $\phi_* I$.

These simple observations easily yield the following corollary.


\hfill


\proposition
\label{_two-embs-in-same-conn-compt_Proposition_}

 Assume that the domains $W_1,\ldots,W_k$ are starshaped.

Then any two embeddings $f, f': \bfW = \bigsqcup_{i=1}^k W_i\to (M,\omega)$
 lie in the same connected component of the space of K\"ahler-type
embeddings $\bfW\to (M,\omega)$ if and only if they lie in the
same orbit of the $\Symp_0 (M,\omega)$-action.
In this case $f,f'$ are holomorphic with respect to complex structures compatible
with $\omega$ that lie in the same orbit of the $\Symp_0 (M,\omega)$-action on $\compK (M)$. In particular,
such $f,f'$ favor a common connected component of $\cmpt (M,\omega)$
and, consequently, a common connected component of $\compK (M)$.

If K\"ahler-type
embeddings $f, f': \bfW\to (M,\omega)$ lie in the same $\Symp (M,\omega)\cap \Diff_0 (M)$-orbit, then they are holomorphic with respect to complex structures compatible
with $\omega$ that lie in the same orbit of the $\Diff_0 (M)$-action on $\compK (M)$. In particular, there exists $[I]\in\cmptteich (M,\omega)$
such that both $f$ and $f'$ are of $[I]$-K\"ahler-type, and $f,f'$ favor a common connected component of $\compK (M)$.
\endproof


\hfill


Let us define a quantity that will be important for the study of the portion of the volume of a target closed symplectic manifold that can be filled out by K\"ahler-type and tame embeddings.


\hfill


\definition
\label{_nu_Definition_}

Assume $M$ is a closed manifold equipped with a K\"ahler-type symplectic form $\omega$.

Define the number $0<\nu_K (M,\omega, \bfW)\leq 1$ by
\[
\nu_K (M,\omega,\bfW) := \sup_\lambda \frac{\Vol (\lambda\bfW,\omega_0)}{\Vol (M,\omega)},
\]
where the supremum is taken over all $\lambda>0$ such that $\lambda \bfW$ admits a K\"ahler-type embedding into $(M,\omega)$.

Given
$\varepsilon>0$, define the number $0<\nu_{T,\varepsilon} (M,\omega, \bfW)\leq 1$ by
\[
\nu_{T,\varepsilon} (M,\omega,\bfW) := \sup_\lambda \frac{\Vol (\lambda\bfW,\omega_0)}{\Vol (M,\omega)},
\]
where the supremum is taken over all $\lambda>0$ such that $\lambda \bfW$ admits an $\varepsilon$-tame embedding into $(M,\omega)$.
(In view of \ref{_existence-of-Kahler-type-embs-of-small-balls_Remark_}, the set of such $\lambda$ is always non-empty.)


\hfill


\remark

\noindent
1. Clearly,
\[
\nu_K (M,\omega,\bfW)\leq \nu_{T,\varepsilon_1} (M,\omega,\bfW)\leq \nu_{T,\varepsilon_2} (M,\omega,\bfW)
\]
for any $0<\varepsilon_1\leq \varepsilon_2$.

\bigskip
\noindent
2. The equality $\nu_K (M,\omega,\bfW) = 1$, respectively $\nu_{T,\varepsilon} (M,\omega,\bfW) =1$, is equivalent
to the claim that the K\"ahler-type, respectively $\varepsilon$-tame, embeddings of $\lambda\bfW$ into $(M,\omega)$, for all $\lambda>0$, are unobstructed.


\hfill


\section{Main results}
\label{_main-results-all_Section_}

In this section we state the main results of the paper.

For $r>0$ we denote by $B^{2n} (r)$ the closed Euclidean $2n$-dimensional ball of radius $r$ in $\R^{2n}$ centered at the origin.


\hfill


\subsection{K\"ahler-type embeddings of balls in general symplectic manifolds 
and K\"ahler classes on complex blow-ups}
\label{_embs-of-balls-kahler-classes-on-blowups_Subsection_}

In this section $M$, $\dim_\R = 2n$, is a closed connected manifold equipped with a K\"ahler-type symplectic form $\omega$.
Let $k\in\Z_{>0}$, $\bfr= (r_1,\ldots,r_k)\in (\R_{>0})^k$.

Denote by $\hM^k\subset M^k$ the set of $k$-tuples of pairwise distinct points of $M$.

For each $\bfx:= (x_1,\ldots,x_k)\in \hM^k$
and each complex structure $I$ on $M$, let $\tM_{I,\bfx}$
denote the complex blow-up of $(M,I)$ at $x_1,\ldots,x_k$.

For any $I$ and $\bfx$ we will use the same notation for the following objects (suppressing the dependence on $I$ and $\bfx$):

\smallskip
\noindent
- The lift of $I$ to $\tM_{I,\bfx}$ will be denoted by $\tI$.

\smallskip
\noindent
- $\Pi: \tM_{I,\bfx}\to M$ will denote the natural projection.

\smallskip
\noindent
- The cohomology classes Poincar\'e-dual to the fundamental homology classes of the exceptional divisors will be denoted by $e_1,\ldots,e_k\in H^2 (\tM_{I,\bfx};\R)$.

For all $I$ and $\bfx$ the groups $H^2 (\tM_{I,\bfx};\R)$ can be canonically identified (preserving the numbering of the classes $e_1,\ldots,e_n$) and we will use these identifications without further mention.


\hfill


\definition
\label{_cK-ttK_Definition}

Define the set
$\cK (\bfr)\subset \compK (M)\times \hM^k$ as
\[
\cK (\bfr) := \Big\{\ (I,\bfx)\ \big|\ I\in\Diff_0 (M)\cdot\cmpt (M,\omega), \bfx\in\hM^k\ \textit{and}
\]
\[
\Pi^* [\omega] - \pi \sum_{i=1}^k r_i^2 e_i\in H^2 (\tM_{I,\bfx};\R)\ \textit{is K\"ahler with respect to}\ \tI\ \Big\},
\]

Let $\pr: \cK (\bfr)\to \compK (M)$ be the natural projection.

Define the set $\ttK (\bfr)\subset \cmptteich (M,\omega)$ as
\[
\ttK (\bfr) := \pr (\cK (\bfr))/\Diff_0 (M) =
\]
\[
= \Big\{\ \iota\in \cmptteich (M,\omega)\ \big|\ \exists I\in \cmpt (M,\omega), [I]=\iota,\ \textit{and}\ \bfx\in \hM^k\ \textit{s.t.}
\]
\[
\Pi^* [\omega] - \pi \sum_{i=1}^k r_i^2 e_i\in H^2 (\tM_{I,\bfx};\R)\ \textit{is K\"ahler with respect to}\ \tI\ \Big\}.
\]

Let $\cC_0$ be a connected component of $\compK (M)$.

Define the set
$\cK_{\cC_0} (\bfr)\subset \cC_0\times \hM^k$ as
\[
\cK_{\cC_0} (\bfr) := (\cC_0\times \hM^k)\cap \cK (\bfr) =
\]
\[
= \Big\{\ (I,\bfx)\ \big|\ I\in\Diff_0 (M)\cdot \big(\cC_0\cap \cmpt (M,\omega)\big), \bfx\in\hM^k\ \textit{and}
\]
\[
\Pi^* [\omega] - \pi \sum_{i=1}^k r_i^2 e_i\in H^2 (\tM_{I,\bfx};\R)\ \textit{is K\"ahler with respect to}\ \tI\ \Big\}.
\]
Define the set $\ttK_{\cC_0} (\bfr)\subset \cmptteich_{\cC_0} (M,\omega)$ as
\[
\ttK_{\cC_0} (\bfr) := \ttK (\bfr)\cap\Teich_{\cC_0} (M) = \pr \left(\cK_{\cC_0} (\bfr)\right)/\Diff_0 (M) =
\]
\[
= \Big\{\ \iota\in \cmptteich_{\cC_0} (M,\omega)\ \big|\ \exists I\in \cC_0\cap \cmpt (M,\omega), [I]=\iota,\ \textit{and}\ \bfx\in \hM^k\ \textit{s.t.}
\]
\[
\Pi^* [\omega] - \pi \sum_{i=1}^k r_i^2 e_i\in H^2 (\tM_{I,\bfx};\R)\ \textit{is K\"ahler with respect to}\ \tI\ \Big\}.
\]


\hfill


As we will see in \ref{_fibers-of-projection-cK-cC-zero-are-connected_Corollary_}, the fibers of $\pr: \cK (\bfr)\to \compK (M)$ are connected and therefore $\ttK (\bfr)$ (respectively, $\ttK_{\cC_0} (\bfr)$)
is connected if and only if so is $\cK (\bfr)$ (respectively, $\cK_{\cC_0} (\bfr)$).


\hfill


\theorem
\label{_existence-connectedness-Kahler-type-arb-mfds_Theorem_}

Let $M$ be a closed manifold equipped with a K\"ahler-type symplectic form $\omega$.

\bigskip
\noindent
(I)
The set $\bigsqcup_{i=1}^k B^{2n}(r_i)$ admits a K\"ahler-type embedding into $(M,\omega)$ if and only if
for some $I\in \cmpt (M,\omega)$ and some $\bfx\in \hM^k$ the cohomology class $\Pi^* [\omega] - \pi \sum_{i=1}^k r_i^2 e_i\in H^2 (\tM_{I,\bfx};\R)$ is K\"ahler with respect to $\tI$.

More precisely, assume $I$ is a complex structure on $M$ compatible with $\omega$ and $\Sigma\subset (M,I)$ is a proper (possibly empty) complex submanifold. If $M = \C P^n\times \C P^1$ (with the product complex structure), we also allow
$\Sigma = (\C P^{n-1}\times \C P^1) \cup (\C P^n\times \textrm{pt})$ (which is not smooth).

Then the following conditions are equivalent:

\begin{itemize}

\item{} There exists $\bfx = (x_1,\ldots,x_k)\in \hM^k$, $x_1,\ldots,x_k\in M\setminus\Sigma$, such that the cohomology class $\Pi^* [\omega] - \pi \sum_{i=1}^k r_i^2 e_i\in H^2 (\tM_{I,\bfx};\R)$ is K\"ahler with respect to $\tI$.

\item{} There exists an $[I]$-K\"ahler-type embedding $f: \bigsqcup_{i=1}^k B^{2n}(r_i)\to (M\setminus\Sigma,\omega)$
which is holomorphic with respect to the complex structure $\phi_1^* I$ compatible with $\omega$, where
$\{ \phi_t\}_{0\leq t\leq 1}\subset \Diff_0 (M)$, $\phi_0 = Id$, is an isotopy such that $\phi_t (\Sigma) = \Sigma$ for all $t\in [0,1]$.

\end{itemize}

\bigskip
\noindent
(II)
Let $\cC_0$ be a connected component of $\compK (M)$.
Assume that
the set $\ttK_{\cC_0} (\bfr)$ (or, equivalently, the set $\cK_{\cC_0} (\bfr)$)
is connected.

Then any two K\"ahler-type embeddings $\bigsqcup_{i=1}^k B^{2n}(r_i)\to (M,\omega)$ favoring $\cC_0$
lie in the same $\Symp (M,\omega)\cap \Diff_0 (M)$-orbit. In particular, there exists $[I]\in\cmptteich_{\cC_0} (M,\omega)$ such that
both embeddings are of $[I]$-K\"ahler type.

If, in addition, $\SympH (M)$ acts transitively on the set of connected components of $\compK (M)$ compatible with $\omega$, then
any two K\"ahler-type embeddings $\bigsqcup_{i=1}^k B^{2n}(r_i)\to (M,\omega)$
lie in the same $\SympH (M,\omega)$-orbit.


\hfill


For the proof of \ref{_existence-connectedness-Kahler-type-arb-mfds_Theorem_} see Section~\ref{_Kahler-type-embs-of-balls-and-blow-ups-pfs_Section_} (\ref{_existence-connectedness-Kahler-type-arb-mfds-COPY_Theorem_}).


\hfill


\remark

\noindent
1. For an analogue of part (I) of \ref{_existence-connectedness-Kahler-type-arb-mfds_Theorem_} for tame embeddings
see \ref{_tame-embs-existence-via-forms-on-blow-up_Proposition_} below.

\bigskip
\noindent
2. \cite[Remark 5.2]{_Fleming2021_} states that the methods of \cite{_Fleming2021_} can be used to prove a result that, in view of
\ref{_I-Kahler-type-embs-under-diff-names_Remark_},
is equivalent to Part (I) \ref{_existence-connectedness-Kahler-type-arb-mfds_Theorem_} in the case where
$\omega$ is compatible with a complex structure $I$ on $M$, so that the complex manifold $(M,I)$ is projective and
the cohomology class $[\omega]\in H^2 (M;\R)$ is the first Chern class of an ample line bundle on $(M,I)$.


\hfill


Recall that the symplectic blow-up construction (see e.g. \cite[Sec. 7]{_McD-Sal-3_}, cf. Section~\ref{_blow-up-forms-Subsection_} below) associates to a complex blow-up of a symplectic manifold at a point (with respect to, say, an almost complex structure compatible with the symplectic form and integrable near the point) a class of symplectic forms on the blow-up. We will call such a symplectic form on the blow-up a {\bf blow-up symplectic form}. A {\bf rational (4-dimensional) symplectic manifold} is either $\C P^2$ with a standard Fubini-Study symplectic form, or $\C P^1\times \C P^1$ with a product symplectic form, or their blow-up at several points equipped with a blow-up
symplectic form.

A {\bf ruled (4-dimensional) symplectic manifold} is a closed symplectic 4-manifold which is the total space of an $\SP^2$-fiber bundle over a closed connected orientable 2-dimensional surface, so that all the fibers are symplectic submanifolds \cite{_Lal-McD-NewtonInst1996_}.

For rational or ruled closed symplectic manifolds of real dimension $4$ we have a stronger connectedness result, which follows directly
from similar connectedness results for {\sl symplectic} (not necessarily K\"ahler-type) embeddings of balls
that were proved in \cite{_McD-connectedness_} and \cite{_Borman-Li-Wu_}
(for previous partial results see \cite{_McD-Topology1991_, _McD-Cambridge1990_, _Lalonde-MathAnn-1994_, _Biran-IMRN1996_}).


\hfill


\vfil\eject

\theorem
\label{_connectedness-rational-ruled-mfds_Theorem_}

Let $(M,\omega)$ be a rational or ruled closed symplectic manifold, $\dim_\R M = 4$. Let $\Sigma\subset (M,\omega)$
be either of the following:

\bigskip
\noindent
(1) the empty set;

\smallskip
\noindent
(2) a finite union of closed compact symplectic submanifolds (without boundary) of real dimension 2
whose pairwise intersections (if they exist) are transverse and $\omega$–orthogonal;

\smallskip
\noindent
(3) a Lagrangian submanifold which is diffeomorphic to either $\SP^2$ or $\R P^2$.

\bigskip
Then for any $k\in\Z_{>0}$ and any $r_1,\ldots,r_k>0$,
 the following claims hold:

\bigskip
\noindent
(I)
Any two K\"ahler-type, respectively tame, embeddings
$\bigsqcup_{i=1}^k B^4 (r_i)\to (M\setminus \Sigma,\omega)$ lie in the same orbit of the $\Symp^c_0 (M\setminus\Sigma,\omega)$-action.
In particular, (for fixed $r_1,\ldots,r_k$) the space of all such K\"ahler-type, respectively tame, embeddings is connected.

\bigskip
\noindent
(II) If there exists an $[I]$-K\"ahler-type, respectively tame, embedding $\bigsqcup_{i=1}^k B^4 (r_i)\to (M\setminus \Sigma,\omega)$,
then any symplectic embedding $\bigsqcup_{i=1}^k B^4 (r_i)\to (M\setminus \Sigma,\omega)$ is of $[I]$-K\"ahler type, respectively tame.


\hfill


\noindent
{\bf Proof of \ref{_connectedness-rational-ruled-mfds_Theorem_}:}

Any two symplectic embeddings $\bigsqcup_{i=1}^k B^{2n}(r_i)\to (M\setminus \Sigma,\omega)$
lie in the same orbit of the $\Symp^c_0 (M\setminus\Sigma,\omega)$-action:

\smallskip
\noindent
- In case (1) this was proved in \cite[Cor. 1.5]{_McD-connectedness_}.

\smallskip
\noindent
- In case (2), where $\Sigma\subset (M,\omega)$ is just a symplectic submanifold, this was essentially proved in \cite[Cor. 1.5]{_McD-connectedness_} -- see \cite[Prop. 2.1]{_Borman-Li-Wu_} where more details of the proof are spelled out in the case where the symplectic submanifold $\Sigma\subset (M,\omega)$ is a sphere; in fact \cite[Cor. 1.2.13]{_McD-Opshtein_}, used in the proof of \cite[Prop. 2.1]{_Borman-Li-Wu_}, yields the same claim for any $\Sigma$ as in (2).

\smallskip
\noindent
- In case (3), this is proved in \cite[Thm. 1.1]{_Borman-Li-Wu_}.

\smallskip
  Since $[I]$-K\"ahler-type, respectively tame, embeddings into $(M\setminus \Sigma,\omega)$ are mapped by
the $\Symp^c_0 (M\setminus\Sigma,\omega)$-action
into $[I]$-K\"ahler-type, respectively tame, embeddings, we get the
results of both parts of the theorem.
\endproof


\hfill


Now let us discuss applications of \ref{_existence-connectedness-Kahler-type-arb-mfds_Theorem_}
and \ref{_connectedness-rational-ruled-mfds_Theorem_} to particular symplectic manifolds.


\hfill


\subsection{The case of projective spaces and their products}
\label{_main-results-CPn_Subsection_}

Denote by $\omega_{FS}$, or by $\omega_{FS,n}$, the Fubini-Study form on $\C P^n$ normalized so that $\int_{\C P^1} \omega_{FS} = \pi$.
Let $I_{st}$ be the standard complex structure on $\C P^n$.

We also set $\C P^0 := \textrm{pt}$.


\hfill


\vfil\eject

\theorem
\label{_existence-products-of-proj-spaces_Theorem_}

Consider the manifold $M:=\C P^{n_1}\times\ldots\times\C P^{n_m}$, $n_1,\ldots, n_m >0$, $n_1 +\ldots + n_m =: n$, endowed with the symplectic form $\omega_\bfc :=
c_1\omega_{FS,n_1}\oplus\ldots\oplus c_m\omega_{FS,n_m}$, $\bfc := (c_1,\ldots,c_m)\in (\R_{>0})^m$. Let $I$ be the complex
structure on $M$ which is the product of the standard complex structures on the factors.

Let $l_1,\ldots,l_m\in\Z_{>0}$ so that $[l_1:\ldots :l_m] = [c_1:\ldots: c_m]\in \R P^m$.
Let
\[
k := \frac{(n_1+\ldots + n_m)!}{n_1!\cdot\ldots\cdot n_m!} l_1^{n_1}\cdot\ldots\cdot l_m^{n_m}.
\]
Let $\Sigma\subset (M,I)$ be a proper (possibly empty) complex submanifold. If $M = \C P^{n_1}\times \C P^1$, we also allow
$\Sigma = (\C P^{n_1 -1}\times \C P^1) \cup (\C P^{n_1}\times \textrm{pt})$.

Then K\"ahler-type embeddings of disjoint unions of $k$ equal balls into $(M\setminus\Sigma,\omega_\bfc)$
are unobstructed.

More precisely, if $k \Vol (B^{2n} (r),\omega_0) < \Vol (M,\omega)$, then for each proper (possibly empty) complex submanifold $\Sigma\subset (M,I)$, there exists an $[I]$-K\"ahler-type embedding of $\bigsqcup_{i=1}^k B^{2n} (r)$
into $(M,\omega)$ which is holomorphic with respect to a complex structure on $M$ that is compatible with $\omega$
and isotopic to $I$ by an isotopy preserving $\Sigma$ (as a set).


\hfill


For the proof of
\ref{_existence-products-of-proj-spaces_Theorem_} see Section~\ref{_complex-strs-on-proj-spaces-their-products-and-blowups_Section_} (\ref{_existence-products-of-proj-spaces-COPY_Theorem_}).

For $n_1 =\ldots = n_m = 1$, $l_1=\ldots = l_m = l$, $c_1 =\ldots = c_m = 1$ and $\Sigma=\emptyset$,
\ref{_existence-products-of-proj-spaces_Theorem_} yields the following corollary.


\hfill


\corollary
\label{_main-results-product-of-m-CPone_Corollary_}

Let $M:=(\C P^1)^m$. Consider a symplectic form $\omega$ on $M$ which is the product of equal Fubini-Study forms on the $\C P^1$-factors of $M$. Let $I$ be the complex structure on $M$ which is the product of the standard complex structures on the $\C P^1$-factors of $M$.

Then for any $l\in\Z_{>0}$, K\"ahler-type embeddings of disjoint unions of $m! l^m$ equal balls into $(M,\omega)$
are unobstructed.

More precisely, if $m! l^m \Vol (B^{2n} (r),\omega_0) < \Vol (M,\omega)$, then there exists a K\"ahler-type embedding of $\bigsqcup_{i=1}^{m! l^m} B^{2n} (r)$
into $(M,\omega)$ which is holomorphic with respect to a complex structure on $M$ that is compatible with $\omega$
and isotopic to $I$.
\endproof


\hfill


Now let us consider the case of a single complex projective space.


\hfill


\theorem\label{_existence-connectedness-CPn-l-power-n-balls_Theorem_}

\noindent
A. For each $l\in \Z_{>0}$,
K\"ahler-type embeddings of disjoint unions of $l^n$ equal balls into $(\C P^n,\omega_{FS})$ are unobstructed.

More precisely,
if $l^n \Vol (B^{2n} (r),\omega_0) < \Vol (\C P^n,\omega_{FS})$, then for any complex structure $J$ on $\C P^n$ compatible with $\omega_{FS}$ (and, in particular, for the standard complex structure $I_{st}$) and for each proper (possibly empty) complex submanifold $\Sigma\subset (\C P^n,J)$, there exists a $[J]$-K\"ahler-type embedding of $\bigsqcup_{i=1}^{l^n} B^{2n} (r)$
into $(\C P^n,\omega_{FS})$ which is holomorphic with respect to a complex structure on $\C P^n$ that is compatible with $\omega_{FS}$
and isotopic to $J$ by an isotopy preserving $\Sigma$ (as a set).

\bigskip
\noindent
B. The group $\SympH (\C P^n,\omega_{FS})$ acts transitively on the set of connected components of $\compK (\C P^n)$ compatible with
$\omega_{FS}$.

\bigskip
\noindent
C.
For any $k\in\Z_{>0}$ and $r_1,\ldots,r_k>0$, any two K\"ahler-type embeddings
$\bigsqcup_{i=1}^k B^{2n} (r_i)\to (\C P^n,\omega_{FS})$
(if they exist!) lie in the same orbit of the  $\SympH (\C P^n,\omega_{FS})$-action.

They lie in the same orbit of the $\Symp (\C P^n,\omega_{FS})\cap \Diff_0 (\C P^n)$-action
if and only if they favor a common connected component of $\compK (\C P^n)$. In the latter case there exists $[I]\in\cmptteich (\C P^n,\omega_{FS})$ such that
both embeddings are of $[I]$-K\"ahler type.


\hfill


Part A of \ref{_existence-connectedness-CPn-l-power-n-balls_Theorem_}
follows from
\ref{_existence-products-of-proj-spaces_Theorem_}.
For the proof of parts B and C of \ref{_existence-connectedness-CPn-l-power-n-balls_Theorem_} see Section~\ref{_complex-strs-on-proj-spaces-their-products-and-blowups_Section_} (\ref{_existence-connectedness-CPn-l-power-n-balls-COPY_Theorem_}).


\hfill


\remark

Let us say that a K\"ahler-type complex structure $I$ on a closed manifold $M$
is {\bf rigid}, if the connected component of $\compK (M)$ containing $I$ coincides with the $\Diff_0 (M)$-orbit of $I$
(or, in other words, if any K\"ahler-type complex structure on $M$ obtained from $I$ by a smooth deformation is, in fact, isotopic to $I$).

As we will see in
\ref{_complex-str-on-CPn-is-unique_Theorem_}, any
K\"ahler-type complex structure on $\C P^n$ is rigid. More
examples of rigid complex structures include
so-called ${\mathcal Z}$-manifolds
 (i.e., Calabi-Yau manifolds with
$h^{n-1,1}=0$, \cite{_Candelas_etc:rigid_}).

The proof of part C of \ref{_existence-connectedness-CPn-l-power-n-balls_Theorem_} allows to generalize the last claim of part C as follows:
Assume $(M,\omega)$ is a closed manifold such that $\omega$ is compatible with a rigid complex structure $I$. Let $\cC_0$ be a connected component
of $\compK (M)$ containing $I$. Then two K\"ahler-type embeddings $\bigsqcup_{i=1}^k B^{2n} (r_i)\to (M,\omega)$ favoring $\cC_0$
lie in the same orbit of the $\Symp (M,\omega)\cap \Diff_0 (M)$-action.


\hfill


Applying part A of \ref{_existence-connectedness-CPn-l-power-n-balls_Theorem_} with
$\Sigma=\emptyset$, with $I_{st} = J$ being fixed, and with $\omega_{FS}$ being varied by an isotopy, we get the following immediate corollary.


\hfill


\corollary\label{_CPn-l-to-n-balls-flat-Kahler-metric_Corollary_}

Let $l\in \Z_{>0}$.

Then for any $r>0$ such that $l^n\Vol (B^{2n} (r),\omega_0) < \Vol (\C P^n,\omega_{FS})$ there exists a K\"ahler form $\omega$ on $(\C P^n, I_{st})$ isotopic to $\omega_{FS}$ and such that the K\"ahler manifold $(\C P^n, I_{st},\omega)$ admits a K\"ahler (i.e., both holomorphic and symplectic) embedding of $\bigsqcup_{i=1}^{l^n} B^{2n} (r)$ with the standard flat K\"ahler metric on it.
\endproof


\hfill


In particular, part A of \ref{_existence-connectedness-CPn-l-power-n-balls_Theorem_} and \ref{_CPn-l-to-n-balls-flat-Kahler-metric_Corollary_}
mean that by changing either $I_{st}$ or $\omega_{FS}$ by an isotopy one can get a K\"ahler metric on $\C P^n$ which is the standard flat K\"ahler metric on $l^n$ disjoint balls in $\C P^n$ that fill out an arbitrarily large portion of the volume of the manifold.


\hfill


\remark
\label{_CPn-l-to-n-balls-flat-Kahler-metric-prev-results_Remark_}

For $l=1$ and $n=2$ \ref{_CPn-l-to-n-balls-flat-Kahler-metric_Corollary_} was proved in \cite{_Eckl2017_}.


\hfill


Applying \ref{_existence-connectedness-CPn-l-power-n-balls_Theorem_} with $J=I_{st}$ and with the hyperplane $\Sigma = \C P^{n-1}\subset (\C P^n, I_{st})$, we get the following corollary.


\hfill


\corollary
\label{_ball-of-arb-even-dimension-filled-by-several-equal-balls_Corollary_}

For any $l\in \Z_{>0}$ and any $r>0$ such that $l^n \Vol (B^{2n} (r),\omega_0) < \Vol (B^{2n} (1),\omega_0)$ there exists a K\"ahler-type
embedding of $\bigsqcup_{i=1}^{l^n} B^{2n} (r)$
into $(B^{2n} (1),\omega_0)$ which is holomorphic with respect to a complex structure on $B^{2n} (1)$ isotopic to $J_0$ and compatible with $\omega_0$.
\endproof


\hfill


Now let us focus on the case of $(\C P^2,\omega_{FS})$.


\hfill


\proposition
\label{_CP2-only-one-conn-cmpt-of-compK-compatible-with-omegaFS_Proposition_}

There is only one connected component of $\compK (\C P^2)$ compatible with $\omega_{FS}$ -- it is the connected component
of $\compK (\C P^2)$ containing $I_{st}$. Any two complex structures in that connected component are isotopic.


\hfill


For the proof of \ref{_CP2-only-one-conn-cmpt-of-compK-compatible-with-omegaFS_Proposition_} see Section~\ref{_complex-strs-on-proj-spaces-their-products-and-blowups_Section_} (\ref{_CP2-only-one-conn-cmpt-of-compK-compatible-with-omegaFS-COPY_Proposition_}).
We do not know whether $\compK (\C P^2)$ has more than one connected component -- see \ref{_CPn-complex-structures_Remark_}.

\ref{_CP2-only-one-conn-cmpt-of-compK-compatible-with-omegaFS_Proposition_} yields the following immediate corollary.


\hfill


\corollary
\label{_CP2-all-K-type-embs-are-holmwrt-c-s-isotopic-to-standard-one_Corollary_}

Any K\"ahler-type embedding (of any disjoint union of domains) into $(\C P^2,\omega_{FS})$ is of $[I_{st}]$-K\"ahler type: it is symplectic with respect to $\omega_{FS}$ and holomorphic with respect to a complex structure on $\C P^2$ isotopic to $I_{st}$ and compatible with $\omega_{FS}$. \endproof


\hfill


Recall the following theorem of McDuff-Polterovich \cite{_McD-Polt_}.


\hfill


\theorem
\label{_McD-Polt-CP2_Theorem_}

Let $\Sigma\subset (\C P^2, I_{st})$ be a proper (possibly empty) complex submanifold.

\noindent
A. Assume that $1\leq k\leq 8$ and $r_1\geq r_2\geq\ldots\geq r_k>0$.

Then
the following conditions are equivalent:

\begin{itemize}

\item{} There exists a symplectic embedding $\bigsqcup_{i=1}^k B^4 (r_i)\to (\C P^2\setminus \Sigma,\omega_{FS})$.

\item{} The radii $r_1,\ldots,r_k$ satisfy the following inequalities (listed in \cite[Cor. 1.3G]{_McD-Polt_}):

\smallskip
\noindent
(v) $\sum_{i=1}^k r_i^4 < 1$ (volume inequality),

\smallskip
\noindent
(c1) $r_1^2 + r_2^2 < 1$, if $2\leq k\leq 8$,

\smallskip
\noindent
(c2) $r_1^2 +\ldots + r_5^2 < 2$, if $5\leq k\leq 8$,

\smallskip
\noindent
(c3) $2r_1^2 + \sum_{i=2}^7 r_i^2 < 3$, if $7\leq k\leq 8$,

\smallskip
\noindent
(c4) $2r_1^2 + 2r_2^2 + 2r_3^2 + r_4^2 +\ldots +r_8^2 < 4$, if $k=8$,

\smallskip
\noindent
(c5) $2\sum_{i=1}^6 r_i^2 + r_7^2 + r_8^2 < 5$, if $k=8$,

\smallskip
\noindent
(c6) $3 r_1^2 + 2\sum_{i=2}^8 r_i^2 < 6$, if $k=8$.

\end{itemize}

\bigskip
\noindent
B. Symplectic embeddings of a disjoint union of $k=l^2$, $l\in \Z_{>0}$, equal balls into $(\C P^2\setminus\Sigma,\omega_{FS})$ are unobstructed:
if $k \Vol (B^4 (r),\omega_0) < \Vol (\C P^2,\omega_{FS})$, then there exists a symplectic embedding $\bigsqcup_{i=1}^k B^4 (r)\to (\C P^2\setminus\Sigma,\omega_{FS})$.


\hfill


Here is an extension of \ref{_McD-Polt-CP2_Theorem_} to K\"ahler-type embeddings of balls into $\C P^2$.


\hfill


\theorem\label{_CP2-existence-connectedness_Theorem_}

\noindent
A. Assume that $1\leq k\leq 8$ and $r_1\geq r_2\geq\ldots\geq r_k>0$.

Then for any proper (possibly empty) complex submanifold $\Sigma\subset (\C P^2, I_{st})$,
any symplectic embedding $\bigsqcup_{i=1}^k B^4 (r_i)\to (\C P^2\setminus \Sigma,\omega_{FS})$ is, in fact, of $[I_{st}]$-K\"ahler type:
it is holomorphic with respect to a complex structure on $\C P^2$ that is compatible with $\omega_{FS}$ and isotopic
to $I_{st}$ by an isotopy fixing $\Sigma$ (as a set).

Accordingly, by part A of \ref{_McD-Polt-CP2_Theorem_}, such a K\"ahler-type embedding $\bigsqcup_{i=1}^k B^4 (r_i)\to (\C P^2\setminus \Sigma,\omega_{FS})$
exists if and only if the radii $r_1,\ldots,r_k$ satisfy the inequalities (v), (c1)-(c6) listed in part A of \ref{_McD-Polt-CP2_Theorem_}.

\bigskip
\noindent
B. For any proper (possibly empty) complex submanifold $\Sigma\subset (\C P^2, I_{st})$ and any $k=l^2$, $l\in \Z_{>0}$,
any symplectic embedding $\bigsqcup_{i=1}^k B^4 (r)\to (\C P^2\setminus\Sigma,\omega_{FS})$ is, in fact, of $[I_{st}]$-K\"ahler type:
it is holomorphic with respect to a complex structure on $\C P^2$ that is compatible with $\omega_{FS}$
and isotopic to $I_{st}$ by an isotopy preserving $\Sigma$ (as a set).

Accordingly, by part B of \ref{_McD-Polt-CP2_Theorem_}, such K\"ahler-type embeddings  $\bigsqcup_{i=1}^k B^4 (r)\to (\C P^2\setminus\Sigma,\omega_{FS})$ are unobstructed: they exist if and only if $k \Vol (B^4 (r),\omega_0) < \Vol (\C P^2,\omega_{FS})$.

\bigskip
\noindent
C. Let $\Sigma\subset (\C P^2,\omega_{FS})$
be either of the following:

\bigskip
\noindent
(1) the empty set;

\smallskip
\noindent
(2) a finite union of closed compact symplectic submanifolds (without boun\-dary) of real dimension 2
whose pairwise intersections (if they exist) are transverse and $\omega$–orthogonal;

\smallskip
\noindent
(3) a Lagrangian submanifold which is diffeomorphic to either $\SP^2$ or $\R P^2$.

For any $k\in\Z_{>0}$ and $r_1,\ldots,r_k>0$, any two K\"ahler-type embeddings $\bigsqcup_{i=1}^k B^{4} (r_i)\to (\C P^2\setminus\Sigma,\omega_{FS})$ (if they exist!) lie in the same orbit of the  $\Symp^c_0 (\C P^2\setminus\Sigma,\omega_{FS})$-action, meaning that
the space of K\"ahler-type embeddings $\bigsqcup_{i=1}^k B^{4} (r_i)\to (\C P^2\setminus\Sigma,\omega_{FS})$ is connected.


\hfill


Part B of \ref{_CP2-existence-connectedness_Theorem_} follows from
part A of \ref{_existence-connectedness-CPn-l-power-n-balls_Theorem_}
and from part (II) of \ref{_connectedness-rational-ruled-mfds_Theorem_}.
Part C of \ref{_CP2-existence-connectedness_Theorem_} follows directly from part (I) of
\ref{_connectedness-rational-ruled-mfds_Theorem_}.
 For the proof of part A of \ref{_CP2-existence-connectedness_Theorem_} see Section~\ref{_complex-strs-on-proj-spaces-their-products-and-blowups_Section_} (\ref{_CP2-existence-connectedness-COPY_Theorem_}).

Applying \ref{_CP2-existence-connectedness_Theorem_} with the hyperplane $\Sigma=\C P^1\subset \C P^2$ and
recalling that $(\C P^2\setminus \C P^1,\omega_{FS})$ is symplectomorphic to $(B^4 (1), \omega_0)$,
we get the following corollary.


\hfill


\corollary\label{_ball-four-dim-existence-connectedness_Corollary_}

\noindent
A. Assume that $1\leq k\leq 8$ and $r_1\geq r_2\geq\ldots\geq r_k>0$.

Then
any symplectic embedding $\bigsqcup_{i=1}^k B^4 (r_i)\to (B^4 (1),\omega_0)$
is, in fact, of K\"ahler-type: it is holomorphic with respect to a complex structure on $B^4 (1)$ that is compatible with $\omega_0$ and isotopic
to $I_0$.

Accordingly, by part A of \ref{_McD-Polt-CP2_Theorem_}, such a K\"ahler-type embedding $\bigsqcup_{i=1}^k B^4 (r_i)\to (B^4 (1),\omega_0)$
exists if and only if the radii $r_1,\ldots,r_k$ satisfy the inequalities (v), (c1)-(c6) listed in part A of \ref{_McD-Polt-CP2_Theorem_}.

\bigskip
\noindent
B. For any $k=l^2$, $l\in \Z_{>0}$,
any symplectic embedding $\bigsqcup_{i=1}^k B^4 (r)\to (B^4 (1),\omega_0)$ is, in fact, of K\"ahler type:
it is holomorphic with respect to a complex structure on $B^4 (1)$ that is compatible with $\omega_0$
and isotopic to $I_0$.

Accordingly, by part B of \ref{_McD-Polt-CP2_Theorem_}, such K\"ahler-type embeddings  $\bigsqcup_{i=1}^k B^4 (r)\to (B^4 (1),\omega_0)$ are unobstructed: they exist if and only if $k \Vol (B^4 (r),\omega_0) < \Vol (B^4 (1),\omega_0)$.

\bigskip
\noindent
C.
For any $k\in\Z_{>0}$ and $r_1,\ldots,r_k>0$, any two K\"ahler-type embeddings $\bigsqcup_{i=1}^k B^{4} (r_i)\to (B^4 (1),\omega_0)$ (if they exist!) lie in the same orbit of the  $\Symp^c_0 (B^4 (1),\omega_0)$-action, meaning that
the space of K\"ahler-type embeddings $\bigsqcup_{i=1}^k B^{4} (r_i)\to (B^4 (1),\omega_0)$ is connected.
\endproof


\hfill


\remark
\label{_CP2-existence-connectedness_Remark_}

\smallskip
\noindent
 1.
 The description of the K\"ahler cone of the blow-up of $\C P^2$ at $k\geq 10$ (generic) points is a difficult open question related, in particular, to the Nagata conjecture
 (which is equivalent to the claim that a certain cohomology class of the blow-up lies in the boundary of the K\"ahler cone,
see e.g. \cite{_Biran-from-sympl-pack-to-alg-geom_}).
 In particular, if there is a counterexample to the Nagata conjecture for $k\geq 10$, it
would yield that the K\"ahler-type embeddings of $k$ equal balls to $\C P^2$ are not unobstructed -- see \cite[Thm. 1.4.B]{_McD-Polt_}. Since by a theorem of Biran \cite{_Biran-GAFA1997_}, {\sl symplectic} embeddings of a disjoint union of $k\geq 10$ equal balls into $(\C P^2, \omega_{FS})$ are unobstructed, this would imply the existence of such
symplectic embeddings into $(\C P^2, \omega_{FS})$ that are {\sl not} of K\"ahler type.

\smallskip
\noindent
2. The connectedness of the space of {\sl symplectic} embeddings of $\bigsqcup_{i=1}^k B^{4} (r_i)$ into $(\C P^2, \omega_{FS})$ and $(B^4 (1),\omega_0)$
was proved in \cite{_McD-connectedness_} (for previous partial results see \cite{_McD-Topology1991_, _McD-Cambridge1990_, _Lalonde-MathAnn-1994_, _Biran-IMRN1996_}).

\smallskip
\noindent
3. By a theorem of Taubes \cite{_Taubes-MRL1995_, _Taubes-JDG1996_, _Taubes-book-2000_}, any symplectic form on $\C P^2$ can be mapped into a positive multiple of $\omega_{FS}$ by a diffeomorphism of $\C P^2$. Therefore, any symplectic form $\omega$ on $\C P^2$ is of K\"ahler type and the results above stated for $(\C P^2, \omega_{FS})$
hold also for $(\C P^2, \omega)$.


\hfill


\subsection{The cases of $\C P^2\sharp \overline{\C P^2}$ and $\C P^1\times \C P^1$}
\label{_main-results-Hirzberuch_Subsection_}

Consider $\C P^2$ with the standard complex structure on it. Define $M = \C P^2\sharp \overline{\C P^2}$ -- as a smooth manifold -- as the complex blow-up of $\C P^2$ at one point. Let $E\subset M$ be the exceptional divisor and $\C P^1 \subset M$ a projective line.
Given $0 <\lambda <\pi$, let $\omega_\lambda$ be a K\"ahler form on $M$ such that $\langle [\omega_\lambda], [\C P^1]\rangle = \pi$,
$\langle [\omega_\lambda], [E]\rangle = \lambda$.

The manifold $M  = \C P^2\sharp \overline{\C P^2}$ admits K\"ahler-type complex structures $J_{2l+1}$, $l\in\Z_{\geq 0}$, such that $(M,J_{2l+1})$ is biholomorphic to the odd Hirzebruch surface $\F_{2l+1}$. In particular, we assume $J_1 = \tI_{st}$, as we identify the Hirzebruch surface  $\F_1$ with $(\widetilde{ \C P}^2_{I_{st},x}, \tI_{st})$. Each complex structure on $M$ can be mapped by a diffeomorphism into exactly one complex structure $J_{2l+1}$ -- this follows from \cite{_Friedman-Qin-InvMath-1995_}. In fact, this diffeomorphism can be chosen to lie in $\Diff_0 (M)$ -- see \cite{_Meersseman-JEP2019_}.


\hfill


\vfil\eject

\theorem
\label{_Hirzebruch-odd-existence-connectedness_Theorem_}

Let $M = \C P^2\sharp \overline{\C P^2}$ and let $\omega_\lambda$ be a K\"ahler-type form on $M$ as above.

Then the following claims hold:

\bigskip
\noindent
A. Let $2\leq k\leq 8$, $R_1\geq R_2\geq\ldots\geq R_{k-1}>0$. Assume that
$r_1:=R_1\geq\ldots\geq r_{i-1}:=R_{i-1}\geq r_i:=\sqrt{\lambda/\pi}\geq r_{i+1}:= R_i\geq\ldots \geq r_k:= R_{k-1}$.

Then for any proper (possibly empty) complex submanifold $\Sigma\subset (M,J_1)$,
the following conditions are equivalent:

\begin{itemize}

\item{} There exists a $[J_1]$-K\"ahler-type embedding
$\bigsqcup_{i=1}^{k-1} B^4 (R_i)\to (M\setminus\Sigma,\omega_\lambda)$ holomorphic with respect to a complex structure on $M$ that is
compatible with $\omega_\lambda$ and isotopic to $J_1$ by an isotopy preserving $\Sigma$ (as a set).

\item{} The numbers $r_1,\ldots,r_k$ satisfy the inequalities (v), (c1)-(c6) in part A of
\ref{_McD-Polt-CP2_Theorem_}.

\end{itemize}

For $r_1,\ldots,r_k$ satisfying the inequalities (v), (c1)-(c6) in part A of \ref{_McD-Polt-CP2_Theorem_}, any symplectic embedding
$\bigsqcup_{i=1}^{k-1} B^4 (R_i)\to (M\setminus\Sigma,\omega_\lambda)$ is, in fact, of $[J_1]$-K\"ahler type.

\bigskip
\noindent
B.
For any $k\in\Z_{>0}$ and $R_1,\ldots,R_k>0$,
any two K\"ahler-type embeddings $\bigsqcup_{i=1}^k B^4 (R_i)\to (M,\omega_\lambda)$ (if they exist!) lie in the same orbit of the  $\Symp_0 (M,\omega_\lambda)$-action, meaning that
the space of K\"ahler-type embeddings $\bigsqcup_{i=1}^k B^4 (R_i)\to (M,\omega_\lambda)$ is connected.


\hfill


Part B of \ref{_Hirzebruch-odd-existence-connectedness_Theorem_} follows directly from
part (I) of \ref{_connectedness-rational-ruled-mfds_Theorem_}.
For the proof of part A of
\ref{_Hirzebruch-odd-existence-connectedness_Theorem_} see Section~\ref{_Hirzebruch-odd_Section_} (\ref{_Hirzebruch-odd-existence-connectedness-COPY_Theorem_}).


\hfill


\remark

\noindent
1. We do not know under which conditions on $r_1,\ldots,r_k$ the complex structure $J_1$ on $M$ can be replaced in part A of \ref{_Hirzebruch-odd-existence-connectedness_Theorem_} by
$J_{2l+1}$, $l>0$. In order to obtain these conditions one needs to get an explicit description of the K\"ahler cones
of the blow-ups of the Hirzebruch surfaces $\F_{2l+1}$, $l>0$, at $k$ points -- see Section~\ref{_Hirzebruch-odd_Section_}.

\bigskip
\noindent
2. For any $\bigsqcup_{i=1}^k B^4 (R_i)$, the connectedness of the space of {\sl symplectic} embeddings of $\bigsqcup_{i=1}^k B^4 (R_i)$ into
$(M,\omega_\lambda)$ was proved in \cite{_McD-connectedness_} (also see \cite{_Lalonde-MathAnn-1994_} for a previous partial result
for $k=1$).

\bigskip
\noindent
3. \ref{_Hirzebruch-odd-existence-connectedness_Theorem_} admits a generalization concerning K\"ahler-type embeddings
of a disjoint union of $k$ balls into a complex blow-up of $\C P^2$ at $j$ points ($1\leq k\leq 8-j$, $j=1,\ldots,7$) equipped with an appropriate symplectic form.
The proof is a straightforward modification of the proof of \ref{_Hirzebruch-odd-existence-connectedness_Theorem_}.


\hfill


Now let $M = \C P^1\times \C P^1$ -- as a smooth manifold. Let
$\omega_0$ be the Fubini-Study form on $\C P^1$ normalized so that $\int_{\C P^1} \omega_0 =\pi$.
Given $\mu\geq 1$, let $\omega_\mu = \mu \omega_0 \oplus \omega_0$.

The manifold $M = \C P^1\times \C P^1$ admits K\"ahler-type complex structures $J_{2l}$, $l\in\Z_{\geq 0}$, such that $(M,J_{2l})$ is biholomorphic to the even Hirzebruch surface $\F_{2l}$. In particular, $J_0$ can be viewed as the standard complex structure on $M= \C P^1\times \C P^1$.
Each complex structure on $M = \C P^1\times \C P^1$ can be mapped by a diffeomorphism into exactly one complex structure $J_{2l}$ -- this follows from \cite{_Friedman-Qin-InvMath-1995_}. In fact, this diffeomorphism can be chosen to lie in $\Diff_0 (M)$ -- see \cite{_Meersseman-JEP2019_}.


\hfill


\theorem
\label{_Hirzebruch-even-existence-connectedness_Theorem_}

Let $M= \C P^1\times \C P^1$ and let $\omega_\mu$ be as above.

Let $\Sigma\subset (M,J_0)$ be either a proper (possibly empty) complex submanifold or $\Sigma = (\textrm{pt}\times \C P^1) \cup (\C P^1\times \textrm{pt})$.

Then the following claims hold:

\bigskip
\noindent
A. Assume $2\leq k\leq 8$,
$R_1, R_2,\ldots , R_{k-2}>0$, $0< R_{k-1} < 1$.
Consider the numbers
\[
\frac{R_1^2}{\mu+1-R_{k-1}^2},\frac{R_2^2}{\mu+1-R_{k-1}^2},\ldots,
\frac{\mu-R_{k-1}^2}{\mu+1-R_{k-1}^2}, \frac{1-R_{k-1}^2}{\mu+1-R_{k-1}^2},
\]
sort them in the non-increasing order and denote the resulting $k$ numbers by $r_1\geq\ldots\geq r_k$.

Then the following conditions are equivalent:

\begin{itemize}

\item{} There exists a $[J_0]$-K\"ahler-type embedding
$\bigsqcup_{i=1}^{k-1} B^4 (R_i)\to (M\setminus\Sigma,\omega_\mu)$ holomorphic with respect to a complex structure on $M$
that is compatible with
$\omega_\mu$ and isotopic to $J_0$ by an isotopy preserving $\Sigma$ (as a set).

\item{} The numbers $r_1,\ldots,r_k$ defined above satisfy the inequalities in part A of
 \ref{_McD-Polt-CP2_Theorem_}.

\end{itemize}

For $r_1,\ldots,r_k$ satisfying the inequalities (v), (c1)-(c6) in part A of \ref{_McD-Polt-CP2_Theorem_}, any symplectic embedding
$\bigsqcup_{i=1}^{k-1} B^4 (R_i)\to (M\setminus\Sigma,\omega_\mu)$ is, in fact, of $[J_0]$-K\"ahler type.

\bigskip
\noindent
B. Assume that $\mu=1$.

Then for any $l\in\Z_{>0}$,
any symplectic embedding of a disjoint union of $2 l^2$ equal balls into $(M\setminus\Sigma,\omega_1)$
is, in fact, of $[J_0]$-K\"ahler type, and such symplectic (or, equivalently, $[J_0]$-K\"ahler-type) embeddings are unobstructed.

More precisely, if $2 l^2 \Vol (B^4 (r),\omega_0) < \Vol (M,\omega_1)$, then there exists a K\"ahler-type embedding of $\bigsqcup_{i=1}^{2 l^2} B^4 (r)$
into $(M\setminus\Sigma,\omega_1)$ which is holomorphic with respect to a complex structure on $M$ that is compatible with $\omega_1$
and isotopic to $J_0$
 by an isotopy preserving $\Sigma$ (as a set).

\bigskip
\noindent
C. For any $k\in\Z_{>0}$ and $R_1,\ldots,R_k>0$,
any two K\"ahler-type embeddings $\bigsqcup_{i=1}^k B^{4} (R_i)\to (M\setminus\Sigma,\omega_\mu)$ (if they exist!) lie in the same orbit of the  $\Symp^c_0 (M\setminus\Sigma,\omega_\mu)$-action, meaning that
the space of K\"ahler-type embeddings $\bigsqcup_{i=1}^k B^{4} (R_i)\to (M\setminus\Sigma,\omega_\mu)$ is connected.


\hfill


Part B of \ref{_Hirzebruch-even-existence-connectedness_Theorem_} is a particular case of \ref{_main-results-product-of-m-CPone_Corollary_} (for $m=2$).
Part C of \ref{_Hirzebruch-even-existence-connectedness_Theorem_} follows directly from
part (I) of \ref{_connectedness-rational-ruled-mfds_Theorem_}.
For the proof of part A of \ref{_Hirzebruch-even-existence-connectedness_Theorem_} see Section~\ref{_Hirzebruch-even_Section_}
(\ref{_Hirzebruch-even-existence-connectedness-COPY_Theorem_}).


\hfill


\remark

\noindent
1. For $\Sigma = (\textrm{pt}\times \C P^1) \cup (\C P^1\times \textrm{pt})$ the manifold $(M\setminus\Sigma, \omega_\mu)$
is symplectomorphic to $(B^2 (1)\times B^2 (\sqrt{\mu}), \omega_0 \oplus \omega_0)$. Therefore \ref{_Hirzebruch-even-existence-connectedness_Theorem_} yields a straightforward corollary concerning K\"ahler-type embeddings of balls into $(B^2 (1)\times B^2 (\sqrt{\mu}), \omega_0 \oplus \omega_0)$.

For {\sl symplectic} embeddings of balls a similar corollary, as well as the result of part B of \ref{_Hirzebruch-even-existence-connectedness_Theorem_},
were proved in \cite{_McD-Polt_}.

\bigskip
\noindent
2. We do not know under which conditions on $r_1,\ldots,r_k$ the complex structure $J_0$ on $M$ can be replaced in part A of \ref{_Hirzebruch-even-existence-connectedness_Theorem_} by
$J_{2l}$, $l>0$.
In order to obtain these conditions one needs to get an explicit description of the K\"ahler cones
of the blow-ups of the Hirzebruch surfaces $\F_{2l}$, $l>0$, at $k$ points -- see Section~\ref{_Hirzebruch-even_Section_}.

\bigskip
\noindent
3. For any $R_1,\ldots, R_k$ the connectedness of the space of {\sl symplectic} embeddings $\bigsqcup_{i=1}^k B^4 (R_i)\to (M,\omega_\lambda)$ was proved in \cite{_McD-connectedness_}.


\hfill


\subsection{The case of tori and K3 surfaces}
\label{_main-results-tori-K3_Subsection_}

Let $M$, $\dim_\R M = 2n$, be either $\T^{2n}$ or a smooth manifold (of real dimension $4$) underlying a complex K3 surface.

As before, let $\omega$ be a K\"ahler-type symplectic form on $M$.

In the case $M=\T^{2n}$ one can assume, without loss of generality, that $\omega$ is linear (see \ref{_Kahler-type-structures-torus_Example_}).

Recall that $\omega$ is called {\bf rational} if the real cohomology class $[\omega]\in H^2 (M;\R)$ is a real multiple of a rational one,
and {\bf irrational} otherwise.

Irrational K\"ahler-type symplectic forms on tori and K3 surfaces are compatible with so-called Campana-simple complex structures. Namely, given a K\"ahler-type complex structure $I$ on a closed manifold $M$, the union of all positive-dimensional proper complex
subvarieties of $(M,I)$ is
either a countable union of proper analytic subvarieties of $M$ (and hence has a dense connected complement) or the whole $M$ (see e.g. \cite[Remark 4.2]{_EV-JTA_}).
In the former case, the complex structure $I$ is called {\bf Campana-simple} \cite{_Campana:isotrivial_}, \cite{_CDV:threefolds_} (for more details see Section~\ref{_Campana-simple-strs-gen-mfds-main-results_Section_}).

 In fact, a K\"ahler-type complex structure on a torus is Campana-simple if and only if it has no proper positive-dimensional complex subvarieties (see \ref{_any-Campana-simple-cs-on-torus-is-Campana-supersimple_Proposition_}). A symplectic form $\omega$ on the torus or a K3 surface is compatible with such a complex structure if and only if $\omega$ is irrational (see \ref{_torus-linear-Campana-simple-cs-compatible-with-omega_Proposition_} and \ref{_irrational-sympl-form-on-K3-compatible-with-Campana-simple-complex-str_Proposition_}).
See Section~\ref{_hyperkahler-case_Subsection_} for other examples of manifolds admitting Campana-simple complex structures.


\hfill


\theorem\label{_existence-connectedness-irrational-forms-tori-K3-surfaces_Theorem_}

Let $M$, $\dim_\R M = 2n$, be either $\T^{2n}$ or a smooth manifold (of real dimension $4$) underlying a complex K3 surface
and let $\omega$ be a K\"ahler-type symplectic form on $M$.

Assume that $\omega$ is irrational.

Let $k\in\Z_{>0}$, $r_1,\ldots,r_k>0$.

Then the following claims hold:

\bigskip
\noindent
A. K\"ahler-type embeddings of $\bigsqcup_{i=1}^k B^{2n} (r_i)$ into $(M,\omega)$ are unobstructed.

More precisely, assume that $\Vol \left(\bigsqcup_{i=1}^k B^{2n} (r_i),\omega_0\right) < \Vol (M,\omega)$ and $\omega$ is
 compatible with a Campana-simple complex structure $I$.

Then
there exists an $[I]$-K\"ahler-type embedding $\bigsqcup_{i=1}^k B^{2n} (r_i)\to (M,\omega)$.

\bigskip
\noindent
B.
The group $\SympH (M,\omega)$ acts transitively on the set of connected components of
$\compK (M)$ compatible with $\omega$.

\bigskip
\noindent
C. Any two K\"ahler-type embeddings $\bigsqcup_{i=1}^k B^{2n} (r_i)\to (M,\omega)$ (if they exist!) lie in the same orbit of the
$\SympH (M,\omega)$-action. They lie in the same orbit of the $\Symp (M,\omega)\cap \Diff_0 (M)$-action
if and only if they favor a common connected component of $\compK (M)$. In the latter case there exists $[I]\in\cmptteich (M,\omega)$
such that both embeddings are of $[I]$-K\"ahler-type.


\hfill


\remark
\label{_torus-of-dim-2_Remark_}

Assume $M=\T^2$ is a torus of real dimension $2$. In this case all symplectic forms and complex structures on $\T^2$ are of K\"ahler-type and
any symplectic (i.e.,
area-preserving) embedding into $\T^2$ is of K\"ahler-type. In particular, for any symplectic form $\omega$ on $\T^2$ (note that any such form is rational!) the result of part A of \ref{_existence-connectedness-irrational-forms-tori-K3-surfaces_Theorem_} is true
and the space of K\"ahler-type embeddings of any disjoint union of domains into $(\T^2,\omega)$ is path-connected
(since the space of symplectic forms on $\T^2$ is path-connected \cite{_Moser_}).


\hfill


For the proof
of \ref{_existence-connectedness-irrational-forms-tori-K3-surfaces_Theorem_}
see Section~\ref{_proofs-of-main-results-on-embs-of-balls-in-tori-K3_Subsection_} (\ref{_existence-connectedness-irrational-forms-tori-K3-surfaces-COPY_Theorem_}).

For a generalization of  \ref{_existence-connectedness-irrational-forms-tori-K3-surfaces_Theorem_} to arbitrary
 closed symplectic manifolds admitting Campana-simple complex structures compatible with the symplectic form see \ref{_Campana-simple-Kahler-type-existence-connectedness-gen-mfd_Theorem_}. In fact, we will use \ref{_Campana-simple-Kahler-type-existence-connectedness-gen-mfd_Theorem_} in the proof of \ref{_existence-connectedness-irrational-forms-tori-K3-surfaces_Theorem_}.

For tori the existence result in \ref{_existence-connectedness-irrational-forms-tori-K3-surfaces_Theorem_} admits the following equivalent formulation.


\hfill


\corollary
\label{existence-linear-tori_Corollary_}

Assume $\Gamma\subset\C^n$ is a lattice of rank $2n$. Let $M := \C^n/\Gamma$ be the corresponding torus. Let $\omega,I$ be the symplectic and the complex structures on $M$ induced by the standard flat symplectic and complex structures on $\C^n=\R^{2n}$.

Assume that $(M,I)$ admits no proper positive-dimensional
complex subvarieties. Let $k\in\Z_{>0}$,
$r_1,\ldots,r_k>0$ such that
$\Vol \left(\bigsqcup_{i=1}^k B^{2n} (r_i),\omega_0\right) < \Vol (M,\omega)$.

Then there exists a K\"ahler form $\omega'$ on $(M,I)$ isotopic to $\omega$ and a K\"ahler-type embedding
$\bigsqcup_{i=1}^k B^{2n} (r_i)\to (M,\omega')$ which is holomorphic with respect to $I$ -- i.e., an embedding which is symplectic with respect
to $\omega'$ and holomorphic with respect to $I$.


\hfill


\noindent
{\bf Proof of \ref{existence-linear-tori_Corollary_}:}

It suffices to prove a similar claim for the case where $\Gamma = \Z^{2n}$ and $\omega$ and $I$ are compatible symplectic and complex structures on $\T^{2n} = \R^{2n}/\Z^{2n}$, with $I$ being Campana-simple (the proof can be reduced to this case by using
a diffeomorphism $M=\C^n/\Gamma\to \T^{2n} = \R^{2n}/\Z^{2n}$
induced by an $\R$-linear isomorphism of $\R^{2n} = \C^{2n}$ sending $\Gamma$ to $\Z^{2n}$). In order to prove the corollary in this case, apply part A of \ref{_existence-connectedness-irrational-forms-tori-K3-surfaces_Theorem_} and
keep the complex structure fixed while changing the symplectic form by an isotopy
(cf. \ref{_I-Kahler-type-embs-under-diff-names_Remark_}).
\endproof


\hfill


\remark
\label{_existence-connectedness-tori-K3-comments_Remark_}

\noindent
1. For an arbitrary K\"ahler-type symplectic form $\omega$ on $M=\T^{2n}$ or on a manifold $M$ underlying a K3 surface,
the unobstructedness of symplectic -- but not necessarily K\"ahler-type or tame -- embeddings of an arbitrary collection of disjoint balls
was proved in \cite[Thm. 3.1]{_EV-JTA_}. In the case $M=\T^4$ a weaker version of this result had been previously proved by Latschev-McDuff-Schlenk \cite{_LMcDS_}.

\bigskip
\noindent
2. It is a well-known open problem (see e.g. \cite[Ch.15, Prob. 15]{_McD-Sal-3_}) whether $\SympH (\T^{2n},\omega)=\Symp_0 (\T^{2n},\omega)$
for a K\"ahler-type symplectic form $\omega$ on $\T^{2n}$, $n>1$.
If the answer to this question is positive, then
Part B of
\ref{_existence-connectedness-irrational-forms-tori-K3-surfaces_Theorem_}
would imply that the space of the K\"ahler-type embeddings $\bigsqcup_{i=1}^k B^{2n} (r_i)\to (\T^{2n},\omega)$
is path-connected for any $k\in\Z_{>0}$ and $r_1,\ldots,r_k>0$.
Conversely, a counter-example to the latter claim would show that
$\SympH (\T^{2n},\omega)\neq\Symp_0 (\T^{2n},\omega)$.

In the case when $M$ is a manifold underlying a K3 surface it is known by the result of Sheridan-Smith \cite{_Sheridan-Smith-JAMS2020_}
that $\SympH (M,\omega)\neq\Symp_0 (M,\omega)$ for certain
K\"ahler-type symplectic structures $\omega$ on $M$,
  including some irrational ones.
It has been known since the work \cite{_Seidel-2000_} of Seidel that for some (not necessarily irrational) K\"ahler-type symplectic forms $\omega$ on $M$ one has $\Symp (M,\omega)\cap\Diff_0 (M)\neq\Symp_0 (M,\omega)$; by a result of Smirnov \cite{_Smirnov-GAFA2022_}, this definitely holds for some irrational K\"ahler-type symplectic forms $\omega$ on $M$.

\bigskip
\noindent
3. \ref{existence-linear-tori_Corollary_} for one ellipsoid -- and, in particular, for one ball -- was previously proved by Luef and Wang  \cite{_Luef-Wang_} using a method similar to the one in this paper.
Let us note that the paper \cite{_Luef-Wang_} by Luef and Wang provides a very interesting connection between K\"ahler-type embeddings of ellipsoids into tori and Gabor frames -- an important notion in signal processing.

We expect that the technique used in \cite{_Luef-Wang_}, as well as the
techniques developed in \cite{_EV-Selecta_} in relation to symplectic embeddings of ellipsoids into tori and K3 surfaces,
can be used to generalize \ref{_existence-connectedness-irrational-forms-tori-K3-surfaces_Theorem_}
to K\"ahler-type embeddings of ellipsoids.


\hfill


\remark
\label{_Seshadri-constants_Remark_}

The assumption that the symplectic form $\omega$ on $M$ is irrational is essential in \ref{_existence-connectedness-irrational-forms-tori-K3-surfaces_Theorem_}: otherwise there may be obstructions
for the existence of K\"ahler-type embeddings of (disjoint unions of) balls into $(M,\omega)$ that
are independent of the symplectic volume -- for instance, the obstructions coming from Seshadri constants.

For simplicity, we will explain this below in the case of one ball and $M=\T^{2n}$,
 expanding on \cite[Sec. 2.2]{_LMcDS_}.

First, recall the definition of Seshadri constants. Let $(X,I)$ be a complex projective manifold and let $\omega$ be a K\"ahler form on $(X,I)$
such that the cohomology class $[\omega]$ is integral and equals the first Chern class of an ample line bundle on $(X,I)$.
Assume that $x\in X$, $\pi: \tX\to X$ is the complex blow-up of $X$ at $x$, and $E:= \pi^{-1} (x)$ is the exceptional divisor.
Define the {\bf Seshadri constant}
\[
\epsilon (X,I,[\omega],x) := \sup \left\{ s>0\ \left|\right.\ \left\langle\pi^* [\omega],[\tC]\right\rangle - s E \cdot [\tC] \geq 0\right.
\]
\[
\left. \textrm{for all (possibly singular) complex curves}\ \tC\subset \tX\ \right\}.
\]
Observe that $E \cdot [\tC]$ is the multiplicity of the singularity of $C$ at $x$.

Then
\begin{align}
\label{_Seshadri-geq-max-ball-admitting-Kahler-embedding_Equation_}
\epsilon (X,I,[\omega],x) = &  \sup \left\{\ \pi r^2\ |\ \textrm{There exists a K\"ahler form}\ \omega',\ [\omega']=[\omega],\right. \nonumber \\
\ & \textrm{and a K\"ahler embedding}\ f: B^{2n} (r)\to (X,\omega',I),  \\
\ & \left.  f^*\omega' =\omega_0,\ f^* I = J_0,\ f(0)=x\ \right\} \nonumber.
\end{align}
For a proof see \cite[Thm. 0.6]{_Eckl2017_}, \cite[Thm. 1.4]{_Fleming2021_}, \cite[Thm. 1.3]{_WittNystrom2_}, cf. \cite[Thm. A]{_Luef-Wang_}. (It had been previously observed in \cite{_Lazarsfeld-MRL_} that the inequality $\epsilon (X,I,[\omega],x) \geq$ [the right-hand side in \eqref{_Seshadri-geq-max-ball-admitting-Kahler-embedding_Equation_}] follows  from the results of \cite{_McD-Polt_}).

Assume that $\omega$ is a rational symplectic form on $\T^{2n}$. For the study of the symplectic and K\"ahler-type embeddings in $(\T^{2n},\omega)$ we may assume then, without loss of generality, that $[\omega]\in H^2 (\T^{2n};\Z)$. Let $I$ be a complex structure on $\T^{2n}$
compatible with $\omega$. Consider an ample holomorphic bundle $L$ on $\T^{2n}$ such that $c_1 (L) = [\omega]$.
The complex manifold $(\T^{2n},I)$ can be viewed as an abelian variety and $L$ as its polarization.

Since the
 identity component of the group of biholomorphisms of a complex torus acts transitively on the torus,
 $\epsilon (\T^{2n},I,[\omega],x)$ depends only on $I$ and $[\omega]$ and not on $x\in\T^{2n}$, so we will denote $\epsilon (\T^{2n},I,[\omega]):= \epsilon (\T^{2n},I,[\omega],x)$.

Since $\Vol (B^{2n} (r),\omega_0) = \pi^n r^{2n}$, \eqref{_Seshadri-geq-max-ball-admitting-Kahler-embedding_Equation_} implies that
\[
\sup \left\{\ \Vol (B^{2n} (r),\omega_0)\ \left|\right.\ \textrm{There exists a K\"ahler-type embedding}\right.
\]
\[
\left. f: (B^{2n} (r),\omega_0)\to (\T^{2n},\omega)\ \right\} = \sup_{I\in \cmpt (M,\omega)} \epsilon^n (\T^{2n},I,[\omega]),
\]
where the supremum in the right-hand side is taken over all $I$ compatible with $\omega$.

There exist examples where
\[
\sup_{I\in \cmpt (M,\omega)} \epsilon^n (\T^{2n},I,[\omega]) < \Vol (\T^{2n},\omega).
\]
For instance, consider the torus $\T^4 = \R^4/\Z^4$ equipped with the symplectic form $\omega = dp_1\wedge dq_1 + dp_2\wedge dq_2$, where $p_1,p_2,q_1,q_2$ are the coordinates on $\R^4$. For any complex structure $I$ on $\T^4$ compatible with $\omega$ one can biholomorphically identify $(\T^4, I)$ with a principally polarized abelian variety. In this case, by a result of Steffens \cite[Prop. 2]{_Steffens-MathZ_}, we have $\epsilon (\T^4,I,[\omega])\leq 4/3$ for all $I$ compatible with $\omega$. At the same time $\Vol (\T^4,\omega) = \int_{\T^4} \omega^2 = 2$. Thus,
\[
\sup_{I\in \cmpt (M,\omega)}\epsilon^2 (\T^4,I,[\omega]) \leq (4/3)^2 < 2 = \Vol (\T^4,\omega).
\]
This yields an obstruction, independent of the symplectic volume, for the existence of K\"ahler-type embeddings of a ball into $(\T^4,\omega)$:
a ball whose volume is between $(4/3)^2$ and $2= \Vol (\T^4,\omega)$ does not admit K\"ahler-type embeddings into $(\T^4,\omega)$.
At the same time, by
 \cite[Thm. 1.1]{_LMcDS_},
 such a ball does admit a
symplectic embedding into $(\T^4,\omega)$. {\it Thus, there exists
a symplectic embedding of a ball which is not of K\"ahler type}.
We do not know whether such an embedding is tame.


\hfill


Let us now consider tame embeddings of balls into tori and K3 surfaces.


\hfill


\theorem
\label{_packing-by-arb-shapes-main-tori-K3_Theorem_}

Let $M$, $\dim_\R M = 2n$, be either $\T^{2n}$ or a smooth manifold (of real dimension $4$) underlying a complex K3 surface.

Let $\omega_1$, $\omega_2$ be
K\"ahler-type forms on $M$ that are irrational and satisfy
$\int_M \omega_1^n = \int_M \omega_2^n >0$.

Let $W_i\subset \R^{2n}$, $i=1,\ldots,k$,  be compact domains with piecewise-smooth boundary whose interiors contain the origin.
Assume that $H^2 (W_i;\R) = 0$ for all $i$. Set $\bfW:= \bigsqcup_{i=1}^k W_i$.

Then, for any $\varepsilon >0$,
\[
\nu_{T,\varepsilon} (M,\omega_1,\bfW) = \nu_{T,\varepsilon} (M,\omega_2,\bfW).
\]


\hfill


For the proof of \ref{_packing-by-arb-shapes-main-tori-K3_Theorem_} see Section~\ref{_arb-shapes-pfs_Section_} (\ref{_packing-by-arb-shapes-main-tori-K3-COPY_Theorem_}).


\hfill


\remark

\ref{_packing-by-arb-shapes-main-tori-K3_Theorem_} generalizes a similar result in \cite{_EV-JTA_} about symplectic (not necessarily
tame) embeddings and its proof is similar to the proof in \cite{_EV-JTA_}.


\hfill


Let us present applications of \ref{_packing-by-arb-shapes-main-tori-K3_Theorem_}.


\hfill


\corollary
\label{_full-packing-of-tori-K3-by-balls_Corollary_}

Assume:

\begin{itemize}

\item{} $M$ is either $\T^{2n}$ or a smooth manifold underlying a K3 surface.

\item{} $\bfW:= \bigsqcup_{i=1}^k B^{2n} (r_i)$ is a disjoint union of $k$ (possibly different) balls.

\item{}
$\varepsilon >0$.

\end{itemize}

Then the following claims hold:

\bigskip
\noindent
A.
For any irrational K\"ahler-type symplectic form $\omega$ on $M$ we have
\[
\nu_K (M,\omega,\bfW) = 1,
\]
meaning that K\"ahler-type embeddings of $\lambda\bfW$ into $(M,\omega)$ are
unobstructed.

\bigskip
\noindent
B. There exists a $\Diff^+ (M)$-invariant open dense set of K\"ahler-type symplectic forms on $M$, depending on $\bfW$ and $\varepsilon$ and containing, in particular, all irrational K\"ahler-type symplectic forms on $M$,
so that for each $\omega$ in this set
$\nu_{T,\varepsilon} (M,\omega,\bfW) = 1$ -- meaning that $\varepsilon$-tame embeddings of $\lambda\bfW$ into $(M,\omega)$ are
unobstructed.


\hfill


\corollary
\label{_full-packing-of-tori-by-polydisks-parallelepipeds_Corollary_}

Assume:

\begin{itemize}

\item{} $M=\T^{2n}$.

\item{} $\bfW:= \bigsqcup_{i=1}^k W_i$ is either a disjoint union of $k$ identical copies of a $2n$-dimensional polydisk
\[
B^{2n_1} (R_1)\times\ldots\times B^{2n_l} (R_l),\ n_1+\ldots+n_l = n,\ R_1,\ldots, R_l>0, l>1,
\]
or a disjoint union of $k$ identical copies of a parallelepiped
\[
P (e_1,\ldots,e_{2n}) := \left\{ \sum_{j=1}^{2n} s_j e_j, 0\leq s_j \leq 1, j=1,\ldots,2n \right\},
\]
where $e_1,\ldots,e_{2n}$ is a basis of the vector space $\R^{2n}$.

\item{}
$\varepsilon >0$.

\end{itemize}

Then the following claims hold:

\bigskip
\noindent
A.
For any positive volume there exists a dense $\Diff^+ (\T^{2n})$-orbit (of an irrational K\"ahler-type symplectic form depending on $\bfW$) in the space of K\"ahler-type symplectic forms of that volume on $\T^{2n}$
such that for any $\omega'$ in this orbit we have
$\nu_K (\T^{2n},\omega',\bfW) = 1$ --
or, in other words, K\"ahler-type embeddings $\lambda\bfW\to (\T^{2n},\omega')$ are unobstructed.

\bigskip
\noindent
B. There exists a $\Diff^+ (\T^{2n})$-invariant open dense set of K\"ahler-type symplectic forms on $\T^{2n}$, depending on $\bfW$ and $\varepsilon$ and containing, in particular, all irrational K\"ahler-type symplectic forms on $\T^{2n}$,
so that for each $\omega'$ in this set
$\nu_{T,\varepsilon} (\T^{2n},\omega',\bfW) = 1$ --
or, in other words, $\varepsilon$-tame embeddings $\lambda\bfW\to (\T^{2n},\omega')$ are unobstructed.


\hfill


\corollary
\label{_full-packing-of-K3-by-parallelepipeds_Corollary_}

Assume:

\begin{itemize}

\item{} $M$ is a smooth manifold underlying a K3 surface.

\item{} $\bfW:= \bigsqcup_{i=1}^k W_i$ is a disjoint union of $k$ identical copies of a parallelepiped
\[
P (e_1,\ldots,e_{2n}) := \left\{ \sum_{j=1}^{2n} s_j e_j, 0\leq s_j \leq 1, j=1,\ldots,2n \right\},
\]
where $e_1,\ldots,e_{2n}$ is a basis of the vector space $\R^{2n}$.

\item{} $\varepsilon >0$.

\end{itemize}

Then there exists a $\Diff^+ (M)$-invariant open dense set of K\"ahler-type symplectic forms on $M$, depending on $\bfW$ and $\varepsilon$ and containing, in particular, all irrational K\"ahler-type symplectic forms on $M$,
so that for each $\omega'$ in this set
$\nu_{T,\varepsilon} (M,\omega',\bfW) = 1$, meaning that $\varepsilon$-tame embeddings of $\lambda\bfW$ into $(M,\omega')$ are unobstructed.


\hfill


For the proof of \ref{_full-packing-of-tori-K3-by-balls_Corollary_}, \ref{_full-packing-of-tori-by-polydisks-parallelepipeds_Corollary_} and \ref{_full-packing-of-K3-by-parallelepipeds_Corollary_} see Section~\ref{_arb-shapes-pfs_Section_}. These corollaries generalize similar results for symplectic -- but not necessarily tame -- embeddings that were proved in \cite{_EV-JTA_}.


\hfill


\remark
\label{_sympl-tame-Kahler-type_Remark_}

We know only one way to show the possible existence of embeddings (of a disjoint union of domains) into $(M,\omega)$ that are symplectic and non-tame, or tame and not K\"ahler-type -- via the volume restrictions: e.g., if K\"ahler-type embeddings cannot fill more than a certain portion of the volume of $(M,\omega)$ while tame embeddings can. This is how we got in
\ref{_Seshadri-constants_Remark_} that there exist symplectic embeddings (of balls) that are not of K\"ahler-type.
It would be interesting to find out whether there are other ways of detecting symplectic and non-tame, or tame and not K\"ahler-type, embeddings.

Together with \ref{_full-packing-of-tori-K3-by-balls_Corollary_}, this leads to the following questions.


\hfill


\question
\label{_existence-of-sympl-not-tame-or-tame-not-K-type-embs_Question_}

 Does there exist a closed manifold equipped with a K\"ahler-type symplectic form that admits a
symplectic embedding which is not tame, or a tame embedding which is not of K\"ahler-type? In particular, does there exist a rational K\"ahler-type symplectic form $\omega$ on $M=\T^{2n}$, $n>1$, or on a smooth manifold $M$ underlying a complex K3 surface, and a disjoint union of balls $\bfW$, $\Vol (\bfW, \omega_0) < \Vol (M,\omega)$,  so that $\bfW$ does not admit a tame embedding into $(M,\omega)$?


\hfill


\subsection{The IHS-hyperk\"ahler case}
\label{_hyperkahler-case_Subsection_}

K3 surfaces are particular examples of so
called IHS-hyperk\"ahler manifolds,
that is, compact hyperk\"ahler manifolds
of maximal holonomy. Most of the results
of Section~\ref{_main-results-tori-K3_Subsection_}
generalize to such manifolds. In this section we present
relevant definitions and statements.

Recall that a {\bf hyperk\"ahler manifold} is a manifold
equipped with three complex structures $I_1,I_2,I_3$ satisfying the
quaternionic relations and three symplectic forms
$\omega_1,\omega_2,\omega_3$ compatible, respectively, with
$I_1,I_2,I_3$, so that the three Riemannian metrics $\omega_i
(\cdot, I_i\cdot)$, $i=1,2,3$, coincide. Such a collection $\goH=\{ I_1, I_2,I_3, \omega_1,\omega_2,\omega_3\}$ of
complex structures and symplectic forms on a manifold is called a
{\bf hyperk\"ahler structure}. All the complex and
symplectic structures appearing in a hyperk\"ahler
structure induce the same orientation on the manifold.

Assume that $M$ is a closed connected and simply connected manifold.

The Levi-Civita connection of a K\"ahler manifold
preserves its complex structure and the K\"ahler form.
Therefore, the Levi-Civita connection of a hyperk\"ahler
manifold preserves the three complex structures $I_1,I_2,I_3$ and the three
symplectic forms $\omega_1,\omega_2,\omega_3$.
The subgroup of $GL(\R^{4n})$ preserving both the quaternionic
Hermitian metric and the quaternionic action is
called {\bf the compact symplectic group}  and is denoted
by $Sp(n)$. This is the group of all quaternionic
unitary matrices; it can also be obtained as
a maximal compact form of a complex symplectic
group $\Sp(\C^{2n})$. Since the Levi-Civita
connection of $(M, g, I_1,I_2,I_3)$
preserves the metric and the quaternionic operators,
its holonomy belongs to $Sp(n)$. However, a priori the
holonomy group  can be smaller; for example,
for a hyperk\"ahler flat torus, the holonomy is trivial.
The de Rham theorem and Berger's classification of
irreducible holonomies imply that any
compact hyperk\"ahler manifold has a finite
covering which is decomposed into a product
of a torus and manifolds with maximal holonomy;
this result is called {\bf the Bogomolov
decomposition theorem}, \cite{_Bogomolov:decompo_}.
This theorem is the reason why the
maximal holonomy hyperk\"ahler manifolds are also
called {\bf irreducible
holomorphically symplectic (IHS)}.

A maximal holonomy compact hyperk\"ahler manifold
can be characterized cohomologically: $M$
has maximal holonomy if and only $h^{2,0}(M)=1$
and $b_1(M)=0$. Here $h^{2,0}(M)$ denotes the
dimension of the group of closed, holomorphic 2-forms on $M$,
for any of its complex structures.

K3 surfaces, as well as the Hilbert schemes of
points for K3 surfaces, are known to admit IHS-hyperk\"ahler
structures.

We say that a symplectic/complex structure on $M$ is
{\bf of IHS-hyperk\"ah\-ler type}, if it appears
in {\it some} maximal holonomy hyperk\"ahler structure. In
particular, all IHS-hyperk\"ahler-type symplectic and
complex structures are of K\"ahler type.


\hfill


\remark\label{_hyperkahler-structures_Remark_}

\noindent
1. Note that even if a symplectic and a complex structure are of IHS-hyperk\"ahler type and compatible, it still does not mean that they can be included in {\it the same} IHS-hyperk\"ahler structure.

\medskip
\noindent
2. All complex structures and all K\"ahler-type symplectic structures on a K3 surface are in fact of IHS-hyperk\"ahler type -- see e.g. \cite[Prop. 3.1]{_EV-K3_}.
As we have already mentioned, it is an open question whether any symplectic form on a K3 surface is of K\"ahler type.

\bigskip
\noindent
3. Manifolds admitting IHS-hyperk\"ahler structures are also
called hyper\-k\"ahler manifolds of maximal holonomy, because
the holonomy group of a hyper\-k\"ahler manifold is $\Sp(n)$ (the
group of invertible quaternionic $n\times n$-matrices) -- and not
its proper subgroup -- if and only if the hyperk\"ahler structure is IHS
\cite{_Besse:Einst_Manifo_}.

\bigskip
\noindent
4.
Note that the property of being IHS is a topological
invariant: a holomorphically symplectic K\"ahler manifold
$M$, $\dim_\C M=2n$, is IHS if and only if $\pi_1(M)=0$ and
 there exists a primitive integral quadratic form $q$ on $H^2(M;\R)$
(called the Bo\-go\-mo\-lov-Beauville-Fujiki form \cite{_Bea1_}, \cite{_Fujiki:HK_})
that satisfies, for some rational $C>0$,
 the
Fujiki relation $\int_M \eta^{2n} = C q(\eta, \eta)^n$.

Therefore, the set of
IHS-hyperk\"ahler type symplectic/complex structures on
a given IHS-hyperk\"ahler manifold coincides with the set of
hyperk\"ahler-type symplectic/complex structures -- i.e., the symplectic/complex structures appearing in some
hyperk\"ahler structure.


\hfill


Now we can go over Section~\ref{_preliminaries_Section_} and define {\bf IHS-hy\-per\-k\"ahler-type}
embeddings by replacing the words ``K\"ahler-type" in the relevant definitions by the words ``IHS-hy\-per\-k\"ahler-type". In particular, in this way we can define {\bf IHS-hy\-per\-k\"ahler-type symplectic embeddings} of domains in $\R^{4n}$ into $M$ equipped with an IHS-hy\-per\-k\"ahler-type symplectic form.

Note that, in view of part 2 of \ref{_hyperkahler-structures_Remark_}, for symplectic embeddings into a K3 surface being of IHS-hyperk\"ahler type is the same as being of K\"ahler type.

The results of Section~\ref{_main-results-tori-K3_Subsection_} on K\"ahler-type
embeddings generalize to the IHS-hyperk\"ahler case as follows:

\smallskip
\noindent
-  Part A of \ref{_existence-connectedness-irrational-forms-tori-K3-surfaces_Theorem_}
generalizes directly to the IHS-hyperk\"ahler case as long as one assumes that the Campana-simple complex structure $I$ is of IHS-hy\-per\-k\"ahler-type.

\smallskip
\noindent
- The claim in part C of
\ref{_existence-connectedness-irrational-forms-tori-K3-surfaces_Theorem_}
concerning the action of $\Symp (M,\omega)\cap\Diff_0 (M)$
generalizes to the IHS-hyperk\"ahler case, if $\compK (M)$ is replaced by the space of the IHS-hy\-per\-k\"ahler-type complex structures on $M$.
For a general IHS-hyperk\"ahler manifold $(M,\omega)$ we cannot say anything about the transitivity of the
$\SympH (M,\omega)$-action.

\smallskip
Let us discuss the results on
$\varepsilon$-tame
embeddings into $(M,\omega)$ for an IHS-hy\-per\-k\"ahler-type form $\omega$.
One can define
tame and $\varepsilon$-tame
embeddings into $(M,\omega)$ exactly as in \ref{_Kahler-type-partially-Kahler-type-embs_Definition_}
using {\sl K\"ahler-type} complex structures. Another possibility is to modify the definition of
tame and $\varepsilon$-tame
embeddings into $(M,\omega)$
by using {\sl IHS-hy\-per\-k\"ahler-type} complex structures in \ref{_Kahler-type-partially-Kahler-type-embs_Definition_} instead of K\"ahler-type ones. The claim $\nu_{T,\varepsilon} (M,\omega_1,\bfW) = \nu_{T,\varepsilon} (M,\omega_2,\bfW)$ of \ref{_packing-by-arb-shapes-main-tori-K3_Theorem_} holds for the original $\nu_{T,\varepsilon}$ defined in \ref{_nu_Definition_}, as well as for a version of $\nu_{T,\varepsilon}$ defined using the
$\varepsilon$-tame
embeddings in the second sense above,
if $\omega_1$ and $\omega_2$ lie in the same connected component of the space of IHS-hy\-per\-k\"ahler-type symplectic forms on $M$.

In Section~\ref{_hyperkahler-case-proofs_Section_} we discuss how to modify
the proofs of the results of Section~\ref{_main-results-tori-K3_Subsection_} in order to obtain the above-mentioned results in the IHS-hyperk\"ahler case.


\hfill


\section{Plan of the paper and an outline of the proofs}
\label{_plan-of-the-paper-outline-of-the-pfs_Section_}

In this section we describe the plan of the paper and outline the proofs.

In order to simplify the exposition, let us fix a connected component $\cC_0$ of $\compK (M)$ and the radii $r_1,\ldots,r_k>0$
of the balls from the beginning. Let $\bfr := (r_1,\ldots,r_k)$.

In Section~\ref{_Kahler-type-embs-of-balls-and-blow-ups-pfs_Section_} we prove \ref{_existence-connectedness-Kahler-type-arb-mfds_Theorem_}
and
\ref{_Campana-simple-Kahler-type-existence-connectedness-gen-mfd_Theorem_}.
Namely, we define the following spaces:

\smallskip
\noindent
- The space $\Pairs_{\cC_0}$ formed by pairs $(\eta,f)$, where $\eta$ is a K\"ahler-type symplectic form on $M$
isotopic to $\omega$ and $f: \bigsqcup_{i=1}^k B^{2n} (r_i)\to (M,\eta)$ is a K\"ahler-type embedding favoring $\cC_0$.

\smallskip
\noindent
- The Teichm\"uller space $\TeichP_{\cC_0} := \Pairs_{\cC_0}/\Diff_0 (M)$.

\smallskip
\noindent
- The space $\Triples_{\cC_0}$ formed by triples $(I,\bfx,h)$, where $I\in \Diff_0 (M)\cdot\cmpt (M,\omega)$, $\bfx = (x_1,\ldots,x_k)\in \hM^k$, so that the pair $(I,\bfx)$ lies in the set $\cK_{\cC_0} (\bfr)$ defined in \ref{_cK-ttK_Definition}, and
$h:=\bigsqcup_{i=1}^k h_i: \bigsqcup_{i=1}^k B^{2n} (r_i)\to (M,I)$
is a holomorphic embedding such that $h_i (0) = x_i$, $i=1,\ldots, k$.

\smallskip
\noindent
- The Teichm\"uller space $\TeichT_{\cC_0} := \Triples_{\cC_0}/\Diff_0 (M)$.

We then construct a continuous surjective map $\Phi: \TeichT_{\cC_0}\to \TeichP_{\cC_0}$. The construction uses
the symplectic blow-up and blow-down operations similar to the ones used in \cite{_McD-Polt_}. Let us note that
for the local constructions on which these operations are based we use
the regularized maximum of plurisubharmonic functions while in \cite{_McD-Polt_} a different analytic technique is used.
We also need the regularized maximum tool to construct a K\"ahler form on $(\tM_{I,\bfx},\tI)$ in the class
$\Pi^* [\omega] - \pi \sum_{i=1}^k r_i^2 e_i\in H^2 (\tM_{I,\bfx};\R)$ which is standard {\sl on a neighborhood of each
exceptional divisor} -- this is needed to obtain a K\"ahler form on $M$ by the symplectic blow-down construction (unlike in
\cite{_McD-Polt_}, where the goal was to produce a symplectic, not necessarily K\"ahler, form on $M$, and this was done much
easier, using Moser's method for forms on $\tM_{I,\bfx}$).

The symplectic blow-down construction, enhanced as we have just described, and a relative version of Moser's method are used to construct $\Phi$, and the symplectic blow-up is used to verify that
$\Phi$ is surjective. The continuity of $\Phi$
is verified using the observation that the set $\TeichP_{\cC_0}$ is discrete and therefore it suffices
to check that the preimage of a point is open, which is then proved using again a relative version of Moser's method.

We then prove \ref{_existence-connectedness-Kahler-type-arb-mfds_Theorem_} as follows.

The mere existence of $\Phi$ implies that if $\cK_{\cC_0} (\bfr)$ (and hence $\TeichT_{\cC_0}$) is non-empty, there exists a K\"ahler-type embedding $\bigsqcup_{i=1}^k B^{2n} (r_i)\to (M,\omega)$ favoring $\cC_0$.

The connectivity of the set $\ttK_{\cC_0} (\bfr)$ (see \ref{_cK-ttK_Definition}) yields the transitivity of the
$\Symp (M)\cap\Diff_0 (M)$-action on the space of K\"ahler-type embeddings $\bigsqcup_{i=1}^k B^{2n} (r_i)\to (M,\omega)$ favoring $\cC_0$ as follows. Using methods of complex geometry (the Kodaira-Spencer stability and the Demailly-Paun description of the K\"ahler cones of fibers of an analytic deformation family)
we show that the connectivity of $\ttK_{\cC_0} (\bfr)$ is equivalent to the connectivity of $\cK_{\cC_0} (\bfr)$. The latter
easily implies the connectivity of $\TeichT_{\cC_0}$. Since $\Phi$ is continuous and surjective, we get that $\TeichP_{\cC_0}$ is
connected, and since it is discrete, it is just a point. This easily implies the transitivity of the $\Symp (M)\cap\Diff_0 (M)$-action.

In Section~\ref{_Campana-simple-strs-gen-mfds-main-results_Section_} we consider the case of a closed symplectic manifold for which
the symplectic form is compatible with a so-called Campana-simple complex structure -- i.e., a K\"ahler-type complex structure admitting ``few" complex subvarieties. Examples of such manifolds include tori and smooth manifolds underlying complex K3 surfaces equipped with irrational K\"ahler-type symplectic forms. The general results on K\"ahler-type embeddings of balls into such manifolds, proved in Section~\ref{_Campana-simple-strs-gen-mfds-main-results_Section_}, are used further in the paper to deduce the results on such embeddings into tori
and K3 surfaces stated in Section~\ref{_main-results-tori-K3_Subsection_}.
The key existence result on such embeddings is based on the observation, made already in \cite{_LMcDS_}, \cite{_EV-JTA_},
that for any Campana-simple complex structure $I$ compatible
with $\omega$ and any Campana-generic $x_1,\ldots,x_k\in (M,I)$, $\bfx := (x_1,\ldots,x_k)$, we have $(I,\bfx)\in \cK (\bfr)$, as long as
$\Vol (\bigsqcup_{i=1}^k B^{2n} (r_i))< \Vol (M,\omega)$. This follows from the Demailly-Paun description of the K\"ahler cone
of a closed K\"ahler manifold \cite{_Dem-Paun_}.

In Sections~\ref{_complex-strs-on-proj-spaces-their-products-and-blowups_Section_}, \ref{_Hirzebruch-odd_Section_}, \ref{_Hirzebruch-even_Section_},
\ref{_Campana-simple-structures-tori-K3_Section_} we
discuss the cases
of the specific $M$: complex projective spaces and their products, Hirzebruch surfaces, tori and K3 surfaces. In each of these cases we discuss the structure of the Teichm\"uller
space of K\"ahler-type complex structures on $M$ and of its subset formed by equivalence classes of complex structures
compatible with a given symplectic form on $M$, as well as the K\"ahler cones of the complex blow-ups of $M$.
We use this information in order to apply our general results about K\"ahler-type embeddings of balls
and deduce the results in
Section~\ref{_main-results-all_Section_}
about these specific $M$.

In Section~\ref{_arb-shapes-pfs_Section_} we discuss the proofs of the results in Section~\ref{_main-results-tori-K3_Subsection_}.
concerning tame embeddings of arbitrary domains.
The strategy of the proof of
\ref{_packing-by-arb-shapes-main-tori-K3_Theorem_}
is the same as for a similar claim in \cite{_EV-JTA_}.
Namely, with $M$ being either $\T^{2n}$ or a smooth manifold underlying a K3 surface,
the group $\Diff^+ (M)$ of orientation-preserving diffeomorphisms of $M$ acts on the space
$\sympK_1 (M)$
of K\"ahler-type forms on $M$ of
total volume $1$. The function
 $\omega\mapsto \nu_{T,\varepsilon} (M,\omega,\bfW)$ is clearly invariant under the action. We will show
that this function is lower semicontinuous (with respect to the
$C^1$-topology on $\sympK (M)$) -- the proof is a modification of the proof of a similar claim in \cite{_EV-JTA_}.
Then we use the result from \cite{_EV-JTA_} saying that
the orbit of a form $\omega\in\sympK_1 (M)$ under the
action of $\Diff^+ (M)$ is dense in $\sympK_1 (M)$, as long as $\omega$ is irrational. Since the orbits
of both $\omega_1$ and $\omega_2$ are dense in $\sympK_1 (M)$, we get, by the
lower semicontinuity, that
$\nu_{T,\varepsilon} (M,\omega_1, \bfW) \leq \nu_{T,\varepsilon} (M,\omega_2, \bfW)$ and
$\nu_{T,\varepsilon} (M,\omega_1, \bfW) \geq \nu_{T,\varepsilon} (M,\omega_2, \bfW)$,
which means that
$\nu_{T,\varepsilon} (M,\omega_1, \bfW) = \nu_{T,\varepsilon} (M,\omega_2, \bfW)$.
Then the proofs of \ref{_full-packing-of-tori-K3-by-balls_Corollary_}, \ref{_full-packing-of-tori-by-polydisks-parallelepipeds_Corollary_} and \ref{_full-packing-of-K3-by-parallelepipeds_Corollary_} follow the same lines as the proof of \cite[Cor. 3.3]{_EV-JTA_}
where the existence of symplectic embeddings of polydisks into tori was proved.

In Section~\ref{_hyperkahler-case-proofs_Section_} we discuss how the proofs in the K3 case can be modified (or just carried over without any changes) in order to prove the results in the hyperk\"ahler case stated in Section~\ref{_hyperkahler-case_Subsection_}.

Finally, in Section~\ref{_Appendix_Section_} (the Appendix) we discuss various well-known facts concerning the dependence of the Hodge decomposition on the
complex structure, the deformations of complex structures, Moser's method, and Alexander's trick that we use in the paper and
for which we could not find a direct reference in the literature.


\hfill


\section{K\"ahler-type embeddings of balls and blow-ups}
\label{_Kahler-type-embs-of-balls-and-blow-ups-pfs_Section_}

As before, let $M$, $\dim_\R M = 2n$, be a closed connected manifold admitting K\"ahler structures. Let $\omega$ be a K\"ahler-type symplectic form on $M$.

Let $k\in\Z_{>0}$, $r_1,\ldots,r_k>0$, $\bfr = (r_1,\ldots,r_k)$.

Let $\cC_0$ be a connected component of $\compK (M)$.

\bigskip
\noindent
{\bf
From now on,
all the spaces of maps between manifolds and of differential forms and complex structures on manifolds will be endowed with the $C^\infty$-topologies. All the product/quotient spaces will be endowed with the product/quotient topologies.
}


\hfill


\subsection{The $\Symp (M,\omega)\cap \Diff_0 (M)$-action on the space of
K\"ah\-ler-ty\-pe em\-bed\-dings}
\label{_Diff0-cap-Symp-action-on-embeddings_Subsection_}

Denote:

\smallskip
\noindent
$\EmbS_{M,\omega} \left(\bigsqcup_{i=1}^k B^{2n} (r_i)\right)$ -- the space of the symplectic embeddings
$\bigsqcup_{i=1}^k B^{2n} (r_i)\to (M,\omega)$;

\smallskip
\noindent
$\EmbK_{M,\omega} \left(\bigsqcup_{i=1}^k B^{2n} (r_i)\right)\subset \EmbS_{M,\omega} \left(\bigsqcup_{i=1}^k B^{2n} (r_i)\right)$ --
the space of the K\"ahler-type embeddings
$\bigsqcup_{i=1}^k B^{2n} (r_i)\to (M,\omega)$;

\smallskip
\noindent
$\EmbK_{M,\omega,\cC_0} \left(\bigsqcup_{i=1}^k B^{2n} (r_i)\right)$ --
the space of the K\"ahler-type
embeddings
$\bigsqcup_{i=1}^k B^{2n} (r_i)\to (M,\omega)$
favoring $\cC_0$.

\smallskip
Given a symplectic embedding $f:\bigsqcup_{i=1}^k B^{2n} (r_i)\to (M,\omega)$, any symplectic embedding $f':\bigsqcup_{i=1}^k B^{2n} (r_i)\to (M,\omega)$ sufficiently $C^\infty$-close to $f$ can be mapped in $f$ by an element of $\Symp_0 (M,\omega)$, and moreover, this element can be chosen arbitrarily close to the identity, provided $f'$ is sufficiently $C^\infty$-close to $f$. This fact is standard and follows from the well-known Alexander trick and the symplectic isotopy extension results -- see \ref{_Kahler-type-embs-form-an-open-set_Corollary_} for details.

Since symplectomorphisms of $(M,\omega)$ map K\"ahler-type embeddings into K\"ahler-type embeddings, this yields that
$\EmbK_{M,\omega} \left(\bigsqcup_{i=1}^k B^{2n} (r_i)\right)$
is an open $\Symp (M,\omega)$-invariant subset of $\EmbS_{M,\omega} \left(\bigsqcup_{i=1}^k B^{2n} (r_i)\right)$
(in the $C^\infty$-topology). This yields the following corollary.


\hfill


\vfil\eject

\proposition
\label{_symplectic-Teich-is-discrete_Proposition_}

The space
\[
\frac{\EmbK_{M,\omega} \left(\bigsqcup_{i=1}^k B^{2n} (r_i)\right)}{\Symp_0 (M,\omega)},
\]
and consequently the space
\[
\frac{\EmbK_{M,\omega} \left(\bigsqcup_{i=1}^k B^{2n} (r_i)\right)}{\Symp (M,\omega)\cap \Diff_0 (M)},
\]
and its subspace
\[
\frac{\EmbK_{M,\omega,\cC_0} \left(\bigsqcup_{i=1}^k B^{2n} (r_i)\right)}{\Symp (M,\omega)\cap \Diff_0 (M)},
\]
are discrete.
\endproof


\hfill


\definition
\label{_Pairs-TeichP_Definition_}

Define $\Pairs_{\cC_0}$ as the set of pairs $(\eta,f)$, where $\eta\in \Diff_0 (M)\cdot\omega$ and
$f: \bigsqcup_{i=1}^k B^{2n} (r_i)\to (M,\eta)$ is K\"ahler-type
embedding into $(M,\eta)$ favoring $\cC_0$. (Note that the embedding here is into $(M,\eta)$, not $(M,\omega)$!)

Define
\[
\TeichP_{\cC_0} := \Pairs_{\cC_0}/\Diff_0 (M).
\]
For each $(\eta,f)\in \Pairs_{\cC_0}$, denote by $\{\eta,f\}\in \TeichP_{\cC_0}$ the element of
$\TeichP_{\cC_0}$ represented by $(\eta,f)$ -- i.e., the $\Diff_0 (M)$-orbit of $(\eta,f)$.


\hfill


The following proposition follows easily from the definitions.


\hfill


\proposition
\label{_TeichP-homeom-to-symplectic-Teich_Proposition_}

The natural map
\[
\frac{\EmbK_{M,\omega,\cC_0} \left(\bigsqcup_{i=1}^k B^{2n} (r_i)\right)}{\Symp (M,\omega)\cap \Diff_0 (M)}\to \TeichP_{\cC_0}
\]
is a homeomorphism. In particular, in both cases $\TeichP_{\cC_0}$ is discrete.
\endproof


\hfill


A discrete topological space is connected if and only if it is a point. Together with \ref{_TeichP-homeom-to-symplectic-Teich_Proposition_}, this immediately yields the following corollary.


\hfill


\corollary
\label{_transitivity-equiv-to-connectedness-of-TeichP_Corollary_}

The $\Symp (M,\omega)\cap \Diff_0 (M)$-action on $\EmbK_{M,\omega,\cC_0} \left(\bigsqcup_{i=1}^k B^{2n} (r_i)\right)$
is transitive if and only if
$\TeichP_{\cC_0}$ is connected.
\endproof


\hfill


\proposition
\label{_smooth-family-of-pairs-lies-in-the-same-Diff0-orbit_Proposition_}

Assume $(\eta_t,f_t)$, $0\leq t\leq 1$, is a smooth family in  $\Pairs_{\cC_0}$.

Then all $(\eta_t,f_t)$, $0\leq t\leq 1$, lie in the same orbit of the $\Diff_0 (M)$-action on $\Pairs_{\cC_0}$.


\hfill


\noindent
{\bf Proof of \ref{_smooth-family-of-pairs-lies-in-the-same-Diff0-orbit_Proposition_}:}

The forms $\eta_t$, $0\leq t\leq 1$, form a smooth family of cohomologous symplectic forms. Hence, by Moser's theorem \cite{_Moser_}, all these forms can be identified by the $\Diff_0 (M)$-action with $\eta_0$. Under such an identification,
$\{ f_t\}_{0\leq t\leq 1}$ becomes
a smooth family of symplectic embeddings
$\{ f_t: \bigsqcup_{i=1}^k B^{2n} (r_i)\to (M,\eta_0)\}_{0\leq t\leq 1}$. By the symplectic isotopy extension result (see e.g. \cite[Thm. 3.3.1]{_McD-Sal-3_}), all such $f_t$ lie in the same $\Symp_0 (M,\eta_0)$-orbit. This proves the proposition.
\endproof


\hfill


\subsection{From holomorphic to K\"ahler-type embeddings of balls}
\label{_holomorphic-to-Kahler-type-embeddings_Subsection_}

As before, for each $\bfx:= (x_1,\ldots,x_k)\in \hM^k$
and each complex structure $I$ on $M$, let $\tM_{I,\bfx}$
denote the complex blow-up of $(M,I)$ at $x_1,\ldots,x_k$.

For any $I$ and $\bfx$ we will use the same notation for the following objects (suppressing the dependence on $I$ and $\bfx$):

\smallskip
\noindent
- The lift of $I$ to $\tM_{I,\bfx}$ will be denoted by $\tI$. This complex structure is
of K\"ahler-type (see e.g. \cite[Vol I, Prop. 3.24]{_Voisin-book_}). If $I'$ is another complex structure on $M$ coinciding with
$I$ near each $x_i$, $i=1,\ldots,k$, then it also lifts to $\tM_{I,\bfx}$ and its lift will be denoted by $\tI'$.

\smallskip
\noindent
- $\Pi: \tM_{I,\bfx}\to M$ will denote the natural projection.

\smallskip
\noindent
- $E_i := \Pi^{-1} (x_i)$, $i=1,\ldots,k$, is an exceptional divisor.

\smallskip
\noindent
- The cohomology classes Poincar\'e-dual to the fundamental homology classes of $E_i$, $i=1,\ldots,k$,
will be denoted by $e_1,\ldots,e_k\in H^2 (\tM_{I,\bfx};\R)$.

For all $I$ and $\bfx$ the groups $H^2 (\tM_{I,\bfx};\R)$ can be canonically identified (preserving the numbering of the classes $e_1,\ldots,e_n$) and we will use these identifications without further mention.


\hfill


\definition
\label{_triples_Definition_}

Define $\Triples_{\cC_0}$ as the set of triples $(I,\bfx,h)$, where

\begin{itemize}

\item{} $(I,\bfx)\in \cK_{\cC_0} (\bfr)$, meaning that
\[
I\in \Diff_0 (M)\cdot \left(\cmpt (M,\omega)\cap \cC_0\right),
\]
\[
\bfx:= (x_1,\ldots,x_k)\in \hM^k,
\]
and
the cohomology class $\Pi^* [\omega] - \pi \sum_{i=1}^k r_i^2 e_i\in H^2 (\tM_{I,\bfx};\R)$ is K\"ahler with respect to $\tI$.

\item{} $h = \bigsqcup_{i=1}^k h_i : \bigsqcup_{i=1}^k B^{2n} (r_i)\to (M,I)$
is a holomorphic embedding into $(M,I)$ so that $h_i (0) = x_i$, $i=1,\ldots,x_k$.
\end{itemize}

\bigskip
Set $\TeichT_{\cC_0} := \Triples_{\cC_0}/\Diff_0 (M)$.

For each $(I,\bfx,h)\in \Triples_{\cC_0}$, denote by $\{ I,\bfx,h\}\in \TeichT_{\cC_0}$ the element of
$\TeichT_{\cC_0}$ represented by $(I,\bfx,h)$ -- i.e., the $\Diff_0 (M)$-orbit of $(I,\bfx,h)$.


\hfill


For each $\kappa>0$ define
\[
R_\kappa: \R^{2n} \to \R^{2n},\ \cR_\kappa (x):= \kappa x.
\]

The proof of \ref{_existence-connectedness-Kahler-type-arb-mfds_Theorem_} is based on the following key propositions.


\hfill


\proposition
\label{_from-TeichT-to-pairs_Proposition_}

There exists a continuous surjective map $\Phi: \TeichT_{\cC_0}\to \TeichP_{\cC_0}$
with the following properties:

\bigskip
\noindent
1. For each $(I,\bfx,h)\in \Triples_{\cC_0}$,
\[
\Phi \left( \left\{ I,\bfx,h\right\}\right) = \left\{ \eta,h\circ R_\kappa \right\}
\]
for a symplectic form $\eta$ compatible with $I$ and cohomologous to $\omega$ and for a sufficiently small $0<\kappa\leq 1$, so that
the embedding $h\circ R_\kappa: \bigsqcup_{i=1}^k B^{2n} (r_i)\to M$
is symplectic with respect to $\eta$ and holomorphic with respect to $I$.

\bigskip
\noindent
2. If $\eta$ is a symplectic form on $M$ compatible with a complex structure $I$ and cohomologous to $\omega$ and
$f = \bigsqcup_{i=1}^k f_i: \bigsqcup_{i=1}^k B^{2n} (r_i)\to (M,\eta)$ is a symplectic
embedding into $(M,\eta)$ holomorphic with respect to $I$,
then, for $x_i := f_i (0)$, $i=1,\ldots,k$, and $\bfx:= (x_1,\ldots,x_k)$, we have
$(I, \bfx, f)\in \Triples_{\cC_0}$ and $\Phi \left(\left\{ I, \bfx, f \right\}\right) = \{ \eta, f\}$.


\hfill


\proposition
\label{_connectedness-of-TeichT-cC-0_Proposition_}

The space $\TeichT_{\cC_0}$ is connected
if only if the set $\ttK_{\cC_0} (\bfr)$ (see \ref{_cK-ttK_Definition}) is connected.


\hfill


For the proof of \ref{_from-TeichT-to-pairs_Proposition_} see Section~\ref{_blow-up-forms-Subsection_}. The proof  of \ref{_connectedness-of-TeichT-cC-0_Proposition_} will be given in Section~\ref{_connectedness-TeichT-C-zero_Subsection_}.

Let us now prove
\ref{_existence-connectedness-Kahler-type-arb-mfds_Theorem_}. For convenience, we restate it here.


\hfill


\theorem  {\bf (=\ref{_existence-connectedness-Kahler-type-arb-mfds_Theorem_})}
\label{_existence-connectedness-Kahler-type-arb-mfds-COPY_Theorem_}

Let $M$ be a closed manifold equipped with a K\"ahler-type symplectic form $\omega$.

\bigskip
\noindent
(I)
The set $\bigsqcup_{i=1}^k B^{2n}(r_i)$ admits a K\"ahler-type embedding into $(M,\omega)$ if and only if
for some $I\in \cmpt (M,\omega)$ and some $\bfx\in \hM^k$ the cohomology class $\Pi^* [\omega] - \pi \sum_{i=1}^k r_i^2 e_i\in H^2 (\tM_{I,\bfx};\R)$, is K\"ahler with respect to $\tI$.

More precisely, assume $I$ is a complex structure on $M$ compatible with $\omega$ and $\Sigma\subset (M,I)$ is a proper (possibly empty) complex submanifold. If $M = \C P^{n_1}\times \C P^1$ (with the product complex structure), we also allow
$\Sigma = (\C P^{n_1 -1}\times \C P^1) \cup (\C P^{n_1}\times \textrm{pt})$.

Then the following conditions are equivalent:

\begin{enumerate}

\item{} There exists $\bfx = (x_1,\ldots,x_k)\in \hM^k$, $x_1,\ldots,x_k\in M\setminus\Sigma$, so that the cohomology class $\Pi^* [\omega] - \pi \sum_{i=1}^k r_i^2 e_i\in H^2 (\tM_{I,\bfx};\R)$ is K\"ahler with respect to $\tI$.

\item{} There exists an $[I]$-K\"ahler-type embedding $f: \bigsqcup_{i=1}^k B^{2n}(r_i)\to (M\setminus\Sigma,\omega)$
which is holomorphic with respect to a complex structure $\phi_1^* I$ compatible with $\omega$, where
$\{ \phi_t\}_{0\leq t\leq 1}\subset \Diff_0 (M)$, $\phi_0 = Id$, is an isotopy such that $\phi_t (\Sigma) = \Sigma$ for all $t\in [0,1]$.

\end{enumerate}

\bigskip
\noindent
(II)
Let $\cC_0$ be a connected component of $\compK (M)$.
Assume that
the set $\ttK_{\cC_0} (\bfr)$ (or, equivalently, the set $\cK_{\cC_0} (\bfr)$)
is connected.

Then any two K\"ahler-type embeddings $\bigsqcup_{i=1}^k B^{2n}(r_i)\to (M,\omega)$ favoring $\cC_0$
lie in the same $\Symp (M,\omega)\cap \Diff_0 (M)$-orbit. In particular, there exists $[I]\in\cmptteich_{\cC_0} (M,\omega)$ such that
both embeddings are of $[I]$-K\"ahler type.

If, in addition, $\SympH (M)$ acts transitively on the set of connected components of $\compK (M)$ compatible with $\omega$, then
any two K\"ahler-type embeddings $\bigsqcup_{i=1}^k B^{2n}(r_i)\to (M,\omega)$
lie in the same $\SympH (M,\omega)$-orbit.


\hfill


\noindent
{\bf Proof of \ref{_existence-connectedness-Kahler-type-arb-mfds-COPY_Theorem_} (=\ref{_existence-connectedness-Kahler-type-arb-mfds_Theorem_}):}

Let us prove part (I) of \ref{_existence-connectedness-Kahler-type-arb-mfds-COPY_Theorem_}.

Let $\cC_0$ be a connected component of $\compK (M)$ containing the complex structure $I$.

Define the set $Y\supset\Sigma$ as follows: If $\Sigma$ is a complex submanifold of $(M,I)$, set $Y = \Sigma$.
If $M = \C P^{n_1}\times \C P^1$ and
$\Sigma = (\C P^{n_1 -1}\times \C P^1) \cup (\C P^{n_1}\times \textrm{pt})$, let $Y$ be the union of
$\Sigma$ and of the closure of an open neighborhood of $(\C P^{n_1 -1}\times \C P^1) \cap (\C P^{n_1}\times \textrm{pt})$.

Let us prove $1\Longrightarrow 2$. Assume that $\bfx = (x_1,\ldots,x_k)\in \hM^k$, $x_1,\ldots,x_k\in M\setminus\Sigma$,
so that the cohomology class $\Pi^* [\omega] - \pi \sum_{i=1}^k r_i^2 e_i\in H^2 (\tM_{I,\bfx};\R)$ is K\"ahler with respect to $\tI$. Choose a holomorphic embedding $h = \bigsqcup_{i=1}^k h_i : \bigsqcup_{i=1}^k B^{2n} (r_i+\epsilon)\to (M,I)$,
for some $\epsilon>0$,
so that $h_i (0) = x_i$, $i=1,\ldots, k$, and  $\Image h\cap Y = \emptyset$.
Then $(I,\bfx,h)\in \Triples_{\cC_0}$.

By part 1 of \ref{_from-TeichT-to-pairs_Proposition_},
there exist a symplectic form $\eta$ compatible with $I$ and cohomologous to $\omega$ and a sufficiently small $0<\kappa \leq 1$, so that the embedding $h\circ R_\kappa: \bigsqcup_{i=1}^k B^{2n} (r_i)\to M$ is symplectic with respect to $\eta$ and holomorphic with respect to $I$. In particular, $\Image (h\circ R_\kappa)\cap Y = \emptyset$. The linear family of 2-forms connecting $\eta$
and $\omega$ is formed by cohomologous symplectic forms compatible with $I$ and, in particular, their restrictions to $\Sigma$ are non-degenerate
(in the case  $M = \C P^{n_1}\times \C P^1$ and
$\Sigma = (\C P^{n_1 -1}\times \C P^1) \cup (\C P^{n_1}\times \textrm{pt})$, the restrictions of the forms to both
$\C P^{n_1 -1}\times \C P^1$ and $\C P^{n_1}\times \textrm{pt}$ are non-degenerate).
By a relative version of Moser's theorem (see \ref{_relative-version-of-Moser-method_Proposition_}, parts II and III),
there exists an isotopy $\{ \phi_t\}_{0\leq t\leq 1}\subset \Diff_0 (M)$, $\phi_0 = Id$, $\phi_t (\Sigma) = \Sigma$ for all $t\in [0,1]$, so that $(\phi_1^{-1})^* \omega = \eta$ outside $Y$.
Then $\phi^{-1}\circ h\circ R_\kappa$ is an embedding of $\bigsqcup_{i=1}^k B^{2n} (r_i)$ into $M\setminus \Sigma$
which is symplectic with respect to $\omega$ and holomorphic with respect to the complex structure $\phi^* I$ on $M$ compatible with $\omega$.

Now let us prove $2\Longrightarrow 1$. Assume $f= \bigsqcup_{i=1}^k f_i: \bigsqcup_{i=1}^k B^{2n}(r_i)\to (M\setminus\Sigma,\omega)$ is a K\"ahler-type embedding holomorphic with
respect to a complex structure $J:=\phi_1^* I$, which is compatible with $\omega$, where $\{ \phi_t\}_{0\leq t\leq 1}\subset \Diff_0 (M)$
is an isotopy such that $\phi_0 = Id$, $\phi_t (\Sigma) = \Sigma$ for all $t\in [0,1]$.
Let $x_i := f_i (0)$, $i=1,\ldots,k$. Set $\bfx := (x_1,\ldots,x_k)$.
Then, by property 2 of $\Phi$ listed in \ref{_from-TeichT-to-pairs_Proposition_},
$(J,\bfx,f)\in \Triples_{\cC_0}$,
meaning, in particular, that the cohomology class $\Pi^* [\omega] - \pi \sum_{i=1}^k r_i^2 e_i$ in $H^2 (\tM_{J,\bfx};\R)$ is K\"ahler with respect to $\tJ$. Since $x_1,\ldots,x_k\in M\setminus\Sigma$ and $\phi_1 (\Sigma) = \Sigma$,
we get that the cohomology class $\Pi^* [\omega] - \pi \sum_{i=1}^k r_i^2 e_i$ in $H^2 (\tM_{I,\bfx};\R)$ is K\"ahler with respect to $\tI$.

This finishes the proof of part (I) of \ref{_existence-connectedness-Kahler-type-arb-mfds-COPY_Theorem_}.

Let us prove part (II) of \ref{_existence-connectedness-Kahler-type-arb-mfds-COPY_Theorem_}.

Assume that the set $\ttK_{\cC_0} (\bfr)$
is connected.
By \ref{_connectedness-of-TeichT-cC-0_Proposition_}, this means that
$\TeichT_{\cC_0}$ is connected. Since, by \ref{_from-TeichT-to-pairs_Proposition_},
$\Phi: \TeichT_{\cC_0}\to \TeichP_{\cC_0}$ is a surjective continuous map, we get that
$\TeichP_{\cC_0}$ is connected. Hence, \ref{_transitivity-equiv-to-connectedness-of-TeichP_Corollary_}
implies that $\Symp (M,\omega)\cap\Diff_0 (M)$ acts transitively on $\EmbK_{M,\omega,\cC_0} \left(\bigsqcup_{i=1}^k B^{2n} (r_i)\right)$.
In particular, by \ref{_two-embs-in-same-conn-compt_Proposition_}, there exists $[I]\in\cmptteich_{\cC_0} (M,\omega)$ such that
both embeddings are of $[I]$-K\"ahler type.

Assume, in addition, that $\SympH (M)$ acts transitively on the set of connected components of $\compK (M)$ compatible with $\omega$.
Then any element of $\EmbK_{M,\omega} \left(\bigsqcup_{i=1}^k B^{2n} (r_i)\right)$
can be mapped by $\SympH (M)$ into an element of the space $\EmbK_{M,\omega,\cC_0} \left(\bigsqcup_{i=1}^k B^{2n} (r_i)\right)$.
The transitivity of the $\Symp (M,\omega)\cap\Diff_0 (M)$-action on the latter space yields that
any two K\"ahler-type embeddings
of $\bigsqcup_{i=1}^k B^{2n}(r_i)$ into $(M,\omega)$
can be mapped into each other by an element of $\SympH (M,\omega)$.

This finishes the proof of part (II) of \ref{_existence-connectedness-Kahler-type-arb-mfds-COPY_Theorem_}.
\endproof


\hfill


\subsection{Regularized maximum}
\label{_regularized maximum_Subsection_}

Given a $(1,1)$-form $\theta$ on a complex manifold with a complex structure $I$, we write $\theta>0$ if $\theta (v,I v) >0$ for any non-zero tangent vector $v$ of the manifold. If $\theta'$ is another $(1,1)$-form on the same complex manifold, we write $\theta>\theta'$ if $\theta-\theta'>0$.

Recall from \cite[Lemma 5.18]{_Demailly_Compl_Analytic_Diff_Geom_} that
for any smooth functions $Q_1,Q_2$ on a complex manifold and any $\delta>0$, one can define their regularized maximum $\max_\delta \{ Q_1,Q_2\}$. It is a smooth function on the same manifold, depending smoothly on $\delta$, and satisfying the following conditions:

\bigskip
\noindent
1. $\max \{ Q_1,Q_2\}
\leq
\max_\delta \{ Q_1, Q_2\}
\leq
\max\{ Q_1, Q_2\} + \delta$;

\bigskip
\noindent
2. $\max_\delta \{ Q_1, Q_2\} (y) = Q_j (y)$, if $Q_i (y) +2\delta \leq Q_j (y)$, $i,j=1,2$, $i\neq j$.

\bigskip
\noindent
3. If $\vartheta$ is a
Hermitian form on the manifold and $\sqrt{-1}\partial\bar{\partial} Q_i > -\vartheta$, $i=1,2$, then
$$\sqrt{-1}\partial\bar{\partial} \max_\delta \{ Q_1, Q_2\} > - \vartheta.$$


\hfill


\subsection{Blow-up -- local constructions}
\label{_blow-up-local-constructions-Subsection_}

First, we recall some standard local constructions concerning the blow-up of $\C^n$ at the origin.

Denote by $z_1,\ldots,z_n$ the complex coordinates on $\C^n$.
For $z=(z_1,\ldots,z_n)\in\C^n$ denote
\[
|z| :=\sqrt{|z_1|^2 +\ldots + |z_n|^2}.
\]
As before, let $\omega_0$ be the standard symplectic form on $\C^n\cong\R^{2n}$. It can be written as
\begin{equation}
\label{_omega-zero-ddbar_Equation_}
\omega_0 =
\sqrt{-1} \partial\bar{\partial} |z|^2/2.
\end{equation}

Let $\sigma$ be the Fubini-Study symplectic form on $\C P^{n-1}$ normalized by $\int_{\C P^1} \sigma = \pi$.

Let $\cL \subset \C^n \times \C P^{n-1}$ be the incidence relation, i.e. $\cL := \{ (z,\gol)\ |\ z\in \gol\}$. It is the complex blow-up of $\C^n$ at the origin. Let $pr_1: \cL \to \C^n$, $pr_2: \cL\to \C P^{n-1}$ be the natural projections.

For $r>0$, denote $\cL (r):= pr_1^{-1} (B^{2n} (r))$.
We will also write $\cL (0):= pr_1^{-1} (0)$.

Given $c,r>0$, define a differential 2-form $\rho (c,r)$ on $\cL$ as
\[
\rho (c,r) := c^2 pr_1^* \omega_0 + r^2 pr_2^* \sigma.
\]

The restriction of $pr_1$ to $\cL\setminus pr_1^{-1} (0)$ is a diffeomorphism $\cL\setminus pr_1^{-1} (0)\to \C^n\setminus 0$. One can easily see that
it identifies the form $\rho (c,r)$ on $\cL\setminus \cL (0)$ with the form $(pr_1)_* \rho (c,r)$ on $\C^n\setminus 0$ that can be written as
\begin{equation}
\label{_rho-c-r-ddbar_Equation_}
(pr_1)_* \rho (c,r) =\frac{\sqrt{-1}}{2} \partial\bar{\partial} \left( c^2 |z|^2 + r^2 \log |z|^2 + d\right)
\end{equation}
for any $d\in\R$.

The following two propositions and their proofs describe local constructions needed further on to define symplectic blow-up and blow-down operations. Here we use the regularized maximum for these constructions -- cf. e.g. \cite{_McD-Polt_}, \cite{_McD-Sal-3_},
where different tools were used.


\hfill


\proposition
\label{_reg-max-for-sympl-blow-up_Proposition_}

Let $r,\epsilon>0$.

Then for any sufficiently small $c >0$ there exists a K\"ahler form $\varrho^{c}$ on $\C^n\setminus 0$, depending smoothly on $c$, so that

\begin{itemize}

\item{} $\varrho^{c} = \omega_0$ on a neighborhood of $\C^n\setminus \Int B^{2n} (r+\epsilon)$,

\item{} $\varrho^{c} = (pr_1)_* \rho (c,r)$ on a neighborhood of $B^{2n} (r)\setminus 0$.

\end{itemize}


\hfill


\noindent
{\bf Proof of \ref{_reg-max-for-sympl-blow-up_Proposition_}:}

The functions $s^2 - r^2$ and $r^2 \log \left( (s/r)^2\right)$ (as functions of $s$) coincide with their first derivatives at $s=r$. The first function is convex while the second one is concave and $s^2 - r^2 > r^2 \log \left( (s/r)^2\right)$ for all $s>0$, $s\neq r$.
Therefore for any sufficiently small $c >0$ and a sufficiently small $\sigma >0$ (depending on $c$) the functions $F(s):= s^2 - r^2 - \sigma$ and $G(s):= c s^2 + r^2 \log s^2$ have exactly two points
$s_1, s_2$, $0< s_1< r< s_2< r+\epsilon$, where they coincide, so that  $F(s) < G(s)$ on $[s_1, r)$, while
$F(s) > G(s)$ on $(s_2, r+\epsilon]$.

Choose $\delta>0$ so that $F(r)< G(r) - 2\delta$
and $F(r+\epsilon) > G(r+\epsilon) + 2\delta$.

Define $H: \C^n\setminus 0\to\R$ by $H(z):= \max_\delta \{ F(|z^2|), G(|z^2|)\}$, where $\max_\delta$ is the regularized maximum mentioned in Section~\ref{_regularized maximum_Subsection_}.
It follows from \eqref{_omega-zero-ddbar_Equation_}
and \eqref{_rho-c-r-ddbar_Equation_} and the properties of the regularized maximum (see Section~\ref{_regularized maximum_Subsection_}) that
 $(\sqrt{-1}/2)\cdot \partial\bar{\partial} H$ is a K\"ahler form on $\C^n\setminus 0$ equal to
$(\sqrt{-1}/2)\cdot \partial\bar{\partial} F = \omega_0$ on a neighborhood $\C^n\setminus \Int  B^{2n} (r+\epsilon)$,
and to $(\sqrt{-1}/2)\cdot \partial\bar{\partial} G = (pr_1)_* \rho (c,r)$ on a neighborhood of the sphere $\partial B^{2n} (r)$.

Define the form $\varrho^{c}$ on $\C^n\setminus 0$ as $(pr_1)_* \rho (c,r)$ on $B^{2n} (r)\setminus 0$, $(\sqrt{-1}/2)\cdot \partial\bar{\partial} G$ on the neighborhood of $\partial B^{2n} (r)$ and $(\sqrt{-1}/2)\cdot \partial\bar{\partial} H$ outside
$\partial B^{2n} (r)$. It is a well-defined form. The parameters $\sigma,s_1,s_2,\delta$ in the construction can be chosen to depend smoothly on $c$
and thus $\varrho^{c}$ can be assumed to depend smoothly on $c$.
One easily sees that $\varrho^{c}$ satisfies the bulleted properties listed in the proposition.
This finishes the proof.
\endproof


\hfill


\proposition
\label{_reg-max-for-sympl-blow-down_Proposition_}

Let $r,\epsilon>0$.

Then for any $c>0$ there exists a K\"ahler form $\varrho_{c}$ on $\C^n$, depending smoothly on $c$, so that

\begin{itemize}

\item{} $\varrho_{c} = \omega_0$ on a neighborhood of $B^{2n} (r)$,

\item{} $\varrho_{c} = (pr_1)_* \rho (c,r)$ on a neighborhood of $\C^n\setminus \Int B^{2n} (r+\epsilon)$.

\end{itemize}


\hfill


\noindent
{\bf Proof of \ref{_reg-max-for-sympl-blow-down_Proposition_}:}

For an appropriate constant $d\in\R$ the functions $F(s):= s^2 +d$ and $G(s):= c^2 s^2 + r^2 \log \left( (s/r)^2\right)$ (as functions of $s$) are equal at some point $s_1\in (r,r+\epsilon)$, so that $F(s) > G(s)$ on $[r, s_1)$, $F(s) < G(s)$ on $(s_1, s_2]$ for some $s_2\in (s_1, r+\epsilon)$.

Having chosen such $s_1, s_2$,
pick $\delta>0$ so that $F(r)> G(r) + 2\delta$
and $F(s_2) < G(s_2) - 2\delta$.

Define $H: \C^n\setminus 0\to\R$ by $H(z):= \max_\delta \{ F(|z^2|), G(|z^2|)\}$, where $\max_\delta$ is the regularized maximum mentioned in Section~\ref{_regularized maximum_Subsection_}.
It follows from \eqref{_omega-zero-ddbar_Equation_}
and \eqref{_rho-c-r-ddbar_Equation_} and the properties of the regularized maximum (see
Section~\ref{_regularized maximum_Subsection_}) that $(\sqrt{-1}/2)\cdot \partial\bar{\partial} H$ is a K\"ahler form on $\C^n\setminus 0$ equal to $(\sqrt{-1}/2)\cdot \partial\bar{\partial} F = \omega_0$ on
$B^{2n} (r + \sigma)\setminus \Int B^{2n} (r)$ for some sufficiently small $\sigma>0$, and to
$(\sqrt{-1}/2)\cdot \partial\bar{\partial} G = (pr_1)_* \rho (c,r)$ on a neighborhood of $\partial  B^{2n} (s_2)$.

Define a form $\varrho_{c}$ on $\C^n$ as $\omega_0$ on $B^{2n} (r)$, as $(\sqrt{-1}/2)\cdot \partial\bar{\partial} H$ on $B^{2n} (s_2)\setminus B^{2n} (r)$ and as $(\sqrt{-1}/2)\cdot \partial\bar{\partial} G = (pr_1)_* \rho (c,r)$ on $\C^n \setminus B^{2n} (s_2)$.
It is a well-defined smooth form. The parameters $d,s_1,s_2,\delta$ in the construction can be chosen to depend smoothly on $c$
and thus $\varrho_{c}$ can be assumed to depend smoothly on $c$.
One easily sees that $\varrho_{c}$ satisfies the bulleted properties listed in the proposition.
This finishes the proof.
\endproof


\hfill


\lemma{\bf ($\partial\bar{\partial}$-lemma for a complex blow-up of a ball)}
\label{_ddbar-lemma-on-nbhd-of-exc-divisor_Lemma_}

Let $B\subset \C^n$ be a ball centered at the origin $0\in\C^n$. Let $\tB$ be the complex blow-up of $B$ at $0$.

Then any exact $(1,1)$-form $\upsilon$ on $\tB$ can be written as $\upsilon = \sqrt{-1} \partial \bar{\partial} Q$ for a smooth real-valued function $Q$ on $\tB$.


\hfill


\noindent
{\bf Proof of \ref{_ddbar-lemma-on-nbhd-of-exc-divisor_Lemma_}:}

Let $\pi: \tB\to B$ be the natural projection. Let $\cO_\tB$ be the structure sheaf (i.e., the sheaf of germs of holomorphic functions) on $\tB$.

Note that
$R^i \pi_* \cO_\tB = 0$ for all $i>1$
(see \cite[Prop. 2.14]{_Ueno-LNM1975_},
cf. \cite[(2), p.144]{_Hironaka-AnnMath1964-I_}).
Consequently, $H^1 (\tB, \cO_\tB) = H^1 (B, \pi_* \cO_\tB) = 0$.

By the Dolbeault theorem, this implies that any $\bar{\partial}$-closed $(0,1)$-form on $\tB$ is $\bar{\partial}$-exact. The latter fact yields, by the complex conjugation symmetry, that any $\partial$-closed $(1,0)$-form on $\tB$ is $\partial$-exact.

Since $\upsilon$ is exact, it can be written as $\upsilon = d \sigma$ for some 1-form $\sigma$ on $\tB$. Write $\sigma$ as a sum of its $(1,0)$ and $(0,1)$-parts: $\sigma = \sigma^{1,0} + \sigma^{0,1}$. Then the 1-form $\sigma^{0,1}$ is $\bar{\partial}$-closed and hence $\bar{\partial}$-exact, meaning that $\sigma^{0,1} = \bar{\partial} S$ for a smooth function $S$ on $\tB$. One readily checks that the form $\sigma - d S$ is of type $(1,0)$ and $\partial$-closed, hence $\partial$-exact,
meaning that $\sigma - d S = \partial (-\sqrt{-1} Q)$ for a smooth function $Q$ on $\tB$. Finally, a direct computation yields
$\upsilon = \bar{\partial} \partial (-\sqrt{-1} Q) = \sqrt{-1}\partial \bar{\partial} Q$. Since the form $\upsilon$ is real, the function $Q$ can be also chosen to be real-valued. This finishes the proof of the lemma.
\endproof


\hfill


\subsection{Symplectic forms on complex blow-ups of $M$}
\label{_blow-up-forms-Subsection_}

Let $I$ be a complex structure on $M$ and let $\bfx=(x_1,\ldots,x_k)\in\hM^k$.

Fix an
$I$-holomorphic embedding $h=\bigsqcup_{i=1}^k h_i:
\bigsqcup_{i=1}^k B^{2n} (r_i)\to M$, such that $h_i (0)=x_i$, $i=1,\ldots,k$, and $h$ extends to a holomorphic embedding
$h=\bigsqcup_{i=1}^k h_i:
\bigsqcup_{i=1}^k B^{2n} (r_i+\epsilon)\to M$ for some $\epsilon>0$.
Each map $h_i : B^{2n} (r_i+\epsilon)\to M$, $i=1,\ldots, k$,
lifts to the unique map
$\tildeh_i :\cL (r_i+\epsilon)\to \tM_{I,\bfx}$ such that
$h_i\circ pr_1 = \Pi\circ \tildeh_i$,
where $pr_1: \cL \to \C^n$ is the projection defined in Section~\ref{_blow-up-local-constructions-Subsection_}.
For any $c,r>0$ the map $\tildeh_i$, $i=1,\ldots,k$, identifies the form $\rho (c,r)$ on $\cL (r_i+\epsilon)$ with a form on a neighborhood of the exceptional divisor $E_i\subset \tM_{I,\bfx}$ that will be denoted by $\rho_{x_i} (c,r)$. Each $\rho_{x_i} (c,r)$ is a symplectic form compatible with $\tI$.

Given a cohomology class $\alpha\in H^2 (M;\R)$, define the cohomology class $\talpha_\bfr\in H^2 (\tM_{I,\bfx};\R)$ by
\begin{equation}
\label{_definition-of-talpha-bfr_Equation_}
\talpha_\bfr := \Pi^* \alpha - \pi\sum_{i=1}^k r_i^2 e_i.
\end{equation}


\hfill


\proposition\label{_from-arbitr-Kahler-form-to-blowup-form_Proposition_}

Let $\alpha\in H^2 (M;\R)$.
Assume that the cohomology class $\talpha_\bfr\in H^2 (\tM_{I,\bfx};\R)$
is K\"ahler with respect to $\tI$.

Then there exists a symplectic form $\ttheta$ on $\tM_{I,\bfx}$ which is compatible with $\tI$, represents $\talpha_\bfr$, and equals to
$\rho_{x_i} (1,r_i)$ on a neighborhood of $E_i$ for each $i=1,\ldots,k$.


\hfill


\remark

If one drops the requirement that $\ttheta$ is compatible with $\tI$, the existence of a symplectic form $\ttheta$, $[\ttheta] = \talpha_\bfr$, coinciding with $\rho_{x_i} (1,r_i)$ on a neighborhood of each $E_i$ can be proved much easier using Moser's method, as in \cite[Pf. of Prop. 2.1.C]{_McD-Polt_}. The need to make $\ttheta$ compatible with $\tI$ comes from the fact that otherwise the blow-down construction (see \cite[Sec. 5]{_McD-Polt_}) applied to $\ttheta$ will produce only a symplectic, but not necessarily a holomorphic, embedding of balls. Unlike Moser's method -- which can be applied, as in \cite[Pf. of Prop. 2.1.C]{_McD-Polt_}, also when $I$ is not necessarily integrable outside neighborhoods of the base points $x_1,\ldots,x_k$ -- here we do rely on the integrability of $I$ on the whole $M$, as we use methods of complex geometry.


\hfill


\noindent
{\bf Proof of \ref{_from-arbitr-Kahler-form-to-blowup-form_Proposition_}:}

The argument below uses a result
based on an idea of Demailly and
published in \cite[Thm. 4.1]{_Ornea-Verb_} and, moreover, mimics the idea of its proof.

Since the cohomology class $\talpha_\bfr\in H^2 (\tM_{I,\bfx};\R)$
is K\"ahler, we can pick a symplectic form $\teta$ that represents it and that is compatible with $\tI$.

For each $i=1,\ldots,k$ denote:
\[
U_{x_i} :=  h\left( B^{2n} (r_i + \epsilon)\right),
\]
\[
U'_{x_i} :=  h\left( B^{2n} (r_i + \epsilon/2)\right),
\]
\[
\tU_{x_i}:= \tildeh_i \left(\cL (r_i+\epsilon)\right) = \Pi^{-1} (U_{x_i}),
\]
\[
\tU'_{x_i}:= \tildeh_i \left(\cL (r_i)\right) = \Pi^{-1} (U'_{x_i}).
\]
Thus, $U'_{x_i}\subset U_{x_i}$, $\tU'_{x_i}\subset \tU_{x_i}$.

Let $\delta>0$.

Let $E := \bigcup_{i=1}^k E_i$ be the union of the exceptional divisors of $\tM_{I,\bfx}$.
Since $E$ is a complex submanifold of $(\tM_{I,\bfx},\tI)$, by a result of Demailly-Paun \cite[Lem. 2.1]{_Dem-Paun_}, there exists a smooth real-valued function $G$ on $\tM_{I,\bfx}\setminus E$ so that
\begin{equation}
\label{_d-dbar-G-C-teta_Equation_}
\sqrt{-1}\partial\bar{\partial} G > -\teta
\end{equation}
and $G$ has logarithmic poles along $E_i$, $i=1,\ldots,k$.

The set $\tU_{x_i}$ is biholomorphic to the blow-up of the ball at the origin and $\rho_{x_i} (1,r_i) - \teta$
is an exact $(1,1)$-form on $\tU_{x_i}$. Therefore, by \ref{_ddbar-lemma-on-nbhd-of-exc-divisor_Lemma_},
\begin{equation}
\label{_rho-x-i-1-r-i-minus-ttheta-equals-d-dbar-F-i_Equation_}
\rho_{x_i} (1,r_i) - \teta|_{\tU_{x_i}} = \sqrt{-1}\partial\bar{\partial} F_i
\end{equation}
for
a smooth function $F_i: \tU_{x_i}\to\R$. Adding, if necessary, a constant to $F_i$, we may assume
without loss of generality that
\begin{equation}
\label{_F-i-smaller-than-G_Equation_}
F_i < G - 2\delta\ \textrm{on an open neigborhood of}\ \tU_{x_i}\setminus\tU'_{x_i}, i=1,\ldots,k.
\end{equation}

Since
the function $G$ tends to $-\infty$ as its argument converges to $E$, for any $i=1,\ldots,k$
there exists a neighborhood $V_{x_i}\subset U'_{x_i}$ of $x_i$, whose preimage $\Pi^{-1} (V_{x_i})\subset \tU'_{x_i}$ will be denoted
by $\tV_{x_i} := \Pi^{-1} (V_{x_i})$,  so that
\begin{equation}
\label{_G-smaller-than-F-i_Equation_}
G + 2\delta < F_i\ \textrm{on}\ \tV_{x_i},\ i=1,\ldots,k.
\end{equation}

Define the function $K: \tM_{I,\bfx}\to\R$ as follows
\[
K:= \max_\delta \{ F_i, G\},\ \textrm{on}\ \tU_{x_i}, i=1,\ldots,k,
\]
\[
K:= G,\ \textrm{on}\ \tM_{I,\bfx}\setminus \bigcup_{i=1}^k \tU'_{x_i}.
\]
By \eqref{_F-i-smaller-than-G_Equation_}, \eqref{_G-smaller-than-F-i_Equation_}
and the properties of the regularized maximum (see Section~\ref{_regularized maximum_Subsection_}),
$K$ is a well-defined smooth function on $\tM_{I,\bfx}$, equal to
$F_i$ on $\tV_{x_i}$, to $G$ on $\tM_{I,\bfx}\setminus \bigcup_{i=1}^k \tU'_i$, and to
$\max_\delta \{ F_i, G\}$ on $\tU'_{x_i}\setminus \tV_{x_i}$, $i=1,\ldots,k$.

Define the smooth 2-form $\ttheta$ on $\tM_{I,\bfx}$ by
\[
\ttheta:= \teta + \sqrt{-1}\partial\bar{\partial} K.
\]
The form $\ttheta$ is of type $(1,1)$ with respect to $\tI$.
Let us verify that $\ttheta$ is positive -- i.e. symplectic. We need to verify
that
\begin{equation}
\label{_d-dbar-K-greater-than-minus-ttheta_Equation_}
\sqrt{-1}\partial\bar{\partial} K > -\teta\ \textrm{everywhere on}\ \tM_{I,\bfx}.
\end{equation}
Let us verify the latter inequality case by case
for the points of $\tV_{x_i}$, $\tU'_{x_i}\setminus \tV_{x_i}$, $i=1,\ldots,k$, and $\tM_{I,\bfx}\setminus\bigcup_{i=1}^k \tU'_{x_i}$.

Consider first the case of $\tV_{x_i}$, $i=1,\ldots,k$. By
\eqref{_rho-x-i-1-r-i-minus-ttheta-equals-d-dbar-F-i_Equation_}, and since $\rho_{x_i} (1,r_i) >0$, we have
\begin{equation}
\label{_d-dbar-F-i-greater-than-ttheta_Equation_}
\sqrt{-1}\partial\bar{\partial} F_i = - \teta + \rho_{x_i} (1,r_i) > - \teta,\ \textrm{on}\ \tU_{x_i}, i=1,\ldots,k,
\end{equation}
Hence, for the points of $\tV_{x_i}\subset \tU_{x_i}$, $i=1,\ldots,k$,
\[
\sqrt{-1}\partial\bar{\partial} K = \sqrt{-1}\partial\bar{\partial} F_i  > - \teta,
\]
meaning that the inequality \eqref{_d-dbar-K-greater-than-minus-ttheta_Equation_} is satisfied in this case.

For the points of $\tU'_{x_i}\setminus \tV_{x_i}$, $i=1,\ldots,k$, we have
\[
\sqrt{-1}\partial\bar{\partial} F_i > - \teta
\]
by \eqref{_d-dbar-F-i-greater-than-ttheta_Equation_}, and
\[
\sqrt{-1}\partial\bar{\partial} G > - \teta
\]
by \eqref{_d-dbar-G-C-teta_Equation_}.
Hence, for such points
\[
\sqrt{-1}\partial\bar{\partial} K = \sqrt{-1}\partial\bar{\partial} \max_\delta \{F_i,G\}  > - \teta,
\]
by the third property of the regularized maximum (see Section~\ref{_regularized maximum_Subsection_}).
Thus, the inequality \eqref{_d-dbar-K-greater-than-minus-ttheta_Equation_} is satisfied in this case too.

Finally, for the points of $\tM_{I,\bfx}\setminus\bigcup_{i=1}^k \tU'_{x_i}$,
\[
\sqrt{-1}\partial\bar{\partial} G > - \teta
\]
by \eqref{_d-dbar-G-C-teta_Equation_}.
Hence,
\[
\sqrt{-1}\partial\bar{\partial} K = \sqrt{-1}\partial\bar{\partial} G  > - \teta,
\]
meaning that the inequality \eqref{_d-dbar-K-greater-than-minus-ttheta_Equation_} is satisfied also in this case.

This finishes the verification of \eqref{_d-dbar-K-greater-than-minus-ttheta_Equation_}.

We have shown that $\ttheta$ is a positive $(1,1)$-form with respect to $\tI$, which means that it is compatible with $\tI$. It is cohomologous to $\teta$, meaning that it represents the cohomology class $\talpha_\bfr$. By
\eqref{_rho-x-i-1-r-i-minus-ttheta-equals-d-dbar-F-i_Equation_},
\[
\ttheta = \teta + \sqrt{-1}\partial\bar{\partial} K = \teta + \sqrt{-1}\partial\bar{\partial} F_i = \rho_{x_i} (1,r_i)\ \textrm{on}\ \tV_{x_i}, i=1,\ldots,k.
\]
This finishes the proof
of the proposition.
\endproof


\hfill


Recall that for each $\kappa>0$ we define
\[
R_\kappa: \R^{2n} \to \R^{2n},\ R_\kappa (x):= \kappa x.
\]
The map $R_\kappa$ induces a diffeomorphism of $\cL \subset \C^n \times \C P^{n-1}$ that will be denoted by $\tR_\kappa$:
\[
\tR_\kappa: \cL \to \cL.
\]
For $\kappa\in (0,1]$, denote by $h\circ R_\kappa$ the holomorphic embedding
\[
\textstyle h \circ R_\kappa :=\bigsqcup_{i=1}^k h_i \circ R_\kappa: \bigsqcup_{i=1}^k B^{2n} (r_i)\to (M,I).
\]
For each $i=1,\ldots,k$ and each $0<\kappa\leq 1$ denote
\[
U_{x_i,\kappa} := h_i \circ R_\kappa \left( B^{2n}  \left( r_i+\epsilon\right)\right) =
h_i \left( B^{2n}  \left( \kappa \left(r_i+\epsilon\right)\right)\right).
\]

The following proposition describes the exact version of the symplectic blow-down operation
that we need for our purposes.


\hfill


\proposition
\label{_sympl-blow-down-Case-I_Proposition_}

Let $\alpha\in H^2 (M;\R)$. Let $0<\kappa\leq 1$, $c_1,\ldots,c_k>0$.

Assume that the cohomology class $\talpha_\bfr\in H^2 (\tM_{I,\bfx};\R)$ is represented by
a symplectic form $\ttheta$ coincides with
$\rho_{x_i} (c_i,r_i)$ on the neighborhood $\Pi^{-1} (U_{x_i,\kappa})$
of $E_i$, $i=1,\ldots,k$.

Then there exists a symplectic form $\eta_{\kappa,\ttheta}$ on $M$ with the
following properties:

\begin{itemize}

\item{} $[\eta_{\kappa,\ttheta}] = \alpha$,

\item{} If $\ttheta$ is compatible with $\tI$, then $\eta_{\kappa,\ttheta}$ is compatible with $I$; if $I'$ is a complex structure
on $M$ coinciding with $I$ on $\bigcup_{i=1}^k U_{x_i,\kappa}$
and such that $\ttheta$ tames $\tI'$, then $\eta_{\kappa,\ttheta}$ tames $I'$.

\item{} The embedding $h\circ R_\kappa: \bigsqcup_{i=1}^k B^{2n} (r_i)\to M$ is symplectic with respect to $\eta_{\kappa,\ttheta}$ (and holomorphic with respect to $I$).

\item{} $\eta_{\kappa,\ttheta}$ depends smoothly on $\kappa$.

\end{itemize}


\hfill


\noindent
{\bf Proof of \ref{_sympl-blow-down-Case-I_Proposition_}:}

Assume $\kappa>0$ is sufficiently small so that $\ttheta = \rho_{x_i} (c_i,r_i)$ on $\Pi^{-1} ( U_{x_i,\kappa})$ for all $i=1,\ldots,k$. Then $\tR_\kappa^*\circ \tildeh_i^* \ttheta = \rho (\kappa c_i,r_i)$ on $\cL (r_i+\epsilon)$ for all $i=1,\ldots,k$.

For each $i=1,\ldots,k$ we can use \ref{_reg-max-for-sympl-blow-down_Proposition_} in order to construct a form $\varrho_{\kappa c_i}$ on $B^{2n} (r_i+\epsilon)$, depending smoothly on $\kappa$, that equals $\omega_0$ on a neighborhood of $B^{2n} (r_i)$,
and $(pr_1)_* \rho (\kappa c_i,r_i)$ on a neighborhood of $\partial B^{2n} (r_i+\epsilon)$. Then the form
$(h_i)_* \varrho_{\kappa c_i}$ is defined on $U_{x_i,\kappa}$
and coincides with $\Pi_* \ttheta = \Pi_* \left( \rho_{x_i} (\kappa c_i,r_i)\right)$ near $\partial U_{x_i,\kappa}$.

Consequently, we can define a 2-form $\eta_{\kappa,\ttheta}$ on $M$ by setting it equal to
$\Pi_* \ttheta$ outside $\bigcup_{i=1}^k U_{x_i,\kappa}$ and to $(h_i)_* \varrho_{\kappa c_i}$ on each
 $U_{x_i,\kappa}$, $i=1,\ldots,k$.
One easily checks that $\eta_{\kappa,\ttheta}$ is symplectic and that the bulleted properties in the proposition are satisfied.
\endproof


\hfill


The following proposition is a version of the symplectic blow-up construction -- cf. e.g. \cite{_McD-Polt_}, \cite{_McD-Sal-3_}.


\hfill


\proposition
\label{_sympl-blow-up-regularized-max_Proposition_}

Assume that $\eta$ is a symplectic form on $M$ such that the embedding $h$ is symplectic with respect to $\eta$ and $\eta$ is compatible with $I$ on the image of $h$.

Then for any sufficiently small $c>0$ there exists a symplectic form $\ttheta_c$ on $\tM_{I,\bfx}$ with the following properties:

\smallskip
\noindent
1. $\ttheta_c = \rho_{x_i} (c,r_i)$ on $\Pi^{-1} \left( h_i \left(B (r_i)\right)\right)$, $i=1,\ldots,k$,

\smallskip
\noindent
2. $\Pi_* \ttheta_c = \eta$ outside $\bigsqcup_{i=1}^k h \left( B^{2n} (r_i+\epsilon) \right)$,

\smallskip
\noindent
3. $\ttheta_c$ depends smoothly on $c$,

\smallskip
\noindent
4. If $\eta$ is compatible with, respectively tames, $I$ on the whole $M$, then $\ttheta_c$ is compatible with,
respectively tames, $\tI$ on the whole $\tM_{I,\bfx}$.

\smallskip
\noindent
5. $[\ttheta_c] = \Pi^* [\eta] - \pi \sum_{i=1}^k r_i^2 e_i$.


\hfill


\noindent
{\bf Proof of \ref{_sympl-blow-up-regularized-max_Proposition_}:}

Using \ref{_reg-max-for-sympl-blow-up_Proposition_} for each $B^{2n} (r_i)$, $i=1,\ldots,k$,
for any sufficiently small $c >0$
we can construct a form $\ttheta_c$ on $\tM_{I,\bfx}$ with properties 1-4.
Properties 1 and 2 imply that $\ttheta_c$ is a symplectic form satisfying also property 5.
This finishes the proof.
\endproof


\hfill


\noindent
{\bf Proof of \ref{_from-TeichT-to-pairs_Proposition_}:}

Let $\alpha := [\omega]\in H^2 (M;\R)$.

Consider a triple $(I,\bfx,h)\in \Triples_{\cC_0}$, where $\bfx = (x_1,\ldots,x_k)$. By \ref{_triples_Definition_},
the cohomology class $\talpha_\bfr\in H^2 (\tM_{I,\bfx};\R)$
is K\"ahler with respect to $I$.
Hence, by \ref{_from-arbitr-Kahler-form-to-blowup-form_Proposition_},
there exists a symplectic form $\ttheta$ on $\tM_{I,\bfx}$ which is compatible with $\tI$, represents $\talpha_\bfr$, and equals to
$\rho_{x_i} (1,r_i)$ on a neighborhood of $E_i$ for each $i=1,\ldots,k$.
Consequently, \ref{_sympl-blow-down-Case-I_Proposition_} implies that
for any sufficiently small $\kappa>0$ there exists a symplectic form $\eta_{\kappa,\ttheta}$ on $M$ with the
following properties:

\begin{itemize}

\item{} $[\eta_{\kappa,\ttheta}] = \alpha$,

\item{} $\eta_{\kappa,\ttheta}$ is compatible with $I$,

\item{} $h\circ R_\kappa: \bigsqcup_{i=1}^k B^{2n} (r_i)\to M$ is symplectic with respect to $\eta_{\kappa,\ttheta}$ (and holomorphic with respect to $I$),

\item{} $\eta_{\kappa,\ttheta}$ depends smoothly on $\kappa$.

\end{itemize}

By \ref{_triples_Definition_}, there exists $\phi\in \Diff_0 (M)$ such that $\phi^* I$ is compatible with $\omega$.
Then $\phi^* \eta_{\kappa,\ttheta}$ and $\omega$ are cohomologous symplectic forms compatible with $\phi^* I$. The linear family of 2-forms
connecting $\phi^* \eta_{\kappa,\ttheta}$ and $\omega$ is then a family of cohomologous symplectic forms.
Hence, by Moser's theorem \cite{_Moser_}, $\phi^* \eta_{\kappa,\ttheta}$ lies in the $\Diff_0 (M)$-orbit of $\omega$, and consequently so does $\eta_{\kappa,\ttheta}$.
Thus, for any sufficiently small $\kappa>0$ we have
\[
(\eta_{\kappa,\ttheta}, h\circ R_\kappa)\in \Pairs_{\cC_0}.
\]

Define the map
\[
\Phi: \TeichT_{\cC_0}\to \TeichP_{\cC_0}
\]
by
\[
\Phi \left(\left\{ I,\bfx,h\right\}\right) := \left\{\eta_{\kappa,\ttheta}, h\circ R_\kappa\right\}.
\]
We need to show that the map $\Phi$ is well-defined and has the required properties.

\bigskip
\noindent
{\sl $\Phi$ is well-defined:}

Let us first show that $\left\{ \eta_{\kappa,\ttheta}, h\circ R_\kappa\right\}$
does not depend on the choice of a sufficiently small $\kappa>0$ (for fixed $\ttheta$). Given such $0<\kappa_1 <\kappa_2$,
we have a smooth family of pairs $(\eta_{\kappa,\ttheta}, h\circ R_\kappa)$, $\kappa\in [\kappa_1,\kappa_2]$, lying in $\Pairs_{\cC_0}$. By \ref{_smooth-family-of-pairs-lies-in-the-same-Diff0-orbit_Proposition_}, this family lies in the same
$\Diff_0 (M)$-orbit on $\Pairs_{\cC_0}$. Thus, $\left\{ \eta_{\kappa,\ttheta}, h\circ R_\kappa\right\}$  does not depend on
the choice of $\kappa$.

Let us verify that $\left\{ \eta_{\kappa,\ttheta}, h\circ R_\kappa \right\}$
does not depend on the choice of $\ttheta$. Assume that $\ttheta'$ is another symplectic form
which is compatible with $\tI$, represents $\talpha_\bfr$, and equals
$\rho_{x_i} (1,r_i)$ on a neighborhood of $E_i$ for each $i=1,\ldots,k$.
Pick a sufficiently small $\kappa>0$ so that $\ttheta = \ttheta' = \rho_{x_i} (c_i,r_i)$ on $\Pi^{-1} \left( h_i \left( B^{2n} \left(\kappa \left(r_i+\epsilon\right)\right)\right)\right)$.
Then $\eta_{\kappa,\ttheta}$ and $\eta_{\kappa,\ttheta'}$ are cohomologous symplectic forms compatible with $\omega$ and
coinciding on $h_i \left( B^{2n} \left(\kappa \left(r_i+\epsilon\right)\right)\right)$. By a relative version of Moser's theorem (see \ref{_relative-version-of-Moser-method_Proposition_}, part I), there exists $\phi\in\Diff_0 (M)$ that fixes $h_i \left( B^{2n} \left(\kappa \left(r_i+\epsilon\right)\right)\right)$ pointwise (thus mapping $h\circ R_\kappa$ to itself) and maps $\eta_{\kappa,\ttheta}$ to
$\eta_{\kappa,\ttheta'}$.
Thus, $\left\{ \eta_{\kappa,\ttheta}, h\circ R_\kappa\right\}$ does not depend
on the choice of $\ttheta$.

Finally, a direct check shows that $\Phi \left(\left\{ I,\bfx,h\right\}\right)$ does not depend on a representative
$(I,\bfx,h)$ of $\left\{ I,\bfx,h\right\}\in \TeichT_{\cC_0}$: acting by $\phi\in \Diff_0 (M)$ on
$(I,\bfx,h)$ and applying the previous construction to the resulting element of $\Triples_{\cC_0}$ we get a pair in $\Pairs_{\cC_0}$ obtained from $\left(\eta_{\kappa,\ttheta}, h\circ R_\kappa\right)$
by the action of the same $\phi$. (Here we identify the complex blow-ups of $(M,I)$ at $x_1,\ldots,x_k$
and of $(M,\phi_* I)$ at $\phi(x_1),\ldots,\phi(x_k)$ by the diffeomorphism induced by $\phi$.)

Thus, we have shown that $\Phi$ is well-defined.

\bigskip
\noindent
{\sl $\Phi$ is continuous:}

By \ref{_TeichP-homeom-to-symplectic-Teich_Proposition_}, the space $\TeichP_{\cC_0}$ is discrete. Thus, showing that
$\Phi: \TeichT_{\cC_0}\to \TeichP_{\cC_0}$ is continuous is equivalent to showing that the preimage of a point in $\TeichP_{\cC_0}$
under $\Phi$ is open in $\TeichT_{\cC_0}$, or, in other words, that for any $(I',\bfx',h')\in \Triples_{\cC_0}$ sufficiently close to
$(I,\bfx,h)$ its $\Diff_0 (M)$-orbit $\{ (I',\bfx',h')\}\in\TeichT_{\cC_0}$ is mapped by $\Phi$ to the same point in $\TeichP_{\cC_0}$ as $\{ (I,\bfx,h)\}$.

Given $(I',\bfx',h')\in \Triples_{\cC_0}$, one easily constructs $\phi\in\Diff_0 (M)$ so that

\begin{itemize}

\item{} $\phi$ sends $\bfx'$ into $\bfx$,

\item{} $\phi\circ h'$ coincides with $h$ near each $x_i$, $i=1,\ldots,k$.

\end{itemize}

If $h$ and $h'$, and $\bfx$ and $\bfx'$, are sufficiently $C^\infty$-close, then $\phi$ can be chosen to be close to the identity and therefore it will send
$(I',\bfx',h')$ into a triple close to $(I,\bfx,h)$. Thus, without loss of generality we may replace
$(I',\bfx',h')$ by a triple $(I', \bfx, h')$ that is close to $(I,\bfx,h)$ and satisfies $h\circ R_\kappa = h'\circ R_\kappa$
for any sufficiently small $\kappa>0$. In particular, for any such $\kappa$ the complex structures $I$ and $I'$ coincide on each $h_i \left( B^{2n} \left(\kappa \left(r_i+\epsilon\right)\right)\right)$, $i=1,\ldots,k$.

Arguing as above, we can construct a symplectic form $\ttheta$ on $\tM_{I,\bfx}= \tM_{I',\bfx}$ which is compatible with $\tI'$, represents $\talpha_\bfr$, and for any sufficiently small $\kappa>0$ equals to
$\rho_{x_i} (1,r_i)$ on the neighborhood $\Pi^{-1} \left(h_i \left( B^{2n} \left(\kappa \left(r_i+\epsilon\right)\right)\right)\right)$ of $E_i$ for each $i=1,\ldots,k$.
Consequently,
for any sufficiently small $\kappa>0$ we have a well-defined symplectic form $\eta_{\kappa,\ttheta'}$ on $M$.
This form is compatible with $I'$, represents the class $\alpha$ and coincides with $\eta_{\kappa,\ttheta}$
on each $h_i \left( B^{2n} \left(\kappa \left(r_i+\epsilon\right)\right)\right)$, $i=1,\ldots,k$.

Assume that $I'$ is sufficiently close to $I$ so that it is tamed by the form $\eta_{\kappa,\ttheta}$ (which is compatible with $I$).
Then $\eta_{\kappa,\ttheta}$ and $\eta_{\kappa,\ttheta'}$ are cohomologous symplectic forms on $M$ that both tame $I'$ and coincide
on each $h_i \left( B^{2n} \left(\kappa \left(r_i+\epsilon\right)\right)\right)$, $i=1,\ldots,k$. A relative version of
Moser's theorem (see \ref{_relative-version-of-Moser-method_Proposition_}, part I) implies then that there exists $\varphi\in\Diff_0 (M)$ that maps $\eta_{\kappa,\ttheta'}$
into $\eta_{\kappa,\ttheta}$ and
fixes pointwise
each $h_i \left( B^{2n} \left(\kappa \left(r_i+\epsilon\right)\right)\right)$, $i=1,\ldots,k$. Hence, $\varphi$
maps $\left(\eta_{\kappa,\ttheta'}, h'\circ R_\kappa\right)\in \Pairs_{\cC_0}$ into
$\left(\eta_{\kappa,\ttheta}, h\circ R_\kappa\right)$.

Thus, for any $(I',\bfx',h')\in \Triples_{\cC_0}$ sufficiently close to
$(I,\bfx,h)$, its $\Diff_0 (M)$-orbit
\break
$\{ I',\bfx',h'\}\in\TeichT_{\cC_0}$ is mapped by $\Phi$ to the same point in $\TeichP_{\cC_0}$ as $\{ I,\bfx,h\}\in\TeichT_{\cC_0}$. This finishes the proof of the claim that $\Phi$ is continuous.

\bigskip
\noindent
{\sl $\Phi$ is surjective:}

The proof uses the well-known fact (see e.g. \cite{_McD-Polt_}) that the symplectic blow-down is an ``inverse" of the symplectic blow-up, up to an isotopy.

Namely, let
$(\eta, f)\in \Pairs_{\cC_0}$,
where $f = \bigsqcup_{i=1}^k f_i: \bigsqcup_{i=1}^k B^{2n} (r_i)\to (M,\eta)$ is a symplectic embedding holomorphic with respect to
a complex structure $I$ on $M$ compatible with $\eta$. Let us show that $\{ \eta, f\}\in \TeichP_{\cC_0}$ lies in the image of $\Phi$.

Set
\[
x_i := f_i (0), i=1,\ldots,k,
\]
\[
\bfx := (x_1,\ldots,x_k).
\]
By \ref{_sympl-blow-up-regularized-max_Proposition_},
the cohomology class $\talpha_\bfr\in H^2 (\tM_{I,\bfx};\R)$ is K\"ahler with respect to $\tI$.
Thus, by \ref{_triples_Definition_}, $(I,\bfx,f)\in \Triples_{\cC_0}$.
Moreover,
by \ref{_sympl-blow-up-regularized-max_Proposition_}, the class $\talpha_\bfr$
can be represented, for a sufficiently small $c>0$, by a symplectic form $\ttheta_c$ on $\tM_{I,\bfx}$ with the following properties:

\begin{itemize}

\item{} $\ttheta_c = \rho_{x_i} (c,r_i)$ on $\Pi^{-1} \left( f_i \left(B (r_i)\right)\right)$, $i=1,\ldots,k$,

\item{} $[\ttheta_c] = \Pi^* [\eta] - \pi \sum_{i=1}^k r_i^2 e_i$,

\item{} $\ttheta_c$ is compatible with $\tI$.

\end{itemize}

Consequently, \ref{_sympl-blow-down-Case-I_Proposition_} implies that
for $\kappa=1$ there exists a symplectic form $\eta_1 (\ttheta_c)$ on $M$ with the
following properties:

\begin{itemize}

\item{} $[\eta_1 (\ttheta_c)] = \alpha$,

\item{} $\eta_1 (\ttheta_c)$ is compatible with $I$,

\item{} The $f\circ R_1 = f: \bigsqcup_{i=1}^k B^{2n} (r_i)\to M$ is symplectic with respect to $\eta_1 (\ttheta_c)$ (and holomorphic with respect to $I$).

\end{itemize}

This yields that

\smallskip
\noindent
1. $\talpha_\bfr\in H^2 (\tM_{I,\bfx};\R)$ is K\"ahler with respect to $\tI$.
Thus, by \ref{_triples_Definition_}, $(I,\bfx,f)\in \Triples_{\cC_0}$.

\smallskip
\noindent
2. $\eta_1 (\ttheta_c)$ coincides with $\eta$ on each $f_i \left( B^{2n} \left(\kappa r_i\right)\right)$, $i=1,\ldots,k$.

\smallskip
A relative version of
Moser's theorem (see \ref{_relative-version-of-Moser-method_Proposition_}, part I) implies then that there exists $\psi\in\Diff_0 (M)$ that maps $\eta_1 (\ttheta_c)$
into $\eta$ and
fixes pointwise each
each $f_i \left( B^{2n} \left(\kappa r_i\right)\right)$, $i=1,\ldots,k$. Hence,
\[
\psi \left(\eta_1 \left(\ttheta_c\right), f\right) = (\eta, f).
\]
Thus,
\[
\Phi \left(\left\{ I,\bfx,f\right\}\right) =
\left\{ \eta_1 \left(\ttheta_c\right), f\right\} = \left\{ \eta, f\right\}\in \TeichP_{\cC_0}.
\]
This shows that $\{ \eta, f\}\in \TeichP_{\cC_0}$ lies in the image of $\Phi$. Hence, $\Phi$ is surjective.

We have shown that $\Phi$ is well-defined, continuous and surjective. By construction, it clearly satisfies the
properties 1 and 2 listed in the statement of
\ref{_from-TeichT-to-pairs_Proposition_}.
This finishes the proof of \ref{_from-TeichT-to-pairs_Proposition_}.
\endproof


\hfill


The tools developed above allow us also to prove the following existence result for tame embeddings of balls, which is an
analogue of part (I) of \ref{_existence-connectedness-Kahler-type-arb-mfds_Theorem_}.


\hfill


\proposition
\label{_tame-embs-existence-via-forms-on-blow-up_Proposition_}

The set $\bigsqcup_{i=1}^k B^{2n}(r_i)$ admits a tame embedding into $(M,\omega)$ if and only if
for some $I\in \compK (M,\omega)$, tamed by $\omega$, and some $\bfx\in \hM^k$, the class $\Pi^* [\omega] - \pi \sum_{i=1}^k r_i^2 e_i\in H^2 (\tM_{I,\bfx};\R)$, can be represented by a form taming $\tI$ and equal to $\rho_{x_i} (c_i,r_i)$, for some $c_i>0$, near each $E_i$,
$i=1,\ldots,k$.

More precisely, given $I\in \compK (M,\omega)$, tamed by $\omega$, the following conditions are equivalent:

\begin{itemize}

\item{} There exists $\bfx = (x_1,\ldots,x_k)\in \hM^k$ so that the cohomology class $\Pi^* [\omega] - \pi \sum_{i=1}^k r_i^2 e_i\in H^2 (\tM_{I,\bfx};\R)$ can be represented by a form taming $\tI$ and equal to $\rho_{x_i} (c_i,r_i)$, for some $c_i>0$, near each $E_i$,
$i=1,\ldots,k$.

\item{} There exist a tame embedding $f: \bigsqcup_{i=1}^k B^{2n}(r_i)\to (M,\omega)$
and $\phi\in\Diff_0 (M)$, so that $f$ is holomorphic with respect to the complex structure $\phi^* I$ on $M$ tamed by $\omega$.

\end{itemize}


\hfill


\noindent
{\bf Proof of \ref{_tame-embs-existence-via-forms-on-blow-up_Proposition_}:}

Assume that for some $I\in \compK (M,\omega)$, tamed by $\omega$, and some $\bfx\in \hM^k$ the cohomology class
$\Pi^* [\omega] - \pi \sum_{i=1}^k r_i^2 e_i\in H^2 (\tM_{I,\bfx};\R)$ can be represented by a form taming $\tI$ and equal to $\rho_{x_i} (c_i,r_i)$, for some $c_i>0$, near each $E_i = \Pi^{-1} (x_i)$,
$i=1,\ldots,k$. Then one can deduce from \ref{_sympl-blow-down-Case-I_Proposition_}
that there exist a symplectic form $\eta$ on $M$, cohomologous to $\omega$ and taming $I$, and an embedding $g: \bigsqcup_{i=1}^k B^{2n} (r_i)\to M$ symplectic with respect to $\eta$ and holomorphic with respect to $I$.
The linear family of closed 2-forms connecting $\omega$ and $\eta$ is formed by cohomologous forms each of which tames $I$
and is therefore symplectic. Therefore, by Moser's theorem \cite{_Moser_}, $\omega = \phi^* \eta$ for some $\phi\in\Diff_0 (M)$.
The embedding $f:= \phi^{-1}\circ g : \bigsqcup_{i=1}^k B^{2n} (r_i)\to (M,\omega)$ is symplectic with respect to $\omega$
and holomorphic with respect to the complex structure $\phi^* I$ on $M$ tamed by $\omega$, meaning that $f$ is tame.
This proves one of the implications in the required equivalence.

Let us now prove the opposite implication.

Assume that $f = \bigsqcup_{i=1}^k f_i: \bigsqcup_{i=1}^k B^{2n}(r_i)\to (M,\omega)$ is a symplectic embedding
holomorphic with respect to the complex structure $J:=\phi^* I$, $\phi\in\Diff_0 (M)$, tamed by $\omega$.
Set $y_i:= f_i (0)$, $i=1,\ldots,k$, and $\bfy:= (y_1,\ldots,y_k)$.

By \ref{_sympl-blow-up-regularized-max_Proposition_},
there exists a symplectic form $\ttheta$ on $\tM_{J,\bfy}$ with the following properties:

\smallskip
\noindent
1. $\ttheta = \rho_{y_i} (c,r_i)$, for a sufficiently small $c>0$,
on a neighborhood of the exceptional divisor $\Pi^{-1} (y_i)\subset \tM_{J,\bfy}$, $i=1,\ldots,k$,

\smallskip
\noindent
2. $\ttheta$ tames $\tJ$,

\smallskip
\noindent
3. $[\ttheta] = \Pi^* [\omega] - \pi \sum_{i=1}^k r_i^2 e_i\in H^2 (\tM_{J,\bfy};\R)$.

\smallskip
Set $x_i: = \phi (y_i)$, $i=1,\ldots,k$, and $\bfx := (x_1,\ldots,x_k)$.
The diffeomorphism $\phi: M\to M$ induces a diffeomorphism $\tM_{J,\bfy}\to \tM_{I,\bfx}$
identifying the form $\ttheta$ on $\tM_{J,\bfy}$ with a form on $\tM_{I,\bfx}$
taming $\tI$, equal to $\rho_{x_i} (c,r_i)$ near each $E_i = \Pi^{-1} (x_i) \subset \tM_{I,\bfx}$,
$i=1,\ldots,k$, and lying in the cohomology class $\Pi^* [\omega] - \pi \sum_{i=1}^k r_i^2 e_i\in H^2 (\tM_{I,\bfx};\R)$.

This finishes the proof of the opposite implication and of the proposition.
\endproof


\hfill


\subsection{Connectedness of $\TeichT_{\cC_0}$}
\label{_connectedness-TeichT-C-zero_Subsection_}


Let $I$ be a K\"ahler-type complex structure on $M$. Equip $\hM^k$ and $\hM^k\times M$ with the product complex structures induced by $I$.

Denote by $\Xi\subset \hM^k\times M$
the following incidence variety:
\[
\Xi:= \left\{ \ \left( \left(x_1,\ldots,x_k\right),x \right) \in \hM^k\times M\ \big|\ x =x_i\ \textrm{for some}\ i=1,\ldots,k\ \right\}.
\]
It is a smooth complex submanifold of $\hM^k\times M$. Denote by $\cM$ the complex blow-up of $\hM^k\times M$ along $\Xi$. It is a complex manifold that projects holomorphically to $\hM^k\times M$ and consequently to $\hM^k$.
Denote the latter projection by $\pi: \cM\to\hM^k$. This projection defines a complex analytic family of
complex manifolds equipped with K\"ahler-type complex structures.

The fiber of $\pi$ over each $\bfx=(x_1,\ldots,x_k)\in \hM^k$
can be canonically identified with the complex blow-up $\tM_{I,\bfx}$ of $(M,I)$ at $x_1,\ldots,x_k$. The holonomy of the Gauss-Manin connection on $\pi$ is trivial and the homology/cohomology groups of all the fibers of $\pi$ are then identified with the homology/cohomology groups of the complex blow-up of $M$ at $k$ distinct ordered points (the identification does not depend on the choice of the points, since the homology/cohomology groups of the complex blow-ups for different choices can be canonically identified). The identifications preserve the Hodge decompositions and, in particular, cohomology classes of type $(1,1)$ of one fiber are identified with the $(1,1)$-classes of another fiber.


\hfill


\proposition
\label{_kahler-locus_Proposition_}

Assume that $\beta$ is a $(1,1)$-cohomology class of the complex blow-up of $(M,I)$ at $k$ points.

Denote by $\cB_\beta\subset \hM^k$ the set of all $\bfx\in \hM^k$
such that $\beta\in H_\tI^{1,1}(\tM_{I,\bfx};\R)$ is K\"ahler.

Then $\cB_\beta$ is connected.


\hfill


\noindent
{\bf Proof of \ref{_kahler-locus_Proposition_}:}

If $\cB_\beta$ is empty, the claim is trivial. Assume it is non-empty. Since $\beta\in H_\tI^{1,1}(\tM_{I,\bfx};\R)$ for all $\bfx\in \hM^k$, the Kodaira-Spencer stability theorem implies (see \ref{_compK-path-connected-Kahler-coh-classes_Proposition_}) that
the set $\cB_\beta$ is open.

A theorem of Demailly-Paun \cite[Thm. 0.9]{_Dem-Paun_} about the K\"ahler cones of the fibers of a complex analytic deformation family implies in our situation that there exists a countable union $\bigcup_\nu S_\nu$ of proper
analytic subsets $S_\nu \subsetneqq\hM^k$ such that the K\"ahler cones
of the fibers $\pi^{-1} (\bfx)$ coincide for all $\bfx\in \hM^k\setminus \bigcup_\nu S_\nu$ under our identifications of the cohomologies of the fibers.

Since the set $\cB_\beta\subset \hM^k$ is open, it intersects $\hM^k\setminus \bigcup_\nu S_\nu$.
Therefore $\beta\in H_\tI^{1,1}(\tM_{I,\bfx};\R)$ is K\"ahler for all $\bfx\in \hM^k\setminus \bigcup_\nu S_\nu$,
meaning that  $\hM^k\setminus \bigcup_\nu S_\nu\subset \cB_\beta$. Being the complement of a countable union of proper analytic subsets of a connected complex manifold, the set $\hM^k\setminus \bigcup_\nu S_\nu$ is connected and dense in $\hM^k$. Hence,
the set $\cB_\beta$ is connected too.
\endproof


\hfill


\corollary
\label{_fibers-of-projection-cK-cC-zero-are-connected_Corollary_}

The fibers of the projection $\pr: \cK (\bfr)\to \compK (M)$ are connected.


\hfill


\noindent
{\bf Proof of \ref{_fibers-of-projection-cK-cC-zero-are-connected_Corollary_}:}

For each $I\in \Diff_0 (M,\omega)\cdot \cmpt (M,\omega)$ and each $\bfx\in\hM^k$ the cohomology class
$\talpha_\bfr = \Pi^* [\omega] - \pi \sum_{i=1}^k r_i^2 e_i\in H^2 (\tM_{I,\bfx};\R)$
is of type $(1,1)$. The image of $\pr$ lies in $\Diff_0 (M,\omega)\cdot\cmpt (M,\omega)$ and
the fiber of $\pr$ over $I\in\Image\pr$ is the set of all $\bfx\in \hM^k$
such that the $(1,1)$-class $\talpha_\bfr$ is K\"ahler.
By \ref{_kahler-locus_Proposition_}, applied to $\beta = \talpha_\bfr$, this fiber is connected.
\endproof


\hfill


\proposition
\label{_the-space-of-holom-embs-of-ball-is-connected_Proposition_}

Assume $I$ is a complex structure on $M$, $x\in M$, $r>0$.

Then the space of holomorphic embeddings $h: B^{2n} (r)\to (M,I)$, such that $h(0)=x$, is path-connected.


\hfill


\noindent
{\bf Proof of \ref{_the-space-of-holom-embs-of-ball-is-connected_Proposition_}:}

Assume $h: B^{2n} (r)\to (M,I)$, $h(0)=x$, is a holomorphic embedding. Let $g: B^{2n} (r)\to (M,I)$, $g(0)=x$, be another holomorphic
embedding. We want to show that $g$ can be connected to $h$ by a smooth path of holomorphic embeddings $B^{2n} (r)\to (M,I)$ mapping $0$ to $x$.

We can assume without loss of generality that $\Image g\subset \Image h$. Indeed,
for a sufficiently small $\kappa>0$ we have $\Image (g\circ R_\kappa)\subset \Image h$ and $g=g\circ R_1$ can be
connected to $g\circ R_\kappa$ by the family $\{ g\circ R_t\}_{t\in [\kappa,1]}$.

Assuming $\Image g\subset \Image h$, in order to connect $g$ to $h$ by a smooth path of holomorphic embeddings $B^{2n} (r)\to (M,I)$ mapping $0$ to $x$ it suffices to show that the holomorphic embedding
$f:= h^{-1} \circ g: B^{2n} (r)\to B^{2n} (r)$, $h^{-1} \circ g (0) = 0$, can be connected to the identity by a smooth path of holomorphic embeddings $B^{2n} (r)\to B^{2n} (r)$ mapping $0$ to $0$.
The latter result is standard and is proved by a version of the Alexander trick -- see \ref{_Alexander-trick_Proposition_}, part II.

This finishes the proof of the proposition.
\endproof


\hfill


\noindent
{\bf Proof of \ref{_connectedness-of-TeichT-cC-0_Proposition_}:}

The space $\TeichT_{\cC_0} = \Triples_{\cC_0}/\Diff_0 (M)$ is connected if only if the space
$\Triples_{\cC_0}$ is connected.

The projection $(I,\bfx,h)\mapsto (I,\bfx)$ defines a surjective continuous map $\Triples_{\cC_0}\to \cK_{\cC_0} (\bfr)$.
In view of \ref{_the-space-of-holom-embs-of-ball-is-connected_Proposition_}, the fibers of this map are connected.
Hence, $\Triples_{\cC_0}$ is connected if and only $\cK_{\cC_0} (\bfr)$
is connected.
By \ref{_fibers-of-projection-cK-cC-zero-are-connected_Corollary_}, the fibers of the projection $\pr: \cK_{\cC_0} (\bfr)\to \cC_0$ are connected. Thus,
$\cK_{\cC_0} (\bfr)$ is connected if and only if $\pr (\cK_{\cC_0} (\bfr))$
is connected.
The set $\pr (\cK_{\cC_0} (\bfr))$ is connected if and only if the set $\ttK_{\cC_0} (\bfr) = \pr (\cK_{\cC_0} (\bfr)/\Diff_0 (M)$ is connected.
Summing up, we see that $\TeichT_{\cC_0}$ is connected if and only if $\ttK_{\cC_0} (\bfr)$ is connected.
This finishes the proof of the proposition.
\endproof


\hfill


\section{The case of complex projective spaces and their products}
\label{_complex-strs-on-proj-spaces-their-products-and-blowups_Section_}

As above, we denote by $I_{st}$ the standard complex structure on $\C P^n$ and by $\omega_{FS}$, or by $\omega_{FS,n}$, the standard Fubini-Study symplectic form on $\C P^n$ normalized so that $\int_{\C P^1} \omega_{FS} = \pi$. Let $h\in H^2 (\C P^n;\Z)$ be the generator of $H^2 (\C P^n;\Z)$ so that $[\omega_{FS}] = \pi h$.

Let us prove \ref{_existence-products-of-proj-spaces_Theorem_}. For convenience, we restate it here.


\hfill


\theorem {\bf (= \ref{_existence-products-of-proj-spaces_Theorem_})}
\label{_existence-products-of-proj-spaces-COPY_Theorem_}

Consider the manifold $M:=\C P^{n_1}\times\ldots\times\C P^{n_m}$, $n_1,\ldots, n_m >0$, $n_1 +\ldots +n_m =: n$, endowed with the symplectic form $\omega_\bfc :=
c_1\omega_{FS,n_1}\oplus\ldots\oplus c_m\omega_{FS,n_m}$, $\bfc := (c_1,\ldots,c_m)\in (\R_{>0})^m$. Let $I$ be the complex
structure on $M$ which is the product of the standard complex structures on the factors.

Let $l_1,\ldots,l_m\in\Z_{>0}$ so that $[l_1:\ldots :l_m] = [c_1:\ldots: c_m]\in \R P^m$.
Let
\[
k := \frac{(n_1+\ldots n_m)!}{n_1!\cdot\ldots\cdot n_m!} l_1^{n_1}\cdot\ldots\cdot l_m^{n_m}.
\]
Let $\Sigma\subset (M,I)$ be a proper (possibly empty) complex submanifold. If $M = \C P^{n_1}\times \C P^1$, we also allow
$\Sigma = (\C P^{n_1 -1}\times \C P^1) \cup (\C P^{n_1}\times \textrm{pt})$.

Then K\"ahler-type embeddings of disjoint unions of $k$ equal balls into $(M\setminus\Sigma,\omega_\bfc)$
are unobstructed.

More precisely, if $k \Vol (B^{2n} (r),\omega_0) < \Vol (M,\omega)$, then for each proper (possibly empty) complex submanifold $\Sigma\subset (M,I)$, there exists an $[I]$-K\"ahler-type embedding of $\bigsqcup_{i=1}^k B^{2n} (r)$
into $(M,\omega)$ which is holomorphic with respect to a complex structure on $M$ that is compatible with $\omega$
and isotopic to $I$ by an isotopy preserving $\Sigma$ (as a set).


\hfill


\noindent
{\bf Proof of \ref{_existence-products-of-proj-spaces-COPY_Theorem_} (= \ref{_existence-products-of-proj-spaces_Theorem_}):}

The proof mimics the proof of \cite[Thm. 1.5.A]{_McD-Polt_}.

Namely, as it is shown in \cite[Pf. of Thm. 1.5.A]{_McD-Polt_}, there exist $x_1,\ldots,x_k\in M$, $\bfx := (x_1,\ldots,x_k)\in\hM^k$, and a holomorphic map $F: \tM_{I,\bfx}\to (\C P^{n-1}, I_{st})$ so that

\begin{itemize}

\item{} $F$ induces a biholomorphism on each exceptional divisor,

\item{} $\Pi^* [\omega] - F^* [\omega_{FS,n-1}]\in H^2 (\tM_{I,\bfx};\R)$ is a linear combination of the cohomology classes $e_1,\ldots,e_k$ Poincar\'e-dual to the fundamental homology classes of the exceptional divisors of $\tM_{I,\bfx}$. Here $\Pi: \tM_{I,\bfx}\to M$ is the blow-up projection, as before.

\end{itemize}

Let $c>0$. A direct check using the above-mentioned properties of $F$
shows that the form $\tomega_c := (1+c)^{-1} (\Pi^* \omega + c F^* \omega_{FS,n-1})$
is a symplectic form on $\tM_{I,\bfx}$ compatible with $\tI$ and representing the cohomology class
$\Pi^* [\omega] - \pi \sum_{i=1}^k (1+c)^{-1} c e_i\in H^2 (\tM_{I,\bfx};\R)$.
Therefore, by part (I) of \ref{_existence-connectedness-Kahler-type-arb-mfds_Theorem_} (=\ref{_existence-connectedness-Kahler-type-arb-mfds-COPY_Theorem_}), for $r_c:=\sqrt{(1+c)^{-1} c}$,
there exists a symplectic embedding $\bigcup_{i=1}^k B(r_c)\to (M\setminus\Sigma,\omega)$
and an isotopy $\{ \phi_t\}_{0\leq t\leq 1}\subset \Diff_0 (M)$, $\phi_0 = Id$, $\phi_t (\Sigma) = \Sigma$ for all $t\in [0,1]$, so that $f$ is holomorphic with respect to the complex structure $\phi_1^* I$ which is compatible with $\omega$.
As it is shown in \cite[Pf. of Thm. 1.5.A]{_McD-Polt_}, as $c\to +\infty$, the symplectic volume
of $\bigcup_{i=1}^k B(r_c)$ tends from below to the symplectic volume of $(M,\omega)$.
This readily yields the claim of the theorem.
\endproof


\hfill


Let us know consider the case $M=\C P^n$.

The following result is classical -- see \ref{_CPn-complex-structures_Remark_} for the references.


\hfill


\theorem
\label{_complex-str-on-CPn-is-unique_Theorem_}

\noindent
(I) Any K\"ahler-type complex structure on $\C P^n$ can be identified with $I_{st}$ by a diffeomorphism lying in $\Diff_H (\C P^n)$.

\bigskip
\noindent
(II) Any two K\"ahler-type complex structures on $\C P^n$ lying in the same component of $\compK (\C P^n)$ can be
mapped into each other by an element of $\Diff_0 (\C P^n)$.
\endproof


\hfill


\remark
\label{_CPn-complex-structures_Remark_}

\noindent
1. Part (I) of \ref{_complex-str-on-CPn-is-unique_Theorem_} is due to Hirzebruch and Kodaira \cite{_Hirzebruch-Kodaira-JMPA1957_} for odd $n$ and to Yau \cite{_Yau-PNAS1977_} for even $n$ -- see \cite{_Tosatti-CPn_} for an exposition  and a proof. (The result in \cite{_Tosatti-CPn_} does not state explicitly that the diffeomorphism in part (I) of the theorem acts trivially on homology but this can be easily deduced from the proof there.)

\bigskip
\noindent
2.
It is not known whether each complex structure on $\C P^n$, $n\geq 3$, is of K\"ahler type. For $n=1$ this is obviously true.
For $n=2$ this is true by a deep result of Yau \cite{_Yau-PNAS1977_}.

\bigskip
\noindent
3. Part (II) of the theorem implies that $\Diff_H (\C P^n)=\Diff_0 (\C P^n)$ if and only if $\compK (\C P^n)$
is connected.
It is known that $\Diff_H (\C P^n)\neq \Diff_0 (\C P^n)$ for $n=3$ \cite[Thm. 1.2]{_Kreck-Su-MCG2020_} and $n=4$ \cite[Rem. II.11]{_Brumfiel_}. Thus, $\compK (\C P^n)$, $n=3,4$, has more than one connected component. We do not know whether the same facts hold for $n=2$ and $n>4$.

\bigskip
\noindent
4.
Part (II) of \ref{_complex-str-on-CPn-is-unique_Theorem_} follows from a theorem of Fr\"ohlicher and Nijenhuis \cite{_Frolicher-Nijenhuis-PNAS_}, cf. \cite[Thm. 6.3]{_Kod-Spen-AnnMath-I-II-1958_}.


\hfill


\vfil\eject

\corollary
\label{_CPn-Teich_Corollary_}

Each connected component of the space $\Teich (\C P^n)$ is a point. The $\Diff_H (\C P^n)$-action
on $\Teich (\C P^n)$ is transitive. \endproof


\hfill


\corollary
\label{_CPn-SymplH-action-on-Teich_Corollary_}

The group $\SympH (\C P^n,\omega_{FS})\subset \Diff_H (\C P^n)$ acts transitively on the set of connected components of
$\Teich (\C P^n)$ compatible with $\omega_{FS}$.


\hfill


\noindent
{\bf Proof of \ref{_CPn-SymplH-action-on-Teich_Corollary_}:}

Assume $J$ is a complex structure on $\C P^n$ compatible with $\omega_{FS}$. Then, by part (I) of
\ref{_complex-str-on-CPn-is-unique_Theorem_}, there exists $\phi\in\Diff_H (\C P^n)$ such that $\phi^* J = I_{st}$.
Then $\phi^* \omega_{FS}$ is a symplectic form compatible with $I_{st}$ and cohomologous to $\omega_{FS}$.
Therefore the straight path in the space of 2-forms connecting these forms is formed by cohomologous symplectic forms, and consequently, by Moser's theorem \cite{_Moser_}, there exists $\psi\in\Diff_0 (\C P^n)$ such that $\psi^* \phi^* \omega_{FS} = \omega_{FS}$. Then $\phi\psi\in\SympH (\C P^n,\omega_{FS})$ and $(\phi\psi)^* J = \psi^* I_{st}$ lies in the $\Diff_0 (\C P^n)$-orbit
of $I_{st}$. This means that the connected component of $\Teich (\C P^n)$ containing $[J]$ is mapped under the action of $\phi\psi\in\SympH (\C P^n,\omega_{FS})$ into the connected component of $\Teich (\C P^n)$ containing $[I_{st}]$. This yields the proposition. \endproof


\hfill


We are ready to prove \ref{_existence-connectedness-CPn-l-power-n-balls_Theorem_} -- let us first recall it here.


\hfill


\theorem {\bf (=\ref{_existence-connectedness-CPn-l-power-n-balls_Theorem_})}
\label{_existence-connectedness-CPn-l-power-n-balls-COPY_Theorem_}

\noindent
A. For each $l\in \Z_{>0}$,
K\"ahler-type embeddings of disjoint unions of $l^n$ equal balls into $(\C P^n,\omega_{FS})$ are unobstructed.

More precisely,
if $l^n \Vol (B^{2n} (r),\omega_0) < \Vol (\C P^n,\omega_{FS})$, then for any complex structure $J$ on $\C P^n$ compatible with $\omega_{FS}$ (and, in particular, for the standard complex structure $I_{st}$) and for each proper (possibly empty) complex submanifold $\Sigma\subset (\C P^n,J)$, there exists a $[J]$-K\"ahler-type embedding of $\bigsqcup_{i=1}^{l^n} B^{2n} (r)$
into $(\C P^n,\omega_{FS})$ which is holomorphic with respect to a complex structure on $\C P^n$ that is compatible with $\omega_{FS}$
and isotopic to $J$ by an isotopy preserving $\Sigma$ (as a set).

\bigskip
\noindent
B. The group $\SympH (\C P^n,\omega_{FS})$ acts transitively on the set of connected components of $\compK (\C P^n)$ compatible with
$\omega_{FS}$.

\bigskip
\noindent
C.
For any $k\in\Z_{>0}$ and $r_1,\ldots,r_k>0$, any two K\"ahler-type embeddings $\bigsqcup_{i=1}^k B^{2n} (r_i)\to (\C P^n,\omega_{FS})$ (if they exist!) lie in the same orbit of the  $\SympH (\C P^n,\omega_{FS})$-action.

They lie in the same orbit of the $\Symp (\C P^n,\omega_{FS})\cap \Diff_0 (\C P^n)$-action
if and only if they favor a common connected component of $\compK (\C P^n)$. In the latter case there exists $[I]\in\cmptteich (\C P^n,\omega_{FS})$ such that both embeddings are of $[I]$-K\"ahler type.


\hfill


\noindent
{\bf Proof of \ref{_existence-connectedness-CPn-l-power-n-balls-COPY_Theorem_} (=\ref{_existence-connectedness-CPn-l-power-n-balls_Theorem_}):}

As we already explained in Section~\ref{_main-results-CPn_Subsection_}, we only need to prove parts B and C.

Part B is exactly
\ref{_CPn-SymplH-action-on-Teich_Corollary_} that we have just proved.

Let us prove part C. By \ref{_CPn-Teich_Corollary_}, for each connected component $\cC_0$ of $\compK (\C P^n)$
the corresponding connected component $\Teich_{\cC_0} (\C P^n)$ of $\Teich (\C P^n)$ is just a point. Consequently, for any $\bfr$ the set $\ttK_{\cC_0} (\bfr)\subset \Teich_{\cC_0} (\C P^n)$ is connected. Now applying \ref{_existence-connectedness-Kahler-type-arb-mfds_Theorem_}
along with part B of this theorem we get part C.

This finishes the proof.
\endproof


\hfill


Now consider the case of $\C P^2$. Similarly to our previous notation, we write  $(\widehat{\C P^2})^k$ for the set of all ordered $k$-tuples of pairwise distinct points in $\C P^2$.
As before, given $\bfx = (x_1,\ldots,x_k)\in (\widehat{\C P^2})^k$, the manifold $\widetilde{\C P^2}_{I_{st},\bfx}$
is the complex blow-up of $(\C P^2,I_{st})$ at $x_1,\ldots,x_k$.
Given $\bfr = (r_1,\ldots,r_k)\in (\R_{>0})^k$ and $J\in\compK (\C P^2)$, set $\alpha := [\omega]$ and
\[
\talpha_\bfr := \Pi^* [\omega] - \pi \sum_{i=1}^k r_i^2 e_i \in H^2 (\widetilde{\C P^2}_{J,\bfx};\R).
\]


\hfill


\definition
\label{_CP2-points-in-gen-position_Definition_}

Assume that $k\in\Z$, $1\leq k\leq 8$.

We say that $k$ points in $\C P^2$ are in {\bf general position} if

\smallskip
\noindent
(i) no three points are on a line;

\smallskip
\noindent
(ii) no six points are on a conic;

\smallskip
\noindent
(iii) no cubic passes through the points with one of the points being a singular
point.


\hfill


\proposition
\label{_tuples-of-pts-in-gen-position-in-cp2-form-dens-conn-set_Proposition_}

Let $1\leq k\leq 8$.

Then the set of tuples $\bfx = (x_1,\ldots,x_k)\in (\widehat{\C P^2})^k$ such that $x_1,\ldots,x_k$ are in general position is connected
and dense in $(\widehat{\C P^2})^k$.


\hfill


\noindent
{\bf Proof of \ref{_tuples-of-pts-in-gen-position-in-cp2-form-dens-conn-set_Proposition_}:}

Consider the product complex structure on $(\widehat{\C P^2})^k\subset (\C P^2)^k$.
In view of \ref{_CP2-points-in-gen-position_Definition_},
the $k$-tuples $(x_1,\ldots,x_k)\in(\widehat{\C P^2})^k$ such that $x_1,\ldots,x_k$
are {\sl not} in general position lie in a proper analytic subset of
$(\widehat{\C P^2})^k$. This means that the set of $k$-tuples $(x_1,\ldots,x_k)\in(\widehat{\C P^2})^k$
such that $x_1,\ldots,x_k$
{\sl are} in general position
contains the complement of a proper analytic subset of $(\widehat{\C P^2})^k$. Thus, it is connected (because $(\widehat{\C P^2})^k$ itself is connected) and
 dense (and open) in $(\widehat{\C P^2})^k$.
\endproof


\hfill


\proposition
\label{_CP2-blown-up-at-1-8-points-Kahler-cone_Proposition_}

Assume that $1\leq k\leq 8$.
Let $\bfr = (r_1,\ldots,r_k)\in (\R_{>0})^k$, $r_1\geq r_2\geq\ldots\geq r_k$. Let $\cC_0$ be the connected component of $\compK (\C P^2)$ containing $I_{st}$.

Then for any $\bfx = (x_1,\ldots,x_k)\in (\widehat{\C P^2})^k$ such that $x_1,\ldots,x_k\in\C P^2$ are in general position,
the class $\talpha_\bfr\in H^2 (\widetilde{\C P^2}_{I_{st},\bfx};\R)$ is K\"ahler with respect to $\tI_{st}$ if and only if the numbers $r_1,\ldots,r_k$ satisfy the following inequalities (listed in \cite[Cor. 1.3G]{_McD-Polt_}):

\smallskip
\noindent
(v) $\sum_{i=1}^k r_i^4 < 1$,

\smallskip
\noindent
(c1) $r_1^2 + r_2^2 < 1$, if $2\leq k\leq 8$,

\smallskip
\noindent
(c2) $r_1^2 +\ldots + r_5^2 < 2$, if $5\leq k\leq 8$,

\smallskip
\noindent
(c3) $2r_1^2 + \sum_{i=2}^7 r_i^2 < 3$, if $7\leq k\leq 8$,

\smallskip
\noindent
(c4) $2r_1^2 + 2r_2^2 + 2r_3^2 + r_4^2 +\ldots +r_8^2 < 4$, if $k=8$,

\smallskip
\noindent
(c5) $2\sum_{i=1}^6 r_i^2 + r_7^2 + r_8^2 < 5$, if $k=8$,

\smallskip
\noindent
(c6) $3 r_1^2 + 2\sum_{i=2}^8 r_i^2 < 6$, if $k=8$.


\hfill


\noindent
{\bf Proof of \ref{_CP2-blown-up-at-1-8-points-Kahler-cone_Proposition_}:}

If $1\leq k\leq 8$ and the points $x_1,\ldots,x_k\in\C P^2$ are in general position (such points exist by \ref{_tuples-of-pts-in-gen-position-in-cp2-form-dens-conn-set_Proposition_}), then $(\widetilde{\C P^2}_{I_{st},\bfx},\tI_{st})$,
for $\bfx = (x_1,\ldots,x_k)$, is a del Pezzo surface (see e.g. \cite[Thm. 8.1.25]{_Dolgachev-book_}).
As such, it is a good generic surface in the sense of
\cite[Def. 2.1]{_Friedman-Morgan-JDG88_} (see  \cite[Rem. 2.9]{_Friedman-Morgan-JDG88_}). The K\"ahler cone of $(\widetilde{\C P^2}_{I_{st},\bfx},\tI_{st})$
is the cone of classes in $H^2 (\widetilde{\C P^2}_{I_{st},\bfx};\R)$ that have a positive square and
evaluate positively on the homology classes of holomorphically
embedded 2-spheres with self-intersection index equal to $-1$ and on the class
$3[\C P^1] - \sum_{i=1}^k [E_i] \in H_2 (\widetilde{\C P^2}_{I_{st},\bfx};\Z)$ \cite{_Friedman-Morgan-JDG88_}.
Together with the information on the homology classes of holomorphically
embedded 2-spheres with self-intersection index equal to $-1$ given in \cite{_Demazure-LNM777-1980_},
this yields that $\talpha_\bfr$ is K\"ahler
with respect to $\tI_{st}$ if and only if the numbers $r_1,\ldots,r_k$ satisfy the inequalities (v), (c1)-(c6) -- see \cite[Cor. 1.3G]{_McD-Polt_}.

This finishes the proof of the proposition.
\endproof


\hfill


Let us recall the following result of Gromov.


\hfill


\proposition {\bf \cite{_Gromov_}}
\label{_Gromov-Symp-CP2_Proposition_}

The group $\Symp (\C P^2,\omega_{FS})$ is homotopy equivalent to $PU (3)$ and, in particular, connected.


\hfill


Recall \ref{_CP2-only-one-conn-cmpt-of-compK-compatible-with-omegaFS_Proposition_}:


\hfill


\proposition {\bf (=\ref{_CP2-only-one-conn-cmpt-of-compK-compatible-with-omegaFS_Proposition_})}
\label{_CP2-only-one-conn-cmpt-of-compK-compatible-with-omegaFS-COPY_Proposition_}

There is only one connected component of $\compK (\C P^2)$ compatible with $\omega_{FS}$ -- it is the connected component
of $\compK (\C P^2)$ containing $I_{st}$. Any two complex structures in that connected component are isotopic.


\hfill


\noindent
{\bf Proof of \ref{_CP2-only-one-conn-cmpt-of-compK-compatible-with-omegaFS-COPY_Proposition_} (= \ref{_CP2-only-one-conn-cmpt-of-compK-compatible-with-omegaFS_Proposition_}):}

The proposition follows immediately from \ref{_CPn-SymplH-action-on-Teich_Corollary_}, \ref{_Gromov-Symp-CP2_Proposition_}.
and part (II) of \ref{_complex-str-on-CPn-is-unique_Theorem_}.
\endproof


\hfill


Now we are ready to prove \ref{_CP2-existence-connectedness_Theorem_} -- let us recall it here.


\hfill


\theorem {\bf (=\ref{_CP2-existence-connectedness_Theorem_})}
\label{_CP2-existence-connectedness-COPY_Theorem_}

\noindent
A. Assume that $1\leq k\leq 8$ and $r_1\geq r_2\geq\ldots\geq r_k>0$.

Then for any proper (possibly empty) complex submanifold $\Sigma\subset (\C P^2, I_{st})$,
any symplectic embedding $\bigsqcup_{i=1}^k B^4 (r_i)\to (\C P^2\setminus \Sigma,\omega_{FS})$ is, in fact, of $[I_{st}]$-K\"ahler type:
it is holomorphic with respect to a complex structure on $\C P^2$ that is compatible with $\omega_{FS}$ and isotopic
to $I_{st}$ by an isotopy fixing $\Sigma$ (as a set).

Accordingly, by part A of \ref{_McD-Polt-CP2_Theorem_}, such a K\"ahler-type embedding $\bigsqcup_{i=1}^k B^4 (r_i)\to (\C P^2\setminus \Sigma,\omega_{FS})$
exists if and only if the radii $r_1,\ldots,r_k$ satisfy the inequalities (v), (c1)-(c6) listed in part A of \ref{_McD-Polt-CP2_Theorem_}.

\bigskip
\noindent
B. For any proper (possibly empty) complex submanifold $\Sigma\subset (\C P^2, I_{st})$ and any $k=l^2$, $l\in \Z_{>0}$,
any symplectic embedding $\bigsqcup_{i=1}^k B^4 (r)\to (\C P^2\setminus\Sigma,\omega_{FS})$ is, in fact, of $[I_{st}]$-K\"ahler type:
it is holomorphic with respect to a complex structure on $\C P^2$ that is compatible with $\omega_{FS}$
and isotopic to $I_{st}$ by an isotopy preserving $\Sigma$ (as a set).

Accordingly, by part B of \ref{_McD-Polt-CP2_Theorem_}, such K\"ahler-type embeddings  $\bigsqcup_{i=1}^k B^4 (r)\to (\C P^2\setminus\Sigma,\omega_{FS})$ are unobstructed: they exist if and only if $k \Vol (B^4 (r),\omega_0) < \Vol (\C P^2,\omega_{FS})$.

\bigskip
\noindent
C. Let $\Sigma\subset (\C P^2,\omega_{FS})$
be either of the following:

\bigskip
\noindent
(1) the empty set;

\smallskip
\noindent
(2) a finite union of closed compact symplectic submanifolds (without boun\-dary) of real dimension 2
whose pairwise intersections (if they exist) are transverse and $\omega$–orthogonal;

\smallskip
\noindent
(3) a Lagrangian submanifold which is diffeomorphic to either $\SP^2$ or $\R P^2$.

For any $k\in\Z_{>0}$ and $r_1,\ldots,r_k>0$, any two K\"ahler-type embeddings $\bigsqcup_{i=1}^k B^{4} (r_i)\to (\C P^2\setminus\Sigma,\omega_{FS})$ (if they exist!) lie in the same orbit of the  $\Symp^c_0 (\C P^2\setminus\Sigma,\omega_{FS})$-action, meaning that
the space of K\"ahler-type embeddings $\bigsqcup_{i=1}^k B^{4} (r_i)\to (\C P^2\setminus\Sigma,\omega_{FS})$ is connected.


\hfill


\noindent
{\bf Proof of \ref{_CP2-existence-connectedness-COPY_Theorem_} (=\ref{_CP2-existence-connectedness_Theorem_}):}

As we already explained in Section~\ref{_main-results-CPn_Subsection_},
 part B of \ref{_CP2-existence-connectedness_Theorem_} follows from
\ref{_existence-connectedness-CPn-l-power-n-balls_Theorem_}
and from part (II) of \ref{_connectedness-rational-ruled-mfds_Theorem_},
while part C of \ref{_CP2-existence-connectedness_Theorem_} follows directly from part (I) of
\ref{_connectedness-rational-ruled-mfds_Theorem_}.

Let us prove part A.

Assume $f: \bigsqcup_{i=1}^k B^4 (r_i)\to (\C P^2,\omega_{FS})$ is a symplectic embedding. Then, by \ref{_McD-Polt-CP2_Theorem_},
the radii $r_1,\ldots, r_k$ satisfy the conditions (v), (c1)-(c6).
Consequently, \ref{_CP2-blown-up-at-1-8-points-Kahler-cone_Proposition_} together with part (I) of \ref{_existence-connectedness-Kahler-type-arb-mfds_Theorem_} (=\ref{_existence-connectedness-Kahler-type-arb-mfds-COPY_Theorem_}) yields that there exists an $[I_{st}]$-K\"ahler-type embedding $f': \bigsqcup_{i=1}^k B^4 (r_i)\to (\C P^2,\omega_{FS})$. By part (II) of \ref{_connectedness-rational-ruled-mfds_Theorem_}, this implies that $f$ is of $[I_{st}]$-K\"ahler type too.
This finishes the proof of part A and of the theorem.
\endproof


\hfill


\section{The case of $\C P^2\sharp \overline{\C P^2}$}
\label{_Hirzebruch-odd_Section_}

Let $M = \C P^2\sharp \overline{\C P^2}$. As a smooth manifold, $M$ will be identified with the complex blow-up $\widetilde{ \C P}^2_{I_{st},x}$
of $(\C P^2, I_{st})$ at a point $x\in\C P^2$. Let $J_1$ be the complex structure on $M$ identified in this way with the blow-up complex structure
$\tI_{st}$ on $\widetilde{ \C P}^2_{I_{st},x}$.

Denote by $E\subset M = \widetilde{ \C P}^2_{I_{st},x}$ the exceptional divisor.
Let $[\C P^1]\in H_2 (M;\Z)$ be the homology class of a projective line $\C P^1\subset M$. Denote by $h,e\in H^2 (M;\Z)$ the cohomology classes Poincar\'e-dual to the homology classes of $[\C P^1]$ and $[E]$.

Given $0 <\lambda <\pi$, assume that $\omega_\lambda$ is a K\"ahler form on the complex manifold $(M,J_1)$ such that $\langle [\omega_\lambda], [\C P^1]\rangle = \pi$,
$\langle [\omega_\lambda], [E]\rangle = \lambda$.
Thus, $[\omega_\lambda] = \pi h - \lambda e$.

{\bf Further on in this section we fix
\[
2\leq k\leq 8,
\]
\[
R_1\geq R_2\geq\ldots\geq R_{k-1}>0,
\]
and set
\[
r_1:=R_1\geq\ldots\geq r_{i-1}:=R_{i-1}\geq r_i:=\sqrt{\lambda/\pi}\geq
\]
\[
\geq r_{i+1}:= R_i\geq\ldots \geq r_k:= R_{k-1}.
\]}

For each $(x_1,\ldots,x_{k-1})\in \hM^{k-1}$
and each complex structure $I$ on $M$,

\begin{itemize}

\item{} $\tM_{I,x_1,\ldots,x_{k-1}}$
will denote the complex blow-up of $(M,I)$ at $x_1,\ldots,x_{k-1}$.

\item{} $\tI$ will denote the lift of $I$ to $\tM_{I,x_1,\ldots,x_{k-1}}$.

\item{} $\Pi_M: \tM_{I,x_1,\ldots,x_{k-1}}\to M$ denotes the natural projection.

\item{} $e_{1,M},\ldots,e_{k-1,M}\in H^2 (\tM_{I,x_1,\ldots,x_{k-1}};\R)$ are the cohomology classes Poincar\'e-dual to the fundamental homology classes of the exceptional divisors.

\end{itemize}

We say that $k-1$ points $x_1,\ldots, x_{k-1}\in M$ are in {\bf general position} if $x_1,\ldots, x_{k-1}\in M\setminus E\cong \C P^2\setminus x$ and the $k$ points $x_1,\ldots, x_{k-1}, x$ are in general position as points of $\C P^2$ (see \ref{_CP2-points-in-gen-position_Definition_}).


\hfill


\proposition
\label{_Hizebruch-odd-tuples-of-pts-in-gen-position-form-dens-conn-set_Proposition_}

The set of tuples $(x_1,\ldots,x_{k-1})\in \hM^{k-1}$ such that $x_1,\ldots,x_{k-1}$ are in general position is connected
and dense in $\hM^{k-1}$.


\hfill


\noindent
{\bf Proof of \ref{_Hizebruch-odd-tuples-of-pts-in-gen-position-form-dens-conn-set_Proposition_}:}

Arguing as in \ref{_tuples-of-pts-in-gen-position-in-cp2-form-dens-conn-set_Proposition_}, we get that
the set of $(k-1)$-tuples $(x_1,\ldots,x_{k-1})\in\hM^{k-1}$
such that $x_1,\ldots,x_k$
are in general position
contains the complement of a proper analytic subset of $\hM^{k-1}$. Thus, it is connected (since $\hM^{k-1}$ is connected) and
 dense (and open) in $\hM^{k-1}$.
\endproof


\hfill


\proposition
\label{_Hirzebruch-odd-blown-up-at-k-points-Kahler-cone_Proposition_}

Let $x_1,\ldots, x_{k-1}\in M$ be in general position.

Then the cohomology class
\[
\Pi_M^* [\omega_\lambda] - \pi \sum_{i=1}^{k-1} R_i^2 e_{i,M} \in H^2 (\tM_{J_1,x_1,\ldots, x_{k-1}};\R)
\]
is K\"ahler with respect to $\tJ_1$ if and only if the numbers $r_1,\ldots,r_k$ satisfy the inequalities (v), (c1)-(c6) listed in \ref{_CP2-blown-up-at-1-8-points-Kahler-cone_Proposition_}.


\hfill


\noindent
{\bf Proof of \ref{_Hirzebruch-odd-blown-up-at-k-points-Kahler-cone_Proposition_}:}

Since the points $x_1,\ldots, x_{k-1}$ lie in $M\setminus E\cong \C P^2\setminus x$, we can view them as points of $\C P^2$.
Denote by $\Pi_{\C P^2}: \widetilde{\C P^2}_{I_{st},x_1,\ldots, x_{k-1}, x}\to \C P^2$ the blow-up projection and by
$e_{1,\C P^2},\ldots,e_{k-1,\C P^2}\in H^2 (\widetilde{\C P^2}_{I_{st},x_1,\ldots, x_{k-1}, x};\R)$ the cohomology classes Poincar\'e-dual to the fundamental homology clas\-ses of the exceptional divisors.

There exists a natural biholomorphism between $(\tM_{J_1,x_1,\ldots, x_{k-1}}, \tJ_1)$ and
$(\widetilde{\C P^2}_{I_{st},x_1,\ldots, x_{k-1}, x},\tI_{st})$.
Since $[\omega_\lambda] = \pi h - \lambda e$,
this biholomorphism identifies the cohomology class
\[
\Pi_M^* [\omega_\lambda] - \pi \sum_{i=1}^{k-1} R_i^2 e_{i,M} \in H^2 (\tM_{J_1,x_1,\ldots, x_{k-1}};\R)
\]
with the cohomology class
\[
\Pi_{\C P^2}^* (\pi h) - \lambda e - \pi\sum_{i=1}^{k-1} R_i^2 e_{i,\C P^2}\in H^2 (\widetilde{\C P^2}_{I_{st},x_1,\ldots, x_{k-1}, x};\R).
\]
(Here $e\in H^2 (\tM_{J_1,x_1,\ldots, x_{k-1}};\R)$ is identified with the cohomology class in $\widetilde{\C P^2}_{I_{st},x_1,\ldots, x_{k-1}, x}$ associated to the exceptional fiber over $x$). Now the proposition follows from \ref{_CP2-blown-up-at-1-8-points-Kahler-cone_Proposition_}.
\endproof


\hfill


Now we are ready to prove \ref{_Hirzebruch-odd-existence-connectedness_Theorem_} -- let us recall it here.


\hfill


\theorem {\bf (=\ref{_Hirzebruch-odd-existence-connectedness_Theorem_})}
\label{_Hirzebruch-odd-existence-connectedness-COPY_Theorem_}

Let $M = \C P^2\sharp \overline{\C P^2}$ and let $\omega_\lambda$ be a K\"ahler-type form on $M$ as above.

Then the following claims hold:

\bigskip
\noindent
A. Let $2\leq k\leq 8$, $R_1\geq R_2\geq\ldots\geq R_{k-1}>0$. Assume that
$r_1:=R_1\geq\ldots\geq r_{i-1}:=R_{i-1}\geq r_i:=\sqrt{\lambda/\pi}\geq r_{i+1}:= R_i\geq\ldots \geq r_k:= R_{k-1}$.

Then for any proper (possibly empty) complex submanifold $\Sigma\subset (M,J_1)$,
the following conditions are equivalent:

\begin{itemize}

\item{} There exists a $[J_1]$-K\"ahler-type embedding
$\bigsqcup_{i=1}^{k-1} B^4 (R_i)\to (M\setminus\Sigma,\omega_\lambda)$ holomorphic with respect to a complex structure on $M$ that is
compatible with $\omega_\lambda$ and isotopic to $J_1$ by an isotopy preserving $\Sigma$ (as a set).

\item{} The numbers $r_1,\ldots,r_k$ satisfy the inequalities (v), (c1)-(c6) in part A of
 \ref{_McD-Polt-CP2_Theorem_}.

\end{itemize}

For $r_1,\ldots,r_k$ satisfying the inequalities (v), (c1)-(c6) in part A of \ref{_McD-Polt-CP2_Theorem_}, any symplectic embedding
$\bigsqcup_{i=1}^{k-1} B^4 (R_i)\to (M\setminus\Sigma,\omega_\lambda)$ is, in fact, of $[J_1]$-K\"ahler type.

\bigskip
\noindent
B.
For any $k\in\Z_{>0}$ and $R_1,\ldots,R_k>0$,
any two K\"ahler-type embeddings $\bigsqcup_{i=1}^k B^4 (R_i)\to (M,\omega_\lambda)$ (if they exist!) lie in the same orbit of the  $\Symp_0 (M,\omega_\lambda)$-action, meaning that
the space of K\"ahler-type embeddings $\bigsqcup_{i=1}^k B^4 (R_i)\to (M,\omega_\lambda)$ is connected.


\hfill


\noindent
{\bf Proof of \ref{_Hirzebruch-odd-existence-connectedness-COPY_Theorem_} (=\ref{_Hirzebruch-odd-existence-connectedness_Theorem_}):}

As we already explained in Section~\ref{_main-results-Hirzberuch_Subsection_}, part B follows from part (I) of \ref{_connectedness-rational-ruled-mfds_Theorem_}.

Let us prove part A.

Assume that the numbers $r_1,\ldots,r_k$ satisfy the conditions (v), (c1)-(c6). Then, by \ref{_Hirzebruch-odd-blown-up-at-k-points-Kahler-cone_Proposition_}, combined with part (I) of \ref{_existence-connectedness-Kahler-type-arb-mfds_Theorem_}, the union $\bigsqcup_{i=1}^{k-1} B^4 (r_i)$ admits a
K\"ahler-type embedding into $(M,\omega_\lambda)$ which is holomorphic with respect to a complex structure on $M$ isotopic to $J_1$ and compatible with $\omega_\lambda$, meaning that the embedding is of $[J_1]$-K\"ahler type.
 By part (II) of \ref{_connectedness-rational-ruled-mfds_Theorem_}, this implies that
any symplectic embedding $\bigsqcup_{i=1}^{k-1} B^4 (r_i)\to (M,\omega_\lambda)$ is of $[J_1]$-K\"ahler type.

Conversely, assume
there exists a K\"ahler-type embedding $f=\bigsqcup_{i=1}^{k-1} f_i: \bigsqcup_{i=1}^{k-1} B^4 (r_i)\to (M,\omega_\lambda)$ which is holomorphic with respect to a complex structure $I$ compatible with $\omega_\lambda$ and
such that $I = \phi^* J_1$ for some $\phi\in\Diff_0 (M)$.
Set $x_i := f_i (0)$, $i=1,\ldots,k-1$, and $\bfx := (x_1,\ldots,x_{k-1})$, $\phi (\bfx) := \left( \phi(x_1),\ldots,\phi(x_{k-1})\right)$.

Recall that a K\"ahler-type embedding of a closed ball is defined as a K\"ahler-type embedding of a larger open ball. Therefore using
\ref{_tuples-of-pts-in-gen-position-in-cp2-form-dens-conn-set_Proposition_} and composing $f$, if necessary, with a small parallel translation of the balls in $\R^{2n}$,
we may assume without loss of generality that the points $\phi (x_1),\ldots,\phi (x_{k-1})\in M$ are in general position.

By part (I) of \ref{_existence-connectedness-Kahler-type-arb-mfds_Theorem_},
the cohomology class $\talpha_\bfr\in H^2 (\tM_{I,\bfx};\R)$ is K\"ahler with respect to $\tI$.
The diffeomorphism $\phi: M\to M$ induces a diffeomorphism $\tM_{I,\bfx}\to \tM_{J_1,\phi(\bfx)}$ that identifies $\tI$ with $\tJ_1$ as well as the classes $\talpha_\bfr$ in $H^2 (\tM_{I,\bfx};\R)$
and in $H^2 (\tM_{J_1,\phi(\bfx)};\R)$ (abusing the notation, we use the same letter for both classes).
Therefore $\talpha_\bfr\in H^2 (\tM_{J_1,\phi(\bfx)};\R)$ is K\"ahler with respect to $\tJ_1$.
Since $\phi (x_1),\ldots,\phi (x_{k-1})$ are in general position,
\ref{_Hirzebruch-odd-blown-up-at-k-points-Kahler-cone_Proposition_} yields that
the conditions (v), (c1)-(c6) are satisfied.
This proves part A and finishes the proof of the theorem.
\endproof


\hfill


\section{The case of $\C P^1\times \C P^1$}
\label{_Hirzebruch-even_Section_}

Let $M=\C P^1\times \C P^1$ and $J_0$ the standard product complex structure on $M$.

Let $J_{2l}$, $l\in\Z_{\geq 0}$, be a complex structure on $M$ such that $(M,J_{2l})$ is biholomorphic to the Hirzebruch surface
$\F_{2l}$. In particular, we assume that $J_0$ is the standard product complex structure on $M$.
Let $b,f\in H^2 (M;\Z)$ be the cohomology classes Poincar\'e-dual, respectively, to $[\C P^1\times pt]$ and $[F]= [pt \times \C P^1]$.

Let $\omega_0$ be the Fubini-Study form on $\C P^1$ normalized so that $\int_{\C P^1} \omega_0 =\pi$. For each $\mu\geq 1$ define the symplectic form
$\omega_\mu$ on $M$ as
\[
\omega_\mu := \mu \omega_0 \oplus \omega_0.
\]
Thus,
$[\omega_\mu] = \pi (\mu b + f)$.

{\bf Further on in this section, we fix
\[
2\leq k\leq 8,
\]
\[
R_1, R_2,\ldots , R_{k-2}>0, 0< R_{k-1} < 1,
\]
and define $r_1\geq r_2\geq\ldots\geq r_k$ as the numbers
\[
\frac{R_1^2}{\mu+1-R_1^2},\frac{R_2^2}{\mu+1-R_1^2},\ldots, \frac{R_{k-2}^2}{\mu+1-R_1^2},
\frac{\mu-R_1^2}{\mu+1-R_1^2}, \frac{1-R_1^2}{\mu+1-R_1^2},
\]
sorted in the non-increasing order.
}

Abusing the notation, for a blow-up of $M$ (respectively, of $\C P^2$) at any finite collection of points $x_1,x_2,\ldots$ (respectively, $y_1,y_2,\ldots$)
we will use the same letter $\Pi_M$ (respectively, $\Pi_{\C P^2}$) for the blow-up projection
and the same letters $e_{1,M}, e_{2,M},\ldots$ (respectively, $e_{1,\C P^2}, e_{2,\C P^2},\ldots$) for the cohomology classes Poincar\'e-dual to the homology classes of the exceptional
divisors.
For a complex structure $I$ on $M$ we will denote by $\tI$ the lift of $I$ to any complex blow-up of $(M,I)$, and similarly we will denote  by
$\tI_{st}$ the lift of the standard complex structure $I_{st}$ on $\C P^2$ to any complex blow-up of $\C P^2$.

Let $x_1,\ldots, x_{k-1}\in M$ be distinct points in $M$.

Pick $y_{k-1}, y_k\in\C P^2$. Consider the complex blow-ups $\tM_{J_0,x_{k-1}}$ of $(M,J_0)$ at $x_{k-1}$
and $(\widetilde{\C P^2}_{I_{st},y_{k-1},y_k},\tI_{st})$ of $(\C P^2,I_{st})$ at $y_{k-1},y_k$.

There exists a biholomorphism $\Psi:  (\tM_{J_0,x_{k-1}},\tJ_0)\to (\widetilde{\C P^2}_{I_{st},y_{k-1},y_k}, \tI_{st})$ such that
the induced isomorphism
\[
\Psi^*: H^2 (\widetilde{\C P^2}_{I_{st},y_{k-1},y_k};\Z)\to H^2 (\tM_{J_0,x_{k-1}};\Z)
\]
identifies
\[
\Pi_M^* b, \Pi_M^* f, e_{k-1,M}\in H^2 (\tM_{J_0,x_{k-1}};\Z)
\]
respectively with
\[
\Pi_{\C P^2}^* h - e_{k-1,\C P^2}, \Pi_{\C P^2}^* h - e_{k,\C P^2}, \Pi_{\C P^2}^* h - e_{k-1,\C P^2} - e_{k,\C P^2} \in  H^2 (\widetilde{\C P^2}_{I_{st},y_{k-1},y_k};\Z).
\]

Identify the complement of the exceptional divisor in $\tM_{J_0,x_{k-1}}$ with $M\setminus x_{k-1}$ and the complement of the union of the exceptional divisors in $\widetilde{\C P^2}_{I_{st},y_{k-1},y_k}$ with $\C P^2\setminus \{ y_{k-1},y_k\}$. The restriction of
$\Psi$ to the complements defines then a map $\psi: M\setminus x_{k-1}\to \C P^2\setminus \{ y_{k-1},y_k\}$.

We say that the points $x_1,\ldots, x_{k-1}\in M\setminus x_{k-2}$ are {\bf in general position} if the points
$\psi (x_1),\ldots, \psi (x_{k-2}), y_{k-1}, y_k\in \C P^2$
are distinct and in general position as $k$ points of $\C P^2$. It is easy to see that this definition does not depend on the choices of $y_1,y_2$ and $\Psi$.

Arguing as in the proofs of \ref{_tuples-of-pts-in-gen-position-in-cp2-form-dens-conn-set_Proposition_} and
\ref{_Hizebruch-odd-tuples-of-pts-in-gen-position-form-dens-conn-set_Proposition_}, and using the fact that any two points of $M$
can be mapped into each other by a biholomorphism of $(M,J_0)$,
we get the following claim.


\hfill


\proposition
\label{_Hizebruch-even-tuples-of-pts-in-gen-position-form-dens-conn-set_Proposition_}

The set of tuples $(x_1,\ldots,x_{k-2}, x_{k-1})\in \hM^{k-1}$ that are in general position is connected
and dense in $\hM^{k-1}$.
\endproof


\hfill


\proposition
\label{_Hirzebruch-even-blown-up-at-k-points-Kahler-cone_Proposition_}

Assume $x_1,\ldots, x_{k-1}\in M$ are in general position.

Then the cohomology class
\[
\Pi_M^* [\omega_\mu] - \pi \sum_{i=1}^{k-1} R_i^2 e_{i,M} \in H^2 (\tM_{J_0,x_1,\ldots, x_{k-1}};\R)
\]
is K\"ahler with respect to $\tJ_0$ if and only if the numbers $r_1,\ldots,r_k$ satisfy the inequalities (v), (c1)-(c6) listed in \ref{_CP2-blown-up-at-1-8-points-Kahler-cone_Proposition_}.


\hfill


\noindent
{\bf Proof of \ref{_Hirzebruch-even-blown-up-at-k-points-Kahler-cone_Proposition_}:}

The biholomorphism $\Psi:  (\tM_{J_0,x_{k-1}},\tJ_0)\to (\widetilde{\C P^2}_{I_{st},y_{k-1},y_k}, \tI_{st})$
induces a biholomorphism
\[
\Upsilon: \left(\tM_{J_0,x_1,\ldots,x_{k-1}},\tJ_0\right)\to
\left(\widetilde{\C P^2}_{I_{st},\psi(x_1),\ldots,\psi (x_{k-2}), y_{k-1}, y_k}, \tI_{st}\right)
\]
that identifies the cohomology class
\[
\Pi_M^* [\omega_\mu] - \pi \sum_{i=1}^{k-1} R_i^2 e_{i,M} = \pi
\left(\mu \Pi_M^* b + \Pi_M^* f - \sum_{i=1}^{k-1} R_i^2 e_{i,M}\right)\in H^2 (\tM_{J_0,x_1,\ldots,x_{k-1}};\R)
\]
with
\[
\pi \cdot \Big( \mu \left(\Pi_{\C P^2}^* h - e_{k-1,\C P^2}\right) + \left(\Pi_{\C P^2}^* h - e_{k,\C P^2}\right) -
\]
\[
- R_{k-1}^2 \left(\Pi_{\C P^2}^* h - e_{k-1,\C P^2} - e_{k,\C P^2}\right) - \sum_{i=1}^{k-2} R_i^2 e_{i,\C P^2} \Big)
\in H^2 \left(\widetilde{\C P^2}_{I_{st},\psi(x_1),\ldots,\psi (x_{k-2}), y_{k-1}, y_k};\R\right).
\]
Recalling the description of $\Psi^*: H^2 (\widetilde{\C P^2}_{I_{st},y_{k-1},y_k};\Z)\to H^2 (\tM_{J_0,x_{k-1}};\Z)$
above and the fact that $[\omega_{FS}] = \pi h\in H^2 (\C P^2;\R)$, we get that
\[
\pi \cdot \Big( \mu \left(\Pi_{\C P^2}^* h - e_{k-1,\C P^2}\right) + \left(\Pi_{\C P^2}^* h - e_{k,\C P^2}\right)  -
\]
\[
- R_{k-1}^2 \left(\Pi_{\C P^2}^* h - e_{k-1,\C P^2} - e_{k,\C P^2}\right) -
 \sum_{i=1}^{k-2} R_i^2 e_{i,\C P^2} \Big) =
\]
\[
\pi \cdot\Big( \left(\mu +1 - R_{k-1}^2\right) \Pi_{\C P^2}^* h  - \sum_{i=1}^{k-2} R_i^2 e_{i,\C P^2} -
\]
\[
- \left(\mu - R_{k-1}^2\right) e_{k-1,\C P^2} -
(1- R_{k-1}^2) e_{k,\C P^2} \Big) =
\]
\[
= \left(\mu+ 1 - R_{k-1}^2\right) \Bigg( \Pi_{\C P^2}^* (\pi h) - \pi\cdot \sum_{i=1}^{k-2} \frac{R_i^2}{\mu +1 - R_{k-1}^2} e_{i,\C P^2} -
\]
\[
- \frac{\mu - R_{k-1}^2}{\mu +1 - R_{k-1}^2} e_{k-1, \C P^2} - \frac{1 - R_{k-1}^2}{\mu +1 - R_{k-1}^2} e_{k,\C P^2}
\Bigg) =
\]
\[
= \left(\mu+ 1 - R_{k-1}^2\right) \Bigg( \Pi_{\C P^2}^* [\omega_{FS}] - \pi \sum_{i=1}^{k-2} \frac{R_i^2}{\mu +1 - R_{k-1}^2} e_{i,\C P^2} -
\]
\[
- \frac{\mu - R_{k-1}^2}{\mu +1 - R_{k-1}^2} e_{k-1,\C P^2} - \frac{1 - R_{k-1}^2}{\mu +1 - R_{k-1}^2} e_{k,\C P^2}
\Bigg)
\]

Thus, the cohomology class
\[
\Pi_M^* [\omega_\mu] - \pi \sum_{i=1}^{k-1} R_i^2 e_{i,M} = \pi \cdot
\left(\mu \Pi_M^* b + \Pi_M^* f - \sum_{i=1}^{k-1} R_i^2 e_{i,M}\right)\in H^2 (\tM_{J_0,x_1,\ldots,x_{k-1}};\R)
\]
is K\"ahler with respect to $\tJ_0$ if and only if the cohomology class
\[
\Pi_{\C P^2}^* [\omega_{FS}] - \pi \sum_{i=1}^{k-2} \frac{R_i^2}{\mu +1 - R_{k-1}^2} e_{i,\C P^2} -
\frac{\mu - R_{k-1}^2}{\mu +1 - R_{k-1}^2} e_{k-1,\C P^2} -
\]
\[
-
\frac{1 - R_{k-1}^2}{\mu +1
- R_{k-1}^2} e_{k,\C P^2}
\in H^2 \left(\widetilde{\C P^2}_{I_{st},\psi(x_1),\ldots,\psi (x_{k-2}), y_{k-1}, y_k};\R\right)
\]
is K\"ahler with respect to $\tI_{st}$.
Since the numbers
\[
\frac{R_1^2}{\mu +1 - R_{k-1}^2},\ldots ,  \frac{R_{k-2}^2}{\mu +1 - R_{k-1}^2},
\frac{\mu - R_{k-1}^2}{\mu +1 - R_{k-1}^2}, \frac{1 - R_{k-1}^2}{\mu +1 - R_{k-1}^2},
\]
sorted in non-increasing order,
are exactly $r_1,\ldots, r_k$, and the points $\psi(x_1),\ldots,\psi (x_{k-2}), y_{k-1}, y_k\in \C P^2$
are in general position, the proposition follows from
\ref{_CP2-blown-up-at-1-8-points-Kahler-cone_Proposition_}.
\endproof


\hfill


Now we are ready to prove \ref{_Hirzebruch-even-existence-connectedness_Theorem_} -- let us recall it here.


\hfill


\theorem {\bf (=\ref{_Hirzebruch-even-existence-connectedness_Theorem_})}
\label{_Hirzebruch-even-existence-connectedness-COPY_Theorem_}

Let $M= \C P^1\times \C P^1$ and let $\omega_\mu$ be as above.

Let $\Sigma\subset (M,J_0)$ be either a proper (possibly empty) complex submanifold or $\Sigma = (\textrm{pt}\times \C P^1) \cup (\C P^1\times \textrm{pt})$.

Then the following claims hold:

\bigskip
\noindent
A. Assume $2\leq k\leq 8$,
$R_1, R_2,\ldots , R_{k-2}>0$, $0< R_{k-1} < 1$.
Consider the numbers
\[
\frac{R_1^2}{\mu+1-R_{k-1}^2},\frac{R_2^2}{\mu+1-R_{k-1}^2},\ldots,
\frac{\mu-R_{k-1}^2}{\mu+1-R_{k-1}^2}, \frac{1-R_{k-1}^2}{\mu+1-R_{k-1}^2},
\]
sort them in the non-increasing order and denote the resulting $k$ numbers by $r_1\geq\ldots\geq r_k$.

Then the following conditions are equivalent:

\begin{itemize}

\item{} There exists a $[J_0]$-K\"ahler-type embedding
$\bigsqcup_{i=1}^{k-1} B^4 (R_i)\to (M\setminus\Sigma,\omega_\mu)$ holomorphic with respect to a complex structure on $M$
that is compatible with
$\omega_\mu$ and isotopic to $J_0$ by an isotopy preserving $\Sigma$ (as a set).

\item{} The numbers $r_1,\ldots,r_k$ defined above satisfy the inequalities in part A of
 \ref{_McD-Polt-CP2_Theorem_}.

\end{itemize}

For $r_1,\ldots,r_k$ satisfying the inequalities (v), (c1)-(c6) in part A of \ref{_McD-Polt-CP2_Theorem_}, any symplectic embedding
$\bigsqcup_{i=1}^{k-1} B^4 (R_i)\to (M\setminus\Sigma,\omega_\mu)$ is, in fact, of $[J_0]$-K\"ahler type.

\bigskip
\noindent
B. Assume that $\mu=1$.

Then for any $l\in\Z_{>0}$,
any symplectic embedding of a disjoint union of $2 l^2$ equal balls into $(M\setminus\Sigma,\omega_1)$
is, in fact, of $[J_0]$-K\"ahler type, and such symplectic (or, equivalently, $[J_0]$-K\"ahler-type) embeddings are unobstructed.

More precisely, if $2 l^2 \Vol (B^4 (r),\omega_0) < \Vol (M,\omega_1)$, then there exists a K\"ahler-type embedding of $\bigsqcup_{i=1}^{2 l^2} B^4 (r)$
into $(M\setminus\Sigma,\omega_1)$ which is holomorphic with respect to a complex structure on $M$ that is compatible with $\omega_1$
and isotopic to $J_0$
 by an isotopy preserving $\Sigma$ (as a set).

\bigskip
\noindent
C. For any $k\in\Z_{>0}$ and $R_1,\ldots,R_k>0$,
any two K\"ahler-type embeddings $\bigsqcup_{i=1}^k B^{4} (R_i)\to (M\setminus\Sigma,\omega_\mu)$ (if they exist!) lie in the same orbit of the  $\Symp^c_0 (M\setminus\Sigma,\omega_\mu)$-action, meaning that
the space of K\"ahler-type embeddings $\bigsqcup_{i=1}^k B^{4} (R_i)\to (M\setminus\Sigma,\omega_\mu)$ is connected.


\hfill


\noindent
{\bf Proof of \ref{_Hirzebruch-even-existence-connectedness-COPY_Theorem_} (=\ref{_Hirzebruch-even-existence-connectedness_Theorem_}):}

As we already explained in Section~\ref{_main-results-Hirzberuch_Subsection_}, part B of the theorem follows from \ref{_existence-connectedness-CPn-l-power-n-balls_Theorem_} and part C follows from part (I) of \ref{_connectedness-rational-ruled-mfds_Theorem_}.

The proof of part A
follows the proof of part A of \ref{_Hirzebruch-odd-existence-connectedness-COPY_Theorem_}, with
\ref{_Hirzebruch-even-blown-up-at-k-points-Kahler-cone_Proposition_}
used instead of
\ref{_Hirzebruch-odd-blown-up-at-k-points-Kahler-cone_Proposition_}.
\endproof


\hfill


\section{Existence of K\"ahler-type embeddings of balls and Cam\-pa\-na-simple complex structures}
\label{_Campana-simple-strs-gen-mfds-main-results_Section_}

In this section we discuss K\"ahler-type embeddings of balls and other domains into symplectic manifolds admitting K\"ahler-type complex structures with ``few" complex subvarieties. The relevant claims about such complex structures on tori, K3 surfaces and IHS-hyperk\"ahler manifolds will be discussed in Sections~\ref{_Campana-simple-structures-tori-K3_Section_} and \ref{_hyperkahler-case-proofs_Section_}.
The results of Sections~\ref{_main-results-tori-K3_Subsection_}, \ref{_hyperkahler-case_Subsection_} on the existence of
K\"ahler-type embeddings of disjoint unions of balls into tori, K3 surfaces and IHS-hyperk\"ahler manifolds will then follow as particular cases of the general
result (\ref{_Campana-simple-Kahler-type-existence-connectedness-gen-mfd_Theorem_}) stated further on in this section.

As before, let $M$, $\dim_\R = 2n$, be a closed connected manifold equipped with a K\"ahler-type symplectic form $\omega$.

Let us recall the definition of a Campana-simple complex structure.

Assume $I$ is a K\"ahler-type complex structure on $M$. Let $\goU$ be the union of all positive-dimensional proper complex
subvarieties of $(M,I)$. The set $\goU$ is
either a countable union of proper analytic subvarieties of $M$ (and hence has a dense connected complement) or the whole $M$ (see e.g. \cite[Remark 4.2]{_EV-JTA_}).

If $\goU$ is a countable union of proper analytic subvarieties of $M$, the complex structure $I$ is called {\bf Campana-simple} and
the points of $M\setminus \goU$ are called {\bf Campana-generic with respect to $I$}.

Examples of manifolds with
a Campana-simple  complex structure
include the generic complex tori of complex dimension $>1$
(see \ref{_irrational-sympl-form-on-torus-compatible-with-Campana-simple-complex-str_Proposition_}), as well as generic K3 surfaces (see \ref{_irrational-sympl-form-on-K3-compatible-with-Campana-simple-complex-str_Proposition_}),
and generic deformations of hyperk\"ahler manifolds \cite{_Verbitsky:Deforma_}.
One can show that a torus equipped with a Campana-simple  complex structure does not admit {\sl any} positive-dimensional
proper complex subvarieties (see \ref{_any-Campana-simple-cs-on-torus-is-Campana-supersimple_Proposition_}), and the same holds for a generic complex K3 surface \cite{_Verb-trianalytic-GAFA1998_}.


\hfill


\remark
\label{_Campana-simple-vs-projectivity_Remark_}

\noindent
1. If $I$ is Campana-simple, then the complex manifold $(M,I)$ is not projective, since otherwise it would have admitted a globally defined meromorphic function
$f$ and could have been represented as the union of zero divisors of the functions $f-a$, for all
$a\in \C$ and the zero divisor of
$f^{-1}$. Consequently, the symplectic form $\omega$ compatible with a Campana-simple complex structure cannot be rational -- otherwise,
for a Campana-simple complex structure $I$ compatible with $\omega$ the complex manifold $(M,I)$ would have been projective,
by the Kodaira embedding theorem \cite{_Kod-emb_}.

\bigskip
\noindent
2. Conjecturally (see \cite[Question 1.4]{_Campana:isotrivial_},
\cite[Conj. 1.1]{_CDV:threefolds_}), any closed connected manifold equipped with a Campana-simple complex structure is bimeromorphic to a
hyperk\"ahler orbifold or a finite quotient of a torus.
For complex threefolds this was proved in \cite{_Horing-Peternell-Sal-Lake-City_}.


\hfill


Campana-simple
complex struc\-tures on $M$ form a $\Diff (M)$-invariant subset of the set of all K\"ahler-type complex structures on $M$.


\hfill


\definition
\label{_Teich-s_Definition_}

Denote by $\Teich^s (M)\subset \Teich (M)$
the set of points of
$\Teich (M)$ that can be represented by Campana-simple complex structures.

If $\cC_0$ is a connected component of $\compK (M)$, denote by $\Teich^s_{\cC_0} (M)\subset \Teich_{\cC_0} (M)$
the set of points of
$\Teich_{\cC_0} (M)$ that can be represented by Campana-simple
complex structures.
In other words,
\[
\Teich^s_{\cC_0} (M) = \Teich^s (M)\cap \Teich_{\cC_0} (M).
\]


\hfill


\theorem
\label{_Campana-simple-Kahler-type-existence-connectedness-gen-mfd_Theorem_}

 Let $\cC_0$ be
a connected component of $\compK (M)$ compatible with $\omega$. Let $k\in\Z_{>0}$, $\bfr= (r_1,\ldots,r_k)\in (\R_{>0})^k$.

Assume $\Vol \left(\bigsqcup_{i=1}^k B^{2n} (r_i),\omega_0\right)<\Vol (M,\omega)$.

Then the following claims hold:

\smallskip
\noindent
(I)
Let $M$, $\dim_\R M = 2n$, be a closed manifold.
Assume that $\omega$ is a symplectic form on $M$ compatible with a Campana-simple complex structure $I$ on $M$.

Then K\"ahler-type embeddings $\bigsqcup_{i=1}^k B^{2n} (r_i)\to (M,\omega)$ are unobstructed.

More precisely, if $\Vol \left(\bigsqcup_{i=1}^k B^{2n} (r_i),\omega_0\right)<\Vol (M,\omega)$,
then $\bigsqcup_{i=1}^k B^{2n} (r_i)$ admits an $[I]$-K\"ahler-type embedding into $(M,\omega)$.

\bigskip
\noindent
(II) Assume that
$\cmptteich_{\cC_0} (M,\omega)\cap \Teich^s_{\cC_0} (M)$ is a dense connected subset of $\cmptteich_{\cC_0} (M,\omega)$.

\smallskip
Then any two K\"ahler-type embeddings $\bigsqcup_{i=1}^k B^{2n}(r_i)\to (M,\omega)$ favoring $\cC_0$
lie in the same orbit of the $\Symp (M,\omega)\cap \Diff_0 (M)$-action.
In particular, there exists $[I]\in\cmptteich_{\cC_0} (M,\omega)$
such that both embeddings are of $[I]$-K\"ahler-type.

If, in addition, $\SympH (M)$ acts transitively on the set of connected components of $\compK (M)$ compatible with $\omega$, then
any two K\"ahler-type embeddings $\bigsqcup_{i=1}^k B^{2n}(r_i)\to (M,\omega)$
lie in the same orbit of the $\SympH (M,\omega)$-action.


\hfill


For the proof of \ref{_Campana-simple-Kahler-type-existence-connectedness-gen-mfd_Theorem_}
we need the following claim.


\hfill


\proposition
\label{_Campana-simple-Kahler-cone_Proposition_}

Assume $I$ is a Campana-simple complex structure on $M$, a cohomology class $\alpha\in H^2 (M;\R)$ is K\"ahler with respect to $I$,
and the points $x_1,\ldots,x_k\in M$ are Campana-generic with respect to $I$. Set $\bfx := (x_1,\ldots,x_k)$.

Then the cohomology class
\[
\talpha_\bfr := \Pi^* \alpha - \pi\sum_{i=1}^k r_i^2 e_i \in H^2 (\tM_{I,\bfx};\R)
\]
is K\"ahler with respect to $\tI$ if and only if $\talpha_\bfr^n >0$,
or, equivalently,
if and only if
\[
\Vol \left(\bigsqcup_{i=1}^k B^{2n} (r_i),\omega_0\right) < \langle \alpha^n, [M]\rangle.
\]


\hfill


\noindent
{\bf Proof of \ref{_Campana-simple-Kahler-cone_Proposition_}:}

This is proved in \cite[Thm. 8.6]{_EV-JTA_} -- the proof is a straightforward application of the Demailly-Paun theorem \cite[Thm. 0.1]{_Dem-Paun_} describing the K\"ahler cone of an arbitrary closed manifold equipped with a K\"ahler-type complex structure.
\endproof


\hfill


Now we can prove \ref{_Campana-simple-Kahler-type-existence-connectedness-gen-mfd_Theorem_}.


\hfill


\noindent
{\bf Proof of \ref{_Campana-simple-Kahler-type-existence-connectedness-gen-mfd_Theorem_}:}

The complex structure $I$ appearing in the statement of part (I) of the theorem is compatible with $\omega$. Hence,
the cohomology class $\alpha :=[\omega]$ is K\"ahler with respect to $I$.

Pick points $x_1,\ldots,x_k\in M\setminus\Sigma$ that are Campana-generic with respect to $I$.
It is possible to pick these points in the complement of $\Sigma$, because the set of Campana-generic points of $(M,I)$ is dense in $M$
(see above)
and $\Sigma$ is a proper complex submanifold of $(M,I)$.
Set $\bfx := (x_1,\ldots,x_k)$.

Assume that
\[
\Vol \left(\bigsqcup_{i=1}^k B^{2n} (r_i),\omega_0\right) < \Vol (M,\omega) = \langle \alpha^n, [M]\rangle.
\]
Then \ref{_Campana-simple-Kahler-cone_Proposition_} yields that $\talpha_\bfr := \Pi^* \alpha - \pi\sum_{i=1}^k r_i^2 e_i\in H^2 (\tM_{I,\bfx};\R)$ is K\"ahler with respect to $\tI$.
Now
 the claim
of part (I)
follows from part (I) of \ref{_existence-connectedness-Kahler-type-arb-mfds_Theorem_}.

Let us prove part (II).

Assume that $\cmptteich_{\cC_0} (M,\omega)\cap \Teich^s_{\cC_0} (M)$ is a dense connected subset of $\cmptteich_{\cC_0} (M,\omega)$.
The argument we have just used in the proof of part (I) shows that $I\in\ttK_{\cC_0} (\bfr)$.
Therefore
\[
\cmptteich_{\cC_0} (M,\omega)\cap \Teich^s_{\cC_0} (M)\subset \ttK_{\cC_0} (\bfr).
\]
Hence, $\ttK_{\cC_0} (\bfr)$ is also a dense connected
subset of $\cmptteich_{\cC_0} (M,\omega)$.
Now part (II) of \ref{_existence-connectedness-Kahler-type-arb-mfds_Theorem_} yields the needed claim.

This finishes the proof of \ref{_Campana-simple-Kahler-type-existence-connectedness-gen-mfd_Theorem_}.
\endproof


\hfill


\section{The case of tori and K3 surfaces}
\label{_Campana-simple-structures-tori-K3_Section_}

In this section we discuss Campana-simple complex structures on tori and K3 surfaces and show how to deduce \ref{_existence-connectedness-irrational-forms-tori-K3-surfaces_Theorem_}
from
\ref{_Campana-simple-Kahler-type-existence-connectedness-gen-mfd_Theorem_}.


\hfill


\subsection{Campana-simple complex structures on tori}
\label{_Campana-simple-structures-torus_Subsection_}

Denote by $\Lin\subset \compK (\T^{2n})$ the set of linear complex structures on $\T^{2n}$ that are compatible with the standard orientation.
It is naturally identified with the space of linear complex structures on $\R^{2n}$.

It is well-known (see e.g. \cite[Prop. 2.5.2]{_McD-Sal-3_}) that the space $\Lin$ can be identified with
$GL^+ (2n,\R)/GL(n,\C)$ which is a connected smooth manifold.
Here the group $GL (n,\C)$ is embedded in $GL (2n,\R)$ by the map
\[
C\mapsto \left(
           \begin{array}{cc}
             \textrm{Re}\ C & -\textrm{Im}\ C \\
             \textrm{Im}\ C & \textrm{Re}\ C \\
           \end{array}
         \right).
\]
The tangent space $T_J \Lin$ of $\Lin$ at a point
$J\in\Lin$ is the space of all real $2n\times 2n$ matrices $R$ satisfying $RJ+JR=0$.
The complex structure on $T_J \Lin$ is given by $R\mapsto JR$ (since $GL^+ (2n,\R)/GL(n,\C)$ is symmetric space, this
 almost complex structure on $\Lin$ is indeed integrable
\cite{_Besse:Einst_Manifo_}).
Thus, further on $\Lin$ will be viewed as a connected complex manifold.

The following claim is an easy exercise (cf. \cite[Sec. 6]{_EV-JTA_}).


\hfill


\proposition
\label{_Lin-smooth-deformation_Proposition_}

The space $T^{2n}\times \Lin$
can be equipped with a complex structure so that
$T^{2n}\times \Lin\to \Lin$
is a complex-analytic deformation family such that for each $J\in\Lin$ the complex structure on the fiber over $J$ is $J$ itself.
\endproof


\hfill


For $l=1,\ldots,n$ denote by $\Lambda^{2l}$ be the space of exterior $2l$-forms on $\R^{2n}$ and by $\Lambda^{2l}_\Q\subset \Lambda^{2l}$
the space of exterior $2l$-forms on $\R^{2n}$ with rational coefficients.

Let $\omega$ be a linear symplectic form on $\T^{2n}$. Denote by $\bomega\in \Lambda^2$ its lift to $\R^{2n}$ -- it is a linear symplectic form on $\R^{2n}$.

For a linear complex structure $J$ on $\T^{2n}$ denote by $\bJ$ the linear complex structure on $\R^{2n}$ which is the lift of $J$.


\hfill


\proposition
\label{_Kahler-classes-wrt-linear-cs-linear-sympl-forms_Proposition_}

The cohomology class $[\omega]\in H^2 (\T^{2n};\R)$ of the linear symplectic form $\omega$ is K\"ahler with respect to a linear complex structure $J$ on $\T^{2n}$ if and only if $\omega$ is compatible with $J$.


\hfill


\noindent
{\bf Proof of \ref{_Kahler-classes-wrt-linear-cs-linear-sympl-forms_Proposition_}:}

If $\omega$ is compatible with $J$ then, obviously, $[\omega]\in H^2 (\T^{2n};\R)$ is K\"ahler with respect to $J$.

Assume $[\omega]\in H^2 (\T^{2n};\R)$ is K\"ahler with respect to $J$ and let us prove that $\omega$ is compatible with $J$ -- i.e., $\omega$ is K\"ahler on $(\T^{2n},J)$. Since $[\omega]\in H^2 (\T^{2n};\R)$ is K\"ahler with respect to $J$, there exists a K\"ahler form $\eta$ on $(\T^{2n},J)$ cohomologous to $\omega$. Viewing $\T^{2n}$ as a compact Lie group and averaging $\eta$ over $\T^{2n}$, we get a $\T^{2n}$-invariant -- hence, linear -- form $\zeta$ which is K\"ahler on $(\T^{2n},J)$ and cohomologous to $\omega$. Since there is only one linear
form representing a given cohomology class, we get that $\zeta=\omega$, which implies that $\omega$ is K\"ahler on $(\T^{2n},J)$.

This finishes the proof of the proposition.
\endproof


\hfill


Denote by $\Lincmpt (\omega)$
the subset of $\Lin$ formed by the linear complex structures on $\T^{2n}$ compatible with the linear symplectic form $\omega$.
It is naturally identified with the space of the linear complex structures on $\R^{2n}$ compatible with $\bomega$.

For each $l=1,\ldots,n$ and each $a\in H^{2l} (\T^{2n};\R)$, $a\neq 0$, denote by $\cT_a\subset \Lin$ the set of $J\in\Lin$ such that
$a\in H_J^{l,l} (\T^{2n})$. It can also be viewed as the Hodge locus of $a$ for the complex-analytic deformation
family $T^{2n}\times \Lin\to \Lin$ given by \ref{_Lin-smooth-deformation_Proposition_}.
Consequently, $\cT_a$ is a complex subvariety of $\Lin$ -- see e.g. \cite[Vol. 2, Lem. 5.13]{_Voisin-book_}.


\hfill


\example
\label{_Lincmpt-omega-is-analytic-subvariety_Example_}

Let $\omega$ be a linear symplectic form on $\T^{2n}$.

Let $a=[\omega]\in H^2 (\T^{2n};\R)$. Then the set $\cT_a$ is the set of linear complex structures on $\T^{2n}$ preserving $\omega$ (because a real cohomology class of degree 2 is of type $(1,1)$ with respect to a complex structure if and only if it is preserved by it). In this case $\cT_a$ is a smooth complex submanifold of $\Lin$ and this can be verified by elementary means.

Indeed, first of all, $\cT_a$ is smooth (see e.g. \cite[Lem. 2.5.5]{_McD-Sal-3_}). If $J\in\cT_a$, then $T_J \cT_a$
is formed by $2n\times 2n$ matrices $R$ satisfying the two conditions: $RJ+JR=0$ and $R^t A R =0$, where $A$ is the matrix of the linear symplectic form $\bomega$ on $\R^{2n}$. Since $J^t A J = A$ (because  $J\in\cT_a$), $JR$ satisfies the same two conditions meaning that $T_J \cT_a$ is a complex vector subspace of $T_J \Lin$. This means that $\cT_a$ is a complex submanifold of $\Lin$.

The set $\Lincmpt (\omega)$ is then an open contractible submanifold of $\cT_a$ (see e.g. \cite[Lem. 2.5.5]{_McD-Sal-3_}).


\hfill


\proposition
\label{_any-Campana-simple-cs-on-torus-is-Campana-supersimple_Proposition_}

Any Campana-simple complex structure on $\T^{2n}$ admits no proper positive-dimensional complex subvarieties.


\hfill


\noindent
{\bf Proof of \ref{_any-Campana-simple-cs-on-torus-is-Campana-supersimple_Proposition_}:}

A Campana-simple complex structure on $\T^{2n}$ is of K\"ahler type and therefore can be mapped by a diffeomorphism
of $\T^{2n}$ into a linear complex structure (see \cite[Prop. 6.1]{_EV-JTA_}),
which is also Campana-simple. Thus, it suffices to prove the proposition for a Campana-simple linear complex structure
$J$.

Assume, by contradiction, that such a $J$
admits a positive-dimensional proper complex subvariety.
The parallel translations of $\T^{2n}$ preserve $J$ and therefore any point of $\T^{2n}$ is contained in the image of that subvariety under a parallel translation, which contradicts the assumption that $J$ is Campana-simple.
This finishes the proof of the proposition.
\endproof


\hfill


Denote by $\Lin^s\subset\Lin$ the space of linear Campana-simple complex structures on $\T^{2n}$.


\hfill


\proposition
\label{_torus-linear-Campana-simple-cs-compatible-with-omega_Proposition_}

Let $\omega$ be a linear symplectic form on $\T^{2n}$.

Then
$\omega$ is compatible with a Campana-simple linear complex structure (i.e., $\Lin^s\cap \Lincmpt (\omega)\neq\emptyset$)  if and only if $\omega$ is irrational. Moreover, in this case $\Lin^s\cap \Lincmpt (\omega)$ is a dense connected subset of $\Lincmpt (\omega)$.


\hfill


\noindent
{\bf Proof of \ref{_torus-linear-Campana-simple-cs-compatible-with-omega_Proposition_}:}

Assume $\omega$ is rational.
Then, by part 1 of
\ref{_Campana-simple-vs-projectivity_Remark_}, there are no Campana-simple complex structures compatible with $\omega$.

Now assume $\omega$ is irrational, or, equivalently, that $\bomega$ is not a real multiple of a form in $\Lambda^2_\Q$.
Let us prove that $\Lin^s\cap \Lincmpt (\omega)$ is a dense connected subset of $\Lincmpt (\omega)$.

We start with the following lemma.


\hfill


\lemma
\label{_powers-of-irrational-form-are-irrational_Lemma_}

Assume that a linear symplectic form $\eta$ is not a real multiple of a form in $\Lambda^2_\Q$.

Then for all $l=1,\ldots,n-1$ the form $\eta^l$ is not a real multiple of a form in $\Lambda^{2l}_\Q$.


\hfill


\noindent
{\bf Proof of \ref{_powers-of-irrational-form-are-irrational_Lemma_}:}

Assume by contradiction that $\eta^l$ is a real multiple of a form in $\Lambda^2_\Q$. Since $\eta$ is symplectic, $\eta^l\neq 0$ and we can assume without loss of generality that $\eta^l\in\Lambda^{2l}_\Q$.

Consider the Galois group $\Gal (\C/\Q)$ of $\C$ over $\Q$ (i.e., the group of field automorphisms of $\C$) and its action on the coefficients of $\eta$ with respect to the standard basis of $\R^{2n}$. It is proved in \cite{_Lee_} that if $\eta^l=\xi^l$ for some $\xi\in\Lambda^2$ and $l=1,\ldots,n-1$, then $\xi = \pm \eta$. This implies that the action of each element of $\Gal (\C/\Q)$ either fixes all the coefficients of $\eta$ or multiplies each of them by $-1$. The complex numbers preserved by the action of the whole group $\Gal (\C/\Q)$ are exactly the rationals (a proof of this well-known fact can be easily extracted e.g. from \cite{_Yale Automorphisms of the Complex Numbers_}). Therefore the square of each coefficient of $\eta$ lies in $\Q$. This easily implies that there exists a positive $r\in \Q$ such that the field $\Q [\sqrt{r}]$ contains all the coefficients of $\eta$.

If all elements $\Gal (\C/\Q)$ fix all the coefficients of $\eta$, we get that $\eta\in\Lambda^2_\Q$, in contradiction to the hypothesis of the lemma. Therefore there exists an element of $\Gal (\C/\Q)$ whose action multiplies all the coefficients of $\eta$ by $-1$, meaning that all these coefficients are rational multiples of $\sqrt{r}$. Therefore, $\eta$ is a positive real multiple of a form in $\Lambda^2_\Q$, which contradicts again the hypothesis of the lemma. Thus, the assumption that $\eta^l$ is a real multiple of a form in $\Lambda^2_\Q$ has led us to a contradiction, which means that it is false. This proves the lemma.
\endproof


\hfill


Since
$\bomega$ is not a real multiple of a form in $\Lambda^2_\Q$, \ref{_powers-of-irrational-form-are-irrational_Lemma_} yields that for each $l=1,\ldots,n-1$ the form $\bomega^l\in\Lambda^{2l}$ is not a real multiple of a form in $\Lambda^{2l}_\Q$.

Let $a\in H^{2l} (\T^{2n};\Q)$, $a\neq 0$, $l=1,\ldots,n-1$, and let us view $a\in H^{2l} (\T^{2n};\Q)$ as an exterior $2l$-form on $\R^{2n}$ with rational coefficients: $a\in \Lambda^{2l}_\Q$. We also view $\Lincmpt (\omega)$ as the set of linear complex structures on $\R^{2n}$ compatible with $\bomega$.

The set $\cT_a \cap \Lincmpt (\omega)$ is a complex subvariety of
the complex manifold $\Lincmpt (\omega)$.

We claim that this complex subvariety is proper.

Indeed, assume, by contradiction, that it is not. Then, since $\Lincmpt (\omega)$ is connected, $\cT_a \cap \Lincmpt (\omega) = \Lincmpt (\omega)$,
meaning that any linear complex structure on $\R^{2n}$ compatible with $\bomega$ also preserves the exterior $2l$-form $a$.

Denote by $G\subset SL (2n,\R)$ the group of linear automorphisms of $\R^{2n}$ preserving $\bomega$. It is isomorphic, as a real Lie group, to $Sp (2n,\R)$. The set $\Lincmpt (\omega)$ is invariant under conjugations in $G$ and thus generates a normal subgroup $H$ of $G$. But $G\cong Sp (2n,\R)$ is a simple Lie group and therefore each normal
subgroup of $G$ has to either lie in its center (equal to $\{Id,-Id\}$) or coincide with the whole $G$ (see e.g. \cite{_Rag_}). The former option is clearly not satisfied for $H$ and therefore we get that $H = G$, meaning that $a$ is preserved by the whole group $G$. By a theorem of Weyl \cite[Ch. VI, Sec. 1]{_Weyl_}  (see e.g. \cite[Sect. 5.3.2]{_Goodman_Wallach_} for a modern presentation), the only exterior $2l$-forms preserved by the whole $G$ are the real multiples of $\bomega^{2l}$. Thus, $a\in \Lambda^{2l}_\Q$ is a real multiple of $\bomega^{2l}$ which, as we have seen above, is a not a real multiple of any form with rational coefficients. We have obtained a contradiction, and this proves the claim that
$\cT_a \cap \Lincmpt (\omega)$ is a proper complex subvariety of
of the complex manifold $\Lincmpt (\omega)$.

The claim we have proved implies that the
complement in the connected complex manifold $\Lincmpt (\omega)$ of the countable union of proper complex subvarieties $\cT_a \cap \Lincmpt (\omega)$, $a\in H^{2l} (\T^{2n};\Q)$, $a\neq 0$, $l=1,\ldots,n-1$, is connected and dense in $\Lincmpt (\omega)$.
Any complex structure $J$ in this complement is Campana-simple. Indeed, the $2l$-dimensional complex subvarieties of $(\T^{2n}, J)$, $l=1,\ldots,n-1$, all have non-zero rational fundamental classes; the corresponding Poincar\'e-dual non-zero rational cohomology classes are of type $(l,l)$ (see e.g. \cite[Vol. 1, Sec. 11]{_Voisin-book_}), which is impossible since $J\notin \cT_a$ for all $a\in H^{2l} (\T^{2n};\Q)$, $a\neq 0$. Therefore $(\T^{2n}, J)$ has no complex subvarieties of dimension $2l$, $l=1,\ldots,n-1$, meaning that $J$ is Campana-simple.

Thus, the set $\Lin^s\cap \Lincmpt (\omega)$ of Campana-simple complex structures lying in $\Lincmpt (\omega)$ contains a subset which is connected and dense in $\Lincmpt (\omega)$. Hence, $\Lin^s\cap \Lincmpt (\omega)$ is also connected and dense in $\Lincmpt (\omega)$.

This finishes the proof.
\endproof


\hfill


Since $\Lin$ is connected, there is a connected component of $\compK (\T^{2n})$ containing it -- denote this connected component by $\cC_l$.

Let $x=0\in \T^{2n}=\R^{2n}/\Z^{2n}$. Let $\Diff (\T^{2n},x)$ be a subgroup of $\Diff (\T^{2n})$ formed by the diffeomorphisms of $\T^{2n}$ fixing $x$. Let $\Diff_0 (\T^{2n},x)$ be the identity component of $\Diff (\T^{2n},x)$. The group $\Diff_0 (\T^{2n},x)$ lies in $\Diff_0 (\T^{2n})$ and is, in fact, the identity component of the subgroup of
$\Diff_0 (\T^{2n})$ formed by the elements of $\Diff_0 (\T^{2n})$ fixing $x$.

Let $\Diff_H (\T^{2n},x) := \Diff_H (\T^{2n})\cap \Diff (\T^{2n},x)$. Clearly, $\Diff_0 (\T^{2n},x)\subset \Diff_H (\T^{2n},x)$.


\hfill


\proposition
\label{_any-complex-structure-in-cC-l-can-be-mapped-to-a-linear-one-by-unique-diffm-in-DiffH_Proposition_}

Let $I\in \compK (\T^{2n})$.

Then the following claims hold:

\smallskip
\noindent
A. There exist unique $J_I\in\Lin$ and $\phi_I\in\Diff_H (\T^{2n},x)$ so that $J_I = \phi_I^* I$.

\smallskip
\noindent
B. The complex structure $J_I$ and the diffeomorphism $\phi_I$ depend
continuously on $I$.

\smallskip
\noindent
C. If $I\in\cC_l$, then $\phi_I\in\Diff_0 (\T^{2n},x)$.

\smallskip
\noindent
D. If $I\in\cC_l$ and $\psi^* I \in \Lin$ for some $\psi\in\Diff_0 (\T^{2n})$, then $\psi^* I = I$.


\hfill


\noindent
{\bf Proof of \ref{_any-complex-structure-in-cC-l-can-be-mapped-to-a-linear-one-by-unique-diffm-in-DiffH_Proposition_}:}

Let us prove part A of the proposition.

The proof of the existence of $\phi_I$ and $J_I$ goes as in the proof of \cite[Prop. 6.1]{_EV-JTA_}.

Namely, recall the construction of the Albanese map: Let $\Omega_I^1 (\T^{2n})$ be the space of holomorphic complex-valued holomorphic 1-forms on the complex
manifold $(\T^{2n},I)$. It is a complex vector space of dimension $n$ that can be identified with $H_I^{1,0} (\T^{2n};\C)$. Pick a basis
$\varsigma_1,\ldots,\varsigma_n$ of $\Omega_I^1 (\T^{2n})$. The Albanese map $\Alb_I: \T^{2n}\to \C^n/\Gamma_I$
is defined by the map $y\mapsto \left(\int_\gamma \varsigma_1,\ldots, \int_\gamma \varsigma_n\right)$, where $\gamma$ is an arbitrary smooth path from $x=0$ to $y$ and $\Gamma_I$ is a lattice in $\C^n$ formed by the vectors $\left(\int_\gamma \varsigma_1,\ldots, \int_\gamma \varsigma_n\right)$
for all the loops $\gamma$. The Albanese map $\Alb_I$ is a biholomorphism between $(\T^{2n},I)$ and the complex torus $\C^n/\Gamma_I$ with
the complex structure induced by the standard complex structure on $\C^n$ -- see e.g. \cite[Vol. 1, Def. 12.10 and Thm.
12.15]{_Voisin-book_}.

Denote by $\pi_I: \C^n\to \C^n/\Gamma$ the natural projection
and let $F_I:\R^{2n}\to\C^n$ be an $\R$-linear isomorphism of
vector spaces mapping $\Z^{2n}$ to $\Gamma_I$. Then $F_I$ covers a
diffeomorphism $\barf_I: T^{2n} = \R^{2n}/\Z^{2n}\to \C^n/\Gamma_I$,
i.e., $\barf_I\circ \pi_I = \pi_I\circ F_I$.
One easily sees that $F_I$ can be chosen in such a way that $\barf_I$ induces the same map $H^* (\R^{2n}/\Z^{2n})\to H^* (\C^n/\Gamma)$
as $\Alb_I$.

Clearly, $\barf_I^* I$ is a linear complex structure on
$T^{2n}=\R^{2n}/\Z^{2n}$. Therefore $\phi_I:=\barf_I^{-1}\circ \Alb_I$ is a
diffeomorphism mapping $I$ to a linear complex structure $J_I$ on
$T^{2n}$, acting as identity on $H^* (T^{2n})$ and sending $x=0$ to itself.
Thus, $\phi_I\in\Diff_H (\T^{2n},x)$. This proves the existence of $J_I$ and $\phi_I$.

Let us prove the uniqueness of $J_I$ and $\phi_I$. It suffices to show that if $J_0,J_1\in \Lin$ and $\phi^* J_0 = J_1$ for some $\phi\in\Diff_H (\T^{2n},x)$, then $J_0 = J_1$ and $\phi=Id$.

Denote the lifts of $J_0,J_1$ to $\R^{2n}$ respectively by $\bJ_0$, $\bJ_1$.
The diffeomorphism $\phi: \T^{2n}\to \T^{2n}$ lifts to a diffeomorphism $\bphi: \R^{2n}\to \R^{2n}$ such that $\bphi^* \bJ_0 = \bJ_1$ and for each $v\in\Z^{2n}$ there exists $v'\in\Z^{2n}$ so that $\bphi (x + v) = \bphi (x) + v'$ for all $x\in\Z^{2n}$.

The group $GL^+ (2n,\R)$ acts transitively on $\Lin$ (see e.g. \cite{_McD-Sal-3_}, Prop. 2.5.2). Therefore there exists
$F\in GL^+ (2n,\R)$
so that $F^* \bJ_1 = \bJ_0$. Then $\Gamma:= F(\Z^{2n})$ is a lattice in $\R^{2n}$ and $F\bphi: \R^{2n}\to\R^{2n}$ is a diffeomorphism preserving $\hJ_0$ so that for each $v\in\Z^{2n}$ there exists $v''\in\Gamma$ so that $F\bphi (x + v) = F\bphi (x) + v''$ for all $x\in\Z^{2n}$.

Using an appropriate linear change of coordinates on $\R^{2n}$, we can assume without loss of generality that $\bJ_0$ is the standard linear complex structure on $\R^{2n}$ coming from the identification $\R^{2n}\cong\C^n$
and thus $F\bphi: \C^n\to\C^n$ is a holomorphic map. The partial derivatives of this map with respect to the complex coordinates on $\C^n$ are $\Z^{2n}$-invariant and holomorphic, hence descend to holomorphic maps on the standard complex torus $\T^{2n}$, hence are constant. Therefore $F\bphi: \C^n\to\C^n$ is a complex-linear affine map and consequently $\bphi$ is of the form $\bphi (x) = Ax + c$ for some $A\in SL (2n,\Z)$, $c\in\R^{2n}$.
Since $\phi$ acts trivially on $H_* (\T^{2n})$, the matrix $A$ is just the identity and $\bphi$ is just a parallel translation by $c$. Hence, $\phi$ is a parallel translation too, and since $\phi$ preserves $x$, we have $\phi=Id$. Consequently $J_0 = J_1$.

This proves the uniqueness of $J_I$ and $\phi_I$ and finishes the proof of part A.

Part B follows from \ref{_dependence-of-space-of-holomorphic-1-forms-on-the-complex-structure_Proposition_}
and from the construction of $J_I$ and $\phi_I$ in part A.

Let us prove part C. Assume $I\in\cC_l$. Then it can be connected to some $J\in\Lin$ by a continuous path $\{ I_t\}\subset\cC_l$. Considering the path $\{ \phi_{I_t}\}\subset \Diff_H (\T^{2n}, x)$, which is
 continuous
 by part B, and noting that $\phi_J= Id$, we immediately obtain that $\phi_I\in\Diff_0 (\T^{2n}, x)$. This proves part C.

Let us prove part D. If $I\in\cC_l$ and $\psi^* I \in \Lin$ for some $\psi\in\Diff_0 (M)$, then precomposing $\psi$ with a parallel translation, which is a biholomorphism of $(\T^{2n}, I)$, we can assume without loss
of generality that $\psi (x) = x$. Thus, $\psi\in\Diff_H (\T^{2n},x)$. By part A, this means that $\psi^* I = J_I = I$.
This proves part D and finishes the proof of the proposition.
\endproof


\hfill


\proposition\label{_Compl-homeom-to-Teich-cCl_Proposition_}

The restriction of the projection $\cC_l\to\cC_l/\Diff_0 (\T^{2n}) =: \Teich_{\cC_l} (\T^{2n})$ to $\Lin$ is a homeomorphism $\Lin\to\cC_l/\Diff_0 (\T^{2n}) = \Teich_{\cC_l} (\T^{2n})$
that identifies $\Lincmpt (\omega)\subset \cC_l$ with $\cmptteich_{\cC_l} (\T^{2n},\omega)\subset \Teich_{\cC_l} (\T^{2n})$.


\hfill


\noindent
{\bf Proof of \ref{_Compl-homeom-to-Teich-cCl_Proposition_}:}

Parts A and C of \ref{_any-complex-structure-in-cC-l-can-be-mapped-to-a-linear-one-by-unique-diffm-in-DiffH_Proposition_}
imply that every orbit of the $\Diff_0 (\T^{2n})$-action on $\cC_l$ intersects $\Lin$. By part D of \ref{_any-complex-structure-in-cC-l-can-be-mapped-to-a-linear-one-by-unique-diffm-in-DiffH_Proposition_}, it intersects $\Lin$ at exactly one point.

Thus, the map $\Lin\to\cC_l/\Diff_0 (\T^{2n})$, associating to each $I\in\Lin$ its $\Diff_0 (\T^{2n})$-orbit, is well-defined and bijective. Clearly, it is also continuous.

By part B of \ref{_any-complex-structure-in-cC-l-can-be-mapped-to-a-linear-one-by-unique-diffm-in-DiffH_Proposition_},
the map $\cC_l\to \Lin$, $I\mapsto J_I$, is continuous. This means that the inverse of the map $\Lin\to\cC_l/\Diff_0 (\T^{2n})$ is continuous too.

Thus, the projection $\Lin\to\cC_l/\Diff_0 (\T^{2n})$ is a homeomorphism.
Clearly, the image of $\Lincmpt (\omega)$ under this homeomorphism lies in $\cmptteich_{\cC_l} (\T^{2n},\omega)$. Let us show that it coincides with $\cmptteich_{\cC_l} (\T^{2n},\omega)$, or, equivalently, that for any complex structure $I\in\cC_l$ compatible with $\omega$ we have $J_I\in \Lincmpt (\omega)$.
Indeed,
since $\omega$ is compatible with $I$, the cohomology class $[\omega]$ is K\"ahler with respect to $I$, hence also with respect to $J_I$, because $J_I$ obtained from $I$ by the action of $\Diff_0 (\T^{2n})$.
Now
\ref{_Kahler-classes-wrt-linear-cs-linear-sympl-forms_Proposition_}
implies
that $J_I\in \Lincmpt (\omega)$. Consequently, the homeomorphism $\Lin\to\cC_l/\Diff_0 (\T^{2n})$
identifies $\Lincmpt (\omega)$ with $\cmptteich_{\cC_l} (\T^{2n},\omega)$.
\endproof


\hfill


\proposition\label{_SympH-acts-transitively-on-set-of-compon-of-space-of-cplx-str-torus_Proposition_}

Let $\omega$ be a K\"ahler-type form on $\T^{2n}$.

Then $\SympH (\T^{2n},\omega)$ acts transitively on the set of connected components of $\compK (\T^{2n})$ compatible with $\omega$.


\hfill


\noindent
{\bf Proof of \ref{_SympH-acts-transitively-on-set-of-compon-of-space-of-cplx-str-torus_Proposition_}:}

By \cite[Prop. 6.1]{_EV-JTA_}, any K\"ahler-type form on $\T^{2n}$ can be mapped by a diffeomorphism of $\T^{2n}$ acting trivially on $H_* (\T^{2n})$ to a linear symplectic form. Thus, we can assume without loss of generality that our K\"ahler-type form $\omega$ is linear.

Let $I\in\compK (\T^{2n})$ be an arbitrary complex structure compatible with $\omega$. By \cite[Prop. 6.1]{_EV-JTA_}, there exists a diffeomorphism $\phi$ of $\T^{2n}$ acting trivially on $H_* (\T^{2n})$ so that the complex structure $\phi^* I$ is linear. The symplectic form $\phi^* \omega$ is then compatible with $\phi^* I$ and has the same cohomology class as $\omega$.
Thus, the cohomology class
$[\omega]$ is K\"ahler with respect to the linear complex structure $\phi^* I$.
Therefore, by \ref{_Kahler-classes-wrt-linear-cs-linear-sympl-forms_Proposition_}, the linear form
$\omega$ itself is K\"ahler with respect to $\phi^* I$.

We have obtained that $\phi^* \omega$ and $\omega$ are cohomologous and compatible with the same complex structure $\phi^* I$. Therefore the straight path in the space of 2-forms connecting these forms is formed by cohomologous symplectic forms. By Moser's theorem \cite{_Moser_}, this means that there exists $\psi\in\Diff_0 (\T^{2n})$ such that
$\psi^* \phi^* \omega = \omega$. Thus, $\phi\psi\in\SympH (\T^{2n},\omega)$
and moreover $(\phi\psi)^* I = \psi^* \phi^* I$ lies in the same connected component of $\compK (\T^{2n})$ as $\phi^* I$ -- i.e., in $\cC_l$.

We have proved that any complex structure $I\in\compK (\T^{2n})$ compatible with $\omega$ can be mapped by an element of $\SympH (\T^{2n},\omega)$ into a complex structure lying in $\cC_l$. This readily implies the proposition.
\endproof


\hfill


\proposition
\label{_irrational-sympl-form-on-torus-compatible-with-Campana-simple-complex-str_Proposition_}

Let $\omega$ be a K\"ahler-type symplectic form on $\T^{2n}$.

Then $\omega$ is compatible with a Campana-simple complex structure if and only if
$\omega$ is irrational.

In this case for any connected component $\cC_0$ of
$\compK (\T^{2n})$ compatible with $\omega$ the set
\break
$\cmptteich_{\cC_0} (\T^{2n},\omega)\cap \Teich^{s} (\T^{2n})$
is a dense connected subset of $\cmptteich_{\cC_0} (\T^{2n},\omega)$.


\hfill


\noindent
{\bf Proof of \ref{_irrational-sympl-form-on-torus-compatible-with-Campana-simple-complex-str_Proposition_}:}

Assume that $\omega$ is rational.
Then, by part 1 of
\ref{_Campana-simple-vs-projectivity_Remark_}, there are no Campana-simple complex structures compatible with $\omega$.

By \cite[Prop. 6.1]{_EV-JTA_}, any K\"ahler-type form on $\T^{2n}$ can be mapped by a diffeomorphism of $\T^{2n}$ to a linear symplectic form.
Thus, we can assume without loss of generality that our K\"ahler-type form $\omega$ is linear.
In view of \ref{_SympH-acts-transitively-on-set-of-compon-of-space-of-cplx-str-torus_Proposition_}, we can also assume without loss of generality
that $\cC_0 = \cC_l$. Now the proposition follows from
\ref{_torus-linear-Campana-simple-cs-compatible-with-omega_Proposition_},
\ref{_Compl-homeom-to-Teich-cCl_Proposition_}.
\endproof


\hfill


\subsection{Campana-simple complex structures on K3 surfaces}
\label{_Campana-simple-structures-K3_Subsection_}

Let $M$ be a smooth manifold underlying a K3 surface.

The space $\Teich (M)$ has the structure of a complex manifold (see e.g. \cite{_Catanese:moduli_}).
More
precisely, a combination of the fundamental theorems of
Kuranishi \cite{_Kuranishi-AnnMath-1962_} and Bogomolov-Tian-Todorov
\cite{_Bogomolov_,_Tian_,_Todorov_} shows that for any $J_0\in\compK (M)$ the Kuranishi
space of $(M,J_0)$ (the base of a certain universal complex-analytic local
deformation family of a complex structure defined in \cite{_Kuranishi-AnnMath-1962_}) is
a smooth complex manifold as long as the canonical bundle of
$(M,J_0)$ is trivial, which, of course, is true in the hyperk\"ahler
case, and in particular, for K3 surfaces. Moreover (see e.g. \cite{_Catanese:moduli_}), there exists a
homeomorphism between a neighborhood
$\ttU_{[J_0]}$
of $[J_0]$ in $\Teich (M)$
and an open subset of the Kuranishi space of $(M,J_0)$. This homeomorphism maps
each $[J]\in \ttU_{[J_0]}$ to the point
$x$ of the Kuranishi space
corresponding to a complex structure on $M$ representing $[J]$.
The pullback of
the universal complex-analytic local deformation family over the Kuranishi space
under the homeomorphism, together with a smooth trivialization of the universal fibration over a neighborhood of $x$ in the Kuranishi space, induces a complex-analytic deformation family $M\times
\ttU_{[J_0]}\to \ttU_{[J_0]}$.
The
smooth and the complex structures on $\ttU_{[J_0]}$ are the pullbacks of the
corresponding structures on the Kuranishi space.


\hfill


\remark
\label{_Teich-K3-not-Hausdorff_Remark_}

Note that the complex manifold $\Teich (M)$ is not Hausdorff (see \cite{_Atiyah-analyt-surf-ProcRoySoc1958_}, \cite[p.238]{_Burns-Rapoport_}, cf. \cite{_Salamon-K3_}).


\hfill


The description above also implies that for any $J\in\compK (M)$
any smooth path in $\Teich (M)$ starting at $[J]$ can be lifted to a smooth path in $\compK (M)$ starting at $J$.

Let $q: H^2 (M;\C)\times H^2 (M;\C)\to\C$ be the intersection form. Its restriction to $H^2 (M;\R)$ is a real-valued symmetric form of signature $(3,19)$.
A subspace of $H^2 (M;\R)$ is called positive if the restriction of $q$ to it is positive. We use $\bot$ to denote the orthogonal complement with respect to $q$ in $H^2 (M;\R)$.

Let $I\in \compK (M)$. Then
$$\dim_\C H_I^{2,0} (M;\C) = \dim_\C H_I^{0,2} (M;\C) = 1,\ \dim_\C H_I^{1,1} (M;\C) =20$$
(see e.g. \cite{_Huybrechts-K3-book_}). The spaces $H_I^{2,0} (M;\R)\oplus H_I^{0,2} (M;\R)$ and $H^{1,1}_I (M;\R)$ are orthogonal with respect to $q$.
If $\Omega$ is a holomorphic $2$-form spanning $H^{2,0}_I (M;\C)$, then $\left(H_I^{2,0} (M;\C)\oplus H_I^{0,2} (M;\C)\right)\cap H^2 (M;\R)$
is the $2$-dimensional plane spanned by $\Re \Omega$ and $\Im \Omega$. Denote this plane by $W_I$. The basis $\{ \Re \Omega, \Im \Omega \}$ defines an orientation of the plane $W_I$. We have $W_I^\bot = H_I^{1,1} (M;\C)\cap H^2 (M;\R)$.
The restriction of $q$ to $W_I$ is positive and the restriction of $q$ to its orthogonal complement $H_I^{1,1} (M;\R)$ in $H^2 (M;\R)$ has signature $(1,19)$. Consequently, the set $\{ a\in H_I^{1,1} (M;\R)\ |\ q(a,a)>0 \}$ has two connected components: the one that contains the K\"ahler cone of $(M,I)$ is called the {\bf positive cone of $(M,I)$} and will be denoted $\Pos (M,I)$, and the other one
is $-\Pos (M,I)$.

Let $Gr_{++} (H^2 (M;\R))$ be the Grassmanian of positive oriented $2$-dimen\-sio\-nal planes in
\break
$H^2 (M;\R)$.
There is a canonical diffeomorphism identifying $Gr_{++} (H^2 (M;\R))$
with an open subset $\{ l\in \mP H^2 (M;\C)\ |\ q(l,l)=0,\ q(l,\bl) >0 \}$
of the quadric $\{ l\in \mP H^2 (M;\C)\ |\ q(l,l)=0 \}$
(see e.g. \cite[Ch.6, Prop. 1.5]{_Huybrechts-K3-book_} and \cite[Claim 2.9]{_V-Duke_}).
Thus $Gr_{++} (H^2 (M;\R))$ has the structure of a connected complex manifold.

Define the {\bf period map} $\Per: \Teich (M)\to Gr_{++} (H^2 (M;\R))$ by
\[
\Per ([I]) := W_I,
\]
where $[I]\in\Teich (M)$ is the element of $\Teich (M)$ represented by $I$.

The global Torelli theorem for K3 surfaces (see \cite{_Huybrechts-K3-book_} and the references therein to the original papers)
implies the following claim.


\hfill


\theorem
\label{_K3-Torelli_Theorem_}

\noindent
A. There exists a subset $\cS\subset \Teich (M)$ so that for each connected component $\cC_0$ of $\compK (M)$
the intersection $\cS\cap \Teich_{\cC_0} (M)$ is dense and connected in $\Teich_{\cC_0} (M)$
and the map
\[
\Per: \Teich_{\cC_0} (M)\to Gr_{++} (H^2 (M;\R))
\]
is a surjective local biholomorphism which is injective on $\cS\cap \Teich_{\cC_0} (M)$.

\bigskip
\noindent
B. Assume that $I_0, I_1\in\compK (M)$, $[I_0], [I_1]\in\cS$ and $\Per ([I_0]) = \Per ([I_1])$.

Then there exists a unique diffeomorphism $\phi\in\Diff_H (M)$ such that $\phi^* I_1 = I_0$.
\endproof


\hfill


By a {\bf $(-2)$-class} we will mean an integral cohomology class $b\in H^2 (M;\Z)$ such that $q(b,b)=-2$.

For each $I\in \compK (M)$ denote by $S(M,I)$ the set of all $(-2)$-classes $b\in H^2 (M;\Z)$ that are of type $(1,1)$ with respect to $I$.

The K\"ahler cone of $(M,I)$ admits the following description.


\hfill


\theorem {\bf (see e.g. \cite[Sec. 8, Thm. 5.2]{_Huybrechts-K3-book_}, cf. \cite[Thm. 6.2]{_Am-Ver-IMRN2015_})}
\label{_K3-Kahler-cone_Theorem_}

The K\"ahler cone of $(M,I)$ is a connected component of the set
\[
\Pos (M,I)\setminus \cup_{b\in S(M,I)} \left( b^\bot\cap \Pos (M,I)\right).
\]
In particular, if $b^\bot\cap \Pos (M,I)=\emptyset$ for all $b\in S(M,I)$, then the K\"ahler cone of $(M,I)$ is just $\Pos (M,I)$.
\endproof


\hfill


The next result follows from the Calabi-Yau theorem (formerly the Calabi conjecture) \cite{_Yau-CPAM-1978_} combined with Bochner's technique -- see \cite{_Bea1_} (cf. e.g. \cite[Thm. 5.11]{_Joyce-lectures-book_}) for the proof of the existence and uniqueness; see \cite[Thm. 6.3]{_Fujiki-Schumacher-PublRIMS1990_} for the proof of the smooth dependence.


\hfill


\theorem
\label{_K3-Calabi-Yau_Theorem_}

Let $I\in \compK (M)$ and let $a\in H^2 (M;\R)$ be a K\"ahler class with respect to $I$.

Then there exists a unique $\eta\in \sympK (M)$ compatible with $I$ such that $[\eta]=a$ and the Riemannian metric $\eta(\cdot, I\cdot)$
is Ricci-flat. The form $\eta$ depends smoothly on the pair $(I,a)$ (as $I$ varies in a complex-analytic deformation family).
\endproof


\hfill


\remark

For any K\"ahler structure $(\eta, I)$ on $M$ such that the Riemannian metric $\eta(\cdot, I\cdot)$
is Ricci-flat,
the form $\eta$ and the complex structure $I$ can be included together in a hyperk\"ahler structure $(\omega_1:=\eta,\omega_2,\omega_3, I_1:=I,I_2,I_3)$.


\hfill


\proposition
\label{_K3-adjusting-cs-to-get-Ricci-flat-metric_Proposition_}

Let $I\in \compK (M)$ be compatible with $\omega$.

Then there exists $J\in \compK (M)$ isotopic to $I$,
so that the Riemannian metric $\omega (\cdot, J\cdot)$ is Ricci-flat.


\hfill


\noindent
{\bf Proof of \ref{_K3-adjusting-cs-to-get-Ricci-flat-metric_Proposition_}:}

By \ref{_K3-Calabi-Yau_Theorem_}, the complex structure $I$ is a part of a (unique) K\"ahler structure $(\eta,I)$ on $M$ for which the corresponding Riemannian metric $\eta (\cdot, I\cdot)$ is Ricci-flat, and moreover $[\eta]=[\omega]$. The symplectic forms $\omega$ and $\eta$ are cohomologous and compatible with the same $I$ and therefore, by Moser's theorem \cite{_Moser_}, there exists $\phi\in\Diff_0 (M)$ such that $\phi^* \eta = \omega$. Then $J:=\phi^* I$ is compatible with $\omega$ and the Riemannian metric $\omega (\cdot, J\cdot)$ is Ricci-flat.
\endproof


\hfill


Let $\omega$ be a K\"ahler-type symplectic form on $M$.

Define $\Dom_\omega\subset Gr_{++} \left(H^2 (M;\R)\right)$ as the set of positive oriented $2$-di\-men\-sio\-nal planes $W\subset H^2 (M;\R)$ such that $[\omega]\in W^\bot$
and $(\R [\omega]\oplus W)^\bot$ does not contain any $(-2)$-classes.


\hfill


\proposition
\label{_K3-Dom-omega_Proposition_}

The set $\Dom_\omega$ is a dense connected subset of the submanifold $Gr_{++} \left([\omega]^\bot\right)$ of $Gr_{++} \left(H^2 (M;\R)\right)$.


\hfill


\noindent
{\bf Proof of \ref{_K3-Dom-omega_Proposition_}:}

We have $q ([\omega],[\omega])>0$. The restriction of $q$ to $[\omega]^\bot$ is a quadratic form of signature $(2,19)$.
The space $Gr_{++} ([\omega]^\bot)$ of positive oriented $2$-dimensional planes in $[\omega]^\bot$ is a smooth connected manifold diffeomorphic to
$SO (2,19)/SO (2)\times SO (19)$.

If $a\in H^2 (M;\Z)$ is a $(-2)$-class, then $a$ and $[\omega]$ are linearly independent and the restriction of $q$ to $a^\bot\cap [\omega]^\bot$ has signature $(2,18)$. The space $Gr_{++} \left(a^\bot\cap [\omega]^\bot\right)$ of positive oriented $2$-dimensional planes in
$a^\bot\cap [\omega]^\bot$
is a closed submanifold of $Gr_{++} \left([\omega]^\bot\right)$ diffeomorphic to
$SO (2,18)/SO (2)\times SO (18)$ and thus has real codimension $2$ in $Gr_{++} \left([\omega]^\bot\right)$.

The set $\Dom_\omega$ is the complement in $Gr_{++} \left([\omega]^\bot\right)$ of the (at most) countable union of the closed codimension-$2$ submanifolds $Gr_{++} \left(a^\bot\cap [\omega]^\bot\right)$ for all $(-2)$-classes $a\in H^2 (M;\Z)$ such that $q(a,[\omega])=0$. Therefore $\Dom_\omega$ is a dense connected subset of the submanifold $Gr_{++} \left([\omega]^\bot\right)$ of $Gr_{++} \left(H^2 (M;\R)\right)$.
 (The density follows from Baire's theorem; for a proof of the connectivity see e.g. \cite[Lem. 4.10]{_V-Duke_}.)
\endproof


\hfill


\proposition
\label{_K3-complex-str-compatible-with-omega-description_Proposition_}

Assume that $\cC_0$ is a connected component of $\compK (M)$ compatible with the symplectic form $\omega$ (i.e., $\cC_0$ contains a complex structure compatible with $\omega$).

Then $\cmptteich_{\cC_0} (M,\omega) = \Per^{-1} (\Dom_\omega)\cap \Teich_{\cC_0} (M)$ and it is
a dense connected subset
of the closed submanifold $\Per^{-1} \left(Gr_{++} \left([\omega]^\bot\right)\right)\cap \Teich_{\cC_0} (M)$ of $\Teich_{\cC_0} (M)$.


\hfill


\noindent
{\bf Proof of \ref{_K3-complex-str-compatible-with-omega-description_Proposition_}:}

Let us show that $\cmptteich_{\cC_0} (M,\omega) \subset \Per^{-1} (\Dom_\omega)\cap \Teich_{\cC_0} (M)$.
Consider a complex structure $I\in \cC_0$ compatible with $\omega$. Then $[\omega]\in H_I^{1,1} (M;\R)$.
Therefore $W_I$ is orthogonal to $[\omega]$.
Moreover, $\R[\omega]\oplus W_I$ cannot contain any $(-2)$-class in its orthogonal complement: indeed, if $a\in H^2 (M;\Z)$ is
orthogonal to $(H_I^{2,0} (M;\C)\oplus H_I^{0,2} (M;\C))\cap H^* (M;\R)$, then $a\in H^{1,1}_I (M;\R)$. However, by \ref{_K3-Kahler-cone_Theorem_}, a $(-2)$-class in $H_I^{1,1} (M;\R)$ cannot be orthogonal to a K\"ahler class with respect to $I$. Thus, we have proved that $\Per ([I]) = W_I$ belongs to $\Dom_\omega$. This yields that $\cmptteich_{\cC_0} (M,\omega) \subset \Per^{-1} (\Dom_\omega)\cap \Teich_{\cC_0} (M)$.

Now let us show that $\cmptteich_{\cC_0} (M,\omega) \supset \Per^{-1} (\Dom_\omega)\cap \Teich_{\cC_0} (M)$.

By the hypothesis of the proposition, there exists a complex structure $I\in\cC_0$ compatible with $\omega$.
In view of \ref{_K3-adjusting-cs-to-get-Ricci-flat-metric_Proposition_}, we may assume, without loss of generality, that the metric $\omega (\cdot, I\cdot)$ is Ricci-flat.

Let $W$ be an arbitrary plane belonging to $\Dom_\omega$. By \ref{_K3-Dom-omega_Proposition_}, it can be connected with $W_I$
by a smooth path of planes in $\Dom_\omega$. Since $\Per$ is a local biholomorphism (by part A of \ref{_K3-Torelli_Theorem_}) and
by the above-mentioned path-lifting property of the $\compK (M)\to\Teich (M)$, this path can be lifted to a smooth path in $\cC_0$ -- i.e.,
there exists a smooth path
of complex structures $\{ I_t\}$, $0\leq t\leq 1$, in $\cC_0$ such that $I_0 = I$ and $W_{I_1} = W$ and $W_{I_t}\in \Dom_\omega$ for each $t\in [0,1]$.
Consequently, for each $t\in [0,1]$ the plane $W_{I_t}$ is orthogonal to $[\omega]$ and $(\R [\omega]\oplus W)^\bot$ does not contain any $(-2)$-classes. Therefore, by \ref{_K3-Kahler-cone_Theorem_}, for each $t\in [0,1]$ either $[\omega]$ or $-[\omega]$ lies in the K\"ahler cone of $(M,I_t)$.

We claim that it is always $[\omega]$, and not $-[\omega]$.

Indeed, consider the set $A$ of $t\in [0,1]$ for which it is $[\omega]$.
The set $A$ is non-empty since $0\in A$: the cohomology class $[\omega]$ is K\"ahler with respect to $I=I_0$.
The set $A$ is also open:
if for some $t_0\in [0,1]$ the class $[\omega]$ is K\"ahler with respect
to $I_{t_0}$, then for any $t\in [0,1]$ sufficiently close to $t_0$, the class $[\omega]_{I_t}^{1,1} = [\omega]$ is K\"ahler with
respect to $I_t$  (by the corollary of the Kodaira-Spencer stability theorem \cite[Thm. 5.6]{_EV-JTA_}, see Section~\ref{_deformations-of-cs-Hodge-decomposition_Subsection_}).
For the same reason the complement of $A$ in $[0,1]$ (i.e., the set of $t\in [0,1]$ for which $-[\omega]$ is K\"ahler with respect to $I_t$) is also open. Thus, $A$ is open, closed and non-empty, meaning that $A=[0,1]$, which proves the claim.

Thus, for each $t\in [0,1]$ the class $[\omega]$ is K\"ahler with respect to a K\"ahler-type complex structure $I_t$ depending smoothly on $t$
and for $t=0$ the metric $\omega (\cdot, I_0\cdot)$ is Ricci-flat. By \ref{_K3-Calabi-Yau_Theorem_}, there exists a smooth family of cohomologous symplectic forms $\omega_t$, $0\leq t\leq 1$, such that $\omega_0 = \omega$ and $\omega_t$ is compatible with $I_t$ (and, moreover, the Riemannian metric $\omega_t (\cdot, I_t\cdot)$ is Ricci-flat) for each $t\in [0,1]$. Now Moser's theorem \cite{_Moser_} implies  that there exists $\phi\in\Diff_0 (M)$ such that $\phi^* \omega_1 = \omega$.
Consequently, the complex structure $\phi^* I_1$ is compatible with $\omega$, meaning that $[I_1]\in \cmptteich_{\cC_0} (M,\omega)$.
Hence, $W=W_{I_1}$ lies in $\Per \left(\cmptteich_{\cC_0} (M,\omega)\right)$. This yields $\cmptteich_{\cC_0} (M,\omega) \supset \Per^{-1} (\Dom_\omega)\cap \Teich_{\cC_0} (M)$.

We have shown that $\cmptteich_{\cC_0} (M,\omega) = \Per^{-1} (\Dom_\omega)\cap \Teich_{\cC_0} (M)$. Since, by part A of \ref{_K3-Torelli_Theorem_}, $\Per: \Teich_{\cC_0} (M)\to Gr_{++} (H^2 (M;\R))$
is a surjective local biholomorphism which is injective on a dense connected subset of $\Teich_{\cC_0} (M)$ and since, by \ref{_K3-Dom-omega_Proposition_}, $\Dom_\omega$ is a dense connected manifold of a closed submanifold of $Gr_{++} \left(H^2 (M;\R)\right) = \Per \left(\Teich_{\cC_0} (M)\right)$,
we get that $\cmptteich_{\cC_0} (M,\omega)$ is a dense connected subset of the closed submanifold $\Per^{-1} \left(Gr_{++} \left([\omega]^\bot\right)\right)\cap \Teich_{\cC_0} (M)$ of $\Teich_{\cC_0} (M)$.
\endproof


\hfill


\proposition\label{_SympH-acts-transitively-on-set-of-compon-of-space-of-cplx-str-K3_Proposition_}

The group $\SympH (M,\omega)$ acts transitively on the set of connected components of $\compK (M)$ compatible with $\omega$.


\hfill


\noindent
{\bf Proof of \ref{_SympH-acts-transitively-on-set-of-compon-of-space-of-cplx-str-K3_Proposition_}:}

Assume $\cC_0$ and $\cC_1$ are connected components of $\compK (M)$ compatible with $\omega$. We need to show that there exists a symplectomorphism
in $\SympH (M,\omega)$ mapping $\cC_1$ into $\cC_0$.

By \ref{_K3-complex-str-compatible-with-omega-description_Proposition_},
\[
\Per \left(\cmptteich_{\cC_0} (M,\omega)\right) =
\Per \left(\cmptteich_{\cC_1} (M,\omega)\right) = \Dom_\omega.
\]
Pick $I_0\in\cmptteich_{\cC_0} (M,\omega)\cap \cS$, $I_1\in\cmptteich_{\cC_1} (M,\omega)\cap\cS$,
so that $\Per ([I_0]) = \Per ([I_1])$. By part B of \ref{_K3-Torelli_Theorem_}, there exists $\phi\in\Diff_H (M)$ such that $\phi^* I_1 = I_0$.
Then $\phi^* \omega$ and $\omega$ are cohomologous and compatible with the same complex structure $I_0$. Therefore the straight path in the space of 2-forms connecting these forms is formed by cohomologous symplectic forms. By Moser's theorem \cite{_Moser_}, this means that there exists $\psi\in\Diff_0 (M)$ such that
$\psi^* \phi^* \omega = \omega$. Thus, $\phi\psi\in\SympH (\T^{2n},\omega)$,
and moreover $(\phi\psi)^* I_1 = \psi^* \phi^* I_1$ lies in the same connected component of $\compK (M)$ as $I_0 = \phi^* I_1$ -- i.e., in $\cC_0$.

Thus, the symplectomorphism $\phi\psi\in\SympH (M,\omega)$ maps $\cC_1$ into $\cC_0$.
\endproof


\hfill


\proposition
\label{_irrational-sympl-form-on-K3-compatible-with-Campana-simple-complex-str_Proposition_}

The K\"ahler-type symplectic form $\omega$ is compatible with a Campana-simple complex structure on $M$ if and only if
$\omega$ is irrational.
Moreover, in this case for any connected component $\cC_0$ of
$\compK (M)$ compatible with $\omega$ the set
$\cmptteich_{\cC_0} (M,\omega)\cap \Teich^{s} (M)$
is a dense connected subset of $\cmptteich_{\cC_0} (M,\omega)$.


\hfill


\noindent
{\bf Proof of \ref{_irrational-sympl-form-on-K3-compatible-with-Campana-simple-complex-str_Proposition_}:}

Assume $\omega$ is rational.
Then, by part 1 of
\ref{_Campana-simple-vs-projectivity_Remark_}, there are no Campana-simple complex structures compatible with $\omega$.

Now assume $\omega$ is irrational. Let us show that
the set
of Campana-simple structures lying in $\cmptteich_{\cC_0} (M,\omega)$
is a dense connected subset of $\cmptteich_{\cC_0} (M,\omega)$.

Let $b\in H^2 (M;\Q)$, $b\neq 0$. Since $[\omega]\in H^2 (M;\R)$ is not a real multiple of a rational cohomology class (and, in particular, not a real multiple of $b$), $\dim_\R b^\bot\cap [\omega]^\bot = 20$.  Consequently, the space $Gr_{++} \left(b^\bot\cap [\omega]^\bot\right)$ of positive oriented $2$-dimensional planes in
$b^\bot\cap [\omega]^\bot$ has real codimension $2$ in $Gr_{++} \left([\omega]^\bot\right)$ if $q(b,b)<0$ and is empty otherwise.

As we have seen in the proof of \ref{_K3-Dom-omega_Proposition_}, $\Dom_\omega$ is the complement in $Gr_{++} \left([\omega]^\bot\right)$ of an (at most) countable union of closed codimension-$2$ submanifolds $Gr_{++} \left(a^\bot\cap [\omega]^\bot\right)$ for all $(-2)$-classes $a\in H^2 (M;\Z)$ such that $q(a,[\omega])=0$. Therefore the complement in $\Dom_\omega$ to the countable union of codimension-2 submanifolds $Gr_{++} \left(b^\bot\cap [\omega]^\bot\right)$, $b\in H^2 (M;\Q)$, $b\neq 0$, is a dense connected subset of $\Dom_\omega$. Denote the latter subset by $A$.

By part A of \ref{_K3-Torelli_Theorem_} and \ref{_K3-complex-str-compatible-with-omega-description_Proposition_}, $\cmptteich_{\cC_0} (M,\omega)\cap \Per^{-1} (A)$ is a dense connected subset of  $\cmptteich_{\cC_0} (M,\omega)$.
We claim that each $I$ such that $[I]\in \Per^{-1} (A)\cap \cmptteich_{\cC_0} (M,\omega)$ is Campana-simple. Indeed, if $(M,I)$ admits a proper positive-dimensional analytic subvariety -- i.e., a (possibly singular) complex curve -- then the cohomology class $b\in H^2(M;\Q)$, $b\neq 0$, Poincar\'e-dual to the fundamental class of the subvariety is of type $(1,1)$ with respect to $I$ (see e.g. \cite[Vol. 1, Sec. 11]{_Voisin-book_}).
Therefore the 2-dimensional plane $\Per ([I])\in A$ is orthogonal to $b$ in $H^2 (M;\R)$ which contradicts the definition of $A$.

Thus, $\cmptteich_{\cC_0} (M,\omega)\cap \Per^{-1} (A)$ is a dense connected subset of the set
$\cmptteich_{\cC_0} (M,\omega)\cap \Teich^{s} (M)$.
Hence, $\cmptteich_{\cC_0} (M,\omega)\cap \Teich^{s} (M)$
is a dense connected subset of $\cmptteich_{\cC_0} (M,\omega)$.
\endproof


\hfill



\subsection{K\"ahler-type embeddings of balls into tori and K3 surfaces -- the proofs}
\label{_proofs-of-main-results-on-embs-of-balls-in-tori-K3_Subsection_}


Now we are ready to prove \ref{_existence-connectedness-irrational-forms-tori-K3-surfaces_Theorem_}. For convenience
we restate it here.


\hfill


\vfil\eject

\theorem {\bf (\ref{_existence-connectedness-irrational-forms-tori-K3-surfaces_Theorem_})}
\label{_existence-connectedness-irrational-forms-tori-K3-surfaces-COPY_Theorem_}

Let $M$, $\dim_\R M = 2n$, be either $\T^{2n}$ or a smooth manifold (of real dimension $4$) underlying a complex K3 surface
and let $\omega$ be a K\"ahler-type symplectic form on $M$.

Assume that $\omega$ is irrational.

Let $k\in\Z_{>0}$, $r_1,\ldots,r_k>0$.

Then the following claims hold:

\bigskip
\noindent
A. K\"ahler-type embeddings of $\bigsqcup_{i=1}^k B^{2n} (r_i)$ into $(M,\omega)$ are unobstructed.

More precisely, assume that $\Vol \left(\bigsqcup_{i=1}^k B^{2n} (r_i),\omega_0\right) < \Vol (M,\omega)$ and $\omega$ is
 compatible with a Campana-simple complex structure $I$.

Then
there exists an $[I]$-K\"ahler-type embedding $\bigsqcup_{i=1}^k B^{2n} (r_i)\to (M,\omega)$.

\bigskip
\noindent
B.
The group $\SympH (M,\omega)$ acts transitively on the set of connected components of
$\compK (M)$ compatible with $\omega$.

\bigskip
\noindent
C.
Any two K\"ahler-type embeddings $\bigsqcup_{i=1}^k B^{2n} (r_i)\to (M,\omega)$ (if they exist!) lie in the same orbit of the
$\SympH (M,\omega)$-action. They lie in the same orbit of the $\Symp (M,\omega)\cap \Diff_0 (M)$-action
if and only if they favor a common connected component of $\compK (M)$. In the latter case there exists $[I]\in\cmptteich (M,\omega)$
such that both embeddings are of $[I]$-K\"ahler-type.


\hfill


\noindent
{\bf Proof of \ref{_existence-connectedness-irrational-forms-tori-K3-surfaces-COPY_Theorem_} (=\ref{_existence-connectedness-irrational-forms-tori-K3-surfaces_Theorem_}):}

Let us first consider the torus case.
By \ref{_SympH-acts-transitively-on-set-of-compon-of-space-of-cplx-str-torus_Proposition_}, in this case we may assume without loss of generality that $\cC_0 =\cC_l$.
Now part A of the theorem follows from part (I) of \ref{_Campana-simple-Kahler-type-existence-connectedness-gen-mfd_Theorem_} and \ref{_irrational-sympl-form-on-torus-compatible-with-Campana-simple-complex-str_Proposition_}.
Part B of the theorem is exactly \ref{_irrational-sympl-form-on-torus-compatible-with-Campana-simple-complex-str_Proposition_}.
Part C of the theorem follows from part B and from part (II) of \ref{_Campana-simple-Kahler-type-existence-connectedness-gen-mfd_Theorem_} and \ref{_irrational-sympl-form-on-torus-compatible-with-Campana-simple-complex-str_Proposition_}.
This finishes the proof in the torus case.

Let us now prove the result in the K3 case.
Here part A of the theorem follows from part (I) of \ref{_Campana-simple-Kahler-type-existence-connectedness-gen-mfd_Theorem_} and
\ref{_irrational-sympl-form-on-K3-compatible-with-Campana-simple-complex-str_Proposition_}.
Part B is exactly \ref{_irrational-sympl-form-on-K3-compatible-with-Campana-simple-complex-str_Proposition_}.
Part C of the theorem follows from part B and from part (II) of \ref{_Campana-simple-Kahler-type-existence-connectedness-gen-mfd_Theorem_}.
This finishes the proof in the K3 case.
\endproof


\hfill



\section{Tame
embeddings of arbitrary domains}
\label{_arb-shapes-pfs_Section_}


In this section we will prove the results on
tame
embeddings into tori and K3 surfaces stated in Section~\ref{_main-results-tori-K3_Subsection_}.

Let $M$, $\dim_\R M = 2n$, be a closed connected oriented manifold admitting K\"ahler structures (compatible with the orientation).

Denote by $\Omega^m (M)$ the space of smooth differential $m$-forms on $M$, $m=0,1,\ldots,2n$.
Let $\symp (M)$ denote the space of symplectic structures on $M$
and $\comp (M)$ the space of complex structures on $M$.
We equip these spaces with the $C^1$-topologies.

Let $W_i\subset \R^{2n}$, $i=1,\ldots,k$,  be compact domains  with piecewise-smooth boundary, so that their interiors contain the origin and $H^2 (W_i;\R) = 0$.
Denote
\[
\bfW:= \bigsqcup_{i=1}^k W_i.
\]


\hfill


\proposition
\label{_tame-embed-perturbed_Proposition_}

Assume $\omega$ is a symplectic form on $M$
and
$f: \bfW\to (M,\omega)$ is a
tame
embedding holomorphic with respect to a complex structure $I$ on $M$ tamed by $\omega$.
Let $\cV$ be a sufficiently small neighborhood of $I$ in $\comp (M)$ so that each $I'\in\cV$ is tamed by $\omega$.

Then there exists a neighborhood $\cU= \cU (\cV)$ of $\omega$ in $\symp (M)$ such that for each $\omega'\in\cU$
there exists a tame
embedding
$f': \bfW\to (M,\omega')$
holomorphic with respect to a complex structure $I'\in\cV$ which is tamed by $\omega'$ and isotopic to $I$.


\hfill


\noindent
{\bf Proof of \ref{_tame-embed-perturbed_Proposition_}:}

The following lemma is standard -- it is proved by an easy modification of the proof of a similar result in \cite[Lem. 9.3]{_EV-JTA_}.


\hfill


\lemma
\label{_small_primitive_for_a_small_exact_form_Lemma_}

Let $V\subset M$, $\dim_\R V = 2n$, be a compact
submanifold with a
piecewise-smooth boundary.
Let $\cX$ be a neighborhood of zero in the space $\Omega^1 (M)$.

Then there exists a neighborhood $\cX'$ of zero in the space $\Omega^2 (M)$,
so that for each $\eta\in\cX'$, which is exact on $V$, one can
choose $\sigma\in \cX$ such that $d\sigma = \eta$ on $V$.
\endproof


\hfill


Assume that the tame embedding $f:=\bigsqcup_{i=1}^k f_i : \bfW\to (M,\omega)$ extends (as a
tame embedding) to an embedding of an open neighborhood of
$\bfU=\bigsqcup_{i=1}^k U_i$, where $U_i$, $i=1,\ldots,k$, is a compact domain with a piecewise-smooth boundary containing $ W_i$ in its interior and  homotopy equivalent to $W_i$, so that
$f$ is holomorphic (on a neighborhood of $\bfU$) with respect to $I$.

Choose a neighborhood $\cU_1$ of $\omega$ in $\symp (M)$ so that for each $\omega'\in \cU_1$ the following conditions are satisfied:

\begin{itemize}

\item{} For each $t\in [0,1]$ the form $\omega + t(\omega-\omega')$ is symplectic.

\item{} $\omega'$ is tamed by any $I'\in\cV$.

\end{itemize}

\smallskip
Choose a neighborhood $\cW$ of $Id$ in the group $\Diff_0 (M)$ so that $\phi_* I \in \cV$
for each $\phi\in \cW$.

Choose a neighborhood $\cY$ of $0$ in the space of time-dependent smooth vector fields on $M$, equipped with the $C^1$-topology, so that for each vector field $v\in \cY$ the time $[0,1]$-flow of $v$
generates a diffeomorphism lying in $\cW$ and such that
for each $t\in [0,1]$ the time-$[0,t]$ flow of $v$ maps each $f(W_i)$ into $f(\Int U_i)$, $i=1,\ldots,k$.

Since for all $i=1,\ldots,k$ we have $H^2 (W_i;\R) = H^2 (U_i;\R) =0$, for each $\omega'\in \cU_1$ the restriction of the form $\omega - \omega'$ to $f(\bfU)$ is exact. Each $\sigma\in \Omega^1 (M)$
satisfying $d\sigma = \omega - \omega'$ defines a smooth time-dependent vector field $v$ on $M$ by the formula
\[
\left(\omega + t(\omega-\omega')\right) (v_t, \cdot) = \sigma (\cdot),\ t\in [0,1].
\]
In view of \ref{_small_primitive_for_a_small_exact_form_Lemma_}, applied to $f(\bfU)$, there exists a neighborhood $\cU_2\subset\cU_1$ of $\omega$ in $\sympK (M)$ so that for each $\omega'\in \cU_3$ one can choose $\sigma$ as above so that the corresponding $v$ lies in $\cY$.

Let $\omega'\in \cU_2$ and let $\{ \phi_t\}\subset \Diff_0 (M)$, $t\in [0,1]$, be the time-$[0,t]$ flow of $v\in\cY$. Then there exists an open neighborhood $Z_i\subset \Int\, U_i$ of each $W_i$, $i=1,\ldots,k$, so that $\phi_t (Z_i)\subset f(\Int\, U_i)$ for all $t\in [0,1]$, $i=1,\ldots,k$.
Moser's argument\footnote{Let us note that in the proof of \cite[Prop. 9.1]{_EV-JTA_} we applied Moser's argument to symplectic forms on the {\sl domain of $f$}, while here we will apply it to symplectic forms on the {\sl target} of $f$ -- i.e., on $M$.} \cite{_Moser_} shows that
\[
\phi_1^* \omega'|_{f(Z_i)} = \omega|_{f(Z_i)},\  i=1,\ldots,k.
\]
We claim that
\[
g:=\phi_1\circ f = \bigsqcup_{i=1}^k \phi_1\circ f_i: \bfW\to (M,\omega')
\]
is a tame embedding.

Indeed, it is clearly a smooth embedding extending to a neighborhood $Z:=\bigsqcup_{i=1}^k Z_i$ of $\bfW$ so that $g(Z)\subset \Int f(\bfU)$. It is also symplectic on $\bfZ$:
\[
g^* \omega' = (\phi_1\circ f)^* \omega' = f^* \phi_1^* \omega' = f^* \omega = \omega_0,
\]
where $\omega_0$ is the standard symplectic form on $\R^{2n}$.
Note that since the embedding $f$ is holomorphic with respect to the complex structure $I$, the embedding $g$ is holomorphic with respect to the complex structure $J:=(\phi_1)_* I$.

By our construction, $\phi_1\in\cW$, hence $J\in \cV$. Since $\omega'\in \cU_2\subset \cU_1$, the complex structure
$J$ is tamed by $\omega'$.
This shows that for $\cU:= \cU_2$ and any $\omega'\in\cU$ there exists a tame embedding
$g: \bfW\to (M,\omega')$
holomorphic with respect to the complex structure
$J$ tamed by $\omega'$ and isotopic to $I$.
\endproof


\hfill


\proposition
\label{_L-varepsilon-tame-embed-perturbed_Proposition_}

Let
$\varepsilon >0$.

Assume $\omega$ is a K\"ahler-type symplectic form on $M$
and
$f: \bfW\to (M,\omega)$ is an
$\varepsilon$-tame
embedding holomorphic with respect to a K\"ahler-type complex structure $I$ on $M$ tamed by $\omega$.

Then there exists a neighborhood $\cU$ of $\omega$ in $\sympK (M)$ such that for each $\omega'\in\cU$
there exists an $\varepsilon$-tame
embedding
$f': \bfW\to (M,\omega')$
holomorphic with respect to a complex structure which is tamed by $\omega'$ and isotopic to $I$.


\hfill


\noindent
{\bf Proof of \ref{_L-varepsilon-tame-embed-perturbed_Proposition_}:}

Choose a neighborhood $\cV$ of $I$ in $\compK (M)$.

Applying \ref{_tame-embed-perturbed_Proposition_} to $\bfW$, we get a neighborhood $\cU$ of $\omega$ in $\symp (M)$
such that for each $\omega'\in\cU$ there exists a tame
embedding $f': \bfW\to (M,\omega')$ holomorphic with respect to a complex structure $I'\in\cV$
tamed by $\omega'$ so that $I' = \phi_* I\in\cV$ for some $\phi\in\Diff_0 (M)$. In particular, for all $p,q\in\Z_{\geq 0}$,
\[
[\omega']_{I'}^{p,q} = [\omega']_{\phi_* I}^{p,q} = [\omega']_I^{p,q}.
\]

The cohomology class $[\omega']$ depends continuously on $\omega'$ and the Hodge decomposition (with respect to $I$) of a cohomology class
depends continuously on that class. Also, the K\"ahler cone of $(M,I)$ is an open subset of $H^{1,1} (M;\R)$.
In view of these facts, we can choose, if necessary, a smaller neighborhood $\cU$ and assume without loss of generality that for any $\omega'\in\cU$:

\begin{itemize}

\item{} The class
$[\omega']_{I'}^{1,1} = [\omega']_I^{1,1}$ is K\"ahler,

\item{} For $a:= [\omega]_{I'}^{2,0} + [\omega]_{I'}^{0,2} = [\omega']_I^{2,0} + [\omega']_I^{0,2}$, we have $\left| \langle a^n, [M]\rangle \right| < \varepsilon$.

\end{itemize}

This yields that $f': \bfW\to (M,\omega')$ is an $\varepsilon$-tame symplectic embedding holomorphic with respect to the complex structure $I'$ which is tamed by $\omega'$ and isotopic to $I$.
\endproof


\hfill


\corollary

Assume $\omega$ is a K\"ahler-type symplectic form on $M$
and
$f: \bfW\to (M,\omega)$ is a K\"ahler-type embedding
holomorphic with respect to a K\"ahler-type complex structure $I$ on $M$ compatible $\omega$.

Then, for any $\varepsilon >0$,
there exists a neighborhood $\cU$ of $\omega$ in $\sympK (M)$ such that for each $\omega'\in\cU$
there exists an $\varepsilon$-tame
embedding
$f': \bfW\to (M,\omega')$
holomorphic with respect to a complex structure which is tamed by $\omega'$ and isotopic to $I$.
\endproof


\hfill


\proposition
\label{_semicontinuity_Proposition_}

Let $\varepsilon >0$.

Then the function $\omega\mapsto \nu_{T,\varepsilon} (M,\omega, \bfW)$ on $\symp (M)$ is lower
semicontinuous with respect to the $C^1$-topology on $\symp (M)$.


\hfill


\noindent
{\bf Proof of \ref{_semicontinuity_Proposition_}:}

Consider an arbitrary
$\omega\in \symp (M)$.
Let us prove that $\nu_{T,\varepsilon}$ is lower semicontinuous
at $\omega$.

Let $\delta,\lambda >0$.
Assume $f: \lambda\bfW\to (M,\omega)$ is an
$\varepsilon$-tame
embedding holomorphic with respect to a K\"ahler-type complex structure $I$ on $M$ tamed by $\omega$.

By \ref{_L-varepsilon-tame-embed-perturbed_Proposition_}, there exists a neighborhood $\cU$ of $\omega$ in $\sympK (M)$
such that for each $\omega'\in\cU$
there exists an $\varepsilon$-tame
embedding
$f': \bfW\to (M,\omega')$.
Choosing, if necessary, a smaller neighborhood $\cU$, we may assume without loss of generality that
\[
\frac{\Vol (\lambda\bfW,\omega_0)}{\Vol
(M,\omega')} > \frac{\Vol (\lambda\bfW,\omega_0)}{\Vol (M,\omega)} -
\delta,\ \forall \omega'\in \cU.
\]
This shows that $\nu_{T,\varepsilon}$ is lower semicontinuous
at $\omega$.

Thus, $\nu_{T,\varepsilon}$ is lower semicontinuous
at all $\omega\in\symp (M)$, meaning that it is a lower semicontinuous function.
\endproof


\hfill


Let us now prove \ref{_packing-by-arb-shapes-main-tori-K3_Theorem_}. For convenience, we restate it here.


\hfill


\theorem {\bf (= \ref{_packing-by-arb-shapes-main-tori-K3_Theorem_})}
\label{_packing-by-arb-shapes-main-tori-K3-COPY_Theorem_}

Let $M$ be either a torus $\T^{2n}$ or a smooth manifold underlying a K3 surface.

Let $\omega_1$, $\omega_2$ be
K\"ahler-type forms on $M$ so that
$\int_M \omega_1^n = \int_M \omega_2^n >0$ and the
forms $\omega_1$, $\omega_2$ are irrational.

Let $W_i\subset \R^{2n}$, $i=1,\ldots,k$,  be compact domains with piecewise-smooth boundary whose interiors contain the origin.
Assume that $H^2 (W_i;\R) = 0$ for all $i$. Set $\bfW:= \bigsqcup_{i=1}^k W_i$.

Then, for
any $\varepsilon >0$,
\[
\nu_{T,\varepsilon} (M,\omega_1,\bfW) = \nu_{T,\varepsilon} (M,\omega_2,\bfW).
\]


\hfill


\noindent
{\bf Proof of \ref{_packing-by-arb-shapes-main-tori-K3-COPY_Theorem_} (= \ref{_packing-by-arb-shapes-main-tori-K3_Theorem_}):}

Without loss of generality, it is enough to proof the theorem in case $\Vol (M,\omega_1) = \Vol (M,\omega_2) = 1$.

Denote by $\sympK_1 (M)$ the space of K\"ahler forms on $M$ of total volume $1$. Equip
$\sympK_1 (M)$ with the $C^\infty$-topology. The group $\Diff^+ (M)$ of
orientation-preserving diffeomorphisms of $M$ acts on $\sympK_1 (M)$ and the
function $\omega\mapsto \nu_{T,\varepsilon} (M,\omega, \bfW)$ is clearly invariant under the action.

By \cite[Thm. 9.2]{_EV-JTA_} in the torus case and \cite[Thm. 1.7]{_EV-K3_} in the K3 case, the $\Diff^+ (M)$-orbits of $\omega_1$ and $\omega_2$ are dense in $\sympK_1 (M)$. The function $\nu_{T,\varepsilon} (M,\cdot, \bfW)$ is constant
on both orbits. Together with the lower semicontinuity of $\nu_{T,\varepsilon} (M,\cdot, \bfW)$ proved in \ref{_semicontinuity_Proposition_},
this implies that
\[
\nu_{T,\varepsilon} (M,\omega_1, \bfW) \leq \nu_{T,\varepsilon} (M,\omega_2, \bfW),
\]
and, by the same token,
\[
\nu_{T,\varepsilon} (M,\omega_1, \bfW) \geq \nu_{T,\varepsilon} (M,\omega_2, \bfW).
\]
Thus,
\[
\nu_{T,\varepsilon} (M,\omega_1, \bfW) = \nu_{T,\varepsilon} (M,\omega_2, \bfW),
\]
as required. This proves the theorem.
\endproof


\hfill


\corollary
\label{_packing-by-arb-shapes-tori-K3-from-Kahler-type-to-tame_Corollary_}

Let $M$ be either a torus $\T^{2n}$ or a smooth manifold underlying a K3 surface.
Let $W_i\subset \R^{2n}$, $i=1,\ldots,k$,  be compact domains with piecewise-smooth boundary whose interiors contain the origin.

Assume that $H^2 (W_i;\R) = 0$ for all $i$. Set $\bfW:= \bigsqcup_{i=1}^k W_i$.
Let $\varepsilon >0$.

Finally, assume that for any $\delta>0$ there exists an irrational K\"ahler-type form
$\omega_\delta$ on $M$ and $\lambda>0$
so that $\lambda\bfW$ admits a K\"ahler-type embedding into $(M,\omega_\delta)$
and $\Vol (\lambda\bfW,\omega_0)/\Vol (M,\omega_\delta) \geq 1-\delta$.

Then the following claims hold:

\bigskip
\noindent
A. If $\omega_\delta = \omega$ for all $\delta>0$ and $\omega$ is irrational, then there exists a dense $\Diff^+ (M)$-orbit -- the $\Diff^+ (M)$-orbit of $\omega$ -- in the space of K\"ahler-type symplectic forms of volume $\Vol (M,\omega)$ on $M$
such that for any $\omega'$ in this orbit we have
$\nu_K (M,\omega',\bfW) = 1$.

\bigskip
\noindent
B. There exists a $\Diff^+ (M)$-invariant open dense set of K\"ahler-type symplectic forms on $M$, depending on $\bfW$ and $\varepsilon$ and containing, in particular, all irrational K\"ahler-type symplectic forms on $M$,
so that for each $\omega'$ in this set
$\nu_{T,\varepsilon} (M,\omega',\bfW) = 1$.


\hfill


\noindent
{\bf Proof of \ref{_packing-by-arb-shapes-tori-K3-from-Kahler-type-to-tame_Corollary_}:}

Let us prove part A.

The hypothesis of part A implies that $\nu_K (M,\omega,\bfW) = 1$.
Since $\omega$ is irrational, \cite[Thm. 9.2]{_EV-JTA_} in the torus case and \cite[Thm. 1.7]{_EV-K3_} in the K3 case yield that the $\Diff^+ (M)$-orbit of $\omega$ is dense in the space of K\"ahler-type symplectic forms of volume $\Vol (M,\omega)$ on $M$.
Since the function $\omega'\mapsto \nu_K (M,\omega',\bfW)$ is $\Diff^+ (M)$-invariant and $\nu_K (M,\omega,\bfW) = 1$,
we have
$\nu_K (M,\omega',\bfW) = 1$
for any $\omega'$ in the $\Diff^+ (M)$-orbit of $\omega$.

This yields part A of the corollary.

Let us prove part B.

Rescaling, if necessary, the forms $\omega_\delta$, we may assume without loss of generality that they all have the same volume -- say, volume $1$. Since the forms $\omega_\delta$ are all irrational, \ref{_packing-by-arb-shapes-main-tori-K3-COPY_Theorem_} implies then that $\nu_{T,\varepsilon} (M,\omega_\delta,\bfW)$ does not depend on $\delta$. On the other hand, the hypothesis of the corollary implies that $\nu_{T,\varepsilon} (M,\omega_\delta,\bfW) \geq 1 - \delta$ for all $\delta>0$. Hence, $\nu_{T,\varepsilon} (M,\omega_\delta,\bfW) = 1$ for all $\delta>0$. Using again \ref{_packing-by-arb-shapes-main-tori-K3-COPY_Theorem_} and rescaling the forms, we get that
$\nu_{T,\varepsilon} (M,\omega', \bfW) = 1$ for any irrational K\"ahler-type form $\omega'$ on $M$, of any volume.
By \ref{_L-varepsilon-tame-embed-perturbed_Proposition_}, we get $\nu_{T,\varepsilon} (M,\omega'', \bfW) = 1$
also for any K\"ahler-type form $\omega''$ lying in a sufficiently small neighborhood $U (\omega')$ of any irrational form $\omega'$ in $\sympK (M)$. Let $U$ be the $\Diff^+ (M)$-orbit of $\bigcup_{\omega'} U (\omega')$, where the union is taken over all irrational $\omega'\in \sympK (M)$.
Clearly, $U$ is a $\Diff^+ (M)$-invariant open dense subset of $\sympK (M)$, containing all the irrational forms in $\sympK (M)$,
and $\nu_{T,\varepsilon} (M,\omega', \bfW) = 1$ holds for any $\omega'\in U$.

This finishes the proof of part B and of the corollary.
\endproof


\hfill


\corollary {\bf (= \ref{_full-packing-of-tori-K3-by-balls_Corollary_})}
\label{_full-packing-of-tori-K3-by-balls-COPY_Corollary_}

Assume:

\begin{itemize}

\item{} $M$ is either $\T^{2n}$ or a smooth manifold underlying a K3 surface.

\item{} $\bfW:= \bigsqcup_{i=1}^k B^{2n} (r_i)$ is a disjoint union of $k$ (possibly different) balls.

\item{}
$\varepsilon >0$.

\end{itemize}

Then the following claims hold:

\bigskip
\noindent
A.
For any irrational K\"ahler-type symplectic form $\omega$ on $M$ we have
\[
\nu_K (M,\omega,\bfW) = 1,
\]
meaning that K\"ahler-type embeddings of $\lambda\bfW$ into $(M,\omega)$ are
unobstructed.

\bigskip
\noindent
B. There exists a $\Diff^+ (M)$-invariant open dense set of K\"ahler-type symplectic forms on $M$, depending on $\bfW$ and $\varepsilon$ and containing, in particular, all irrational K\"ahler-type symplectic forms on $M$,
so that for each $\omega$ in this set
$\nu_{T,\varepsilon} (M,\omega,\bfW) = 1$ -- meaning that $\varepsilon$-tame embeddings of $\lambda\bfW$ into $(M,\omega)$ are
unobstructed.


\hfill


\noindent
{\bf Proof of \ref{_full-packing-of-tori-K3-by-balls-COPY_Corollary_}
(=\ref{_full-packing-of-tori-K3-by-balls_Corollary_}):}

Assume without loss of generality that
$n>1$. (In the case $n=1$ the result is known -- see \ref{_torus-of-dim-2_Remark_}).

Part A of the corollary is equivalent to the existence claim proved in \ref{_existence-connectedness-irrational-forms-tori-K3-surfaces_Theorem_}.

Part B of the corollary follows then from part B of \ref{_packing-by-arb-shapes-tori-K3-from-Kahler-type-to-tame_Corollary_}.
\endproof


\hfill


Now let us prove \ref{_full-packing-of-tori-by-polydisks-parallelepipeds_Corollary_}. We recall it here.


\hfill


\corollary {\bf (= \ref{_full-packing-of-tori-by-polydisks-parallelepipeds_Corollary_})}
\label{_full-packing-of-tori-by-polydisks-parallelepipeds-COPY_Corollary_}

Assume:

\begin{itemize}

\item{} $M=\T^{2n}$.

\item{} $\bfW:= \bigsqcup_{i=1}^k W_i$ is either a disjoint union of $k$ identical copies of a $2n$-dimensional polydisk
\[
B^{2n_1} (R_1)\times\ldots\times B^{2n_l} (R_l),\ n_1+\ldots+n_l = n,\ R_1,\ldots, R_l>0, l>1,
\]
or a disjoint union of $k$ identical copies of a parallelepiped
\[
P (e_1,\ldots,e_{2n}) := \left\{ \sum_{j=1}^{2n} s_j e_j, 0\leq s_j \leq 1, j=1,\ldots,2n \right\},
\]
where $e_1,\ldots,e_{2n}$ is a basis of the vector space $\R^{2n}$.

\item{} $\varepsilon >0$.

\end{itemize}

Then the following claims hold:

\bigskip
\noindent
A.
For any positive volume there exists a dense $\Diff^+ (\T^{2n})$-orbit (of an irrational K\"ahler-type symplectic form depending on $\bfW$) in the space of K\"ahler-type symplectic forms of that volume on $\T^{2n}$
such that for any $\omega'$ in this orbit we have
$\nu_K (\T^{2n},\omega',\bfW) = 1$ --
or, in other words, K\"ahler-type embeddings $\lambda\bfW\to (\T^{2n},\omega')$ are unobstructed.

\bigskip
\noindent
B. There exists a $\Diff^+ (\T^{2n})$-invariant open dense set of K\"ahler-type symplectic forms on $\T^{2n}$, depending on $\bfW$ and $\varepsilon$ and containing, in particular, all irrational K\"ahler-type symplectic forms on $\T^{2n}$,
so that for each $\omega'$ in this set
$\nu_{T,\varepsilon} (\T^{2n},\omega',\bfW) = 1$ --
or, in other words, $\varepsilon$-tame embeddings $\lambda\bfW\to (\T^{2n},\omega')$ are unobstructed.


\hfill


\noindent
{\bf Proof of \ref{_full-packing-of-tori-by-polydisks-parallelepipeds-COPY_Corollary_} (=\ref{_full-packing-of-tori-by-polydisks-parallelepipeds_Corollary_}):}

Assume without loss of generality that
$n>1$. (In the case $n=1$ the result is known -- see \ref{_torus-of-dim-2_Remark_}).

Let us consider the case where $\bfW$ is the disjoint union of $k$ identical copies of the polydisk
$B^{2n_1} (R_1)\times\ldots\times B^{2n_l} (R_l)$.
In this case the proof is the following straightforward modification of the proof of \cite[Cor. 3.3]{_EV-JTA_}.

Assume, without loss of generality, that
$\Vol (\bfW,\omega_0) = 1$.
For each $m\in\Z_{>0}$ denote by $\Omega_m$ the standard symplectic form on $\R^{2m}$.
For each $i=1,\ldots,l$ set
\[
v_i:= \Vol (B^{2n_i} (R_i), \Omega_{n_i}) = \pi^{n_i} n_i! R_i^{2n_i}.
\]
Note that
\[
\Vol (\bfW, \Omega_{2n}) = k\Vol (B^{2n_1} (R_1)\times \ldots B^{2n_l} (R_l),\Omega_{2n})
= kN\prod_{i=1}^l v_i,
\]
where
\[
N:=\frac{n!}{n_1!\cdot\ldots\cdot n_l!}.
\]
Since $\Vol (\bfW,\omega_0) = 1$, we have
\[
kN v_1\cdot \ldots \cdot v_l =1.
\]

For each $m\in\Z_{>0}$ and $w_1,\ldots, w_m \in \R_{>0}$, set
\[
\bbw:= (w_1,\ldots, w_m) \in \R^m,
\]
\[
v_\bbw := m! w_1\cdot\ldots \cdot w_m,
\]
and denote
by $\omega_\bbw$ the symplectic form
\[
\omega_\bbw := \sum_{i=1}^m w_i dp_i\wedge dq_i
\]
on the torus $\T^{2m} = \R^{2m}/\Z^{2m}$.
Note that
$\Vol (\T^{2m}, \omega_\bbw) = v_\bbw$.

Set
\[
k_1:=k, k_2=\ldots=k_l:=1,
\]
so that
\[
k = k_1\cdot\ldots\cdot k_l.
\]

One can choose $\bbw_i\in \R^{2n_i}$,
$i=1,\ldots,l$,  so that the following conditions
hold:

\smallskip
\noindent (i) $v_{{\bbw_i}} = k_i v_i$ for all $i=1,\ldots,l$.

\smallskip
\noindent (ii) The vector $\bbw := (\bbw_1,\ldots,\bbw_k) \in \R^n$ is
not proportional to a vector with rational coordinates.

\smallskip
\noindent (iii) For each $i=1,\ldots,l$ such that $n_i>1$ the vector $\bbw_i\in \R^{n_i}$ is
not proportional to a vector with rational coordinates.

\smallskip
Conditions (i) and (iii) can be achieved since, by our assumption, $n>1$, and since for any $m>1$ and $C>0$
the subset of the set
$$\{ (w_1,\ldots,w_m)\in\R^m\ |\ w_1,\ldots,w_m>0,
w_1\cdot\ldots \cdot w_m = C\}$$
formed by the vectors that are not
proportional to a vector with rational coordinates is dense.

Consider the symplectic form $\omega_{\bbw}$ on $\T^{2n}$ -- it is of K\"ahler type (since it is a linear symplectic form)
and irrational (because of condition (ii)).
Note that, by condition (i),
\[
\int_{\T^{2n}} \omega_{\bbw}^n = N \prod_{i=1}^l v_{{\bbw_i}}  = N \prod_{i=1}^l k_i v_i =
\]
\[
=
N k \prod_{i=1}^l k_i v_i =
\int_{\T^{2n}} \omega^n = 1.
\]
Also note that, by condition (ii), for any $0< \lambda < 1$ for all $i$
\[
\Vol (\T^{2n_i},\omega_{{\bbw_i}}) = v_{{\bbw_i}} = k_i v_i
> \lambda k_i v_i  = \Vol \bigg( \lambda \bigsqcup_{k_i} (B^{2n_i} (R_i), \Omega_{n_i}) \bigg),
\]
where $\bigsqcup_{k_i}$ denotes the disjoint union of $k_i$ equal copies
of the ball. By condition (iii), for each $i=1,\ldots,l$, such that $n_i>1$, the form $\omega_{{\bbw_i}}$ on $\T^{2n_i}$ is irrational. Therefore, by \ref{_existence-connectedness-irrational-forms-tori-K3-surfaces_Theorem_} and \ref{_torus-of-dim-2_Remark_}, there exists a
K\"ahler-type embedding
\[
f_i: \bigsqcup_{k_i} \left(\lambda B^{2n_i} (R_i),
\Omega_{n_i}\right) \to (\T^{2n_i}, \omega_{\bbw_i}).
\]
Accordingly, the direct product of all such embeddings $f_i$,
$i=1,\ldots,l$, is a K\"ahler-type embedding
\[
f: (\lambda\bfW, \Omega_{2n}) \to
(\T^{2n},\omega_{\bbw}),
\]
where $\lambda\bfW$ is the disjoint union
of $k_1\cdot\ldots\cdot k_l=k$
copies of $\lambda \left(B^{2n_1} (R_1)\times \ldots \times
B^{2n_l} (R_l)\right)$.
Thus, we have found an irrational K\"ahler-type form $\omega_{\bbw}$ of volume $1$ on $\T^{2n}$ such that
for any $0<\lambda<1$ there exists a K\"ahler-type embedding of $(\lambda\bfW, \Omega_{2n} = \omega_0)$ into
$(\T^{2n},\omega_{\bbw})$.
Now \ref{_packing-by-arb-shapes-tori-K3-from-Kahler-type-to-tame_Corollary_} yield the claims of parts A and B
of \ref{_full-packing-of-tori-by-polydisks-parallelepipeds-COPY_Corollary_}
in the case where $\bfW$ is the disjoint union of $k$ equal polydisks.

Now let us consider the case where $\bfW$ is the disjoint union of $k$ identical copies of
a parallelepiped
\[
P (e_1,\ldots,e_{2n}) := \left\{ \sum_{j=1}^{2n} s_j e_j, 0\leq s_j \leq 1, j=1,\ldots,2n \right\},
\]
where $e_1,\ldots,e_{2n}$ is a basis of the vector space $\R^{2n}$.

Note that for any $0<\lambda <1$, putting $k$ copies of $\lambda P (e_1,e_2,\ldots,e_{2n})$, alongside each other in the $e_1$-direction we get a K\"ahler-type embedding of the disjoint union of $k$ identical copies of $\lambda P (e_1,e_2,\ldots,e_{2n})$ into $P (k e_1,e_2,\ldots,e_{2n})$. Considering the composition of this K\"ahler-type embedding with K\"ahler-type, or $\varepsilon$-tame, embeddings of $P (k e_1,e_2,\ldots,e_{2n})$
into $\T^{2n}$, we readily get that it suffices to prove the corollary for just one parallelepiped $P (e_1,\ldots,e_{2n})$
(i.e., for the case where $k=1$).

Therefore let us prove the corollary for one parallelepiped $P (e_1,\ldots,e_{2n})$.

Consider the $2n\times 2n$-matrix $Z (e_1,e_2,\ldots,e_{2n})$ whose columns are
the vectors $e_1,e_2,\ldots,e_{2n}$ spanning $P (e_1,\ldots,e_{2n})$.

For each $m\in\Z_{>0}$ let ${\it Mat}_m (\R)$, respectively ${\it Mat}_m (\Q)$, be the space of real, respectively rational, $m\times m$-matrices.

Let $S\subset {\it Mat}_{2n} (\R)$
be the set of matrices $Z\in {\it Mat}_{2n} (\R)$ such that for all $\lambda\in \R_{>0}$
\[
\lambda Z^t
\begin{pmatrix}
0 & I_n\\
-I_n & 0
\end{pmatrix}
Z \notin {\it Mat}_{2n} (\Q).
\]
(Here $I_n$ is the identity $n\times n$-matrix.)


\hfill


\vfil\eject

\lemma
\label{_S-is-dense_Lemma_}

The set $S$ is dense in ${\it Mat}_{2n} (\R)$.


\hfill


\noindent
{\bf Proof of \ref{_S-is-dense_Lemma_}:}

If $Z$ has the form
\[
\begin{pmatrix}
A & B\\
C & D
\end{pmatrix},
\]
where $A,B,C,D\in {\it Mat}_n (\R)$, then
\[
\lambda Z^t
\begin{pmatrix}
0 & I_n\\
-I_n & 0
\end{pmatrix}
Z =
\begin{pmatrix}
A^t C - C^t A & A^t D - C^t B\\
B^t C - D^t A & B^t D - D^t B
\end{pmatrix}.
\]
Since $n>1$, the set of $Z\in {\it Mat}_{2n} (\R)$ such that $A,C\in {\it Mat}_n (\Q)$ and  $B^t D - D^t B\notin {\it Mat}_n (\Q)$
is dense in ${\it Mat}_{2n} (\R)$ and also lies in $S$. Thus, $S$ is dense in ${\it Mat}_{2n} (\R)$.

This proves the lemma.
\endproof


\hfill


Since $S\subset {\it Mat}_{2n} (\R)$ is a dense set,
it suffices to prove the corollary only in the case when $Z (e_1,e_2,\ldots,e_{2n})$ lies in $S$.

Indeed, if $\Vol \left(\lambda P \left(e_1,\ldots,e_{2n}\right),\omega_0\right) < \Vol (\T^{2n},\omega)$ for some $\lambda >0$ and some K\"ahler-type form $\omega$ on $\T^{2n}$,
then there exists a basis $e'_1,\ldots,e'_{2n}$ of $\R^{2n}$ such that $Z (e'_1,e'_2,\ldots,e'_{2n})\in S$,
the parallelepiped $P (e'_1,\ldots,e'_{2n})$
contains the parallelepiped $P (e_1,\ldots,e_{2n})$ and has only slightly bigger volume than it, so that
$\Vol (\lambda P (e'_1,\ldots,e'_{2n}),\omega_0) < \Vol (\T^{2n},\omega)$. If $\lambda P (e'_1,\ldots,e'_{2n})$ admits a K\"ahler-type, or respectively an $\varepsilon$-tame, embedding into $(\T^{2n},\omega)$, then the composition of this embedding with the inclusion
$\lambda P (e_1,\ldots,e_{2n})\to \lambda P (e'_1,\ldots,e'_{2n})$ is a K\"ahler-type, or respectively an $\varepsilon$-tame, embedding of
$\lambda P (e_1,\ldots,e_{2n})$ into $(\T^{2n},\omega)$.

Now, in view of \ref{_packing-by-arb-shapes-tori-K3-from-Kahler-type-to-tame_Corollary_}, it suffices to show that
for any  $P (e_1,\ldots,e_{2n})$, such that $Z (e_1,e_2,\ldots,e_{2n})\in S$,
there exists an irrational
K\"ahler-type symplectic form $\omega$ on $\T^{2n}$ satisfying the following condition:
for any $\lambda>0$, such that $\Vol \left( \lambda P \left(e_1,\ldots,e_{2n}\right)\right) < \Vol (\T^{2n},\omega)$,
there exists a K\"ahler-type embedding $P (e_1,\ldots,e_{2n})\to (\T^{2n},\omega)$.

Pick $P (e_1,\ldots,e_{2n})$, such that $Z (e_1,e_2,\ldots,e_{2n})\in S$ and
let us construct such an $\omega$.

Without loss of generality, assume that $\Vol \left(\lambda P \left(e_1,\ldots,e_{2n}\right),\omega_0\right) = 1$.

Consider the integral lattice $\Lambda:= \Span_\Z \{ e_1,\ldots,e_{2n}\}\subset \R^{2n}$. The fundamental domain of $\Lambda$
is exactly $P (e_1,\ldots,e_{2n})$.
The space $\R^{2n}/\Lambda$
is a torus and the form $\omega_0$ on $\R^{2n}$ induces
a K\"ahler-type symplectic form on $\R^{2n}/\Lambda$ of total volume $1$ that will be also denoted by $\omega_0$. The matrix $Z (e_1,e_2,\ldots,e_{2n})$ induces a linear map $\R^{2n}\to \R^{2n}$ mapping $\Z^{2n}$ to $\Lambda$. This map induces
a diffeomorphism $F: \T^{2n} = \R^{2n}/\Z^{2n}\to \R^{2n}/\Lambda$. Since $Z (e_1,e_2,\ldots,e_{2n})\in S$,
the form $\omega:= F^* \omega_0$ is an irrational K\"ahler-type form on $\T^{2n}$ of total volume $1$.

For any $0<\lambda<1$ the inclusion $\lambda P (e_1,\ldots,e_{2n})\to (\R^{2n}/\Lambda, \omega_0)$ is a K\"ahler-type embedding.
This embedding is holomorphic with respect to the complex structure on $\R^{2n}/\Lambda$ induced by the standard complex structure $J_0$ on $\R^{2n}$. The composition of this embedding with $F^{-1}$
is then a K\"ahler-type embedding  $\lambda P (e_1,\ldots,e_{2n})\to (\T^{2n},\omega)$ holomorphic with respect to the linear complex
structure $F^* J_0$ on $\T^{2n}$.

Thus, for any $0<\lambda<1$ there exists a K\"ahler-type embedding $\lambda P (e_1,\ldots,e_{2n})\to (\R^{2n}/\Lambda, \omega_0)$,
provided $Z (e_1,e_2,\ldots,e_{2n})\in S$.

This finishes the proof of the corollary in the case where $\bfW$ is the disjoint union of $k$ identical copies of
a parallelepiped.
\endproof


\hfill


Let us prove \ref{_full-packing-of-K3-by-parallelepipeds_Corollary_}. We recall it here.


\hfill


\corollary {\bf (= \ref{_full-packing-of-K3-by-parallelepipeds_Corollary_})}
\label{_full-packing-of-K3-by-parallelepipeds-COPY_Corollary_}

Assume:

\begin{itemize}

\item{} $M$ is a smooth manifold underlying a K3 surface.

\item{} $\bfW:= \bigsqcup_{i=1}^k W_i$ is a disjoint union of $k$ identical copies of a parallelepiped
\[
P (e_1,\ldots,e_4) := \left\{ \sum_{j=1}^4 s_j e_j, 0\leq s_j \leq 1, j=1,\ldots,4 \right\},
\]
where $e_1,\ldots,e_4$ is a basis of the vector space $\R^4$.

\item{} $\varepsilon >0$.

\end{itemize}

Then there exists a $\Diff^+ (M)$-invariant open dense set of K\"ahler-type symplectic forms on $M$, depending on $\bfW$ and $\varepsilon$ and containing, in particular, all irrational K\"ahler-type symplectic forms on $M$,
so that for each $\omega'$ in this set
$\nu_{T,\varepsilon} (M,\omega',\bfW) = 1$, meaning that $\varepsilon$-tame embeddings of $\lambda\bfW$ into $(M,\omega')$ are unobstructed.


\hfill


\noindent
{\bf Proof of \ref{_full-packing-of-K3-by-parallelepipeds-COPY_Corollary_} (=\ref{_full-packing-of-K3-by-parallelepipeds_Corollary_}):}

Consider the integral lattice $\Lambda:= \Span_\Z \{ 2e_1, 2 e_2, 2e_3, 2e_4\}\subset \R^4$. The parallelepiped $P (e_1,\ldots,e_4)$ is the fundamental domain of the lattice $\Lambda/2$.

Recall the following construction by Kummer of a K3 surface -- see e.g. \cite[p.55]{_Geom-K3-Asterisque_}.

Consider the torus $\T^4: = \R^4/\Lambda$ and let $\pi: \R^4\to \T^4$ be the projection.
Equip the torus $\T^4$ with the complex structure $J_0$ and the symplectic structure $\omega_0$ induced by the standard complex and symplectic structures on $\R^4$.
The reflection with respect to the origin of $\R^4$ defines an involution of $\R^4$, which in turn defines a holomorphic and symplectic involution $\iota: (\T^4,J_0,\omega_0) \to (\T^4,J_0,\omega_0)$.
The involution $\iota$ has $16$ fixed points -- these are exactly the points of $\T^4$ that are the images under $\pi$ of
the points in $\R^4$ all of whose coordinates with respect to the basis $\{ 2e_1, 2e_2,\ldots,2e_{2n}\}$ are either $0$ or $1/2$. Let $(\tilde{\T}^4,\tJ_0)$ be the complex blow-up of $(\T^4,J_0)$
at these 16 points.
The involution $\iota$ lifts to a holomorphic involution $\tiota: \tilde{\T}^4\to \tilde{\T}^4$. Then $M:=\tilde{\T}^4/\tiota$, equipped with the
complex structure induced by $\tJ_0$, is
a complex manifold which is, in fact, a K3 surface (called a Kummer surface). Denote by $I$ the complex structure on $M$ induced by $\tJ_0$.

The complex manifold $(M,I)$ can also be obtained as follows: The space $\T^4/\iota$ is a complex orbifold with $16$ isolated singular orbifold points $x_1,\ldots,x_{16}$. Denote by $\bJ_0$ the complex structure on $\T^4/\iota$ induced by $J_0$. Then $(M,I)$ is the smooth complex manifold obtained by the complex blow-up of $(\T^4/\iota, \bJ_0)$ at $x_1,\ldots,x_{16}$.
Let $pr: \T^4\to \T^4/\iota$ and $\Pi: M\to \T^4/\iota$ be the natural projections. The form $\omega_0$ on $\T^4$ is $\iota$-invariant and induces an orbifold K\"ahler form $\bomega_0$ on the complex orbifold $\T^4/\iota$ which is a smooth K\"ahler form outside the singularities.

The points $x_1,\ldots,x_{16}\in \T^4/\iota$ lie on the boundary of $(pr\circ \pi) \left(P \left(e_1,\ldots,e_4\right)\right)$. The interior $\Int P (e_1,\ldots,e_4)$ of $P (e_1,\ldots,e_4)$ projects diffeomorphically onto its image under $pr\circ \pi$ and is an open dense subset of $\T^4/\iota$. Consequently, $pr\circ \pi: \Int P (e_1,\ldots,e_4)\to \T^4/\iota$ is an embedding which is
holomorphic with respect to $\bJ_0$ and symplectic with respect to $\bomega_0$.
The map
$F:= \Pi^{-1}\circ pr\circ \pi: \Int P (e_1,\ldots,e_4)\to M$ is then an $I$-holomorphic embedding.

Let $0<\lambda<1$. Move $\lambda P (e_1,\ldots,e_4)$ by a parallel translation into the interior of
\break
$P (e_1,\ldots,e_4)$
and then map it by $F$ into $M$. Denote the resulting map by $f: \lambda P (e_1,\ldots,e_4)\to M$. It is an $I$-holomorphic embedding.

Following the proof of \cite[Thm. 3.6]{_EV-Selecta_} (or, alternatively, using \ref{_sympl-blow-up-regularized-max_Proposition_}), one can find rational cohomology classes $b_i\in H^2 (\Pi^{-1} (x_i);\Q)$, $i=1,\ldots,16$, so that
for any $\epsilon>0$ and any sufficiently small $0<\epsilon_1,\ldots,\epsilon_{16} <\epsilon$, there exist
a K\"ahler form $\omega_\epsilon$ on $(M,I)$ and sufficiently small neighborhoods $U_{x_i,\epsilon,\lambda}$ of $\Pi^{-1} (x_i)$, $i=1,\ldots,16$, in $M$
with the following properties:

\begin{itemize}

\item{} $U_{x_i,\epsilon,\lambda}\cap f \left(\lambda P \left(e_1,\ldots,e_4\right)\right) = \emptyset$, $i=1,\ldots,16$.

\item{} The form $\omega_\epsilon$ coincides with $\Pi^* \bomega_0$ outside $\bigcup_{i=1}^{16} U_{x_i,\epsilon,\lambda}$.

\item{} $[\omega_\epsilon] = \Pi^* [\bomega_0] + \sum_{i=1}^{16} \epsilon_i b_i$.

\item{} $\Vol (\T^4/\iota, \bomega_0) > (1-\epsilon)
\Vol (M,\omega_\epsilon)$.

\end{itemize}

If the vector $(\epsilon_1,\ldots, \epsilon_{16})\in\R^{16}$ is not a real multiple of a vector in $\Q^{16}$, the cohomology class
\[
[\omega_\epsilon] = \Pi^* [\bomega_0] + \sum_{i=1}^{16} \epsilon_i b_i\in H^2 (M;\R)
\]
is a not a real multiple of a rational cohomology class, meaning that $\omega_\epsilon$ is irrational.

Our construction implies that $f: \lambda P (e_1,\ldots,e_4)\to (M,\omega_\epsilon)$ is a K\"ahler-type embedding which is holomorphic with respect to $I$ and
\[
\Vol \left(\lambda P (e_1,\ldots,e_4),\omega_0\right) = \lambda^4 \Vol \left(P (e_1,\ldots,e_4),\omega_0\right) =
\]
\[
= \lambda^4 \Vol (\T^4/\iota, \bomega_0) > \lambda^4 (1-\epsilon)
\Vol (M,\omega_\epsilon).
\]

For any $\delta>0$ we can choose
$\lambda < 1$ sufficiently close to $1$ and $\epsilon>0$ sufficiently small
so that $\lambda^4 (1-\epsilon) > 1-\delta$.
Then
\[
f: \lambda P (e_1,\ldots,e_4)\to (M,\omega_\epsilon)
\]
is a K\"ahler-type embedding
and
\[
\Vol (\lambda P (e_1,\ldots,e_4),\omega_0)/\Vol (M,\omega_\epsilon) \geq 1-\delta.
\]
By \ref{_packing-by-arb-shapes-tori-K3-from-Kahler-type-to-tame_Corollary_}, this implies the corollary.
\endproof


\hfill


\section{Hyperk\"ahler case -- proofs}
\label{_hyperkahler-case-proofs_Section_}

As we mentioned in Section~\ref{_hyperkahler-case_Subsection_}, the results of Section~\ref{_main-results-tori-K3_Subsection_}
generalize to the IHS-hyperk\"ahler case. In this section we
explain how the proofs of the results in the K3 case should be adjusted for general IHS-hyperk\"ahler manifolds.

Let $M$ be a closed connected and simply-connected manifold carrying IHS-hyperk\"ahler structures.

Recall from Section~\ref{_hyperkahler-case_Subsection_} that in the IHS-hyperk\"ahler case the words ``K\"ahler-type" should be replaced everywhere by the words ``IHS-hy\-per\-k\"ahler-type". {\bf Accordingly, in this section we will use the same notation as in Section~\ref{_Campana-simple-structures-K3_Subsection_} for IHS-hy\-per\-k\"ahler-type objects as for K\"ahler-type objects. In particular, $\compK (M)$ will denote the space of IHS-hy\-per\-k\"ahler-type complex structures, $\Teich (M) = \compK (M)/\Diff_0 (M)$ the Teichm\"uller space of IHS-hy\-per\-k\"ahler-type complex structures etc.; all the symplectic forms and complex structures will be assumed to be of IHS-hy\-per\-k\"ahler-type.}

For general IHS-hyperk\"ahler manifolds:

\begin{itemize}

\item{} The Teichm\"uller space $\Teich (M) = \compK (M)/\Diff_0 (M)$ has the structure of a (non-Haus\-dorff) complex manifold for the same reasons as in Section~\ref{_Campana-simple-structures-K3_Subsection_} (see \cite{_Catanese:moduli_}).

\item{} Instead of the intersection form one should consider the Bogomolov-Beauville-Fujiki (BBF) symmetric bilinear form $q$ on $H^2 (M;\R)$
(see \cite{_Bea1_}, \cite{_Fujiki:HK_}) -- for K3 surfaces the intersection form and the BBF form coincide, up to a multiplicative factor. The BBF form is a primitive integral quadratic form of signature
$(3, b_2-3)$, where $b_2 = \dim_\R H^2 (M;\R) >0$.
The space $Gr_{++} (H^2 (M;\R))$ is defined then using the BBF form $q$.

We will use $\perp$ to denote orthogonal complements with respect to $q$.

\item{} Part A of the Torelli theorem (\ref{_K3-Torelli_Theorem_}) remains true for IHS-hyperk\"ahler manifolds -- see \cite{_V-Duke_}.

\item{} In the definition of $S(M,I)$, for $I\in \compK (M)$, $(-2)$-classes should be replaced by the so-called MBM cohomology classes in $H^2 (M;\Z)$ defined in \cite{_Am-Ver-JGF2015_}. With $S(M,I)$ defined in this way,  \ref{_K3-Kahler-cone_Theorem_} remains true for general IHS-hyperk\"ahler manifolds.

\item{} The analogues of \ref{_K3-Calabi-Yau_Theorem_} and \ref{_K3-adjusting-cs-to-get-Ricci-flat-metric_Proposition_} remain true for general IHS-hy\-per\-k\"ah\-ler manifolds.

\end{itemize}

Let $\omega$ be an IHS-hy\-per\-k\"ahler-type symplectic form on $M$.

The set $\Dom_\omega$ is then defined literally as in Section~\ref{_Campana-simple-structures-K3_Subsection_} with $(-2)$-classes being replaced by the MBM classes. \ref{_K3-Dom-omega_Proposition_} and \ref{_K3-complex-str-compatible-with-omega-description_Proposition_} remain true in this case -- in their proofs one should replace
$19$ by $b_2 -3$, $18$ by $b_2 -4$, and $(-2)$-classes by MBM classes \cite{_Am-Ver-JGF2015_}.

Next, let us state the analogue of \ref{_irrational-sympl-form-on-K3-compatible-with-Campana-simple-complex-str_Proposition_} in the IHS-hyperk\"ahler case.


\hfill


\proposition
\label{_irrational-sympl-form-on-IHS-hyperkahler-mfd-compatible-with-Campana-simple-complex-str_Proposition_}

The symplectic form $\omega$ is compatible with a Campana-simple IHS-hyperk\"ahler-type complex structure on $M$ if and only if
$\omega$ is irrational. In this case, for any connected component $\cC_0$ of
$\compK (M)$ compatible with $\omega$ the set
$\cmptteich_{\cC_0} (M,\omega)\cap \Teich^{s} (M)$ of Campana-simple structures lying in $\cmptteich_{\cC_0} (M,\omega)$
is a dense connected subset of $\cmptteich_{\cC_0} (M,\omega)$.


\hfill


The proof of \ref{_irrational-sympl-form-on-IHS-hyperkahler-mfd-compatible-with-Campana-simple-complex-str_Proposition_}
that we are going to give will not be a direct generalization of the proof of \ref{_irrational-sympl-form-on-K3-compatible-with-Campana-simple-complex-str_Proposition_} and requires
a number of preliminaries.

Consider an IHS-hy\-per\-k\"ahler-type structure $\goh=\{ I_1, I_2,I_3, \omega_1,\omega_2,\omega_3\}$ on $M$.
For any $a,b,c\in\R$, $a^2+b^2+c^2=1$, the tensor $aI_1+bI_2+cI_3$
is a complex structure of hyperk\"ahler type on $M$. We will call
any such  $aI_1+bI_2+cI_3$ {\bf a complex structure induced by $\goh$}. Clearly, all complex structures induced by $\goh$ lie in the same connected
component of $\compK (M)$ -- we will say that this connected component is {\bf favored} by $\goh$.

Following \cite{_Verbitsky:Symplectic_II_},
we say that a closed subset $Z\subset M$ is
{\bf trianalytic} (with respect to $\goh$), if $Z$
is a complex
subvariety of $(M,I)$ for each complex structure $I$ induced by
$\goh$.


\hfill

\theorem {\bf \cite[Thm. 4.1]{_Verbitsky:Symplectic_II_}}
\label{_G_M_invariant_implies_trianalytic_Theorem_}

Let $a\in  H^* (M;\Q)$, $a\neq 0$.

Then $a$ is Poincar\'e-dual to the fundamental class $[Z]$ of a trianalytic $Z\subset M$ if and only if $a$ is invariant with respect to the action of any induced complex structure on $H^* (M;\Q)$.
\endproof

\hfill


Let $\cC_0$ be a connected component  of $\compK (M)$.

Let $Z$ be a closed subvariety of $M$ for a certain (IHS-hyperk\"ahler-type) complex structure on $M$ lying in $\cC_0$.
Following \cite{_Kurnosov:tori_,_SV:k-symplectic_}, we say that $Z$ is {\bf absolutely trianalytic}
(with respect to $\cC_0$)
if the cohomology class Poincar\'e-dual to the fundamental class of $Z$ is
of type $(p,p)$ for all complex structures $I\in\cC_0$.
By \ref{_G_M_invariant_implies_trianalytic_Theorem_}, such a
manifold is trianalytic with respect to any IHS-hyperk\"ahler structure favoring $\cC_0$.

Given $a\in H^* (M;\Q)$, $a\neq 0$, denote by
$\THodge_{\cC_0, a}\subset \Teich_{\cC_0} (M)$ the set of all $[I]\in \Teich_{\cC_0} (M)$ such that $a\in H_I^{1,1} (M;\R)$.
It follows e.g. from the Torelli theorem for hyperk\"ahler manifolds \cite{_V-Duke_} that $\THodge_{\cC_0, a}$ is an analytic submanifold of
$\Teich_{\cC_0} (M)$ and $\Per (\THodge_{\cC_0, a}) = Gr_{++} \left(a^\bot\right)$.

For a real vector space $W$ equipped with a non-degenerate indefinite
quadratic form we denote by $SO^+ (W)$ the identity component of $SO (W)$. Similarly, $SO^+ (p,q)$ stands for the identity component of the Lie group $SO (p,q)$.

For the proof of \ref{_irrational-sympl-form-on-IHS-hyperkahler-mfd-compatible-with-Campana-simple-complex-str_Proposition_} we will need the following algebraic lemma.


\hfill


\lemma
\label{_rotations_generate_SOpq_Lemma_}

Let $W$ be a real vector space equipped with a non-degenerate
quadratic form $q$ of signature $(2,l)$.

Consider the set $\cF (W)$ of $F\in SO^+(W)$ fixing pointwise a linear subspace $U\subset W$, $\dim_\R U = \dim_\R W- 2$ whose orthogonal complement
is a positive 2-dimensional plane in $W$.

Then the set $\cF (W)$
generates $SO^+(W)$.


\hfill


\noindent
{\bf Proof of \ref{_rotations_generate_SOpq_Lemma_}:}

Let $N$ be the
subgroup of $SO^+(W)$ generated by $\cF (W)$.
Clearly, $N$ is normal.

Assume $l\neq 2$.
Then the group
$SO^+(W)$ is a connected simple Lie group (in the sense of Lie theory). Therefore it follows from \cite{_Rag_} that any normal subgroup of  $SO^+(W)$ is either
the whole $SO^+(W)$ or is contained in the center of $SO^+(W)$ (which is discrete). Hence, $N=SO^+(W)$, meaning that
$\cF (W)$
generates $SO^+(W)$ for $l\neq 2$.

Let us consider the case $l=2$ -- i.e., the case when $q$ has signature $(2,2)$.
Identify $(W,q)$ with the space ${\it Mat}_2 (\R)$ of real $2\times 2$-matrices equipped with the bilinear form $\langle A,B\rangle = -tr (AJBJ)$, where
$J=\left(
           \begin{array}{cc}
             0 & -1 \\
             1 & 0 \\
           \end{array}
         \right)$.
Consider the following action of $SL(2,\R)\times SL (2,\R)$ on $W\cong {\it Mat}_2 (\R)$:
each $h_1\times h_2\in SL(2,\R)\times SL (2,\R)$ acts on $A\in {\it Mat}_2 (\R)$ by
$A\mapsto h_1 A h_2^{-1}$.
This action preserves the bilinear form and induces a homomorphism $P: SL (2,\R)\times SL (2,\R)\to SO^+(W)$
whose kernel is $\{ \pm Id\}$. The subgroup $P^{-1} (N)\subset SL(2,\R)\times SL (2,\R)$ is normal.
The Lie group $SL(2,\R)\times SL (2,\R)$ is a connected semi-simple Lie group. Therefore it follows from \cite{_Rag_} that any normal subgroup of $SL(2,\R)\times SL (2,\R)$ which is not the whole $SL(2,\R)\times SL (2,\R)$ is contained either in $SL (2,\R)\times \{ \pm Id\}$ or in
$\{ \pm Id\}\times SL (2,\R)$.

Let us show that $P^{-1} \left(\cF (W)\right)\nsubseteq SL (2,\R)\times \{ \pm Id\}$ and consequently $P^{-1} (N)\nsubseteq SL (2,\R)\times \{ \pm Id\}$ (the claim $P^{-1} (N)\nsubseteq \{ \pm Id\}\times SL (2,\R)$ is proved similarly). Indeed, if an element $h\times Id\in SL (2,\R)\times SL (2,\R)$ is sent by $P$ into $\cF (W)$, then its action on ${\it Mat}_2 (\R)$ fixes a two-dimensional subspace of ${\it Mat}_2 (\R)$. This easily implies that $h\in SL (2,\R)$ has a two-dimensional eigenspace for the eigenvalue $1$, meaning that $h=Id$. Hence,  $P^{-1} \left(\cF (W)\right)\nsubseteq SL (2,\R)\times \{ \pm Id\}$
and consequently $P^{-1} (N) =  SL(2,\R)\times SL (2,\R)$ and $N = SO^+(W)$, meaning that
$\cF (W)$
generates $SO^+(W)$ also for $l = 2$.

This finishes the proof of the lemma.
\endproof


\hfill


Now we are ready to prove \ref{_irrational-sympl-form-on-IHS-hyperkahler-mfd-compatible-with-Campana-simple-complex-str_Proposition_}.


\hfill


\noindent
{\bf Proof of \ref{_irrational-sympl-form-on-IHS-hyperkahler-mfd-compatible-with-Campana-simple-complex-str_Proposition_}:}

If the cohomology class $[\omega]$ is rational,
then by part 1 of
\ref{_Campana-simple-vs-projectivity_Remark_}, there are no Campana-simple complex structures compatible with $\omega$.

Assume that $\omega$ is irrational and $\cC_0$ is a connected component of
$\compK (M)$ compatible with $\omega$.

The analogue of \ref{_K3-complex-str-compatible-with-omega-description_Proposition_} in the
IHS-hy\-per\-k\"ahler case implies that $\cmptteich_{\cC_0} (M,\omega)$ is a dense connected subset of the analytic subvariety
$\THodge_{\cC_0, [\omega]}$ of $\Teich_{\cC_0} (M)$.
Consequently, for each $a\in H^* (M;\Q)$, $a\neq 0$, there are only two possible options:

\smallskip
\noindent
(1) $\THodge_{\cC_0, a} \cap\THodge_{\cC_0, [\omega]}=\THodge_{\cC_0, [\omega]}$, and consequently  $\THodge_{\cC_0, a} \cap\cmptteich_{\cC_0} (M,\omega) = \cmptteich_{\cC_0} (M,\omega)$,

\smallskip
\noindent
(2) $\THodge_{\cC_0, a} \cap\THodge_{\cC_0, [\omega]}$ is a proper (possibly empty) complex subvariety of $\THodge_{\cC_0, [\omega]}$ and consequently $\THodge_{\cC_0, a} \cap\cmptteich_{\cC_0} (M,\omega)$ is a proper complex subvariety of  $\cmptteich_{\cC_0} (M,\omega)$.


\hfill


\lemma
\label{_not-trianalytic_Lemma_}

Assume $Z\subset M$ is a complex subvariety of $(M,I)$ (some complex structure $I\in \cC_0$) which is not absolutely trianalytic
with respect to $\cC_0$. Let $a\in H^{even} (M;\Q)$ be the cohomology class Poincar\'e-dual to $[Z]$.

Then $\THodge_{\cC_0, a} \cap\THodge_{\cC_0, [\omega]}\nsubseteqq\THodge_{\cC_0, [\omega]}$ (meaning that option (2) holds for $a$).


\hfill


\noindent
{\bf Proof of \ref{_not-trianalytic_Lemma_}:}

We start with a number of preparations.

Given $I\in\cC_0$ and $t\in\R$, consider the Hodge rotation map $H^* (M;\C)\to H^* (M;\C)$ that acts on each $H_I^{p,q} (M;\C)$
as multiplication by $e^{\1 (p-q) t}$. Such a map preserves $H^* (M;\R)$ and its fixed point set is exactly the union of $H_I^{p,p} (M;\C)$ for all $p$. We will view the action of this map on $H^* (M;\R)$ as an element of
$GL \left(H^* (M;\R)\right)$ and denote it by $\rho_I (t): H^* (M;\R)\to H^* (M;\R)$. The restriction of each $\rho_I (t)$ to $H^2 (M;\R)$
preserves the BBF form $q$ and lies in $SO^+ \left(H^2 (M;\R)\right)$. The maps $\rho_I (t)|_{H^2 (M;\R)}$, $t\in\R$, are exactly the maps that act trivially on the subspace $H_I^{1,1} (M;\R)$
and rotate its orthogonal complement, which is the positive 2-dimensional plane $(H_I^{2,0}(M;\C)\oplus H_I^{0,2}(M;\C))\cap H^2 (M;\R)$.

Denote by $G\subset GL \left(H^{even} (M;\R)\right)$ the
subgroup generated by
$\rho_I (t)|_{H^{even} (M;\R)}$ for all $I\in\cC_0$ and $t\in\R$. The
restriction to $H^2 (M;\R)$ defines a homomorphism $G\to
SO^+ \left(H^2 (M;\R)\right)$. It follows from
\cite[Thm. 2.2]{_V:applications_} and
\cite[Corollary 8.2]{_V:Mirror_}
that this homomorphism
is, in fact, an isomorphism.
(On $H^{odd}(M)$ the group generated by $\rho_I (t)$, for all $I\in\cC_0$ and $t\in\R$, acts as
the identity component of $\Spin(H^2 (M;\R))$,
with the center of the spinor group acting as $-\Id$.)

Denote by $G_{[\omega]}\subset G$
the subgroup generated by all Hodge rotation maps $\rho_I (t)$, $t\in\R$, $I\in \cC_0$, $[I]\in \THodge_{\cC_0, [\omega]}$.
Each such $\rho_I (t)$ fixes $[\omega]$ because $[\omega]\in H_I^{1,1} (M;\R)$. Hence, $[\omega]$ is fixed by $G_{[\omega]}$. Thus, $G_{[\omega]}$
lies in the stabilizer of $[\omega]$ in $G\cong SO^+ \left(H^2 (M;\R)\right)$. The latter stabilizer is isomorphic to $SO^+ ([\omega]^\perp)$. Note that
the restriction of the BBF form $q$ to $[\omega]^\perp\subset H^2 (M;\R)$ has signature $(2, b_2-3)$. Now \ref{_rotations_generate_SOpq_Lemma_}
yields that the maps $\rho_I (t)|_{H^2 (M;\R)}$, $t\in\R$, $[I]\in \THodge_{\cC_0, [\omega]}$, generate the stabilizer. In other words,
$G_{[\omega]}$ coincides with the stabilizer of $[\omega]$ in $G$.

For the class $a\in H^{even} (M;\Q)$ given in the hypothesis of the lemma, denote by $G_a$ the stabilizer of
$a$ in $G$.

Now we are ready to prove the lemma.

Assume, by contradiction, that $\THodge_{\cC_0, [\omega]}\subset\THodge_{\cC_0, a}$, meaning that option (1) above, and not (2), holds for $a$.
Then $a$ is of type $(p,p)$ for all
complex structures $I\in\cC_0$ such that $[I]\in \THodge_{\cC_0, [\omega]}$.
In particular, $a$ is fixed by  $\rho_I (t)$ for all such $I$ and all $t\in\R$.
Thus, $G_{[\omega]}$ fixes $a$ or, in other words,
\[
SO^+ (2, b_2-3)\cong G_{[\omega]}\subset G_a\subset G\cong SO^+ (3, b_2-3).
\]
By \cite[Lem. 9.9]{_EV-JTA_},
any intermediate subgroup
$G_{[\omega]} \subset H \subset G$ is
equal either to $G_{[\omega]}$ or to $G$.
This implies that
$G_a$ equals
either $G_{[\omega]}$ or $G$.
Since $[\omega]$ is not rational, \cite[Lem. 9.10]{_EV-JTA_} yields that $G_a\neq G_{[\omega]}$.
Since $Z$ is not
absolutely trianalytic, there exists $I\in\cC_0$ such that $a\notin H_I^{1,1} (M;\R)$ and consequently
$a$ is not fixed by $\rho_I (t)\in G$ for $t\neq 0$. Therefore $\rho_I (t)\notin G_a$ and $G_a\neq G$.
This yields a contradiction and finishes the proof of the lemma.
\endproof


\hfill


Now let us finish the proof
of \ref{_irrational-sympl-form-on-IHS-hyperkahler-mfd-compatible-with-Campana-simple-complex-str_Proposition_}.

Let $\cA\subset H^* (M;\Q)\setminus \{ 0\}$ be the set of all $a\in H^* (M;\Q)$, $a\neq 0$, for
which $\THodge_{\cC_0, a} \cap\THodge_{\cC_0, [\omega]}$ is a proper (possibly empty) complex subvariety of $\THodge_{\cC_0, [\omega]}$ and consequently $\THodge_{\cC_0, a} \cap\cmptteich_{\cC_0} (M,\omega)$ is a proper complex subvariety of  $\cmptteich_{\cC_0} (M,\omega)$.
By \ref{_not-trianalytic_Lemma_}, $\cA$ includes, in particular, all $a\in H^2 (M;\Q)$ that are Poincar\'e-dual to the fundamental classes of
not absolutely trianalytic subvarieties (with respect to $\cC_0$).

The set
\[
S:= \cmptteich_{\cC_0} (M,\omega)\setminus \bigcup_{a\in \cA} \left(\THodge_{\cC_0, a} \cap\cmptteich_{\cC_0} (M,\omega)\right)
\]
is connected and dense in $\cmptteich_{\cC_0} (M,\omega)$. For any $I\in\cC_0$ such that $[I]\in S$, the complex manifold $(M,I)$ can contain only absolutely trianalytic subvarieties. By \cite[Lem. 7.9]{_EV-JTA_}, the union of all such subvarieties of $(M,I)$ is of measure zero, meaning that $I$ is Campana-simple and $[I]\in \cmptteich_{\cC_0} (M,\omega)\cap \Teich^{s} (M)$.  Hence,
$S\subset \cmptteich_{\cC_0} (M,\omega)\cap \Teich^{s} (M)\subset \cmptteich_{\cC_0} (M,\omega)$.
Since $S$ is connected and dense in $\cmptteich_{\cC_0} (M,\omega)$, so is $\cmptteich_{\cC_0} (M,\omega)\cap \Teich^{s} (M)$.

This finishes the proof of the proposition.
\endproof


\hfill


We are now ready to discuss the analogue of \ref{_existence-connectedness-irrational-forms-tori-K3-surfaces_Theorem_}
in the IHS-hyperk\"ahler case, as described in Section~\ref{_hyperkahler-case_Subsection_}.

The proof of the claims of \ref{_existence-connectedness-irrational-forms-tori-K3-surfaces_Theorem_}
in the IHS-hyperk\"ahler case
(apart from the claim on the transitivity of the $\Symp_H (M,\omega)$-action)
is the same as in the K3 case -- one just needs
to use the generalizations of the results of Section~\ref{_Campana-simple-structures-K3_Subsection_} to the IHS-hyperk\"ahler case discussed above
(instead of the corresponding results for the K3 surfaces) and use \ref{_Campana-simple-Kahler-type-existence-connectedness-gen-mfd_Theorem_}
in the simply-connected case involving Campana-simple complex structures.

As for the transitivity of the action of $\Symp_H (M,\omega)$, we cannot say anything in the IHS-hyperk\"ahler case, as we have no analogue
of part B of \ref{_K3-Torelli_Theorem_} for IHS-hyperk\"ahler manifolds.

Finally, let us discuss the results on tame embeddings into $(M,\omega)$ for an IHS-hy\-per\-k\"ahler-type form $\omega$.

Recall from Section~\ref{_hyperkahler-case_Subsection_} that there are two possibilities to define tame embeddings into such an $(M,\omega)$. One possibility is to use K\"ahler-type complex structures tamed by $\omega$. Another possibility is to use instead
IHS-hy\-per\-k\"ahler-type complex structures tamed by $\omega$.

The claim $\nu_{T,\varepsilon} (M,\omega_1,\bfW) = \nu_{T,\varepsilon} (M,\omega_2,\bfW)$ of \ref{_packing-by-arb-shapes-main-tori-K3_Theorem_} holds for both definitions of $\nu_{T,\varepsilon}$, as long as
$\omega_1$ and $\omega_2$ lie in the same connected component of the space of IHS-hy\-per\-k\"ahler-type symplectic forms on $M$.
Indeed, it is easy to see that \ref{_semicontinuity_Proposition_} holds for both definitions of $\nu_{T,\varepsilon}$.
The relevant claim about the density of the $\Diff^+ (M)$-orbits of the form $\omega_1$, $\omega_2$ in the space of forms of the same volume in the same connected component of the space of IHS-hy\-per\-k\"ahler-type symplectic forms on $M$ follows from \cite[Thm. 9.2]{_EV-JTA_}.


\hfill


\section{Appendix}
\label{_Appendix_Section_}

In this section we discuss various well-known facts concerning
the dependence of the Hodge decomposition on the complex structure,
the deformations of complex structures, and the Moser method, that we used in the paper and
for most of which we could not find a direct reference in the literature.


\hfill


\subsection{Dependence of the Hodge decomposition on the complex structure}
\label{_dependence-of-space-of-holomorphic-1-forms-on-the-complex-structure_Subsection_}

In this section we prove two auxiliary facts about K\"ahler manifolds. These facts must have been known to the experts
but we have not been able to find a reference to them.

Assume $M$, $2n = \dim_\R M$, is a closed manifold admitting K\"ahler structures -- i.e., $\compK (M)\neq\emptyset$.
Recall that we use the $C^\infty$-topology on $\compK (M)$.

For each $l=0,1,\ldots, 2n$ denote by $A^l (M)$ the space of smooth complex-valued differential $l$-forms on $M$.
For each $I\in \compK (M)$ let
\[
A^l (M) = \oplus_{p+q=l} A_I^{p,q} (M)
\]
be the $(p,q)$-decomposition of $A^l (M)$.

Let
$\Omega_I^{1,0} (M)$
be the space of $I$-holomor\-phic 1-forms on $M$, with respect to $I$.

By Hodge theory, $\dim_\C H^1 (M;\C) = 2\dim_\C H_I^{1,0} (M) = 2\dim_\C H^0(\Omega_I^1 (M))$.
Let $m:=\dim_\C H^1 (M;\C)/2$.


\hfill


\proposition
\label{_dependence-of-space-of-holomorphic-1-forms-on-the-complex-structure_Proposition_}

Let $J\in \compK (M)$.

Then there exists a neighborhood $\cW$ of $J$ in $\compK (M)$ with the following property:

\smallskip
\noindent
For each $I\in \cW$ one can choose a basis $\theta_1 (I),\ldots,\theta_m (I)$ of $H^0(\Omega_I^1 (M))$
continuously with respect to $I$.


\hfill


\noindent
{\bf Proof of \ref{_dependence-of-space-of-holomorphic-1-forms-on-the-complex-structure_Proposition_}:}

Pick a $J$-invariant Riemannian metric $g$ on $M$. Then there exists a neighborhood $\cW$ of $J$ in $\compK (M)$ such that
for any $I\in\cW$ sufficiently close to $J$
we have an $I$-invariant Riemannian metric $(g+I^* g)/2$ that depends continuously on $I$ and coincides with $g$ for $I=J$.
Together with $I\in\cW$ such an $I$-invariant Riemannian metric defines a Hermitian metric $h_I$ on $(M,I)$ that depends continuously on $I$.

For each $I\in\cW$ and each $l=0,1,\ldots, 2n$ the Hermitian metric $h_I$ allows to define the Laplace-Beltrami operator on $A^l (M)$ that will be
denoted $\Delta_I^l: A^l (M)\to A^l (M)$. This is an elliptic operator. Denote the kernel of $\Delta_I^l$ (i.e., the space of harmonic $l$-forms with respect to $\Delta_I^l$) by $\cH_I^l (M)$. Then $\dim_\C \cH_I^l (M) = \dim_\C H^l (M;\C)$.
Also, $\cH_I^1 (M)\cap A_I^{1,0} (M) = \Omega_I^{1,0} (M)$.
These results follow from the classical
Hodge theory (see e.g. \cite[Sec. IV.4, IV.5]{_Wells_}).

Thus we have obtained that for each $l$ the operators $\Delta_I^l$ form a family of elliptic operators depending
continuously on $I\in \cW$
so that the kernels of these operators all have the same dimension. Then these kernels -- i.e., the spaces $\cH^l_I (M)$, $I\in \cW$ -- form a continuous vector bundle over $\cW$. This
follows from a standard argument essentially contained in \cite[Lem. 2.1]{_Atiyah-LNM103-1969_}. Namely, one passes to a Hilbert completion of
$A^l (M)$ and extends the continuous family of elliptic operators $\Delta_I^l$, $I\in \cW$, to a continuous family of Fredholm operators that, by elliptic regularity, have the same kernels, all of the same dimension. An easy functional analysis argument implies then that the kernels of these Fredholm operators form a continuous vector bundle (in fact, we only need the operators to be bounded with a closed image and these conditions are satisfied by any Fredholm operator).

Since the spaces $\cH_I^1 (M)$, $I\in \cW$, form a continuous vector bundle over $\cW$,
so do their $I$-invariant parts -- i.e., the spaces $\Omega_I^{1,0} (M)$ -- which are all of the same dimension $m$.
This readily implies the proposition.
\endproof


\hfill


\subsection{Deformation families of complex structures}
\label{_deformations-of-cs-Hodge-decomposition_Subsection_}

Let $(M,I)$ be a complex manifold.

Assume $\cX$, $\cB$ are connected (not necessarily Hausdorff)
smooth manifolds. Let $s_0\in \cB$ be a marked point.
Let $\cX\to \cB$ be a smooth map.

If $\cX\to \cB$ is a proper fiber bundle whose fibers are equipped with complex structures depending
smoothly on the fiber and the fiber over $s_0$ is
$(M,I)$, we say that $\cX\to (\cB,s_0)$ is a {\bf smooth deformation family}, or a
{\bf smooth deformation of $(M,I)$}.

If $\cX$, $\cB$ are complex manifolds and the map $\cX\to \cB$
is a proper holomorphic submersion, then it is a fiber bundle (by Ehresmann's lemma, the fibration admits {\it
smooth} local trivializations). The fibers of $\cX\to \cB$ are then closed
complex submanifolds of $\cX$ that are diffeomorphic to $M$. Denote
by $M_s$ the fiber over $s\in \cB$  and denote by $I_s$ the complex
structure on $M_s$ induced by the complex structure on $\cX$. Assume
that $I_{s_0} = I$. In this case we say that $\cX\to (\cB,s_0)$ is
a {\bf complex-analytic deformation family}, or a {\bf complex-analytic deformation of $(M,I)$}.

Given a complex structure $I$ on $M$,
a smoothly trivial smooth (respectively complex-analytic) deformation family
$M\times U\to U$ of $(M,I)$ over an open set $U\subset \C^m$ (respectively $U\subset \R^m$), $s_0=0\in U$,
is called a
{\bf local smooth (respectively complex-analytic) deformation} of $I$. It can be viewed as a smooth
(respectively complex analytic) family $\{ I_s\}_{s\in U}$, $I_0=I$.

Assume $\{ I_s\}_{s\in U}$, $I_0=I$, is a smooth local deformation of a K\"ahler-type complex structure $I$ on $M$.
The Kodaira-Spencer stability theorem \cite{_Kod-Spen-AnnMath-1960_} says that for any symplectic form $\omega$ on $M$ compatible with $I$ there exists a neighborhood $U'$ of $0$ in $U$ and a smooth family $\{ \omega_s\}_{s\in U'}$, $\omega_0=\omega$, of symplectic forms $\omega_s$ on $M$
compatible with $I_s$ for each $s\in U'$; in particular, each complex structure $I_s$, $s\in U'$, is of K\"ahler type.
It follows from the Kodaira-Spencer stability theorem (see \cite[Thm. 5.6]{_EV-JTA_}) that there exists a neighborhood $U''\subset U'$ of $0$ so that for each $s\in U''$ the cohomology class $[\omega]_{I_s}^{1,1}$ is K\"ahler with respect
to $I_s$.


\hfill


\proposition
\label{_compK-path-connected-Kahler-coh-classes_Proposition_}

The space $\compK (M)$ is locally $C^\infty$-path-connected.


\hfill


\noindent
{\bf Proof of \ref{_compK-path-connected-Kahler-coh-classes_Proposition_}:}

Kuranishi's theorem \cite{_Kuranishi-AnnMath-1962_}, \cite{_Kuranishi-Minneapolis-1964_}, \cite{_Kuranishi-Montreal-1969_}
implies that the $\Diff_0 (M)$-action on the space of complex structures on $M$ admits local slices that are finite-dimensional analytic spaces. This implies that the space of complex structures on $M$ is locally path-connected.
Therefore,
since $\Diff_0 (M)$ acts trivially on the cohomology and preserves $\compK (M)$, the Kodaira-Spencer stability theorem  yields that the space $\compK (M)$ is locally $C^\infty$-path-connected.
This proves the proposition.
\endproof


\hfill


\subsection{Relative version of Moser's method}
\label{_Relative-Moser_Subsection_}

Here we recall the following (well-known) relative version of Moser's method \cite{_Moser_}.


\hfill


\proposition
\label{_relative-version-of-Moser-method_Proposition_}

Let $N$ be a closed manifold.

\bigskip
\noindent
I. Assume that $X\subset N$ is a closed domain with a piecewise-smooth boundary, such that $H^1 (X;\R) = 0$.

Let $\{\Omega_t\}$, $t\in [0,1]$, be a smooth family of cohomologous symplectic forms on $N$, so that
all the forms $\Omega_t$, $t\in [0,1]$, coincide on a neighborhood of $X$.

Then there exists an isotopy $\{\phi_t : N\to N\}_{0\leq t\leq 1}$, $\phi_0 = Id$, such that for each $t\in [0,1]$ the diffeomorphism $\phi_t$ is identity on a neighborhood of $X$
and $\phi_t^* \Omega_0 = \Omega_t$.

\bigskip
\noindent
II. Assume that $X\subset N$ is a proper submanifold. Let $\{\Omega_t\}$, $t\in [0,1]$, be a smooth family of cohomologous symplectic forms on $N$, so that
all the forms $\Omega_t|_X$, $t\in [0,1]$, are non-degenerate.

Then there exists an isotopy $\{\phi_t : N\to N\}_{0\leq t\leq 1}$, $\phi_0 = Id$, such that $\phi_t (X)= X$
and $\phi_t^* \Omega_0 = \Omega_t$
for each $t\in [0,1]$.

\bigskip
\noindent
III. Assume that $N = \C P^n\times \C P^1$, $X_1 := \C P^{n-1}\times \C P^1\subset N$, $X_2 := \C P^n\times \textrm{pt}\subset N$. Let $\{\Omega_t\}$, $t\in [0,1]$, be a smooth family of cohomologous symplectic forms on $N$, so that
all the forms $\Omega_t|_{X_i}$, $i=1,2$, $t\in [0,1]$, are non-degenerate.

Then for every open neighborhood $\cU$ of $X_1\cap X_2$ there exists an isotopy $\{\phi_t : N\to N\}_{0\leq t\leq 1}$, $\phi_0 = Id$, such that
\begin{itemize}

\item{} $\phi_t (X_1\cup X_2)= X_1\cup X_2$,

\item{}  $\phi_t^* \Omega_0 = \Omega_t$ outside $\cU$
for each $t\in [0,1]$.

\end{itemize}


\hfill


\noindent
{\bf Proof of \ref{_relative-version-of-Moser-method_Proposition_}:}

Let us prove part I of the proposition.

The proof
follows the lines of the proof of the standard, absolute, version of the same result, as long as one can find a smooth family $\{ \lambda_t\}$, $t\in [0,1]$, of 1-forms on $N$ so that for all $t\in [0,1]$
\[
d\lambda_t = \frac{d}{dt} \Omega_t,
\]
and all the forms $\lambda_t$, $0\leq t\leq 1$, vanish on a neighborhood of $X$.

In order to find such a family, pick first a smooth family $\{ \mu_t\}$, $t\in [0,1]$, of 1-forms on $N$ satisfying
\[
d\mu_t = \frac{d}{dt} \Omega_t
\]
for all $t\in [0,1]$ (this can be done e.g. using Hodge theory -- see e.g. \cite[Thm. 3.2.4]{_McD-Sal-3_}). Note that all the forms $d\mu_t$, $0\leq t\leq 1$, vanish on a neighborhood of $X$.
Now, since $H^1 (X;\R) = 0$, it is easy to show that there exists a smooth family $\{ F_t\}$, $t\in [0,1]$,  of functions defined on a
neighborhood $U$ of $X$
so that $dF_t = \mu_t$ on $U$ for all $t\in [0,1]$. Extend these functions from possibly a smaller neighborhood $U'\subset U$ of the set $X$
to a smooth family of functions on $N$. By a slight abuse of notation, denote the latter family of functions on $N$ also by $\{ F_t\}$, $t\in [0,1]$. Then $\{ \lambda_t := \mu_t - dF_t\}$,  $t\in [0,1]$, is the wanted family of 1-forms on $N$.

This finishes the proof of part I.

Part II is proved in \cite[Cor. 4.1.B]{_McD-Polt_}.

Part III is proved in \cite[Prop. 4.1.C]{_McD-Polt_}.

This finishes the proof
of the proposition.
\endproof


\hfill


\subsection{Alexander's trick}
\label{_Alexander-trick_Subsection_}

Let $W\subset \R^{2n}$ be a compact domain with piecewise-smooth boundary that is starshaped with respect to the origin.


\hfill


\proposition
\label{_Alexander-trick_Proposition_}

Assume $f: W\to \R^{2n}$ is an embedding. Let $\iota: W\to \R^{2n}$ be the inclusion map.

The following claims hold:

\bigskip
\noindent
I. Let $U\subset \R^{2n}$ be an open neighborhood $W$ in $\R^{2n}$.  Assume that the embedding $f$ is symplectic.

If $f$ is sufficiently $C^1$-close to $\iota$,
then there exists a smooth family of symplectic embeddings $\{ f_t: W\to U\}_{0\leq t\leq 1}$ so that $f_0=\iota$, $f_1 = f$.
Moreover, all $f_t$, $t\in [0,1]$, can be made arbitrarily $C^\infty$-close to $\iota$, provided
$f$ is sufficiently $C^\infty$-close to $\iota$.

\bigskip
\noindent
II.
Assume that the embedding $f$ is holomorphic, $f(0)=0$, and $f(W)\subset W$.
Then there exists a smooth family of holomorphic embeddings $\{ f_t: W\to W\}_{0\leq t\leq 1}$ so that $f_0=\iota$, $f_1 = f$
and $f_t (0) = 0$ for all $t\in [0,1]$.


\hfill


\noindent
{\bf Proof of \ref{_Alexander-trick_Proposition_}:}

Let us prove part I.

If $f$ is sufficiently $C^1$-close to $\iota$,
then, composing $f$, if necessary, with a family of small affine symplectomorphisms staring at $Id$, we may assume without
loss of generality, that $f(0)=0$ and $df (0) = Id$.

Consider the following family of maps $\{ g_t\}_{0\leq t\leq 1}$ defined on $W$:
\[
g_t (x) := g (tx)/t,\ \textrm{if}\ t\in (0,1],
\]
\[
g_0 = \iota.
\]
If $f$ is sufficiently $C^1$-close to $\iota$,
then all the maps $g_t$ are symplectic embeddings taking values in $U$. The family $\{ g_t\}_{0\leq t\leq 1}$ may not be smooth at $t=0$ but for an appropriate smooth function $\chi: [0,1]\to [0,1]$, $\chi (0) = 0$, $\chi (1) = 1$,
the family $\{ g_{\chi (t)}\}_{0\leq t\leq 1}$ is a smooth family of symplectic embeddings $W\to U$ connecting $f$ and $\iota$.
The construction clearly shows that
all $f_t$, $t\in [0,1]$, can be made arbitrarily $C^\infty$-close to $\iota$, provided
$f$ is sufficiently $C^\infty$-close to $\iota$.
This finishes the proof of part I.

Let us prove part II.
For $C>0$ consider the following family of maps $\{ g_t\}_{0\leq t\leq 1}$ defined on $W$:
\[
g_t (z) := f (tz)/Ct,\ \textrm{if}\ t\in (0,1],
\]
\[
g_0 := df (0)/C.
\]
For a sufficiently large $C>0$ we have $\Image g_t\subset W$ for all $t\in [0,1]$.
The family $\{ g_t\}_{0\leq t\leq 1}$ may not be smooth at $t=0$ but for an appropriate smooth function $\chi: [0,1]\to [0,1]$, $\chi (0) = 0$, $\chi (1) = 1$,
the family $\{ g_{\chi (t)}\}_{0\leq t\leq 1}$ is a smooth family of holomorphic embeddings $(W,0)\to (W,0)$ connecting $g_0$ and $g_1 = f/C$.
If $C>0$ is sufficiently large, the linear holomorphic embedding $g_0$ can be easily connected by a smooth family of holomorphic embeddings $(W,0)\to (W,0)$ to $Id/C$ and the latter can be connected to $Id$.
\endproof


\hfill


\corollary
\label{_Kahler-type-embs-form-an-open-set_Corollary_}

Let $M^{2n}$ be a manifold and $\omega$ a K\"ahler-type symplectic form on $M$.
Let $f: W\to (M,\omega)$ be a K\"ahler-type/tame/$\varepsilon$-tame embedding (in the latter case we assume that $M$ is closed).

Then any symplectic embedding $f': W\to (M,\omega)$ sufficiently $C^1$-close to $f$ is of K\"ahler-type/tame/$\varepsilon$-tame,
and there exists $\phi\in\Symp_0 (M)$ such that $\phi\circ f = f'$. The symplectomorphism $\phi$ can be made arbitrarily $C^\infty$-close to $Id$ provided $f'$ is sufficiently $C^\infty$-close to $f$.


\hfill


\noindent
{\bf Proof of \ref{_Kahler-type-embs-form-an-open-set_Corollary_}:}

Assume $f$ extends to an $I$-holomorphic symplectic embedding $f: U\to (M,\omega)$, where $U$ is a neighborhood of $W$ in $\R^{2n}$
and $I$ is a complex structure on $M$ with the relevant properties as needed for a K\"ahler-type/tame/$\varepsilon$-tame embedding.

Let $f': W\to (M,\omega)$ be a symplectic embedding (extending to a symplectic embedding of an open neighborhood of $W$). If $f'$ is sufficiently $C^1$-close to $f$, then $f' (W)\subset f(U)$ and
$f^{-1}\circ f': W\to U$ is a symplectic embedding that can be made arbitrarily $C^1$-close to the identity. By part I of \ref{_Alexander-trick_Proposition_}, it can be connected to the inclusion $\iota: W\to U$ by a smooth family of symplectic
embeddings $W\to U$ arbitrarily $C^\infty$-close to the identity, provided $f'$ is sufficiently $C^\infty$-close to $f$. Composing these embeddings with $f$, we get a smooth family of symplectic embeddings $f_t: W\to (M,\omega)$, $0\leq t\leq 1$,
connecting $f$ and $f'$. These embeddings can be made arbitrarily $C^\infty$-close to $f$, provided $f'$ is sufficiently $C^\infty$-close to $f$.

By a standard result concerning extension of symplectic isotopies (proved as \cite[Thm. 3.3.1]{_McD-Sal-3_}),
there exists $\phi\in \Symp_0 (M,\omega)$ such that $\phi\circ f = f'$.
The symplectomorphism $\phi$ can be chosen to be arbitrarily $C^\infty$-close to the identity, if the symplectic
embeddings $f_t$, $0\leq t\leq 1$, are sufficiently $C^\infty$-close to $f$, which is true if $f'$ is sufficiently $C^\infty$-close to $f$.
The symplectic embedding $f'$ is of K\"ahler-type/tame/$\varepsilon$-tame because it is holomorphic with respect to the
complex structure $\phi_* I$ on $M$ satisfying the relevant properties.

This finishes the proof of the corollary.
\endproof


\hfill


\noindent {\bf Acknowledgements:} We are grateful to Y.Karshon, D.Novikov, H.Nuer, L.Polterovich and M.Prokhorova for useful discussions. We thank L.Polterovich and D.Salamon for comments on the preliminary version of this paper, and
 the anonymous referee for numerous comments and corrections.


\hfill


{\small
\noindent {\sc Michael Entov\\
Department of Mathematics,\\
Technion - Israel Institute of Technology,\\
Haifa 32000, Israel}\\
{\tt  entov@technion.ac.il}
}
\\

{\small
\noindent {\sc Misha Verbitsky\\
            {\sc Instituto Nacional de Matem\'atica Pura e
              Aplicada (IMPA) \\ Estrada Dona Castorina, 110\\
Jardim Bot\^anico, CEP 22460-320\\
Rio de Janeiro, RJ - Brasil}

\smallskip
\noindent
and

\smallskip
\noindent
Laboratory of Algebraic Geometry,\\
National Research University HSE, Faculty of Mathematics, \\
7 Vavilova Str., Moscow, Russia}\\
{\tt  verbit@mccme.ru}
}


\begin{thebibliography}{alpha}




\bibitem[AmV1]{_Am-Ver-JGF2015_}
Amerik, E., Verbitsky, M. {\em Teichm\"uller space for
hyperk\"ahler and symplectic structures},
J. Geom. Phys. \textbf{97} (2015), 44-50.

\bibitem[AmV2]{_Am-Ver-IMRN2015_}
Amerik, E., Verbitsky, M.,
{\em Rational curves on hyperkähler manifolds}, Int. Math. Res. Not. (2015), 13009-13045.

\bibitem[At1]{_Atiyah-analyt-surf-ProcRoySoc1958_} Atiya, M.F., {\em On analytic surfaces with double points}, Proc. Roy. Soc. London Ser. A \textbf{247} (1958), 237-244.

\bibitem[At2]{_Atiyah-LNM103-1969_}
Atiyah, M. F., {\em Algebraic topology and operators in Hilbert space}, in  {\em 1969 Lectures in Modern Analysis and Applications I, p. 101-121}, Springer, Berlin, 1969.


\bibitem[BaHPV]{_Barth-Hulek-Peters-vdVen_}
Barth, W., Hulek, K., Peters, C.A.M., Van de Ven, A.,
{\em Compact complex surfaces},
2nd edition, Springer-Verlag, Berlin, 2004.



\bibitem[Bea1]{_Bea1_}
 Beauville, A., {\em Vari\'et\'es K\"ahleriennes dont la premi\`ere classe de Chern est
nulle},  J. Diff. Geom. \textbf{18} (1983), 755-782.

\bibitem[Bea2]{_Geom-K3-Asterisque_}
Beauville, A., et al,
{\em G\'eom\'etrie des surfaces K3: modules et p\'eriodes.
Papers from the seminar held in Palaiseau, October
1981-January 1982}, Ast\'erisque \textbf{126} (1985).



\bibitem[Bes]{_Besse:Einst_Manifo_}
Besse, A., {\em Einstein Manifolds}, Springer-Verlag, New York,
1987.


\bibitem[Bir1]{_Biran-IMRN1996_} Biran, P.,
{\em Connectedness of spaces of symplectic embeddings}, Internat. Math. Res. Notices \textbf{10} (1996), 487–491.



\bibitem[Bir2]{_Biran-GAFA1997_} Biran, P.,
{\em Symplectic packing in dimension 4}, Geom. Funct. Anal. \textbf{7} (1997), 420-437.


\bibitem[Bir3]{_Biran-from-sympl-pack-to-alg-geom_} Biran, P.,
{\em From symplectic packing to algebraic geometry and back}, in {\em European Congress of Mathematics, Vol. II (Barcelona, 2000), 507-524}, Birkh\"auser, Basel, 2001.

\bibitem[Bo1]{_Bogomolov:decompo_}
Bogomolov, F. A., {\em On the decomposition of
K\"ahler manifolds with trivial canonical class}, Math. USSR-Sb.
{\bf 22} (1974), 580-583.


\bibitem[Bo2]{_Bogomolov_}
Bogomolov, F., {\em Hamiltonian K\"ahler manifolds}, Sov. Math. Dokl.
\textbf{19} (1978), 1462-1465.


\bibitem[BSV]{_Borisov-Salamon-Viaclovsky-Duke2011_} Borisov, L., Salamon, S., Viaclovsky, J., {\em Twistor geometry and warped product orthogonal complex structures}, Duke Math. J. \textbf{156} (2011), 125-166.

\bibitem[BLW]{_Borman-Li-Wu_} Borman, M.S., Li, T.-J., Wu, W., {\em Spherical Lagrangians via ball packings and symplectic cutting}, Selecta Math. (N.S.) \textbf{20} (2014), 261-283.

\bibitem[Bru]{_Brumfiel_} Brumfiel, G., {\em Homotopy equivalences of almost smooth manifolds}, Comment. Math. Helv. \textbf{46} (1971), 381-407.


\bibitem[Buch]{_Buchdahl-AnnInstFourier1999_} Buchdahl, N., {\em On compact Kähler surfaces}, Ann. Inst. Fourier \textbf{49} (1999), 287-302.

\bibitem[BurR]{_Burns-Rapoport_} Burns, D., Rapoport, M.,
{\em On the Torelli problem for k\"ahlerian K3 surfaces},
Ann. Sci. \'Ecole Norm. Sup. (4) \textbf{8} (1975), 235-273.


\bibitem[Cam]{_Campana:isotrivial_}
Campana, F., {\em Isotrivialit\'e de certaines familles
K\"ahl\'eriennes de vari\'et\'es non projectives}, Math. Z.
\textbf{252} (2006), 147-156.


\bibitem[CaDV]{_CDV:threefolds_}
Campana, F., Demailly, J.-P., Verbitsky, M., {\em Compact K\"ahler
3-manifolds without non-trivial
  subvarieties},
 Algebr. Geom. \textbf{1} (2014), 131-139.

\bibitem[CDP]{_Candelas_etc:rigid_}
Candelas, P., Derrick, E., Parkes, L., {\em Generalized
Calabi-Yau manifolds and the mirror of a rigid manifold,}
Nucl. Phys. B407 (1993), 115-154.



\bibitem[Cat]{_Catanese:moduli_}
Catanese, F., {\em A Superficial Working Guide to Deformations and
Mo\-du\-li}, in {\em Handbook of moduli, Vol. I, 161-215}, Int. Press,
Somerville, MA, 2013.

\bibitem[Ch]{_Chen-ProcAMS-2020_} Chen, W.,
{\em A characterization of the standard smooth structure of K3 surface},
Proc. AMS \textbf{148} (2020), 2707-2716.

\bibitem[Cr]{_Crist-Gardiner-JDG2019_}
Cristofaro-Gardiner, D.,
{\em Symplectic embeddings from concave toric domains into convex ones. With an appendix by Cristofaro-Gardiner and Keon Choi}, J. Diff. Geom. \textbf{112} (2019), 199-232.

\bibitem[Dem]{_Demailly_Compl_Analytic_Diff_Geom_} Demailly, J.-P., {\em Complex analytic and differential geometry}, 2012. Avai\-la\-ble at \url{https://www-fourier.ujf-grenoble.fr/~demailly/manuscripts/agbook.pdf}.


\bibitem[DP]{_Dem-Paun_}
Demailly, J.-P., Paun, M., {\em Numerical characterization of the
K\"ah\-ler cone of a compact K\"ahler manifold}, Ann. of Math. \textbf{159} (2004), 1247-1274.


\bibitem[Dmz]{_Demazure-LNM777-1980_}
Demazure, M., {\em Surfaces de del Pezzo, II, III, IV, V}, in
{\em Séminaire sur les Singularités des Surfaces, Palaiseau, 1976–1977, pp.23-69,
edited by M.Demazure, H.C.Pinkham and B.Teissier}, Lecture Notes in Mathematics \textbf{777}. Springer, Berlin, 1980.


\bibitem[Dol]{_Dolgachev-book_}
Dolgachev, I.,
{\em Classical algebraic geometry. A modern view}. Cambridge Univ. Press, Cambridge, 2012.


\bibitem[Don]{_Donaldson:ellipt_}
Donaldson, S. K.,
{\em Two-forms on four-manifolds and elliptic equations}, in {\em
Inspired by S. S. Chern, 153-172},
World Sci. Publ., Hackensack, NJ, 2006.


\bibitem[Ec]{_Eckl2017_}
Eckl, T.,
{\em K\"ahler packings and Seshadri constants on projective complex surfaces}, Differential Geom. Appl. \textbf{52}
(2017), 51-63.

\bibitem[EV1]{_EV-JTA_}
Entov, M., Verbitsky, M.,
{\em Unobstructed symplectic packing for tori and hyper-K\"ahler manifolds}, J. Topol. Anal. \textbf{8} (2016), 589-626. {\em Erratum}, J. Topol. Anal. \textbf{11} (2019), 249-250.


\bibitem[EV2]{_EV-Selecta_}
Entov, M., Verbitsky, M.,
{\em Unobstructed symplectic packing by ellipsoids for tori and hyperk\"ahler manifolds}, Selecta Math. (N.S.) \textbf{24} (2018), 2625-2649.

\bibitem[EV3]{_EV-K3_}
Entov, M., Verbitsky, M.,
{\em Rigidity of Lagrangian embeddings into symplectic tori and K3 surfaces}, preprint, arXiv:2105.05971, 2021.


\bibitem[FiS]{_Fintushel-Stern-InventMath1998_} Fintushel, R., Stern, R., {\em Knots, links, and 4-manifolds},
Invent. Math. \textbf{134} (1998), 363-400.


\bibitem[Fl]{_Fleming2021_}
Fleming, A.,
{\em K\"ahler packings of projective complex manifolds}, European J. of Math. \textbf{7} (2021), 1382-1400.

\bibitem[FriM]{_Friedman-Morgan-JDG88_}
Friedman, R., Morgan, J. W., {\em On the diffeomorphism types of certain algebraic surfaces. I}, J. Diff. Geom. \textbf{27} (1988), 297-369.

\bibitem[FriQ]{_Friedman-Qin-InvMath-1995_}
Friedman, R., Qin, Z., {\em
On complex surfaces diffeomorphic to rational surfaces}, Invent. Math. \textbf{120} (1995), 81-117.

\bibitem[FroN]{_Frolicher-Nijenhuis-PNAS_}
Fr\"olicher, A., Nijenhuis, A., {\em A theorem on stability of complex structures},
Proc. Nat. Acad. Sci. U.S.A. \textbf{43} (1957), 239-241.

\bibitem[Fu]{_Fujiki:HK_}
Fujiki, A., {\em On the de Rham Cohomology Group of a Compact
K\"ahler Symplectic Manifold}, Adv. Stud.
Pure Math. \textbf{10} (1987), 105-165.

\bibitem[FuS]{_Fujiki-Schumacher-PublRIMS1990_} Fujiki, A., Schumacher, G., {\em The moduli space of extremal compact Kähler manifolds and generalized Weil-Petersson metrics}, Publ. RIMS \textbf{26} (1990), 101-183.

\bibitem[GV]{_Goodman_Wallach_}
Goodman, R., Wallach, N.R., {\em Symmetry, representations, and
invariants}, Springer, Dordrecht, 2009.



\bibitem[Gro]{_Gromov_} Gromov, M., {\em Pseudoholomorphic curves in
symplectic manifolds}, Invent. Math. \textbf{82} (1985), 307-347.

\bibitem[GHJ]{_Joyce-lectures-book_}
Gross, M., Huybrechts, D., Joyce, D., {\em Calabi-Yau manifolds and related geometries. Lectures from the Summer School held in Nordfjordeid, June 2001}. Springer-Verlag, Berlin, 2003.

\bibitem[Hi]{_Hironaka-AnnMath1964-I_} Hironaka, H., {\em Resolution of singularities of an algebraic variety over a field of characteristic zero. I}, Ann. of Math. \textbf{79} (1964), 109-203.

\bibitem[HiK]{_Hirzebruch-Kodaira-JMPA1957_} Hirzebruch, F., Kodaira, K.,
{\em On the complex projective spaces},
J. Math. Pures Appl. (9) \textbf{36} (1957), 201-216.

\bibitem[H\"oP]{_Horing-Peternell-Sal-Lake-City_}
H\"oring, A., Peternell, T.,
{\em Bimeromorphic geometry of Kähler threefolds}, in {\rm Algebraic geometry: Salt Lake City 2015, 381-402}, AMS, Providence, RI, 2018.


\bibitem[HuKP]{_Huckleberry-Kebekus-Peternell_} Huckleberry, A.T., Kebekus, S., Peternell, T.,
{\em Group actions on ${\mathbb S}^6$ and complex structures on ${\mathbb P}_3$},
Duke Math. J. \textbf{102} (2000), 101-124.

\bibitem[Huyb]{_Huybrechts-K3-book_}
Huybrechts, D., {\em Lectures on K3 surfaces}. Cambridge U. Press, Cambridge, 2016.


\bibitem[Ko]{_Kod-emb_}
Kodaira, K., {\em On K\"ahler varieties of restricted type (an intrinsic characterization of algebraic varieties)},
Ann. of Math. \textbf{60} (1954), 28-48.

\bibitem[KoS1]{_Kod-Spen-AnnMath-I-II-1958_}
Kodaira, K., Spencer, D.C., {\em On deformations of complex analytic structures. I, II}, Ann. of Math. \textbf{67} (1958), 328-466.


\bibitem[KoS2]{_Kod-Spen-AnnMath-1960_}
Kodaira, K., Spencer, D.C., {\em On deformations of complex analytic structures. III. Stability theorems for complex structures},
Ann. of Math. \textbf{71} (1960), 43-76.


\bibitem[KrS]{_Kreck-Su-MCG2020_} Kreck, M., Su, Y.,
{\em Mapping class group of manifolds which look like 3-dimensional complete intersections}, arXiv:2009.08054, 2020.

\bibitem[Ku1]{_Kuranishi-AnnMath-1962_}
Kuranishi, M., {\em On the locally complete families of complex analytic structures}, Ann. of Math. \textbf{75} (1962), 536-577.

\bibitem[Ku2]{_Kuranishi-Minneapolis-1964_}
Kuranishi, M., {\em
New proof for the existence of locally complete families of complex structures}, in
{\em Proc. Conf. Complex Analysis (Minneapolis, 1964), pp. 142-154}, Springer, Berlin, 1965.

\bibitem[Ku3]{_Kuranishi-Montreal-1969_}
Kuranishi, M.,
{\em Deformations of compact complex manifolds}, in {\em Séminaire de Mathématiques Sup\'erieures, No. 39 (Été 1969)}, Les Presses de l'Université de Montréal, Montreal, Quebec, 1971.

\bibitem[Kur]{_Kurnosov:tori_}
Kurnosov, N.,
{\em Absolutely trianalytic tori in the generalized Kummer variety},
Advances in Mathematics \textbf{298} (2016), 473-483.

\bibitem[Lal]{_Lalonde-MathAnn-1994_}
Lalonde, F., {\em Isotopy of symplectic balls, Gromov's radius and the structure of ruled symplectic 4-manifolds}, Math. Ann. \textbf{300} (1994), 273-296.


\bibitem[LalM]{_Lal-McD-NewtonInst1996_} Lalonde, F., McDuff, D., {\em $J$-curves and the classification of rational and ruled symplectic 4-manifolds}, in {\em Contact and symplectic geometry (Cambridge, 1994), 3-42}, Cambridge Univ. Press, Cambridge, 1996.

\bibitem[Lam]{_Lamari-AnnInstFourier1999_}
Lamari, A., {\em Courrants k\"ahl\'eriens et surfaces compactes}, Ann. Inst. Fourier \textbf{49} (1999), 263-285.

\bibitem[LatMS]{_LMcDS_}
Latschev, J., McDuff, D., Schlenk, F., {\em The Gromov width of
4-dimen\-sional tori}, Geom. and Topol. \textbf{17} (2013), 2813-2853.


\bibitem[Laz]{_Lazarsfeld-MRL_} Lazarsfeld, R., {\em Lengths of periods and Seshadri constants of abelian varieties}, Math. Res. Lett. \textbf{3} (1996), 439-447.


\bibitem[Lee]{_Lee_} Lee, J.-L., {\em Non-Degeneracy of 2-Forms and Pfaffian}, Symmetry {\bf 2020}, 12(2), 280. Available at \url{https://doi.org/10.3390/sym12020280}.


\bibitem[LiWo]{_Libgober-Wood-JDG1990_}
Libgober, A., Wood, J.W.,
{\em Uniqueness of the complex structure on Kähler manifolds of certain homotopy types}, J. Diff. Geom. \textbf{32} (1990),
139-154.

\bibitem[LuWa]{_Luef-Wang_} Luef, F., Wang, X.,
{\em Gaussian Gabor frames, Seshadri constants and generalized Buser-Sarnak invariants}, preprint, arXiv:2107.04988, 2021.



\bibitem[McD1]{_McD-Topology1991_} McDuff, D., {\em Blow ups and symplectic embeddings in dimension 4}, Topology \textbf{30} (1991), 409-421.

\bibitem[McD2]{_McD-Cambridge1990_} McDuff, D.,
{\em Remarks on the uniqueness of symplectic blowing up}, in {\em Symplectic Geometry, ed.
by D. Salamon}, Cambridge Univ. Press, Cambridge,
1993, 157-167.

\bibitem[McD3]{_McD-connectedness_}
McDuff, D., {\em From symplectic deformation to isotopy}, in
{\em Topics in symplectic 4-manifolds (Irvine, CA, 1996), 85-99},
Int. Press, Cambridge, MA, 1998.

\bibitem[McD4]{_McDuff-ellipsoids-JT2009_}
McDuff, D.,
{\em Symplectic embeddings of 4-dimensional ellipsoids}, J. Topol. \textbf{2} (2009), 1-22.

\bibitem[McDO]{_McD-Opshtein_} McDuff, D., Opshtein, E., {\em Nongeneric J-holomorphic curves and singular inflation}, Algebr. Geom. Topol. \textbf{15} (2015), 231-286.

\bibitem[McDP]{_McD-Polt_}
McDuff, D., Polterovich, L., {\em Symplectic packings and algebraic geometry.
With an appendix by Yael Karshon}, Invent. Math. \textbf{115} (1994), 405-434.

\bibitem[McDS]{_McD-Sal-3_} McDuff, D., Salamon, D.,
{\em Introduction to symplectic topology},
3nd edition, Oxford Univ. Press, Oxford, 2017.

\bibitem[Mee]{_Meersseman-JEP2019_} Meersseman, L.,
{\em The Teichm\"uller and Riemann moduli stacks}, J. \'Ec. polytech. Math. \textbf{6} (2019), 879-945.

\bibitem[Mos]{_Moser_} Moser, J., {\em On the volume elements on a
manifold}, Trans. AMS \textbf{120} (1965), 288-294.





\bibitem[OV]{_Ornea-Verb_}
Ornea, L., Verbitsky, M., {\em Embeddings of compact Sasakian manifolds}, Math. Res. Lett. \textbf{14} (2007), 703-710.


\bibitem[Rag]{_Rag_} Ragozin, D.L.,
{\em A normal subgroup of a semisimple Lie group is closed}, Proc. of
AMS \textbf{32} (1972), 632-633.

\bibitem[Sal]{_Salamon-K3_} Salamon, D.,
{\em Notes on the Teichm\"uller space of K3}, preprint,
\url{https://people.math.ethz.ch/~salamon/PREPRINTS/atiyahflop.pdf}, 2018.


\bibitem[Sei]{_Seidel-2000_} Seidel, P.,
{\em Graded Lagrangian submanifolds}, Bull. Soc. Math. France \textbf{128} (2000), 103-149.

\bibitem[ShS]{_Sheridan-Smith-JAMS2020_}
Sheridan, N., Smith, I., {\em Symplectic topology of $K3$ surfaces via mirror symmetry}, J. of the AMS \textbf{33} (2020), 875-915.

\bibitem[Siu]{_Siu_}
Siu, Y.T. {\em Every K3 surface is K\"ahler}, Invent. Math. \textbf{73} (1983), 139-150.

\bibitem[Sm]{_Smirnov-GAFA2022_}
Smirnov, G.,
{\em Symplectic mapping class groups of K3 surfaces and Seiberg-Witten invariants}, Geom. Funct. Anal. \textbf{32} (2022), 280-301.

\bibitem[SolV]{_SV:k-symplectic_}
Soldatenkov, A., Verbitsky, M.,
{\em k-symplectic structures and absolutely trianalytic
subvarieties in hyperk\"ahler manifolds},
J. Geom. Phys. \textbf{92} (2015), 147-156

\bibitem[Som]{_Sommese-Quaternionic-MathAnn-1975_} Sommese, A.J.,
{\em Quaternionic manifolds}, Math. Ann. \textbf{212} (1974/75), 191-214.

\bibitem[Ste]{_Steffens-MathZ_} Steffens, A., {\em Remarks on Seshadri constants}, Math. Z. \textbf{227} (1998), 505-510.

\bibitem[StT]{_Streets-Tian_}
Streets, J., Tian, G.,
{\em K\"ahler stability of symplectic forms}, preprint, arXiv:2202.04564, 2022.

\bibitem[Tau1]{_Taubes-MRL1995_} Taubes, C.,
{\em The Seiberg-Witten and Gromov invariants}, Math. Res. Lett. \textbf{2} (1995), 221-238.


\bibitem[Tau2]{_Taubes-JDG1996_} Taubes, C.,
{\em Counting pseudo-holomorphic submanifolds in dimension 4}, J. Diff. Geom. \textbf{44} (1996), 818-893.

\bibitem[Tau3]{_Taubes-book-2000_} Taubes, C.,
{\em Seiberg Witten and Gromov invariants for symplectic 4-manifolds}. International Press, Somerville, MA, 2000.




\bibitem[Tia]{_Tian_}
Tian, G.,
{\em Smoothness of the universal deformation space of compact Calabi-Yau manifolds and its Petersson-Weil metric}, in
{\em Mathematical aspects of string theory (San Diego, Calif., 1986), 629-646},
World Sci. Publishing, Singapore, 1987.


\bibitem[Tod]{_Todorov_} Todorov, A.N., {\em The Weil-Petersson geometry of the moduli space of $SU(n\geq 3)$ (Calabi-Yau) manifolds. I}, Comm. Math. Phys.
\textbf{126} (1989), 325-346.

\bibitem[Tos]{_Tosatti-CPn_} Tosatti, V., {\em Uniqueness of $\C P^n$}, Expo. Math. \textbf{35} (2017), 1-12.

\bibitem[Tr]{_Trusiani-AIF2021_} Trusiani, A., {\em Multipoint Okounkov bodies}, Ann. Inst. Fourier (Grenoble) \textbf{71} (2021), 2595-2646.

\bibitem[Ue]{_Ueno-LNM1975_}
 Ueno, K., {\em Classification theory of algebraic varieties and compact complex spaces. Notes written in collaboration with P. Cherenack}. Lecture Notes in Mathematics, Vol. 439. Springer-Verlag, Berlin-New York, 1975.

\bibitem[Ver1]{_V:Mirror_}
Verbitsky, M.,
{\em Mirror Symmetry for hyperk\"ahler manifolds}, in {\em Mirror symmetry, III (Montreal, PQ, 1995), 115-156}, AMS, Providence, RI, 1999.

\bibitem[Ver2]{_Verbitsky:Symplectic_II_} Verbitsky, M., {\em Hyperk\"ahler embeddings and holomorphic symplectic geometry II}, Geom.
and Funct. Analysis \textbf{5} (1995), 92-104.

\bibitem[Ver3]{_V:applications_}
Verbitsky, M., {\em Cohomology of compact hyper-Kähler manifolds and its applications}, Geom. Funct. Anal. \textbf{6} (1996), 601-611.


\bibitem[Ver4]{_Verbitsky:Deforma_}
Verbitsky, M., {\em Deformations of
trianalytic subvarieties of
hyperk\"ahler manifolds}, Selecta Math. \textbf{4} (1998), 447-490.

\bibitem[Ver5]{_Verb-trianalytic-GAFA1998_} Verbitsky, M.,
{\em Trianalytic subvarieties of the Hilbert scheme of points on a K3 surface}, Geom. Funct. Anal. \textbf{8} (1998), 732-782.

\bibitem[Ver6]{_V-Duke_}
Verbitsky, M.,
{\em A global Torelli theorem for hyperk\"ahler manifolds},
 Duke Math. J. \textbf{162} (2013), 2929-2986.


\bibitem[Voi]{_Voisin-book_} Voisin, C., {\em Hodge theory and complex algebraic geometry I,II}.  Cambridge Univ. Press, Cambridge, 2002.

\bibitem[Wel]{_Wells_} Wells, R.,
{\em Differential analysis on complex manifolds}.
Springer, New York, 2008.


\bibitem[Weyl]{_Weyl_} Weyl, H.,
{\em The classical groups}.
Princeton University Press, Princeton, NJ, 1946.

\bibitem[WN1]{_WittNystrom1_} Witt Nystr\"om, D., {\em Canonical growth conditions associated to ample line bundles}, Duke Math. J. \textbf{167} (2018), 449-495.


\bibitem[WN2]{_WittNystrom2_} Witt Nystr\"om, D., {\em Okounkov bodies and the K\"ahler geometry of projective manifolds}, preprint, arXiv:1510.00510, 2015.

\bibitem[Yale]{_Yale Automorphisms of the Complex Numbers_} Yale, P.B.,
{\em Automorphisms of the Complex Numbers},
Math. Magazine \textbf{39} (1966), 135–141.


\bibitem[Yau1]{_Yau-PNAS1977_} Yau, S.-T.,
{\em Calabi's conjecture and some new results in algebraic geometry}, Proc. Nat. Acad. Sci. U.S.A. \textbf{74} (1977), 1798-1799.

\bibitem[Yau2]{_Yau-CPAM-1978_}
Yau, S.-T.,
{\em On the Ricci curvature of a compact K\"ahler manifold
  and the complex Monge-Amp\`ere equation. I},
Comm. Pure Appl. Math. \textbf{31} (1978), 339-411.





\end{thebibliography}
\end{document}